\documentclass[a4paper,11pt]{amsart}
\usepackage{amsmath,amssymb,amsfonts,amsthm}
\usepackage[top=25mm,bottom=25mm,left=25mm,right=25mm]{geometry}
\usepackage{booktabs,mathtools,mathrsfs}
\usepackage{needspace}
\usepackage{tikz,float,lscape,afterpage}
\usepackage[hidelinks,unicode]{hyperref}
\expandafter\providecommand\csname subjclassname@2020\endcsname
  {\textup{2020} Mathematics Subject Classification}
\theoremstyle{plain}
\newtheorem{theorem}{Theorem}[section]
\newtheorem{proposition}[theorem]{Proposition}
\newtheorem{lemma}[theorem]{Lemma}
\newtheorem{corollary}[theorem]{Corollary}
\newtheorem{conjecture}[theorem]{Conjecture}
\theoremstyle{definition}
\newtheorem{definition}[theorem]{Definition}
\newtheorem{example}[theorem]{Example}
\newtheorem{question}[theorem]{Question}
\newtheorem{remark}[theorem]{Remark}

\numberwithin{equation}{section}
\allowdisplaybreaks[4]
\newcommand{\tr}{\operatorname{tr}}
\newcommand{\arcosh}{\operatorname{arcosh}}
\newcommand{\SL}{\mathrm{SL}}
\newcommand{\PSL}{\mathrm{PSL}}
\newcommand{\HH}{\mathbb H}
\newcommand{\EE}{\mathcal E}
\newcommand{\OO}{\mathcal O}
\newcommand{\kk}{\mathbf k}
\newcommand{\ZZ}{\mathbb Z}
\newcommand{\RR}{\mathbb R}
\newcommand{\LL}{\mathscr L}
\newcommand{\PP}{\mathfrak P}
\newcommand{\Spec}{\mathscr S}
\tikzset{tree example vertex/.style={font=\small,align=center,inner xsep=2pt,inner ysep=3pt}}
\newcommand{\treeexampleedges}{\draw (treeRoot.east) -- (treeOne.west) node[midway,above] {$1$};
 \draw (treeRoot.east) -- (treeTwo.west) node[midway,above] {$2$};
 \draw (treeRoot.east) -- (treeThree.west) node[midway,above] {$3$};
 \foreach \treeParent/\treeLeft/\treeRight/\treeLeftLabel/\treeRightLabel in
  {treeOne/tree12/tree13/2/3,treeTwo/tree21/tree23/1/3,treeThree/tree31/tree32/1/2} {
  \draw (\treeParent.east) -- (\treeLeft.west) node[midway,above] {$\treeLeftLabel$};
  \draw (\treeParent.east) -- (\treeRight.west) node[midway,below] {$\treeRightLabel$};
 }
 \foreach \treeLeaf in {tree12,tree13,tree21,tree23,tree31,tree32} {
  \draw[densely dotted] (\treeLeaf.east) -- ++(0.4,0.25);
  \draw[densely dotted] (\treeLeaf.east) -- ++(0.4,-0.25);
 }
}
\title[GENERALIZING THE GASTER--LOUSTAU THEOREM]{Generalizing the Gaster--Loustau Theorem via Generalized Discrete Markov Spectra}
\author[Yasuaki Gyoda]{Yasuaki Gyoda}
\address[Yasuaki Gyoda]{Institute for Advanced Research, Nagoya University, Furo-cho, Chikusa-ku, Nagoya-shi, 464-8601, Japan}
\email{ygyoda@math.nagoya-u.ac.jp}
\hypersetup{pdftitle={Generalizing the Gaster–Loustau Theorem via Generalized Discrete Markov Spectra},pdfauthor={Yasuaki Gyoda}}
\date{}
\usepackage{lmodern}
\newbox\sourcestorusbox
\newcommand{\declaresourcestorus}[2]{\begingroup\pdfcompresslevel=0\immediate\pdfobj stream attr{/ShadingType 4 /ColorSpace /DeviceRGB /BitsPerCoordinate 16 /BitsPerComponent 8 /BitsPerFlag 8 /Decode [-1.65 1.65 -1.25 1.25 0 1 0 1 0 1] /Filter [/ASCIIHexDecode /FlateDecode]}{#2 >}\edef\sourcestorusobject{\the\pdflastobj}\setbox\sourcestorusbox=\hbox to330bp{\vbox to250bp{\vfil\pdfliteral{q 100 0 0 100 165 125 cm /Torus sh Q}}\hfil}\immediate\pdfxform resources{/Shading << /Torus \sourcestorusobject\space 0 R >>}\sourcestorusbox\expandafter\xdef\csname sourcestorus#1\endcsname{\the\pdflastxform}\endgroup}
\newcommand{\sourcestorusimage}[3]{\resizebox{#2}{#3}{\pdfrefxform\csname sourcestorus#1\endcsname}}
\declaresourcestorus{plain}{
78DA64DD77BC55D5D53560F29A98C4682489C618357663124D2C31B1448DB1F78222F6DE51B197602422082ABDF7DE7BEFBD
09D2140111A4A32022F6E81BCBF7CE7DD6B3D8F2FDB57FE7DCBDD79E73CC31C69C6BDF73CFAD52F79E4E03860F1F5EA54E95
4E13D3715AF975EDD6959FD77EA2D3D038D65FDCE9C938BAAE6E83CAFBCE7BEA4F959F3FF971A77F56D6D9DE62C68C19559E
9EBDFD91B973E7C6EB46AFBCF24A95178FEF746879BD17F6A89C5FEFEAED4F4E9F3E3D9F57EFA4ED8FCC993327AFFBD8FCCA
758FD5AA1C9FEEBBFDC99933675679ACE9F6C3E7CD9B97EFF3E2C0ED474E9B362D5FFFE283DB0F8F63A3A61D07C575EEFF62
93B45EB5ED47CE9A35AB4AAD5B3EE8F9EAABAFE6F5DCEFE1072BD73DF459C77E716CF2BF1FF49C3A756A5EB749E70F7ACE9E
3DBB4AFBC7DEFF7ACA9429555A9DB66D6BFCBCE5CC6DDB8B757F50F9F97DA76CDB1EEBBBCFBDD3B66D8DEBEF78F1FDAF23DF
3B0FD9B6EBFCF9F3F3F5D6B54E9D2A9DF72BE3AD4E752FEAFC93EFBC4EF5A93DA2F23EFC5C57BF4DA715E5F3D4413CE2106F
978E5B574D9E3C39E7D7EEED6DBB463CF055AFFF3B1EBA60C1822A4FDD5D591F7EEEEF3C757AFADDEDBB2D5CB8B04AC3D59D
5E28C791EB93F0759F7ABFDABE5BC4250EB888078E37DFB67555AC7FCBBE5BDF8F78DC2FD7F5F50F26C57D1F3FB8725FF515
6FA3033BED5D8E036FD4435C2F1EFDC1A4886707FE95FAECD053AA57AA03FCC58107B5667C70FAA2458B76D425ADF7FFD531
E10A4FF9CA53DD1E5E5B89FFC1FD3AEE5DCE2FE39F7057AF8C7FBABF3A94EA5B9C07FFA6353A3E57C6A5E9C6CA7D7A0E79AF
41F0045F3ABFB3F5FDD0A73AC1AFC9DF3E383DF2968F3AC03FD72BD549BCD6556FEB679EA4F8D405BE787CDFE46D8D172F5E
9CEB2B6F757EB04F25AFFB5EEB7043193F7CBAEEB2F71A44DDAE3FF4BDAEDFE14FC23BF3487D535DE591F993E283438BFA1D
1697716CF9C7CAFD33CE295E78CB1FAE2DCFD8D638D611A7F8F022F321E527AF8C4FC245BCF0C2A73B67BE3FF0B5D75EAB52
F3A94A9C77B768FF7C19A75CDF948F78FA0DDCF23F93264DCABCE8B1EABDAE51AF9C6FCA6F67FEA8ABF35BF76BFF4D199736
352BF7FF3F7DAD2EEB83DEEA1EDF79E1775E273DC95F1EE2AF3DBF723E3E5AAFFEED9DBB94AFA70FF1A943BB53DE1F18F5CE
B825BCE08A37D5AB6DF99FF8F9D5676E392470C657FAAA53E583AB5E7FFDF52A4F75A8DC17BFC4E53C7AF8E76E1FFC2DCE6F
D8BCF3A5E5F8F0FA9EEA15DCEEB8BDDD4BE5BAE19FFB3FFFD9077F8B7833AF52FDE48507F874CBCAAD6DE3BEEAA29EEA23BE
ECB3FB6CFB78C99225551EBFBF12273ECA0F3E70C1DF7693DBEF51AE77FB0E953C5E1ED7696D394FFAC03FF935DCBCEDE332
6FE10187124F567FD7BF13AF122FF0413E74F2D0EFB73DF1C61B6FECE0495A0F5F32AF523DD531E39970C49747CEA9E485A7
70A2237AE83B71CB213157E1031EE04FE6438A0B2F76F0AD729EBA5B17AF3B1FB6B56DD4A3C98A4E5794F1A5D7C14DDEBD30
E2119778F0461D1A2FDDF644F03FF339D55FDDD5178FD4197FF0467E779E5139FFD673DA3E5CE675D663E2ABF8F1409EF890
7D2EE9B3E6B9EFBF1A75CCBC4AF8E117DEE3E50DDF7FEFCA880B9FC4259E5A475470C33BF594BF3AAB2F5EE193BC3B4E6FFB
7599F79DC655F2E62397D77AF7C2C8EF8ADAEF3E1679C15D9EF06FBEBD63AF112346E4FAD1ABFACA5B9DE1A89E2DC6BEFF6A
AC0307F9E3AF78C441C7F284071CE00F7779CB172FE4AD4E747257ADAD1F2C5DBAB4CAFD553BF62AEB1ADEE296BFF86F3BB5
B2FE4D27B539A5CC23FCA583EEFF7DEFCAE8A35917894FDE8727FCE0DAEACB0E2BE27DF8E2357EB8BFFBD291FB0FAAFBEE63
155D77FD7EACC33FF851DDFDBB7C3A72E4C81DAF93DFC0171E70E832AACDE2326FBACEAEE45D7B43974FCBBEE43EF5FFDA65
F6A851A3F2BAFC43FE78D1B6C9D60F024F75510FF51C71D5A6DB232FF9C98B7EF0A3C6999BD7068FE8932FD5A9B2EDF065CB
965579EA952EB3CBBA11AFF3F8C23FFFB86D8FE5CB975769787197A7468F1E9DE3A6E37B77EDB0A25C077C910F7EC88B7EC4
F5FC6BDBF68878F00FEF6FEDFBDEFA785FDDF142FDC597FBDF4DEFCF8A38D505DFD5E7F1E15D9E2AFB83BCE1062FFA6BFB59
FB0F830F78463778206FF9BE7C63973DC68C199371A17B3A9177C3AEEFCF8ABCF8CB45FD37DD1EF7BFF8B54DED229F9B2FA8
C47FFD19ADDE29EB69073F2BFCDDD157139F13EFF00D2EF4FF50EBF72F7EF3CD37B34EDCC7FA999FE93E78AA2EEA818F99E7
8947F8F348B72E7B9475036FF9CA939FD127FDF799B8796D1C33DF525C7897F9997889CFAEA7A74E57BFB73E74401FEE63FD
1DBAA85C8F874DEEEEFCDCD8B16373FDF850F756AD5F28EBBCC7E44A7DF05B9D1BDFF3FEC5810FFEE02F1E65DD253EE2E15D
BBB5FFB0EC63F494799EF82D6F79D2957CB3BE523EF886077CA4E6B2ADBD830F99B709EF9D752D2EF1D0279DDCF0F1963EB1
7EAD369D9F2BEB43FD3B6C68B72570A01FBACC7A4C7CC9BAC4DBC45779F33DF51CFEC74DED0ADEA4BAD0BDFAF0413CABD6FE
DDD322EFACC7848B7A35BFBBD35BE3C68DCBF5CE7D081F124E780177F56FF1BBADBD639D31FFDD70C9C4891373BCE2148F38
F84DC633E148675947095778D2ED0DF756F2E613B9BEA9AEF0833F3EAA031ED0F3DD87BCB778C58A1555EE6FD3E9AD72FDD5
F39C95EB0745DF3FB3CEBADD22DFF3BEBFE192C8E7F6FDDB6D29C785BF5917095FF8D01FBD76DF634B9F58C76B3CF67EAE4B
AA83FAE0A5FBBB6FABEB3BED15E7AB1FFDE18DFB0FBCE3DDD32ACF4DBB358875F9193FE5CBF23CFF9C0D4DA25E757FD4EDA9
F1E3C7E7F3F865E7296D67854FE1373E665D245CE159FB3FDD9E0A3EF103F7AFBF67B78B264C9890D7E57F70C1BB366FBFB7
38F0C42FF9C94B7DD5152FE816FFAEB96CF3A9A12FF9CA934FF3113E59A7CAD6D62B57AEACF2D487DD2E0ABFA06379F4BAA9
D54165BFE493AEE75FFFACB6F5D95867D2D4B537C43A7433F1D975BBC5B1E12FBA7D2F8EF2E74FF7F6E8B457B9FE784ACFE2
7D7EE0D667DF7AEBAD1D3A4DFC831BDED3DB6DCF6C793570CAFC49BCC32371F3A3C7FB6D3D29E257773ACBFEF67EB7EFC56B
3A8613FCE14EF76DCFEA786CF00A5FE9957FA9C7A57B6CBABF98CFE934E1018797AB766D13BE032FBEC52FAEED5AA9CF3537
B57CB4984F923EE1D7B0FAD69302B75BFEDE7656D967E85B9DF06EF4B80D4DE2F70A3B7451D1D38E3925E92BF11DCFE9D67D
AC0F5FFEF5F00FDE9B1538675DA4FBD047D651E2231EAAAFBAD2013CE1C88F1FD9DEB54DC44BAFEA96FD2AF951EF219B4F8D
63E673E2315D64DEA678F1B75B8336F5E2C817E8835EDDC7FAEE47F71D9FDB52E879871E2BF7C3F326BB76FD593CDF536F7E
89F7EADCE87FDF9B55E802DF12DFF12EEB3AF1176FEF9AD0F1D872DDE82FEB22E9011EE2A637F1D3371EE14FD66FCA070FF1
807FDD7FEF7B87AE5AB56A07CF13DEFF9F3FA478C549E7F476E396CDED62FD5A5F76FD59D99FD5BFC3FE1D0A7CE890BEF15E
DD867EBCE9FE78EE9FF59E7894759FF88CC770C207FE81177D7FDEA24BD9477B77ADE8945FE3EB9573DEE9157864FD26BCD4
B7D9275D6AC6EF71F0809FFDE3F4B537C41C7EFAC1ABDB879FE8E7993F09573C523F3C6ADEEABD43637D7583071C6E7CA24D
BDB2CFF00FF18B9BFFA98B7AD02B5FB09E75F463719FB56A5D8DEFE837D5531DF1076FD421D737E9409DC57DF5B64A1DAAB7
6A5EF42FFCE32377BFB5A55AAC271E71E89F0FECDEA566D99FF1EBF60D1D8E2DFB1E3DC157FDE0BCB31F74FBE5E676711FAF
E9CAFBB9EEA99EEADFF36FAD7F1F3CE33BF895759DE2124FCBA59D07C43A78C327F05E5C034E7CA75779FF2AFF0B3A6DA855
F1A91EF7441C7C5A3F78EE8B1ED7C6FA5EF3F54E1BDBFF327C93DEE863FA5FDE5E1BFAE61F5357AF6E1F47F586379C9FF961
8F6B83C7FCCCFDEB6DEF7174FCFED7FDF8341CF1BDCDE15BAA459DF15ADEF2CD7C483CC02F38C85F9FE13B787FED3DEF167A
55177EA03E7C90CFD7A9B265E0EAD5ABAB3C5DB5C7D1112FBF919FF3F8EF3F9FDAD22ACE8717BD4E38745D8D383658DF7D63
FC1E5FFE7CF4BEDF771E50E605FEF217713D5FBFB23EBCF0156E7442CFB7D7D85C3DF0C12B3CC52F71F3C7C73FDB72C59A35
6B76F021E9152F9E38A0FBC6784DF7F080337CF94A9BFE9D1646DDF198EEF9CB75335AFFBEEC5B2F2DEAFE4CD197122E7C91
FFAAEB65D76C5C173CA05B3835DC6DCB15810F5EC2159EEEC767DCF7D6FDDAFFB2EC7B7C6276DD9527469DE900FF77E8A8A2
B31DF357E5351DE03F1EE0F5A85B37D42AFA57C29D3F3EFCF8E675813F5F109778E8C8FDE9495C78876FE2C05FBCC513FCA0
A7478FE9FE4CE441DFEA0C7FB8EB27FC880FF66AFDEE4971CCBA487AA023FC97071DECD05DE53C3CEFDAAEED9AF8393FA257
F7E31F1DFA6CAE1EF9F0097189A7F1A86E5F46BFC4273E4D3F78D4A8E1E67565BDE33B9EF77EABE5C8C08F8F9B1BB28E927E
E48DFFF44707775FD56961D9FFE95EDDF0054FB24F269DCB37FB43C20BFFF1893F3E5065F3D911177DA9C3CEFE232EF1F011
7ABE69F3BB47BFFDF6DB551EBAA0DB97C16B3AC397CCE7C46378C011BFF0AAFD431DDF087CF802BFA113F51CD27BE3BAA2AE
097FFEAD0EA7CE7D7B6DFCFCE4F7567E1B7DD09CE5FD936AAF3C31F8E0E778C5E7F04B9FC0EFEADFDB3428F231F758F7EFAD
577F5E9E1BE08E3FCD5A74AB11BF67C133BE6B1DD79B5FF0537DF054BDF1B4D9B79BCF8EF5F1024EF0B9E9CEB66BCA3EC7CF
CC07EE77CEA5EBCE8D3CD5575DF9813AA9B3FAF223F7B13EDCE0A50FE021FEE185F8AE39BCE5C8B21FE32B3FBAA7D5BBFF8C
F3F1849EF0455EF2D1FF1FA8D9AD46D413AFF1D4FDAEAAD67C97F27DB30FA85FC2F98EC61DDF28E74DAF7C850F75DDE7DDA3
034775574FF5771EDD3ABFC7C96D1A85FEF81A5EB7BCA36B87F8DC02DEF00FFA149778F4297D60D6DE2BBF8D639D2ABDBA46
1CFA80FEF4DCD65ECD2A9FBF4BAF533FA137F9F5FB66D3A072DF85EF453F5AFF51C4EF7E7C6BCA3DAB3F2FFACE5EBD9AC591
8FBA7FF68D541F75E9F4DB0E53C38FE89F2EEB2DED75677CBE487CFA069EE0077EC1976E5ADFFBEE3FA38EFD6E6E7158E0C3
AFFB0FACD49F6EE42D5F7EAA1FD5A9F26ECDF5EBD75779FAF7BDEE8CFBF333F9C1072EFAA4EBF5837F7679B77AACC34FE9E4
BA9736BD19FAC307FE84170D26F43AA0F83C58CA9F7FC39D1F8C1BBAEEDC38F223793C7F73E5BE70C16FF8DCD7BE6B87F063
FCA21B71671F3FE3DDBD631DFAE417779CF8CE0BF1F3CCDBA403FC853FFF568727CEEA7540FC9CDFC0151FE81E2FE0062FBE
45B7F26E30B912E74B1D7B8E0D3CE0C55FDBFCA1CBF4E273A849577C85AF5D3FAB4DA3B28F663F4938C167878E2A3AA327FD
030F2E1FB1E1B3CA7E22E931E9832EE8959EE8481C7C4B3CEAC14F1FDEFC4EDFC8F7D6373B4C2DFB33DFCA3A4AF1D2D3BC5E
6FDE11F5E01FE2C027FA18F19FF51FC531EB3AE900FFB34FA6FBBBAF3CF0143F1FBDADE7D8B21FA83F5EE1137D671D25FDD0
235CAA77697158D9D7B35E52DE74B343D795EBE942BDD4491FE383FCBAE7A39BDE8C63E3C77BFE25F6A5F8C4EFB33F245EE3
33BDE169A37D2AF5B22EDF6A3FE79D17CA7E416FF2EEF26AFBA6C163FEC827F89678C5D96B76EB43A29EFA0AFD657F48F9D3
4BD645D2317DC0F78A254D2F29E37CF71B5DA697FB11FF900F7F9117BEF2EB07DABE5335EE8B1778876F59BFA99E749CF99F
780FA7EC27895F798E69D8F32F81031DE31DFFE223376FDD7457AC9BFD30E527AF76EB3B2F0C7CF8057F32179A2FFFF6E4AA
7BA34FAA377CD5DDF9CE330FD29B7A0EBE68C367455D531DF507F5CC7C487E8917E6DDBFAC5E714A5C2F2EFD06BFAF3E6BE3
D5119FB9501CFF386CF5A7A18766E7F4287EDF8DDF7CDDFA275CFBE61DF15CDF7DB22E52DDE8030FE8A0D95D95BABBAFFB99
E3062C68F64CF0549F1E74548577993F8937EA9FFD24D5597DF108AEF05477F5E63737FFB97DD3B2FFF351B8C14BFFC1633E
75EF159B2E8978CC49F23977F9DAE5713D7F741FEBCB579F94B738AE3DA9F521E53E907998F48A8FD96712CFE0FCC0B81E6D
C387F09A0EC4293EF38BFBD1BDFBDE7940E5F329F2A0FBCC93547F7CA953A5EF91F15ADFD0FFF813FFEBB2D7A6BB22EEE7D6
F7FD75ACE33CFDCA797CC4F92D77EFBE3A7E3F8E977CEB9903FAFE3AFC821FBB7FF7EBDA9D1B47BE494F742E3F79E987FAD8
8C1756DD1BC77AD3FB149F8B169F3E937D26E10D677A9677DF091BAF8EA375F9DFA495AB3F8DA3B9413D2EFEF7FABDCAFB78
3CA3037D5CBC730F5E714A1CF9BB7E57A7CAC6959B366DCA38D265AB4195F53ABED8694AE0C6D7F88AEBF5A17FCEDA387BE3
C68D559E3EB5CFF6C89BEFC21B6FAE3CAED933E5BED1A04B9F3EC5E7B4135E7CBDEFC92D3B45BCFA131DF04D793C7FF2C6D9
113F1DC3159E7C9C3EAF1FBE71FFB56BD7663CE1A8CFCB27F78F461BEB475ED9AF927EE0833F7C178FD48FDF8C3D72EDF238
F205FE72E7AF37DE12EBDCF749F7D5C5BC937847974FDCDCA74FF0989FC1158E973DD0B84E19CFCCFFA43F3AA07FF8357878
63FDC00DCEFA09BC5F7AA2CF99F13958F5E0DF3BF458D1315DE2091FC297ACF7A43FBA5317F5E0BB70E7E38F5CB7F18002FF
D9DD7F5CFCDD45D2339FCBBA4B71D15FF6B7842B3CB30F249EE227FFBDA167BB73CB3EAB2FE247B57DD60F0A1ED03DBDD2A9
EBF9A575B25E924EE8EEB6C73A4D29FB3F1FA557FA75BF471BF439338E7C052F76E8BAB22E5D647D259CE80C0FF18F7FE02B
3D0DFFD7FABD624EE083E21567D67BE2297E661D25FDC89B7EF0F1E59E1B0F8838E14737F432F8374D26051FF5FB21CD2A7C
6F7C5EEF19E579509FC9FE20FFC47375543F7D5B9DAE9ED5B253B9BFF1657DA5C7C51BF70FBE7ACD2FDB6DAAE0B7B0EDD2FB
224EBEAB7E78CC971F3C66C34B81037F721FEB67DDA5BAD19FFBF11DF7CDBA48BAA78F2E07756C1DD7F16BFE044FF3055CE5
CF57E070CF49DD7F1CE7A93F7FA253FA54E79ECBDB6C0A7EEAA37C25FB58E26F9ECFC6F49E1171D12B5EAB3F3EE261F6B154
57FCC9F35AF2915B3EDE7041D957E521FE760F77DB37F0E5237CABD99EBD7F1B732C3EE53906DF527DF04E3DF0BAE9A20D2F
95FB913A0FDCBC7E501CCDD3E6FCD30E5AF564F97DAFCDBFFA16BED668BCBE6FF8259CF51F789BEB4FACBB6274CC99EE830F
7C1A2F9C6F5E771DDEE2AB3A9B5FC577E62BAB5B84AEB22FA57AAA63FF562D8ADF27E8DF78671DD79B37F1EFF2754D2695FB
9BFAAA2BDFA22BFE755FD5F58BCA3E831FEA627F72DCB7CB2E0B5CE48D5FEAA65EF087BB7EF5C087BD7F1B3E875FF87BF39A
8EADCB7D84AFD7A9D2FFABE27949EA3FFA62F6FF749DF3CD797039FF81B53F8E38F096EEF1F7B955FD3FA87C1EAAB2AE7E28
DF63EF5A7A5FEC2FE5AD2EFABFFA88FBBA23DB6C2AF7B1AC8BA403FA78E6D0FE1F443CFA80FCEA8DE93F237EAF210EFD475E
F2317FF155BED6798F54E746DDF68D79194E7C68673F76BEB8F999F85B2CE8552BE2A1677EA94FE89B75AAACFBF6DD77DFCD
AFF5AD7F6E5DB735DECF7E95EAAFEE4F9FD77F46ACC7E7E0819F78494FD957537EF2923F5F69B97DFDA2721F13EFF37B55E2
69D0A47FEDA82F7CF58D6EDD3B5C19FEC4F7F9C49F3E7BE38F51C7A37BBE7E41F01C2FF8145CFB3CB4BE6F1CCD1DFAF3B435
AB9E8C235CF8FE132BD60D7EE79D77F29CA4AE97EEBAF6ADA84FC79F76DD337C93CFF2217DD27DE6D45E31BA7C3F7D60E25F
57B788235CE8193EEA7DD5BD2DBE57EE63FA01DDDEB06AFD4F8AB9A456FFDAC15FBEA56E740F0F38F017BE73E757EBDE8BBC
F98BFA343862DDE0CD9B3757E9734AEB221FFD8FCE76E8B4E207F40A377899536A1ED7AB56F4497CE137D94F92FEE82EF33C
E915DF5FBAA5FFBEE16FF8A18F659DA6B8E8555DF9FE23F3D79D1FF5CDBC4AFD04BFD4898F8E6EB6F6C771CCFE90F44117EA
C4E7D5CB1C85170BAA2CBBACF27BABA4C7A4433A555F75D51F86566BF46EAC63DE3137EED07BE57A3AA26778C339F336F92E
FEB6FE57CF4EA13B3EC28FB31E13AE74890FFA365E3CDAB9FFBE81235FC1C7ACAFA42B79679F49FCC5DBCCCFC44B7E432778
FAF2F1EBCE0F7E663FE12389B7E60CFCBBB2EEBA9F846E1A1FD5AFF87D25DEE803FAD68D3774B8B2DCF76EFBA8EB9EB1DFF2
73FECD3FE89BAE5DCFBFADC30FF883EBF1929F3E387BED84C029F33CD501DFF507FEDA6E9775EF45FEB9BFA478C589C7FC60
E8976BDF8A63D645D283FAF04DF7E97E68858FEEC7A7DC77E0B266C5E76ECD43780A7FB89B3BB25F25BEE009FCE88FEE96AC
59FC4D594FAF5FF7FA0571C45B3E81BF7C4E1EE2CFFE96F8813FEA546366EB7BCBFD965EE957FD1E5AD2AF5DD10F93CEF0FD
9E7FF7EC14CF4BD481CFF1854B07367A378E971CF97287E2EF34DB75DE1AFCD33FF8105FA6FF5BBF58FBF6860D1B32CEE61D
788B9B0F89BF67D5765F06BFF44B7ACDFC4FFCA20375A5B3A627AD9D10F864DF4EF9C9ABE9077DBF2D3E079EF4C3A7D50DEF
F02DF33CF11BEFDA1DDC6368F1FD04C977F853F69384379CF53F7CBA66E5DA1E81137EABFF80DAEB7E1247BAE223F7BDBA76
D7729E78205FF8C39D0FE57A251F55377CC4CF1D3E3F78AFC84BDFD0E79E7B6BF0AE95CF91A5D7A95FD987D8379DBE66D591
317778DF6BFB0CF5D5C7D4D93EE9E4B3562C8D79DB7ACF1C3E78D7D00BDF1597F3ED975C675E77BFB39F797B59E0596FC4A0
95517771EB0F0FFEACEFB7E1DFF8421FFD1AB47C33E23237E0AF3AAA9FBE87C7D57ED7EC92725F158F38CCDD998789D7F8A8
1FEC988357AFDDB2654B7EAD0FD5AEBA7A41BC8FA778843F999FC96FF0F4E98B06AD8CFCF9195CEDC3FE3C7ED9AAA81B5CF9
8EB8EA6E5FBDE0BDF7DEDBE1E7C94F3BFD70EDDB71CCBC4D7CA59B5BEEECBCB5DCA7725F4ABEEFE7DE6FD1BBEFF9810F7DF2
CB060D07B58CBAAB9F3EC32FC5251E7330DC2FFCF9EA3531AFD0A7FAAB3B7CF595270F5DDD3CF2878FFD1B9CE4CF875AFE79
EDAE713EFE983BF0E8CED53D86C67E45DE7C37F7DFEDEDBE2CF761F18B3BCFA78967FC19DFB29FA7FB58FF89C707B50C7EF1
27F5CFFE997C44FCE2E0D3E2C9F34BE2C50BAB56370F3EECF0878A7FEC98AF939F241DD32FBF53B7DE8755F494759AD6A3D7
97AE19746CC1C7546F7D803FF2BBBB56ADF95DE4A37EFCFAD153561F13F1663F493AC0FFE1F7BDF4DFC0C9DC3A6C79A5DF75
ABDA6970F8A63EC7B7F439FABF69D735FDC277EDE7FE74D992D3235E7CE930AA7BBF785F9FE0C7744DCFF48EFF7C820ECC5F
E68EA907AD3A328EE645FCB8ACF99A73C3A776F849655DFAAD794DDFF38397F8C2FFF887FAA84BD65DAA0FFD8943FF1CFFDC
DBCBE268EE10EFAC092B9616FBD6A42FBA92379E557FA2E59BE5B9009EF6CB70CD3E96789DE7B3DE838E0D7CE91BDFE10317
F31D1EC3059FE901DF5F9AB9FA98783FF78B841FDCF0010FF4E1DE77B7AD1F799B63F801DEE22BBF1CFCDF26F5C227CCA3F8
D8FA8A3E4F06FE7C90AFE3A73E8FA78D0F1D383CFC894EF41975D237466E5BBDA6F8FB68FC4C7D034FB3EE52FDE98FAEF87B
AD6BDFEE5BF64D7A536775E55BEA6BEEC09757CF5EB62A8E70303FC223FB58AAB3FAAA93FAA9973EC0A7DB1EBAE677715EF6
EDE477DECFF54A738ABA99DBF0E6AA8FD7EC55DE67B84FD72FD7F48BE351072DFE26EAF4FBBD167D3FF482BFB75FD0BD5FEC
6BF53B7DE0A1550387C7BE1FBFF0978FBB6EE77E79E3F64E83CB7D35F324F907BEA8377D0FB97FCDB971743DBFB78E7CE569
DEC9BC4ABCC02F3E010F38C0979E9AF47ABBEFD6AD5BB36FF10FBE8157F88447FC839FB89EEFD3DB6D5FAD3E3CEE7BCFA03E
4FC6EF15E1C7E79A6E18F8ABE2F362490FFCD5FDE9832E064C68714CD4DFFC878FD96752DDD44BDD8F7C7641A127F5AF5365
D8DD8197FEA32F9AC3E960F19025A7579E530DBBB6F279D5CA79FA61F6ED945FCE2BE18BF770A64FFE52B3EEAA2F0357F5BE
A64DDBFAE53EDFF9AFDD1A16DF0B95FA135F79E6D061D7461E7C571EF44ABFEA94E7E0EF37A9579E53EA8D19F6A7C05F3EFA
09FFDD31EFAEAC1FF1B7FDB0F7CB811BBFE34FEAAEDE7CD4F5FA4AED3FAF7C24D6514773997AF678B4E335C1237D9DEEE1C2
B7E0A3DFD3C37547AEBE3878872F74465F4F9F37EC4F9127DF821F9F9277DD57573EF2FEFBEF577970FF81BF2ABEB720D599
1FD02D9DF5DFA3A2CBCCE7A44BBC6ED064E8E6D8FFC0977FC3852F3F79DFCA938BBE9678821FF4870FFA035EC041FEFA9E7A
5FDCF2A5FF96EB9EFD2DE527AF1D3CAFE803DFF9295FEEB8CBEAC32BBC483A4A7CC443E7F19B16B5577D19756FD16EC0ADB1
FFA273BE4827EAF042CB952707FEF4478FF89D7992FA2DBE3C516BE8E688971FA873E66DF263FCB52FB52F3E63E8CADA3197
78DF6BFB453E0A0F38641D25FCE8C9BEF96F07AC6813EBBB0F5CF8027CEC8BDCEFDCBFAF9A14FB8CACF7C45FBCC51BBEFC68
9FB756453CF88EA7F8291EFB7371BD74CBD07A71C44BFD24EB37E9960F88537CF6433B7CA3721E1DF5EDDBBAF8DE5AF316DD
2FAFB3A0F87E65BEBB74EBA2EFC791CEAF3CABC531E53925F7EB57BA358CF3F4DDACAF54073AE3B3FCE8EE31AB0E0B5DDDB5
67EF97E373E8D6D397B21E930EE5ED7EFAB3FBF27DEB59C77EFD2F2D96EF1E71C1DB7E005E174D79FBFAE8EBD97F125EF89F
FD2AF11A9FF90E3EF65CBEFAE238667D251FA2333CA3B3976ABCB56ADBB66DFF37AF0DAD173AE00774A37FF3919B8F5EF555
F8BBF8C56D7E97AF7DB8BCE175FD4B1DAF29CF35F860BEC38B9A770EB835E647B8F2DDCCF3C46FF5C10F3E872799CF891778
3DECFB8DFA848F98E3ED33F1043FF82E5DF1DF5A5BDEAA11AF1B1F3574D7E2F3294927FC5E9EFA957C3B3CDBEB85384FBFE4
DFB9BFA4BCE59BFD30F140BDE858FDD5BD6B9FAE97455F301FF04575543FFDD9DCAC7E571CF2F68D95DF0326FF4CBCC007BA
A453BE0C4FFB83BC3F3AB37FF3C09DBFF35DFBEF630F5AB22A788A2FE64773DFE44B57D62EBE9721C5293EF3AFFEC237DB9C
B7EAB038F2633FF7BE75CD0563A7AC9A1447F391FBCF58BFA24DF1BD9C4B86EE1AFEAA3E591FF89078808FEAABCFA873E649
E223BEC0C1FE191EF2A1D726BF7CABC6071F7C90F974F5D5ADBF2ECF4175AA8CFC2870D61FF4E3E7D68FDC10F7CBAF535FCA
3C497304BEF49AD67E5BE46DDEA24BFE02D72EAFAFFA2A8ECF1C307243C4CB9FDC1FBEFACF88BFBF7D7D1C07BDDEECE1E075
DE3F26BE645F4AFC966FE64FF2333CAA377DE4D8E273D0291F7E9F7D2CD54DBD9A7E30E4BE389F3EF9229FD5A7EA5459FE9F
B8BF3912CFAEBE65D5FCD08FF3F483DA572FDF14D7EB2B7475FBFFAEFC7D79BFC447EE3F61C5D0E2F789A97E7C4B1DCD35F8
38B7E6F2DD2BFBDDA4EF5407F8F30579D4EDBB7C53E0F6F4A923C7469E74AC3EB9DFD5E9F542E0AEEF79ADBF795FBDCCE9EA
061FB89837EF99D8BF797CCECE7DF837BCF8F59343970F89D7F8668EC33B3E7EC97F1BF529CFAD0DBA8CAC15FD52BDF97AEE
17E9FEEE8B7FFC72D0DB15FF90E74D35BA5E167FAFAE3FA9A3FAE9033B7454D1193D797EF087D717FF24F0A6E3ACC7A40FBA
C83C4C7E808F78821F7C457CFA8938E944FD5FB874F990EDDBB7EFD0458A973E1EFCD990FB42DF70E267FC952FF2C3276E1E
592BF2A0273CCA7A4FFCC2ABAC8BA4073A6ABBB65F9F8883EFF3C5B7B6CE2BFE1ED7DCA7AFECD063E57A7CC61B3EFBD84F97
5F19AFF934DF964FA7AF7ADE537E7EC6E7E4C79FF892F9068FAF7F78E5C0E893FA063F6E3E6FC5D08827F785842B3C5F7A62
E4EEC5F3CBC44B7E9FF99F780FA7ECF3C9EFAC8F579E53E057E67FAA331DF45FD4EA86D08D799A5EF567BEB270CD92E2C8C7
E5D7E1072B7F5F9E37E18E572D7A0F1E15F7E59BFC182FAF3DADFDB6F2BC436F7CAB6FE755F34347D99712BFF08ACE2FECD8
F0ABB2DED5453DCC2FEA46DF74CD4FAF38BAD9C3E579905FCA5BBEFC95DFE2937AD0D98B1F2DBFF2C30F3FACD26397CEC5E7
42CD557C2BF33CF52B7CE7FB7CE49E97DF5C5F3C674F7CC663BCC8BC4D7304FE3EDA60E4EEC103BCA34B753747ABBF7AE833
EA92F9997884A774C5771E7A64D9D2789D7D26F1000E990FC9CFF2F38EF346BCBC6EDDBACC7B3E97FD24E5217EF8C085FFA9
271EE10FFCCD2FEA90FBCE81FDFAC4DF5DE83F5EEB6FDED7A7F1E9964BDEDA2D5EE30D3EE24FCD6B068F8ADF47B90F1FAD53
656CEF625F9E7C91EFDBB79FD9F4ADC7036FFBE3E7B68E6D13AF9DC7A79D9F9F3BA4EBE8077F7B345A39308EF68DD63DFF27
2B2F0CDCE579EB2E3DEF89FDA73E47F77C00EF735F48F9C90B1FF0004FDDD7FDECFF3C573D62ECAB1DA3CE7CEB99BDC6B629
9EB725DDC3C9F38FD3767BF39BC857DEE2D6F7C4CF8FF4853A55DE383D78586FE9D8E2EF47E0C8179DC7676B3FF1C6D171BD
BAABB73E2F1ECF41C495F99C7483D77C935FB77EE4CDF565FFE61F8D5F5EB6F4A38F3EAA32FC839746C7F360FBCC91B735FC
2A8EF65D70BCE4C19575239F3E5FB4BBB7D8DFA7F99BBE1F1A33E2E5E0BDBAF103FD93DFBF71F4E29FC491CF89579CE2A367
71D239BCEB3678E3E8789F8F5C754CAB1BCA736BEB2B066D0C3EE807FC4F5EF2B16FEB707FDF1E819F3980CF3DFDFBB135E3
7D3E812F70B3BF849FBAF2CD273F5EF245D49B4FD00DBD64FF4CF8C1CDF39813EB2D5B1F3C55FF3C777DD5F9DBF2730BFEC8
A7EE5FBBECB838EED075C50FE8BBC184B18746FFC24BBE9B7D23E9832E9AED39FCC3E2F739C937F92B7D1DBED7BC63632EA0
33FB13F8561BB1F2F0E27E7491E2A20FBC3277E3171C3CAF81077EE0C50B7B2CF9E2934F3ED9A1AFA42B7AEC7A778FE3A36E
E62D3EBD43BF95F3F02AFB43AABFBAF33BFCED3CECADDDE2A81EFA9BBAE8AFF47CC7C7CB3F8A63E667EA57780A3778D9C7E0
B73E80E799B7A9CFE3EF90854DBA967F2F8AA74F9C35B6F89E1CBCC3EBCC9FC41BFE8AC77CF6B15A4BBA95E744FC8677E655
AA337E65FD267CD52BF32AF9197EF109BEA14FBCD471CCE4E88BF8CAFFE00C5F7D1B0EF68BF0C8BA4E78C0812ED5597DCDE7
E6BBD19FAEBC308EE67A3F9FF8E05B8F177F0794E675F5AB31FFAD13236E3C344FE1A3E74CC7EFF9C6411117BEE7FE3E68D0
C6785EA0CFF373FED8FCC065C77DFCF1C7F93CF385F31FF870F887C55C93D6E3D33BF707EBE4BAA77AABB3FBD1DF8B4D9674
FBF4D34FB37FB98FF5F3F39534370D6B5AF1DD3CFF34EFDB233EEF683E801BBCCCD7E65F784FFBECCD6FE298F7E13F7B6974
794EA7C31A3F6C776F79DECCBE91F8277E71E8F3E2C97C4BBCC6BB5EE7776A13F73767F32D75F4BC4C3D1FBD6DCCE49843E1
447F74456F0F7DB6E488CAEF27263E13F7E3437C99EFA85BFB3D967F14BA7AEE8B890F54FE4EB5721E7F6DB160F849B10EBF
E637FA049DDC7BCDD29AC5E731F8438A579CEAA78F0DBC68E5E1716CFBE1C02FE33EFA2B1F1AB8B245BFE2FBDFD27E8FFE9E
F9E1C40742E7F4273FBEC087EB54597C7BF9357FA9DD7A71B5C8D77C8C1F378C7A734E1C1B3F3EE66FD15FF9025FA9B77DE2
3F8AE7250917FE94E7BA363D8E8FE74AFAABE700E69715E7BEDA318EFC401DE08F677C1ADF727D535DE929F327ED1BF0C8FC
8BEF739E5EB63E8EF80AA7BAF72CAE16BC5547F5E3DFF0818BB9121FF82E5EC097CE9FFAD3E283E2A88EF673EAB9836F159E
E21D1CF5437876DAD2A77ED4D95CC5FF9EAE3AF11FC5E7E513EFF022F33FF1054FF4BFCB8E6DD2B5BC1FCA3C4C7D151FB35F
259D37DE6DC9119F7FFEF90E9E247DE34BDECFA53ED0E7F8B74E8C235FB49E75323F132FF136F32DE184770DD64FF822EA84
7F7C6507FF2BD7ABAFE776478F5E3C3EF2E12BFA8CBEA31F39F25D7E9B7597F08633DF3C7FE60B5BCAFE092FBCABBF78F141
912F9EE087BC731F4938B5EABDB4661C1F6A38E66F11BFBEB0739FA879DCF0938AFF7796FA52F6C9749DF39DA7BF39DF5C85
DFB7DEB56C551CF507FD025E9987696EC2C7ACAF84BB7CF9325D3D50F7F58671E407FCC1BC423F74FED8FC45F3E3A80E7C94
7F663E24BEE0C513074CF82274090FFAE037F08673EEFB7B0E2C9E53E7FE9AF0818BF928CFABA33BB529F7CBACC7747FF7C5
43CFFBF0914EF862F70BDE9C13C766E78CEE1CF7E1BFFC467CFAAD38E10C5F7ED97FEF36AF85FFD967F2A1EE53BA6D2CFFBE
888F66FC13EE786A7EE4AFF33F7CE3A0CAF71B4D68567C3E39E9901FA80F3FE00375AA4CBD3FEA454F7C8A1FD1135F6AFDAF
6177C6CFF55DFE54F747536F8CFBBA9E4FF07BFE8FC7FC862E5F3C7ED1FC2FBEF822BFC657EFE3F18E3E3AFFE138EA4FD50E
6BD1AFBC9F701E3DD71E31FFD6E2F89FA937066FF042DEB9BEC947D5199FF443BCCA3827BEC33BCF5DBFEC533F3E4FAF1F67
BF4D7937ABFD7AC3FFFCE73F795EC1E33BD7BC51378EF5F79C7A6CF1B9F284233DE18FFD251EC10B4E752F9A7F6B1C33AFD2
9C9EF789B95E95FAAB1B9DE3EBC30F2EBAAACC13F8C12DF321F561BC803BBE3E75F7FC634247F40D0F38E023FDD11D1DF289
4E6D97AD2AF72575C587EC33C74C68569E53BDE6C3DEEFF0EC804F230F731E1F5517F5308FECE061E57E70566FF382BAAB23
BCD5F3A90FA71E1B7889171F737F1937BA73CCA9FA8AE7569E075EB072C5F3E13FDEF7DA732075543FFEE139E6D94FADF865
AC633D75159778C4A15F8927D737D5551DACEFB9AAFB786E22BECB0E5A716D1C1BFE624AB1EFC66F7A72BEF33C3F1839AEE1
F0E2F95E7A5EE5F91E9FA083FA6DE61F13BCCABA4875F27EE35113F6289EBB277FA47B3EC2FFF889E7A67FDFF5CDDD222FF9
651F4BF76FD474D155FFFBBFFFBB8357698EC32F7D179FEEDBEFB579C5EFADCD03FF1E7667F13D4DA9FFD5A932B366F17C22
D5870EE1865FF0D307F405FE6BDF0EC7AB4E7FB36B31A7ED3FF386CA73DECABA784AB778F3F8C1AFAE8EFB78DE7C6483F923
A38FEA9F3BEB954EE14E2F75AABCF25059AFD997E7BF726B5C0767CF85E19DEB957C57DDF8873AC09F5FF09B7607BD5137E6
A1FCFC29E56F5F0E77FD1CFEEAA54EFA79CBDD476D8CD7E600FE577BC3CC1B224E3ACABA4A78C2916E9CF7F8FB5336055FF8
C19AAE331705EFECDB56FF6676DBCAFF7F54AF0A2FD4ADCF931D8BCFA1DBB7EB0F993F292EF10C1BD86842F4A9FCBC37E989
7FDDF0F76E1B630E3457898F8F8A33D731D54F7DE94CFDEB1EFFCAAD65BF749EF7EBFF75E69F8AD7897F78DDF58ADE97C57C
6B7FA58FE5FDE97B6D5E2BEF53CDADE6A3D7CF5F3C3E8E79FD849FFBE05FAE4787578E8BFB663F4979345CFDEAEA2FBFFC32
C70D4FF8DAB7E3D3B5C72EBF398EF2DBC1FBCAF9B91F5D30618F983BF9FE8EF32BE78937C799E2DB196FAFC595F999E61A3C
5577CF11D53FF3F695997F8AEBF0D21C80F70FEEB77072F1774155A7D48EB8F8355D7AAE7FF2CBCBAE8F3E40C774095778C2
0D5EF6F7993FF871FB2BC715CFE5125EFA12DCE42F7EF9F0217EC0079C87D7F2759FCCD7847F9EF3EB0CF834FECEC19C071F
B8E8B7F0E1D370A22B73257DF1297ED672E36BF38AB988AED3F315FACE71A6F8F01BAEB9FFA5F3323E893FF2E7E77ECFC1D7
DB1D3CF4A9E27BEFD27CAF8F653DA67D275D36BC78C6BA889F6EE58D0F7ECF8117FC8CBFD97FE675523ED6D387E0FBF0DA57
6B06DEF6C97873E3F6A53D0BBF4D7A7EFCFE39EB8B7E96F2CEEBA4EBF9B1BA9AE3D4173F1ED93EA576E44FB779FD74BDFBE0
E57D9F8CDA58FCBD4AEAEB9E63D8378F3AAC32AF356B31BE665C6FFEA05779C9C77E37F376F88C75F15A7CD6B5CF1BB266C5
B571E48FE61A3E090FEB5817FFE8B5E9C68593FFFBDFFF5669D87CCEFA329FD4879EF0176F07CF6FB62EF4E23929BFC9FC4C
EB5817BEF40DE7CCC3B4CFC047CF3FE03ABEFE8A5F16DFE7E4B9F9390D87979F03A9873EAF2E19FF345FA883FCC5DBE8C057
6B060EF0B55FE93FF5CDAE71CCF8D061E271C639F123D72D1DE5AD8EF681EA7DDB634B9E8CA33EA6AF99F3E0AD6EEA651DD7
DBC7A88BE7A0EA233FF9CAA7C9AE933F8B7CF83B9C5EBE71C663517779D2BBDFF39CB0EF1B3F0EDFA77B7899137B2D5C7E73
1CADCF2FDD471EEA20BE7CDF743FD7F5FA5ED7678BEF8F4ACF4DF8A5E70DF832E5CB37778B63DEAF0DEA7D59FCFDAAB927F7
979AE36BC67AF44CE78F9C33E791C0513D3B7DD57F59D4CDBE912F9AEFE17ED7D4C5EFC5313F4F78A1E3BEE5E70AD6E74B3B
DFE7916E331E2BE6C4C41B3CBE73F5D0A7E2F3E2FA53EBD9A3F68CFB9A47E9463CE230CFBE3C6ECE23E57AAB67E667BA2F9E
E6E7D4D51A4D283F1F5257FDB65BD5A53DE3E8BA5A5F4EFE2CF890F14DF38438EE7B6DE10F222EFEFE609F792B8AD7090FF5
A06375C137BF97C2BBCCE794973CF37DD3FD76EEDB5EBB3F3CCD5370C51BFD1C7F329E09477CC54FCF63F034D737D51BDE99
57691F895FEAAD5FA93BDCF971C717973C19C781FBB63E2C78E23923DDBB8FBAEC5C2778AB17BEEA2778AB8FE93FFA4ED66F
AA43D31AF3567CFDF5D7F9E77C947FEAFFE6F7E54FCC1F19479FC338B4E52B7B056EE614FB5AFA9EF5C8B2EBE3287FBED8F6
F4C5EFC5BCD9F28E71C5EF7FF557F56C72F7F4AD9F7DF659AE133F116FF69D5457EBCAAFE51F17FEE0ABAFBECA79F039FEA6
2FE953F0CAF1A4FE23AE1C4F8A833FE63AA6E77FEA99FBDA49A3F68CCFADF005BAA975C4EC4DC5F745257E898B5EE9D47C70
F9F1CDD6959F7F34FB64F2E5C5E79692EFE778D3FBE2F373F7A9D566FAD6626E4A3C6BB262F6A6328EF0EE306AC8C3C5DF39
A6F912DF333FD33A78292FFA945F9EDF76E9BF2CFEAE92BFAB8F7AA9ABBCE9865EF407FD42FDFC9EEA4F135EFB79BCE6377C
D0BC91FDB06FAFFF29FE7E31ED57E95E9DE858BDC4051F78890BAFF139EF23BEE9FA6CF1FF4A535FC97C4F3C957FF6CFF6E3
86847EF026D725D503FFD4C5FCA53EB9BE09F79DEBF3C0EE932F2FBE173EDD8F9FD67C6ADE6F8BEF914A75A36F7A37DFD0FB
C1EB66B70D7FA77B75E5BBEA6B7D7EEE3E19FFB45F5107F18A53DCFD9AB4AF1FFEE1B9099DA9BBDF0FA93F1CD5139EE632FE
34EF9D377E5C7C3F497AEE76F9FC65EF949F3F79DF734E3FF71CC179D5672E5D1E47FDF2EE160B8E2ABE3723E5ED7CE7D94F
F3453E89DF9EC35E346079CDF273567CE2F778E5BEE607F76FF7F0C803037F7306BED88F89E3BA97DEB8208ECDEF9E7E663C
A7C333FEE17CE7D9F7A813BCF14AFCF93978CA837FD1251F334F5BF7E6FD5FBF28EE93FD21E9AF45FD79BFFDE69B6FF2DC76
FB638B3617DF5392F2B68EEB779EAFBD769DF5AC4F6F795FF64EEBC3CAFBB3DCDF2F18F2703C7FE46B999FE979009E66FCD3
9CA20E9EEB9F7BD7F25BCB3865FC13EE7C9C2EEFAF3AFBC4382FEB3B1DEF6F33FDCCD877D18FF53DDF779FE6DB679F1878CA
539D47DFF642A7F87D88E7F7A34F79614B3CB7B79EF5D5BBCD1FC62E8BFEA44FE7396CE9A4A6F15C46DCFA51D7BBFB758975
CCDF749FFD30D5B94DCD0547451DF2FAA90FBB4F5E3FADEBE7F0E12F70CA73E0D5BDFE27F2E2CBEAC7B7D47144CD9766C673
85FCFB8B8487DF27FC63F76553CAB8E20F3CE19B714EF8AAAB7E21EF0E2F2EDA1C4779C3C1BAB9BE291F75C878A63914AEB9
7F351A7960EC9FF1555E9E2FCA4F5E9EE3CB4FFF156F970DAF5F14477AA063BA30B7B8AEE7A36F5C10C7EC0FBF9FD434FC99
4FF439AACBF78BFF3795F68174C38FF9337F34BFB94FBF53962E8FA3F5EE7E63ECB2F83B93ECF7C977EFA93EF7AAD01B5EF8
7DC8EF9F5DF044C4A56F987BDD67F0F1CBDE8923DEF05DFC510779C9D37333F3F2AAFB5ED9ABF83EC9E4D7D56BB5AF5F9ED3
D51F1FF06068D7A6BBC6DF1FE4E7C0A96EEA8ABFEA6B6E95C7882B97D72CF6910907FD001E19FFB4EF5107CF83FFD663E933
655E74FEEBE02BE27AF31FDE990BF4D5C567BDF6F3E27BEFAF9FF64CCC3778CF0F0EBA69E6A2A8CF6F364F5F1B753347C8D7
F30F79B7EE37F7AA6FBFFD36EB81AEF0113FE1A81EF8801F793E69D3AF4BF1FF45922FE4BCD23C253FD7C31BAFD4473CE293
AF3CF994F9547DC6B6ADF8309C3D6F85379EDEBBEBAC9A85EFD27D3ADEDB63DA33711DFEE97B77DC3EFFF8A8AB3AB7FA7256
CD8833FB51AA03FE7B8E4C07F44ABFF0C55F7C3617BAAF3EEAFEE212A77AD02BFDBA1F3DF015BAE8F8D3117B17DF8395E615
F569D37FE2F0C015BEFCBFE7E61E9D8BEFEB4B733C5EC313EFE13AF87B2D6F8FCF0BDB6FE35DBE6FEADBEECF0FF8837CF089
4FE095F9E6A2162FCD2CEF3F72FC296E7DCE5C72CBC90B67179F5349B8E2973AAA6BAE6FAA2BFE991FF16ED227CBA614DFB7
B97ECCFAE273BCA94FE7FCD2FBE2F0F3CCB784977A8AD37C25DE3C972EE9F2FDF2F3337EC01FF83BFFBAE595C157C4FBF4B9
733F6CDF61FEF171E467FC4D3D3CBFFCEB494BDA867EE829F321ED33F022D73DEDD7D59F2F1D70EDB41FC53E813FB9BF78D4
C39C7AFDF2C5DF14DF2791F0C8FE73D5C4E1C13F3EE27CF3B3EBE0EA7CD7F309BE217F7DF7D29B9AEE5A9EA3F9D46D1F8DD8
BB8C2F7C3C7783139FB8F38C576E8FFAD287FD408D535F5F54FC1D708A37FBE10163D6073EEEA77E74A08ED6B1DFB11E5F51
E7CEB316CE8E23DEF003FC697BD6D4F6F173FCA6FBEE53FA5E5E7CDF4B9A23D57940A3B6AB8B7921ED0FD5D91C411F33AE5F
FA4C1CF33EEFA1378E2AEFE3E0C3D7E0E47CFB40D7B59BFCCAEDC18BDC0F923E739C695E13AFE764873F3AF7B4A8373FE6AF
FAEAD23A0B9E88239CE08697F0810B7FE147D7EDD3A373FCFF78FAA4E7BB264C6D1FBE8F7F7C826FE01F5EDCB5DBCCBA85AF
A7EBF1151FF083DF9DDFF1854E65DF5317FB24F569FBD9CCBA653F53673CC62F7CD637F0A8C7918BBF89233CE0AF1EF6CD97
3CBEF4C7E5FDBBB8DD079FE1011F79E6B9F79B96B77FE7395C5ADFFED67DBAF619B46FF1BD58694E527FBE7EEB39AFEE17E7
E197E755C7AE7A6DCFC0995EDB3F34615CB9AE7C9AFEE1CD07D40B5FD54DBF825F9FE9AF2F8AA338CD41E2CD7549F5C033FB
D4F3AF5A7A43396F7550177CC53FFAC643F55177F7CBFC4C71AB43C621E5AF6FADEF35ED47E5CFDFACDB67FADA38F63BADD3
A3A16FFB0ABC83BFBEA50E744677FA169D7ACE45AFB9EFFFBDEFE5F1FD1AF4061FFB6D38F11BFE436F9DDB0DFBA8F8FBBE34
DFA843C717472F89E7CBF2D4AF9D6F4E715D3E3F9DE7E7FAA4FA0F6CFCC65171D487F9E22BB397B48D231E7A8E868F19CFB4
6F802B7DDCD178C2B8C089AED4415DCC91FCFDCA87DBAE2ECF7DF842C79DC6BDBA5F1CE90DDFE9CECF9D8F4F788DE7788947
E214B7BE7DE3CF171475C30BFBFBB37EB9F437E57A8EFA49C38951EFFC9C243D37C97DB0C6A07DE339027EE8DBEA30ACE1D2
1FC791EE6F3B75CE55D18FE9159FE806AFC4696E10AF78ECEFC5D561FF29EDCBBAE51FF2B2CF975FF6C3C7462F8939882FF6
F979F7C38207E669F5CFFE7CE7B08F62FF297FFCC237F5185EBD71F13DB8F6B1E2ED387DCE77FC9A0FD0FFED1BA6B40FBDA9
9FBCE4C317E049CF7015BF79501EF88EFFF469DEB9A6EFA2A6B10EBCD5E7F6FD67D48E75C4D761C38CDA65FF82B77CED6BE5
8DF7788DFF78097FF5D087F16874FFA53714CF93CD4BD33A3D5A9E9BF00F8FF0900EE8ACDBB605CFC7D1FEFEB4996F9C53E6
051CCC6DF050077511373DAAA3BA7ABE7254DB051BA29FF1153AA00BBC119738ADCB3FCC476FBD34F7B4D867C10B0FF1B2D7
F77ABF5F7CDF6F9A0BD5DF7970C2DB210B9A378B7D817DA0BAC1C7FE1C4EFA031F5D78E86B7BC6B1D36FC7F72CF3927FE30F
5EE251E64FE20D7E9B43AFFAD16BC7471EEA90F34A73A5FC3C4739A8DFCCFB638ED4AFB28ED2FDD44D7DF0561C1987B44F82
075FBCE9A4B91B237E3A93AF3CF521F19B8FE5A1EFE1F5842D15FFC163FAC067F3CD059F369C589E2BF35CB1ADFB61F17FBB
E901DF3CD7C03B3E4407BDAF5ED4348EF4445F78D9639701AFC46B731BBCF10E0FCD0D5D0E1AD53FF296BFFE9BD7497396F5
E0C9C7E13AE8C1368716FF5728ED5BD421AF9FD63507D0DFAD6F8EEF19F3373D74AB3A7464F17752699E707FF577BEEBF3F9
695DD799972F3BE9F59B8BCFF9A77ADA679F74CFEB0B231FFAD0A72FEED778DF72FDE00D7F3AA71B3AE2CB70B06F8047D7D9
7337C67310BE80C7E2347F8B970F5EFB4DEFF7E3FB6EF0A6D3C6C9FF0A9DD11BBD8ACBFAEE475F9E0BD019FEE12F1ED2C7CD
17CC3EB2786E9274CFE7F0B1FF7F5E3B3E8EF0B43F872B7DD00BBEAAF3CD6B46F52FF6E1C9D7F4313A9B76CA1BE714FF4739
F5B16BBF9D7F68AC8B77FAD065C7356F569EB3F8C5ADFB4DFE57F44BBCE19BD77F35E095D8E7C377C0891D16C77DCDD9EA06
67F8F20BFB860B4F5BB27B79DFD465D4EC23E3B9AABAF225F9DDB87DE8C8D8078AA3F394E935A23ECE534F78DFF2F7E9358A
7E95F2E123E7FDA2DE93653FA10778C2577DD40B8F73BE692E96379CF54B78E3039DE345D663C2591CF0B17F80D3A6072615
F3827D8F391BBFF111CFF9937EB5E4AE051BE2E879C36F9F99B777E0C22FE9164EF44B07EAA24EFC12AF87CC7EBD78BECBC7
AF68DAE6D0F2BCD3EFD0AEAD8ACF1BA639127EB98EA97E78D3E5D5718F9579C35FED8FFFFC8BC523034FBAC177FCA71F7C85
13DEAA375EE53922F1479DF048FDF944AF2AF30F8D237E394F3DFDDCF970841FDECADF3EE39C83971C5DAEBFFE5FBDE9C25D
824FF895F14C732E5CE5236F7AA17375A67738C357BFC12FF3D79AEA33EF0F7EF1697E33BBF5EB0B8BEF394DF5C707FAC13B
3CE4FBF2B50F91377EDA6FE3291CCC13F0E067D79FF14A93E27B0192CEF0C43E165FFA1CD5EBA1E27BF9D25C08BF3C87CCE9
B0B83C8FF05D7C1F396DC9EE71ECD67DE4DD65BCF44FEB9BD7DC27FBCC9FC73D16F33DFEE18BF79D473FEA4647EAA02EE692
D12FBED033786D5F3AE6FD7A4FC6519CE2D3EF9D6F1FE9BA5EE7F7BF39F235F7C8C3F3F1FDFF3065F7A8B77E087F7EA70ED6
3107596FC4B697F72FFFBEDDFD7B3C3AE4A6727CE616EFE7D7E6AA7693CE8DDF7FF00B3E617DFB30F7E931F99526F17BFA9D
7568AEB9FC83454F16CF7913BFEC87FEBEF5F5F1659EE6BEBFAA6BABF89C04DED0015DF0233AE01FFD1E5CB84B1CF108AFE8
077F6FB877E6EAB89E0FF113F1CB6758EDA683624EB46F91B7FBB8AFEBF11BBFF09C7FDD74E7A4E2FF6F675EF2C91B46DE5D
FC3FCEE42F70336FC12FD725D583BFC0D33E0AAEE2B7FF9007FFA6BF71AB971C5DFC9FD906D37E1579F34DF5EFDE6AE6EAF0
5375F5F3DC1F97F47A289E8FA81F5C6F7C62DAAFC26FE46F9FFAC7850B3E8F3AF00FFA84C3CEFA855B3E8F0FBF34E4A6F03D
7E9CFBD1E8FE3717DF5399E2A3637DFBCDE7E6ED1DDF57C02FF9EEABEF2F1E19477DA0C62EF37A449CF437E48456D58BFF57
9EF60F70C52BFAC42F3CC60B7CA6F75F379B7458CCBF746F7D73BAFB98E32EFECDE2C7CA73B8FDE8216D670D087CF93AFEAA
97FAE5FDDE632FF42CCF6BEAAEDE788617784257F4447F837FB6E8C93876BF6EEC51655EF26F78EA2B70B5BF39A5D76B83CA
FCC517FC517FFCB1FFC4A3BC2FFAF9CBFB97E72C3C540FEBE119FEE0131EAB1B3EE30BFEE8BB70A50738AA97F952DD06EDD3
FE57C5FF174A73BC3A67DC125EFA0A5CE1CC77F0146FF9635E3FCDD1EEC32FE97ECA5E157FA01379D30B5CFC9CAEF4CB4BEA
361D549E77FA7C35AF47F0860FE29FFADA77A9B33C9CEFFA01FB75FE28E234BFCA43FD6FE839F6A8627F967C00FFED1BE9C0
CF9DAFDE3DAB8ED8A38CAF7E48DFD7DC34E780E00B1DE778D23C282EF3E379BF5DBCBD5C67FDAD5AF3F935226F7CCFF3C6BC
56D5CB73478E27C5A11FEFACB3E1EB163F16C71E274FF83C3EF7469F74A69EEA4BB7F635E6CD0D6F4CD9BDF27B9C1E6B8BEF
014DF3D98EFCD2FB69CECACFA9A70DFAA41C9F39010EE64B78D0219FA1C7BC4EBADE1C03377D0E7E7D9EECBBBE382FCD49E2
9287BCCC6DF9FCB4AEEB7A779D7340FC3D155ED2ADBAF323F5A76FFDE1B563177C1EC7DCDF37B7FF55B9CFDBE7FDAEDFBCD6
C11BFEA46EEA48CFF6557FB97E7195B83F7DF0A7EB674DF83CFA02DED31B5CE90E4FAFED3AE389F07D3E8727F8C15FE8E0BA
ED23F688FD189DF7FCDB9449810FFFC22FBEC23F665EFBDAA062AE4F7A85277C7BDD34E38932BED613D77533A64C2AF6513B
F9BA7CE935F7CD8D9D3F8AFF7B07773F17BF7CE80A4EF4B5EF31E3BF88F37FF5C37147C47DF45B75E04303EF9F5F238EEAAC
EE7CD6CF9D2F7F75B4FF514FF3EC993F5BFC76591FD9574F1BF449CCD1FC157FCD296FDF356B401CF5F94B3E5CD039D6A58F
DCA7A6F5581BCF23E08587F675F88827F88D2F7C82DF8C7973F1F6C20FF58517FAAE8FFD83FBE1097EE0195EE0093FA213F5
A1979ECB47CF8DCFF3E1253FC623FCC1277CC11F7CE1AF571D34F78791171DC3DF1CAC0E6336D7FB57FCFD84FD997DDCA8DB
5FAC1AF5C8FBB5F4737550177380F3ED875C87E7F015AF3A9A03D4537E74025FFBBF03BE983A20D6D1DFF01E7FF049DDD491
EF8CD8BDF1D8383FEFD752BCF087BB3E446FEA4677D6B19FB01EFCF9B53A0C5DD3E2C9D09DF9D875F8EA3E782B1F7CC3BFBC
4E9A8BADA7FED71D397A6EEC13F882EBDD076FED338EDD7FC1B088975EB39FA475E04BF7FC6F6E8FC555E2D87FCDDC1F06EF
E91F5F7BDF3DAC7D19477D521DF8CDB03D17748EA339F7D47716FFB5CCD3C12BDA5E52CC71692E94B7F9E2824317169FE3C0
23BE241EF1E1A93EB6ACFABCD6F1FF23DC5F3CAEA7577E30F183C56F17FD73F6F81743DF744EAFF295A7BE4E7FD55BCD6A1C
3EC46F725E693E951F9EDBFFE0BB7DCC61C366B78EEBF9D9C02F3B9D527C1E2DCD85D6CB737F8717AB96E70E75CFFBCD54FF
BCFFD9A7DEBFE278CEA0BA67C5519DD59DAFA997F955DDFA7C31B056397F73478E33CD8FE2554773927AF2896B4F1AFF62CC
ADF89DD74FEB9A6FD447BDE85F9DD55D3FCBF3C3278DC796E708F55557BED0EFDE598DE3EFF176EEEBFACF6543E615FD84EE
079CD8FD6FE5E7A63BF24EEFA7B932CF994DFA5C5ACEC75C97DFF73A5DD7FBADC9E7C63E813FE2237E5DBD6DDAE7E18FD977
938EAF6933AC7D7CBF64F697C44F3E85A7FAD3A507B578B23C4758EF9AC32717738A3AF5FDF9B4CFCBB8884BDDD491AF6CB9
75F4FBB19E79FCDD2FC61D51F9FFEA959F8B53DC740C7FF5709D75CCC1EAA2AFA80F1DD917D1137FBFFC88B69794FB331EF1
AD51AB16D68C63C6EF87036BC53E90CF9A6BCE3E6EE1A8B82F5EAB973AE19939F7C4A10B3F8BF3E8493FA8B66BA753CAF375
9F53465D1A7CC44BFEAAEEEAA7FEF4CACFA6EF5BD1299CCC59EB7F3875401CE5275FFEA8CEEAAE9EB9BFD4EA7369EC93F519
3AD337166E58302C8EFC88CFF09D2197CD6B127F27CF7FD53BF7C739DDFF16F1C145FDF10DFFE8D97DE85A1ED6773FFB923F
7C326F6AE8802FE2391CE5ED3EF8ADAEF083A739098FF1046F721D53FDF40FFC319FE1111DCA8F1ED5537DF9205E99BFF18B
3F5D71F59C4B821774AF3EEAC5F7F5D58BAABE5AFC1FEC3C77279D888B5E322F668EBAB4D857253DDB67EC3B60C26F024FFD
CDCF9DBF335E788567F4E0FE74016F7D7EE525B35BC7517DCDFDEADCB7EFD09F9571D7AFDC0FDF76E69FFB677D263CF90B5C
FB9E3C6ED7F85E1B3AA75778EA73701DD477CE25C177BCC76F718A4FFF5367FE347EC1C251C5EF59C6BFF0797C2FBE79DCDC
0E0FF8E88BCECFF37FDA0F8C7CE6E5E27B88F33C9FD63397FDFDB70B6B45FC789ACF4F73B1EBF0ABDAC6199B424FFCC9FC7E
E0E7D3BE0ABCF8F1F07B9BBD18F3A7F9D87AF2932F9FA263B8C1D13AE665EBE1355F9C73E9C2CFE248DF57CF1AB76B9C8FC7
FDF7EEBFBC8CBBFEACBEEA4ACFE6E8E3FFB960EFC085FE863EDBFAD1E2FF57A4F9555C79FDB4AE7EDFEFE689638B793AF90B
1E59C7FC683D78F0B911DB5FED5B3E7A1F5E7476F5D5437F163CA0B781FBCDD814F7CBFC4B7E47BF9E1BD0313EE4B93DF142
9DAFAA36B5F87F787C74F0840E8F14FF3721CDA9F28037FCF9C18046BDFA7C079734FF58AF7A978963232EF5B3BE79D37DD4
59DDF9665E3FAD6B3EEB3F70EAB58183FCD561D03E5D9F8FFCCC733BEE93DE4FD7FBB93CE00C77FAD327DFD87DDED438D20B
1FA11B73C4B9E7CF3FBBF89ED7E41FF8668EC73B3E7DF19657FAC77AFC89EEEC2BE88F6ED41D0FEC338ED83EE7E8B80F1FC5
4FBCA4073CBAEABDFECB63BFCA4FF33EE4EC173E2FF77F79E7E74CA90FE973173EF7F2A8723FEFD760C4A4B2DEF8B7FAE303
1EE94F17B76AF662B97FC2537F852BDEC21D7FF185EF4E797361ADE2EFA4F4C13AAD1F2DF743F8E10F3EE5FEF270AFE2FB16
F5153FC733F8AAA33EA19EF82F3E796F1D317C4EDCDF7C6D3E36579EB2DB823665BDC2138EFA0DBEAAA3BAE6FE7B568747CA
7D581E576CEEFA7C3CBF709EFB9B37DFB972C26FE2083F78F2773F77BEB96AF82F5FE91FBEC5BFF0153FCDF178AAFED59F18
31297C9F0EF14B1F98FFFC82BD63AEE593D6773FD759072EE2E5BB6347CF3F3B8E70CEF37CC2DBCF9DEF7E747969FF594F07
EFF25C9670C713BC81877967ED8FA77D154773EB1FAF9AD73CEECF0FFA2F1AFC70E5FFC057EAAA8FD1153DD257FF5663EE8D
F9989EE98C7E5C4747F083279FB5EE39ABEA5E54D6BB79E11FB7CEBFA1F87E83A43F79C8CBBC4257EE4F5FF2928F3E0A2F3C
C44B7CB5DE8EF9B1729E9FE3ED2F3B8D7E3FFAD1DE9347558DFCF43DBE7BFE65735B15FF473EE972E855B39E0E1FC44BFA57
5FBA55677ABDEADE31F7C6FC8BAFF0873B7DAA6FDEA7A43A9BA3F7AB377151E0CED707EEDBEFF7655CF44FFCBAFCD2E90787
7FF03975D1BFD5471DF4AB1555E71C1D4775E68B333FABF07CC08209550307FE8547D94F8E19FC703C27C56F718A4F3FA63F
78C2174E70A3CBC143A71F1CF7833B3F552FFEAE6EE6B213DE997F48C44DAFFCF4CAE326548DFAAB93FB885F3EF0BB62C9E4
5EC5FF874BD7C3096EFA1C9CF2BE2DF9E0A0A326F72A9E57A6B8E137E8C11EEFC7E790E062BE18B5E0C59E314798BBF3F3FA
F47E9EFFD37C3F62599346F17DEDE669D7E5F5D3BAE618E79BA35D37667BBD45D107AC6BBF30FC472D1B46FEE65AD799EB0F
F9CDF4AF02773E9DCF4FF3ACEBAC6FDD7187D6BD287C6FF08ACE13222EF199BB86DEDCAE5DAC6F8EDC311F57CE779EB9309F
9F5EBB8E2EE8933E460F99DBAA3C47E3259EE337DDAAA7FE30A9D3FCE2FB87F9C95E170D9F137CE22B03260CDFABF8BD5BD2
399FA68F6AEFF4FB7DECCBF89AF5ACAFAF388A53DC786E2EC477BE618EE71FF88EFF74467FF4456F78ADAFBED67F5EF3B8CE
7C71EA6DAFBE5B7CFE2AF90ABFBCF01FB3FB17DFD39A7C918EE04E4FF0810BDFA7173A510FFA73BD3CC4254E7D4B1E74675D
3A87FFCEFDF88AA63DDE8FDF93E91FE6E2DF5D38A7F87F187C1AFFF33E37E940DF3BABF9BC5D8AFF579BFC29F7A3233A4F88
E70BFA529E5F7FD7A451B90FE7F9F8B8177B96FB79EE9BFF69D9B07C7DEE835DDAB52B9FA7CE579E357CAFF0057EE7E87DB8
E439BE6ABD45E53EAFEE79CE4CF55777FEABFEEAA09FCCDB77FE21C5BE6AE5C049C5E7A952BDF531F3ADB9F2BD3346552D1F
BD6F8E1C396976FFF23C87C7EAAA8EEA3A70D9A831C5F384A467BAB49EF96BE3D31317C5F77CAAA37EA99EF890F74D8917F2
CBFBD1D40FAC677DF350C621E5AF5F9AEF8E9937EFA7713D7DD3D3C5FBCD68147AE267FC000E70913FBEE19FFBEBFF6FAF9B
FE551CC5CBE7A6777CF5DD38D259B5DF8D1A13FFE7159FE009477A763F3CDC999778BB338FF17F673DE03F3D884B9C7C73C2
FDF37629FB255C067FAFCFBD657CF55BEBE239DEC3033EFABDFB58DFFD866F9CD128FA04FEF19DC1BF197759BCCFBFF02EFB
CF610327459FCD3E947085335D8A5FDCFABAFA5FB26CEA59E18BFC157FCC9178640E3AB1C6BC95651DE355F6EDC42F73F13E
F78DF955BCD64FC4257EF9F0E5CBD78DBB2CFA3DDEFCE25743CF8AD73FFFDBE046F1737D72D8EFA69E15F31FBCF8F3906613
9717CFFBD26BF8F1D17386CFB93BE2E537F2BEEC8189CB635D71A89B3EB66CE49C07E238E484EE07067E70347798AF0F6832
F9CFB12E5FA727F555EFBC4EBADE5CA21F9E3E66EED6E2FBD0927FB87F7E2E92FA169EE0075D0D7DB6E3DAB8DEBAF9F9DAEB
43AF0E5ED31B1FCDE7EF340FA97B9E8753FDF1EEF26FFADC1B7318BF18F66E9BA18183EBCD5DF9FDF4DA9CA7CEEA8EBF231E
6ADE34FEFEC07979DEF4BED7691E1D7548A3E782FF5E9B63F11C2FF1DDF9794E4EF3B139F4B09A33FA076E7C7AF4570DEAC7
E75E9D67CEC52BFE8E5FCECFF3709A5FD54D7F5B74C2BC9FC6112FF019BFD54B9DF83E9EE8A3F882FF788F1FE3AAD4ED147E
240EF3331DAAA3BAD2AF79948EE9835ED46DDCC573EE2ECF97782F3FF9EA7374A71E597F894797CDEB7E60F43DFD20E3B193
CEE92EEF7393FEF0E88AA3875E1D3AE63B7C579CE23697FD61DA9C6B82DFFCCF75D681973CF4BD397D2ABAE02BE74D9C59FC
7F02BE98FB5A9D8E6BCBFDCC7579FF97FA131EE6FD51E2E3D4F3E66E2DEF1FE463BD4B7ED56668B92F0E9EDFFFB3326FF431
F7B18E7507FF77C407F1F9063EC227729F6DDCBC69B9BFE28D7AA95FEEBF6F377AAEDC87CD23C75F38F7F8B2AFE439759706
F5CB7D35C79FE2D6CFD42BCFD5A96E63CE9CD92CE2CFF35DD29539D0BCB6B9E5985F453D1CBD6FDEE2BBFFF8BF0A05FEF49D
E7DA6FEB167FB779F6AC7FFFB4883BD59DFED51F2EEA8E07DB660C2EF031B78A4F5F3D79EEDC33CBFA1697F96B7DADC97F2E
F6BFA95EEAC42FE8F2F2EF8FF820E6747AF8E0AB015F97EF2B0EEBE5E76D698EC1DF0B064EAB1F7AE5A34316F46A537CBF4B
AA479E2B529E78889743AB8D793AFA133FC2AFEC3FC7F7FF2CEE874F74A98EEA2A0F79D1255EE5B937F14B1EFAF9CA1633FA
C7517CF88CDFF2928F3EED3EE214B7F5F2739DD4A747559B563FF2C53FFEE57E744127F896FB7EE21D5FBC74E098A7A3AFA9
AB3A5C547BF225C5E79E938F5A978EE80AAFF46FFC1AB67CFCC7D127F8A9FA0C5DD3F58CF8FE0DF99B23DCE79223C77F1C7D
4C7CF2D6DFDE386DCE35711C5177F225B1BEFCDD477C744457E6C47D2F1DBB3274C4D7F190AEF0519CE233C74CFA764E9FF2
9C8D4FE212A7FE849F78495743160EAE17FF7796AEF8A839F417A3876E8C38F44973DF81B3A65C10F1F0EFEC93C7F56A13FB
49BE30FC47ED8BE7B7E23637F163F1CB870EF00ACFF022F324E9437EFAD5FC91738F8F63BE6FBA9FF9CCCFF37E34F5A71197
B61A1B3F779EB96CD65FE69E59DEEF899B5F9EF5D75907C73A7CD17ACE773D9CE1CBF7CD5F47AC98B94B5CC7EFC5230E73E1
C8F14D5F8CBA786D1EA417BCA11B7A561775CAEBEC345FD24BD649E2CDE88B5FBA2F7CC979E644FAA017F55137BAA093BC8E
3933CD91B98F1CD4F58CF87BDBECB3498FD907D27DCC2F7F3C7FF6B6B22FBA0F1DCA075F2E3B7670BD623F96FC6EEC8FEBDF
587C1E24C5617EF473E7C3917EF3BE2FE9D83A79FE4CF3A8FE76EA09736E8C3AF33F3ACAFBA9A427F3C209AFBFF283B29FE1
5BEED78977746AFF48AFF0BBF83FED3F2FF7ABA15DFBFD3ABE1715FFF4B709AFCC3AB83C2FD2E1B0EF0FAF1A477EC127F036
EF23127FF14C1DD535F7A3A1ADC696FB528E27C5A13FD1D33935A6DF547C0E20F9A2BAE73936D5DF7A179EDDF4C572FFE12B
F292A7BAE43924D527F797E12FDD57EE33EA92E7BA541FF8C085CEE9E092FF0EAF1A3CC73FB8A88B3ABD3374ECCAC0293F0F
4BF38AF9CDFBCE13D7B99FD7BFB1DC675C672E587BF2940B828FC36AF71811B8C257DF7BFFFCA11BCBF399FBB92E3FEF4973
001D5C7A53BF5FC7DF7BABF7B83ED36F0A1CF026CF07F78D2A7EAF9CE7A3C4A3E9F3E6DC589ECBD5094E70A33F739C38C50D
47BCC233F1EADB2B8E98B94BF88EF5C42F9F8C4FC245DFC5BBF36E987277F419FE6ABDFCDC22F55D7E2F2F79E273EE1789D7
AED3C75E1B3D7B5B1CE583CFF8FDE1823E3DE268CE850B5FBCB8E5A8261117FE592FEF53F5BFA32BFCD5AFC4EDE7DE771EDF
3ABDE5ACE3220F3E34A6FB94BBA3AEEACD7F477C3BF68B622E4CFEABFEF2A12F7A1B7E6FE76A8137DCCD0BF48A5778A60E17
9C31F1CCE2FB5C527FC8C72A638BBFF782CBB08103BB842EE993FFE1275ED2F3A8C913CF8CF8F37C97F2C871A6F8CC19F2A1
47FA34DFED77FBF8E2FB4EF8BDB9E99083A7EE19F8F2517ABAA46E8F11B1CFE417E6AFBDEF1BFE75ACA7FFA8B3BAD3CD8887
DA7E1DEB8BCF5CA39F9C78EDECFA65FF3377FCEE9699BBC6F70DF057FCCF3E9C74907DA7EABF6F29FB8FFBBA5F9EAB12FE70
E7FBFAF6316B661F55F657FA84035CF030F7D3C447BCC043BC1CD9B6C588E08938CC4D3FDB65C0D7C18BAA0FF53F3DF2D13F
F3F93BCF5B89CF788C3753EF9B755C797EC557BCC3C3AC4BFEDFAA73B5F89E56FE356A43E327E27CF73367A9E725D5067689
FD18FFC2FFEC4B4907D671BD79CB75D6812F9CF01CEF1DF15EFEF4957597F21B7375C30BC207DDCF3C856F799E49BCC3FFBC
8F483AA097BC3F4ABAE1A370863B1F3A63FDF4E9711D7FCDF1A438CC67789EE7EDC4F7E1D5FB1C1D7CC44B7D123FF39C9978
3AFC83A107C6EF75F8055FC9FDA871DBAFCB7D897EF37E2AE978ECA7CFFF21FE3E477CE646BCC713BC71BEF3C66DFFF72DF1
FB59F18B5B3F9BD36B76FDF23E048EB9BFDCD56244B9CFF07BE7BB9E5EF63CAE4F8FE027DDC041FE743EF980E9D3E398E7CA
A4CF3C9FEDDFF889721FA08F8B7F36F4C0E2F945E2BDFAD2B93A2F3A68F651E5FD9EB8E1057FF5D0079DEF7A738E79E1EDD5
53F72C1FF3F38FD4EF3776187F55F17934CF93D2F5CB3ACFDCB5FC5CC0FD5CE77DE799E3AC63DD11BB779B1CFF1F40FDF443
3E7DF2CA99E3CBFE87EF67FD6AEA2F83577C9ABE2FEED7E7E8F83C093EC1FBBCBE0D2F28F72B7867FC934FBCD772F8D7DF79
8E98E21E79DB88CAE705930FD281EBDC571CE635EB58978EF1046FE0207F7D575DF1168FE571EE4F9EFF43B92F4D7877EA2F
CB7303BFDBDEB8FFE9E5F9497CFACFF1FF9CB54FD9B7F8DC851D473489F9145FC581CFF8AD0EE71C3BE9F8F04B3E3D6259C7
EA85DFA47CF467739D78C4A70FFFE1BF3346C6757C68F48F466F2E7C3BF9B23A8CA8D9BF6F7C5E840EF999F9E2B0CE53EBC4
3AFC0FAFF23C90F8A5BEEA4AE7E69703FA4C782EAEE7DF39AF948F39003FB33F279EC20B0FF33C98FACAF9FF19BD39E61AB8
8F5B38E9F8E239A97933F527B89E77E4F8834357AEC7F78B3EE93639F8CA17C62C1F7F70F97AF899A7F6A9317249C4ADFFCC
3A6CE6F8F21C8FAFB9FEF8907433727CEBCFE33AF9E7B923D5453DF835BEE327BE667D263EE3775E3FAD6B5EC1DB3C4725FE
D2211CE0C2D7E5254F73D0CF3E1CD82BE63EFD8ADF9C7AD9F45E717F7E3C6A6AB3EE81BF38CC2F789BE7C9C45FBCC5573CA3
97EC63493770BDA845FFBE455F4CFE92EF6B5E4A734FF6FFDF75AC5EEEC3F89FE7EDA403EB593FE395788787EAAC3E788897
7491E7FCA48FF9CF57EA90F75B09F7D15B5FBEB6F8FB9514779E9BD2FB5E9B63F0165FB3FEB7F53A26FE4F1B3EE91B7491FD
67277D782D3FFD509CE2961FBD8863E4B8C1FB851FE43928F90D1DE57D4AD2D34FFFDA7345E0B9479F1E27453CE688B1A7BF
F09BC0559EE6A0ECF367B7FEFC3B7E9F78ACBEEAFDC6F7678C2CEF9FE5637DEB9A9FF039F7D3C46B78C2511F9A3E647AAFF2
DC49B7FAB9FB8A43DFF8CB3E338695FD75E52D53EB445DF3F383747DEE23A737EB5EEE27E610E7BB7EDC847F7F14F37C9E07
D31C48BFE6333A562F75E247747CFAF0A917C7797C7A7D8D09CFC55C939F03A538E8E0C27306EF17F7A54FF38BF35DCF57E0
0647F5535FF5CE73E15E2F5F5BEE0F239FE93221EAA01EFA67D6EBCF7B1D13BFD7C38FCD7D462E89BF5FCCCFD5527CFAE71F
7BCED8A3EC97F290173F30A758C7BAA34F1956BBF83DAC3928E9C075E2119FF8C5AD5FE209BFC697ECFF535FF84DB90F4CBD
78EAC5DF99AF931FE9E7BFFDCBD44EB12EDF852F5EE1D9077B0EEC559E7BE4A9FE677C35694CF0871FE349F6FFC417FDFFA0
AA136F8CFBF1D7B99B670C8BE76879FF9178C06FCE9F39AC767CAE1A9FCC3BE211DFA843DA5F18EBC24B5F1DF593BE5D83BF
78CCFFF2FC77D6BF3F2A7E7F7DF6B3FB97F3C72B3C1BF3C1C899513F7EAA9EEAABAEF4833F791E4B3CCA71A6F8F46973C4BE
4F8EBE3570D21FF82B7CE03579974963E23E796E4DFE4B37EA8E07F87EE1735D26C4FE90AED4EDAC23C76F09BFD707F8D089
03A77D1B71F23FFDE0BC9F8D9C19CF0DD4E9E31A3D4E2A3F77542738E227BEAA67AE6FD2937AA9139F1DB770CCD662AE4B7D
451DF0599FC66B799C73EC98ADD127C53D61F9F82DB14E9E13CD791B5A6E0CFCD5C39C429F19B784A379E71783062F8AF3F5
B1D7AE9BB1C777F69DA97E9F5CDDF5F6323EF0CAF74DF733AFD057F6DBA43338E22D1ED3977982CEF419F1880FDE177CDAB7
6BECBBF8119EEFCCEFD1739AD42F3E2F96E2338FB82EAFC35F13BF2E78BBFD8565FFC647BCC5E3BCBE3927CD37F895F996F8
F2E6DCA99DCAFB6479E26DEEA789BFF495E7F6A4B351B7F7382CF88567FA06BF3FF68AE93F2FFBEE9C6AD3BE2DCFC1F4664E
F9E97BBDD6451DF5D531AFBD785CCC33F231BFE8CFE2970F3D665F4DBA1C7DDBC01F149F0349FEC2CFE8821ED4CF7DF3BC94
E61B78E33F3DE0B77AA91FDD641DA53AE47EB17FCB8D65DF5CB37DE28DF11C2D3F174979D2E30E9FA9AC97F14FB8EB7FE2A2
4F79F2A1534E98BA4BF0589F3147B8AF38C67DFFF9ADC1973CDFA539888FAAA3BAD275DED7247DE775CC5169AED23F7FB76E
EA1DDFD9BFA47AA9135FCBFDE5C426F5CBFEF8CE0BA36F8DFE9F9F1BA57CF0FDFC8E037F107FAF44B7F8A9AFE0A9B9C33AD6
5557F552BF510B3AF52FBECFCDDC94FA271D5FD0A1C761F17D58F8E1FEE2A17BFDFC900E13B7C4DCC28FF37CF6C7178F2BF7
8199F3A6EE529EF3F882F5DC571C633E1A52A3F009F354D2C1C24195BCF9ACBABE7FC5E04531EFE7E7520917F9CA533FA65F
FB0E3AC6373EB233EF4EED3F6946EC63F8B479C47DC5017775C7037E2C7EF9E0559ED312BFB2FFFFF7F9ADE53E605ED86FD6
985E719E3E30FAC5DE0D8BBFAB3507259F1B7D71DB13CAF388BEC7E7FEFCBDA9738AEF274C3E4A1FEAA57EFCE0BC9F0EA911
7ACB7C4D79AA3B1EA8A3FAD1C9477BF75A579E1BE0284EF1E997D3AF9A34A33CDFF0B5EC63C775EA5FE09E783F6EF8F026C5
7EDC1C9778B4EC3753EF28EFC7D4039FF3BC9A78ADCE7FFFE184E2F3267C5AFFDFBBFBD07181973E037FB8F33BFD5FBEF287
5FC633E94C7F10A7B8F9E839170F6F123EAF1EA3B7369F5DFCBFE38497BEAD1EF8834FF4063738E279E67DE2FB942F269C15
78E73934F9F584EEA37A54E6AB34AFC15B3C290EFD1E7EFF386BECCBA13FFD211F6F18D523FA983C274D18FB72F1FC2BCF71
95FBE07F9E9F930EF0993FE2359F3EAAFD9415E5F91CAFF10ACFD4E1FCC77A378CCF85679F4AC7FC7EC2F1EDDB276E29EFC3
D50F1FCF1FDEF684B21FBFFACDD439E15B793F94F4366660A39BC38FE0658E50BFB3BB3E3BA4EC03788EDFF4661EA9FA60DF
EA5107FDEDD3CE1DD6C775E653BCD43FC52F9F1C8F79C67C33BEDB6EC5DFFB277EEB27F898F999784D47F8444F63DEEFB736
F0E60B7C8B0FFDA5DF943F85AEF831BF841B1CE931CFE74997635B35D825E649719B0BF03ED725D569E3C9637AC5733FFD18
2ED9CFF76A3EBBECEBF8421774021FB8E843E2118739838EE8475DF5D5C3FF31E9A0F023FD419F17A7B8E926EB28D581EEE0
4E7F74A1DFD0873CE88E0EC7F57A6E68FC3E48DCE6A2FCF3A44FF9D075D67B8A4BDDD59B6FEEDEB7EBED11FFEEB7769E1AF1
EACFB9BF546B7473B9CFE478521CE36F7A7648E8F6BD1B868E8BE7D9F9794FC20B1FCFFB45BFB5C5BC9574BEE48E292B8237
FA019EE187BAE7F937E9FEFCB3BBED16EBE3D3E8AF3ABC1C7C546FFD766EF5297F2ACF9D7C874EF33E22E9553FFFCD2963DB
47DDF525BAC8FD3AE9C3DC215FF9CB2FE79B7C2B1FC59FF2C9F1A7B8F563FD073EF0CA73DEBD0D7629F7197A3DE9AC4917C5
F9FAC9B87A830E0D5EE5F92EE91ECFF55D7CE7D37FBA76F297655F570F7CC00F7E0067B86F6FDAB7FA779E6F25BCF4ED7D5E
1D7E43DC5F3F19B3B9E76DC5DFCB98A792CF8DB9BAF52FE2E7799E4AFD36CF67D73E37B4ECD36F4D9A745079BF97FB78D2E5
394F0F3A34E66BBC34178853DCB3274CBA289E8FE4393AF98E78C4A1EFD2D14F6EE9B03E74474FD9F776E9F072CCABF8AD2E
EA419FF0555FF5C64F7D174FF50DF9CA5FFD4F7975C2E1C5EF8D521F98506DE879856F98D7122FE00F777E8A57788F5FF0CB
78E27BE2B77CE53FE6B5A65D8ABF034E78E9B7F8C607F14EFFFCF9B1FD9B87EEF499CCFB54773CE0D3BFFD6AD2FF042EFCF5
934E9DA796FBAE7AE778521CFA2EBF3C6BE0D0F3A25FE0CDE25E95B8F814DECFFCF384C3CB730C1FA703FCC1A77533C7B60F
DCF23E3ED54F7D4E3B69ECE48843DF9874D988F38B79CCDC97EAC64FCEDBA7E76DB10F85BB3AE4F7F94EE2DD797D5BFFA2EC
AFF498E7EAA44B730BDCE038F6C1978E8BFE022F7D5EBFF9C79011E7C73C02477D52BEF29F367BECE4E03DBCF43379FFBDFE
E80382C7D6C55B7CA59F294F8D3EA07C3D9CF8A57AA9DF98ED9D37C4EF93F23C95FA069F3BEED793E616DF1F9FFC3EE76B9E
49F3047D65BF4D3AC397CC9FA48F7157F41954F887F92BF91C9DE4BAA43A89537CFA84BE74D0EBE34E0B9EF3EF7127D45B10
B88BCFBC907DFB8F4DBB94FD5B9FDFA35FB7C74207FAE1E63F0FBF21F8929F97A4FAE12F9EE33DDDE57929E94F3C796E4973
035DEBBBF48DE7F8ADDEFAAA78C4472F70A71B38C397AFBDB9CBA4FF29EFEBF0034F336F53DDC61F56A756F17C2BC56DCE59
B069D2DCF21C4CE7D9B79BBEF49DFE0C273CC77B7CE377788707EAA57E637FDCEEF6C0150FF4C3ACFBAA9D37C4F736E013DE
9D33A8CFA0F87D245FA0173A813BDD65FFD9497FF975C2491E74675DFD07CE70D7E77F7DC8883191873EC6A772FC299F6D0B
FB378FEF81C8CF99120FE491F34AD7D3EB9F974CBC30EAAA5F65DC125EFAB7FEF0FB59139B95E770FACA7360D259EE3BF3EA
2D28F79FF19FF75F537CAECE9C95F4CD6FF0046F76EBD3765DE873B75BDA5C19389B17CC1DF2953F1EE0D5CEFACD737BD2F1
D8D35BBC1B38E4B929F5DB7187763F26F49BE7ACE47F6B8E1E775AD4531F52BF578F9A7861F80F7ED33F7E9DFB79BBDBCBFD
C47DDD4F7FD5B7F73A7A6093C0439FA1B3B37FDC7F4DEC23F0159E70A467F9C35F3DF413F1CBE7E3EADD1E2BF77BF8C241FE
7C165FFEFAECF8BA11375FC7ABDC1712BFF059BFC46BBCE283F885DF7083E3C445833E2EF4682E4B3C1ADBAAF183E16F70D4
6FF9FDA1F5C6352DCFED4B4F9ED8ACF83F39C9D7F00FAE19E7A41BFD1C3EF0CAF74DF7D3FFE07EF6CA3AB5CAFCC76B3CC19B
4D6F8F18137D43FF50A757EA8CAF1BF8E7F93AF9235DFCB877CB0B638EA60FFC3A776A8B77BFD3C792EECF59D5FD98F83C3A
7CE597DF4FE7F1AD338F19F471F839FEA9F3C9BF1CFB78E8483FE15BF0842F9F38BAE5848565FF1E57B5C1D7C5E768125EFA
B0BEBAE7713DBE0C7EE97B93570C5D53CC2DE63BF5AFD271583CBFC8734AF27B7D0F9EF0CDF7D5F753DFC52B7CA22B7AD477
E9F2B3CE6DAE0C5CF3F3AAC48F595BC63E1ECF0BD489EFD30B3DD08D3E71C61143D7841EE18EF7B98EA9AEE35BF5AA59FC9E
CD7C94FC491DFED662F46F03CFDCDF133E70E1EFD3F61CD1E93BF360C2559DD59D0EB2AFDEDBF8C1B2BFBACF691F8EE8143C
948F7EB57FF791B3E2BEFAC38C9AA37F1B78E47936DD9F4EF50F7A1D37BC6EB3E27343E6A3345FD00D5DD0C9E7D7B4BCB05C
1775CAEB984FD2BCB2F5F5814DE2B9427E1E9378B3EAE9714D8327797F9CF80D7FB8F3C7D7EFABF09B8FD039DEE22B9ED105
BFA38FECDBDB1B7C5DF66FBAC31BFAD3F7C42F9FB19FB63936D6556FFD0D9FD431FB033FF8B6E3B078CE8E7FFAED4FEE687F
6AF04BFFC453BCC517FD046E70C4A3B3EFED55339E4FF31D7AD78FE95E3DE988AEF89438C5CD578EB96DFCB872FF8183FCF5
5D7DE0888BC7352ECFE1E2CA71A6FBD10B9DA8ABFEFCCB1B079F1FD7EB6FF4A84E74C917F1046FE8868EE04A47FA253DE5B9
EAE2BACDCA7D064EF4456F13FFD1B771E098E7AFE41FF49B7D23DD1F3FF0102FAD47B770F970418F2FCB73035E8EFB7EB391
C55C69CE4AFD79FCC95D3E2E3E4763CE4AFEBDE18691B3824779DF9FF894FDF0276D8E2DFBA23EECBEE2705FF7D3A717751C
3FAE3C5FF291ACB784277CE1062FFAD74FC4296E3A387352DFC6F173BC97AF3CF932DED23FFED25D9E7B93FEF0FAF8A3C65D
5DFC9FF3D44FF0996FE275D633DC128EE34E78F9A4E0255CF461FDB6EAFFF63E3778ADEFAD183EAE71797F82C77CFDC0B346
D52DF6DBC9BFE942BFA48FC9470C685CF031F9203EC229E396F497E34CF1E9DF74A38EEA9AF5F0EF7FFDA1AC8B4FDBB73FB5
DCA7F1068FCEF96FB391655FA4FBB36775F9387C158EE2CAEFA7F3E62F197775D4159FF817FF831B1CF373B743C79D54EE03
5BBA0D3E3FFA8F3E94E79646F5A7873FC83FF7C3EDED1F2BFE8F52E213BFE71F67AC18D038F88997FA211CE0824727BC33A6
7FF0419FC9F775BFD447F129E396709C7ACAE03B8BDFC3267F54677D495EF2C4073CA0373CA7077C97AF3CF9FB84913D0E8C
7E9CE7A3E42BB99EF35E3EA9EC67F3F61DD3BF3C37F0CDF1BBFFFBB2E2F75E294F7382BEB7FB908EDF441CFA15FEF3653AA0
A33CC7263DA943AE4BE23BDF3D7DE6E03B8BDF87A5FAE927BF9A34647CACAB0FFC3FB6DE324CCBEA7BFF36BF2A4A0A487729
5D4A770888D2523A3230D4D09D03C39043A720DDDDA1B474977449A748873F417DFE6B5FEBB36633C7F36A0E98FBBEAEB5CF
75C6DAFB8EA13F8586AD5D27759133E8001EC2CBED112BE3BBB94BFD9AFE5CDEB8364AFA8ACFC2D73D1DD7AE133CC005DF07
67F0C5D7C88FE2912BE3CBF910EBB17D42B2D57DE5BA3CEED485F585FD7D027A7E9170E432E90F731C7C452FF81ABAD9796B
755FE18BCD95BABEF51BC7BD70DFE3CE1CA4F9035FE143EC9CABD061E0F6D77240755F21C14F9DE59C01DE3F78775E45F97F
720E7EC37F780FBFD1173C4667E6237A7DEE479ED017FA64EBD2F59073F0A85C95D969E5DC0F9FB2F39A377F19E3E7CC8BA8
61837D9CC11DDE190F95A75CCFAEAF75DBB9EA5B6B1F08EFC81F72187CC00B3EA23B74882FC207F8617EFE24B29AEFEBA667
FAA5FDDBDC625E1EE183CD3BEA07E8D7E658D531FA229FD019FBEA0F168FD922CF2387D7CF1DD9DD9D0F31D7687EA22F7405
0FD029FD46AFE46AA2970BEACAFDC9C98DB9A62E72DF17CD7CA4BE8E9ED017FD013F748A6ECD9FCB8D7BE1FB34F55337F97A
ABCCF20D725E4BDEC03FF44BDDE818DDA137D68FBE791CF552173AE5F9C7FF0D7889FEF1AF27D5A6FCEBE7347CA2DFF0017E
801B789113F0011EE05FA637C507BCD00BFE826ED6AF1AF28F9CF78217398C3ECA4C989747CE59D121BA200FD007FACBF5FB
CF93FC7C230FC19F7E5CFC67ED03B90FBE8C4EC86DF0012FAB933941F3DBCE912BACB822B89057D44BFDE8DB7CB2FEC8EEBE
DEB6D65998C6F19C394E75865EE9237DC527CA1D9BBA48CEA9C11DFEDB9CA93AE0FE3C9EE7DBEB2E45A78F96E791ABE8F4FD
FB2397C97AD1EBFA9BFDFBBBCF5F32BF90E71B26E69175DB1CA1F9806F653AB52E8E9F33F82EF8D38F63E97F9EE4CF9DF89A
DD97FBE93C613A044FC5F7CF7716D495F991FC8007E8BBD4C28569E43C129E513F759307F028EF9A7535FD9C795E6BCC163F
47E18DF967D521AFE5DEA66533F68A9FD83CA53E44DFE937BA82CFE8045E6FEBB2F88A3CDEE641E50B39C07A593FF9439DD4
0D6FF165F80B1EE0037F0F575957D39F7BF047F848BFE81FEB659DF888F5097C14AF0D2163DCF787833B3960E7F5F9166595
EB9227360F0D597C4578405F2DB792F7EFFF9A9FA93ECA979F98C7F76FDBEFF55AB34D7C9B3CB9BE7EC515F169FC1ABE9EFF
745D1CF77D3FEA17E8139EE347F01D7D918BE88CDC8BB362FC0AF93DF9B963C7B20ED24F7C9F3E1B3E8A0B79637A665DBACE
03516BB6F9730039617E5A7DC65E398FC6176C8E2CB6AC839CDF82AFED5B72AFBA2DFE475EC143F8078FB99E5D5FEB4347F0
183DE1C7E00CEEBB0BAF38E85EDFD53C030F3B57FA65ED313F4F1EED9C3EDACF33F440BD8577AF38E8DE8FA2EBDC7B74D56D
3FD7B90FFE9DFAD1CAF5A20772C6F8AABC8267F8167C801F2F1E0F72DFE3C53C880EE129BC45671B320D2FECF685CC419A87
E809DDA1C32D2D67EF95BAF10B7C65E33F93EBC9EF2DB7D58FC931F0012F74507EC698B5BE8F5A3DE4B2E61F7A2457D025BC
85AFF0037DD16F7466FADCDC279EAF53DBD7BD37F3915C9F1C635DAC07BF86A7F0165CD183AD4BD7097EE8029D9CA9B8F698
3F1FE36B1B3E1AB4D9CDA7E4BFE6DEDD438BB2CADC4EAEC063746DF3ADEA9B7ED127FCF1FD5EC3068BAFBCDF67681C99EFC8
5B74075EE80FFEE32FE800BDA013FA61F593D79AA3967FE367EF9573597C017DA017700527F4827ED00BF9816EECFC65E6DA
8A7E8EC1AFF2E787BF9693E8947572FF67DF8C5FF1DA3993F298FBA35BD67F35DECAF5FEBE0E7D9287F48BFEB16E7088ED5B
E5DE9A5C4FCEC1E1D38654FDCAB8D7B71447E6107B9E3E9EE7FF1A326F967BDF85FAA3E9B5D1F8EBD23F9B9B34F7ECBC7BE2
E2958223F98CAFC04FF84A3DD4B1A14C9F78F23E17721EDCC0113D18FEDA0F3B675CBFAAB43F57E3C7E009BED44FDDE413FD
2CFF64D066DF2FD135F987BE4F7E1FF003BDE177E4E187AF26F691E7919F9BFB4D6B2EF384CD53EA733607CC98374BCE3BE1
B7EDCFA3D784FBF9435FE807FA343DC61D7459FA8B2E1FFC35F391F4D5CE27944FF01C1DC2F7ED83175673FB38F541FA0FEE
868FE2B5E1D2A8D65297E5BCE680F9E4B57E655EF34BE563B969E3AFFB7D051F70C137C1037CE0FB8BBE43E3F8F9046FE021
FE0E1FC90D700017F8481FE92BEB623DE4CA6F9DD784FBF322BE6CFA615DBACE4B155695F6F709E8E176B3C52B653DF82F7C
B539A6EBC26AF26FF801CFF129F86EFBA5B87392CABFC937DB4FB65A5D5AEA27C72CDF153770343FED3DADB9F0035FE027FF
6FBEA23CA0BFF47BE7A5C5FBDDF99BFA383CDA1011FD54F2C0F2597309BDD83CA9BA8187F00FDEDBB9D28FABBAF9F32A7A81
C7E8E6C8D8D5A5FDF9833CC047C19F7ED8BEB1F808F7770FC8439BC3D22DDE2F3AA34F9B164C7AC37D1E54F98A2FDAF96CB6
258FA45FE48CEDC74AAD3C20D72727C1075CC83FF3BD74A35AFBFE879FC12B78B6B5FCCC3AEEF3F5EA0BF80D3C8257E8694F
C8B20FE571E4107D7BFAF6C43E7E6EA11FF401FFD1C3C1AD2B0F887EC1979CB3F97BC6B20F053F70B47DCBBF2BE2CA3A79DC
FE3757C495FA793E758133F8E2E3A637F051BC369419D043F0B03942738F3C615DAC13DEC257F8C7752C3735B7C0197CF12F
7407FFD01F3A253FD0EB85E6ABBA090FF005FCC8FC3632FAA9EFBBF00E1ED24F7B1D227AB2FB7E5FF28AFEC07FF4F0E7E339
4945CFE4017AB01CDA30B38EFC1B5F4077F80EFAC30FCA7E3BE90D39D7A51F1B72F57DC39D7790F39AEF1BAF8D1DEB5E9F61
8ED0DCE03EDC17BFB1F397A8550DFDBC425FF83B3A33FFA21EAD0F5D938BE81BFEC37B78804EE91F7ABD717AC92359077E8D
0E8D7FF44BFB07DFE13F7D7BEF56FF25B29F7AEF76D4E7F2FFE4F6B65273524B3DF81D3AA66FE80BBD9163F491BE9A3F6F1E
D0C3F769F00777720EFF8087F0129DB27EF48A4ED12738A2631EC7BA9F6F1F91DECF6974C83AD03BD7B573F02FE7BE92BE92
E7765E797EE9757F6EDF9277CA1BEE7B359923D48FF175FA45FFCEF65AD5D09F47F135F2903AA9DBFA0A9E8AAFE5F4D639A9
2547D00F7DA47FF8C5C6A223D24BBED9DCA479887EC93F748C9ECAA51A3BD6F7393B47E8B2F2819F7BE0C0FAF14DDBD7951D
9D46EE43DEA22FF48CCEAC4EAD8F7CDAD17EDE5AB75F537F84A78F3B4FCEE0E704BC040FF041DFA64F5D17EB345D6C0FFC03
7D18EEC7FABEE1F3199ED247FA0A7FF05978843F943E3CE50D3947A60EF33BFD7F1EB761CDE0A5EE7D5ECC2F9A27A78604B8
C31B7CF65AA6A5D7C5EFF0057445CE800F78FDF549D4E76E0E632E53FEFDF1F3DC57E287F8323C86FFF81D3AB0B964C4BCB5
720E0A2FAD7EAD9B3C8347E576BECE6FDBF7B6589957FA4F8ED9397BBDA963852FE4D8E61C1366B9BF03CE1CA4BEBEEBE8C2
64EE1C4CFD14BEC01F78834EECDC6AF1B2C1FE7CBBF5AF6917C5EFD00FFAFEAB78C00370034774018FD1C7865A018F6C2ED0
5CC22FC119DC8D97E0A378E11FF4977E8303EBC767E9037D4107C727ACCCEBCF4FF8BB9D17AF9E3F55F024DFD02FB9858EF1
D722B9172613BDC0035B2FF9AB794C3F73D75E9149728D7CB37D6383814DC4BFC949F404CFE1FDDE7A8B4BB97348F57DFA6C
3E5965F052DF2F9F6D1A9DC6F77DF4407FE92BFE7574D18A4CFEBC424ED8BC386F7129799D87F5B38E7CF997ED77DF23A979
75E09DA5598537E40AF5DAFCFD72A97B7D92EBC2DBD87C3D7470D97EB93F75703DFCA3CCF109B3E4FC1ADE6C2A35BAB87BFD
96BC55BF371F7C6FDA45396FC6772ED75A3658EAC02FF035F4088FD1A5F1129C1577F32FAD87FAD02FB9828EC937FA429FE0
11BC82BFD44BFDDCCFFC7971A0437C1A5DE353E81B7CC085BC815FE8029D3C9C3775AC7CAF95CD3BAA5BDB378E1AFBB7CC33
E4A79D078D5A76DE9F57B7FD3D7394F40DDFC157D01D3981FEEE7C357FAAF82EFE8BFEE13FBC8737E88B3EA1B32DEB278D71
E7DFCA6FFC12FEC267F0468F361FAA2EF11B78026F8C07E0AFFD2067C00D1CC1157DA1B78D1D861E977A2DFF35F7F0FB84A1
D30F0A3EE42A7C2FBB757471DF2FD11D7A032F740A2EE8F5C59C814DFCFC8397F85FF24A0B5E085EE41E7A45BFE047FD96D7
9A7FF49D7EE33BF83DFDA5DFE081DED13F7DA47FF83E7EC27AA8CB7CFFDD99A3E4BC16FD733D7C04BC2EB65D76DE9FBFF15F
D33338330FD2EF6BBD37FBBAE33EDC17BDDABCC1F999E28E3E320F5876D1DF6F6DD81DF0CEE608CDCF9D79E62472FB26F547
F4841EC92774891ED13FBAC43F4A57983446CE9D59A7F98BFE3F8FB37AC871CD39DBD7251BE2DE5F434EC2C772C3871EF7FD
F269DBB17FFB39844E360F1B5752EAB23C57FF4667F4853EDD5CB7E0859CABE00BE813BDE0B3E8C672FAC89C44F23D60F063
4399BE8FDCF902B9AD797DAEC7B28BFE5C8B3FA20B7C0A7DDC9F3AFDA07C1F15FE08BF8D378A03B8FCDA7FCAF7EEF57CD503
FE447FE92BFA011F70C1C7C90DF004DFBF92F5ACE97E16EFF39F3FA7C1AF6C7D970DF7F36DF78479DFB9F719AA9FC22FF3B7
C2015FE1BDED97D24D5822D723DF6C3E7E6BD1BB820BF903DEE08F3EE03FBC4407F801FDA5DFF82838833B7882237E815F7E
3C63E65EA9879C79F1B85BDAD7F051BCF0B3C22DE67D27FB217873BA4F80173CC0C7D14D998EE34AFA7EBCA9DB88EB924736
8FA8AFDBEB769A2BE8C67402CE8A3B7C8787F092BE65EFB5F49E5C97DCDB9765412B778EA7BE4F3F63F2A9EF23DFB74CFF5A
3FEB611FF87EDC88A1A277F290BE95EC39E57B399FC62F582FEB240F9EDF1A12E5E70A7C856FF00CDE733DAE4F7DF8FAE767
17B412BDD2A713514BEFF9F30AF9011EB9DA2D7925F8927B07CE2E2AE45EC755DF07A76BFF2C7A57FC031FC177E83F7C80BF
E808BEA227CBA12C8B0AC9F3A97B7BF4B457725FCB43D5BDF59BBE689F8E8D5CF24AF8CA7AC82DD691B7CEE2AD920BDC0F3D
E2FBE8F2F0C2C55BFD1C65DDE48F9D0F699FD02FBE868EB7A6F8319DFC4457E6D3FD03FE5AFE6B5E3DBE346189F0087F47B7
F0171DA12B745376D088EBBE5FE21F69CE2CDA25F727C76C3F967658B4AC9BDCA31ECB4DCD2D7C803CC00FEE85CCDC2BDF53
8A5FE33BE808FDC0577060FDF8327C849FF005DF8257F00CFDC21B746C3EA238800BFC8137F89D9DD7DCEA5953AE8B5EC901
F0A71F1BD6F44F297CB43CD7DC83D7C53B4F737FAF0D9DD31FF485DEECBCFBE1A4DCF26FF2D6E6F804B3FB094FC8497487DE
C0D7EAB17382208FF175FA429F5E3C8E18FADAF98DF2189DA13BFAC03A5817BE6279DC33E03D3AC437E807FE7125EBA25D72
5D7C015FC69F4ADDF8319DFBBBC8DA3FF347FD7F1EB773D7CCE6EE7B12C979F5037045EFE81FFDA25B70DCBC6BB4FBDC22FC
2367F003FC013CF00DD6C3EFB90FFEC0FDD059FA9D8B0BFAF98C1EF1057469F3869DAB057D323D2BCEE06EAFF369FEA163D6
C57AC895985CEC9FD2F7CB5FEF4C3AEBBE1F9EFC577FB2FD5B9F48777E43DE3EBB3C2CDACF217462FE5E64667379DD1B3DA0
57FA0E0FE823FD439F771ECCEE27CFC747D03F7AC197D10D7AC18FD0CDEF451717F4E75D7CF351FC49B98547F8293AD99377
F664F7BA80FA14FCDA141DE06AF9AC3E0D6EE085AF991EC05971B7F9A0C8E8B37E5FC90D3B67527CC10F3CD1DF8BA88E5D03
7EE8FCA5BC815F99B62EAEE5BE1F4073D5F633AD467610FCC8377C25D9C239CFA47EF28AF5DA7E5AF3C0F4A9F5B31E780E8F
E13B7E431FE92BBE52E8F0ECC9EE73F8CA1BF397A493CE4A3FD03FBE4BBFE89F9DE7BE9CDC5BEE4B2EF13CAE43DDD68F2711
CB7DFE9F2FB5B8963F17E2D7FB76CF75BCC467E145CCEB5E414EA00F780DDFE0DFC60E83F2B99FE4AAE6067C8367E881BE65
DDB638839F6FEFC7ED9656E6F8F7EB7459E4CE5734AFB6D61EB758D66DB9AA7E66F9D739E013BCA49E98FD71901FDB1F4F69
22FA4787F8C18BBE915BFC9C80AFF48D3EC253CB8FC27347497DF4FB569D39CFE47196E7EA0B674A04EB04777CFC609FF965
DC9CA97E4D1F58AFED77D5EFE03FFC4307E0F7D98E45D3FD5C321D2A3EE0B52122F284DB9F93AB9AA7E8057F4137F6BA94FA
2F3AC3BF63CE87021CC98DFC7DE79791FB8017F8833BFE62F5B0CFFE2862B99C6F9E2CB668BA3F4F90434FC78EECE0FB32FA
39527BC152B99EE5A1E28AEE53369D9746D6432E8153CE2D0B0FB9BFFBABB907AFD10BFA412FE8047EF3BC3C8B16B87992F5
FF567AE121B94FEC9C8EC9B941F97C1F7DF0CEE4DEA26BFC1A7FC107C803FCC0F687C5FB5F171F2327CD4FE34E69227311BE
034FE12DBCC63F4A2D1AB758DE7F056FF033F80CBFB91ED7C787369F097861F9AC7962BAEED5B1AB3C1F7D9B7FE97DA9037F
80DFF804B961E743CA03D3BFE209BEE4C047FF8C29E0E651CD559BEF8B4D1D297C205777E599B6C17D8F81FA143E41FFD12F
7A46BFE8167E980F47469EF0FD185FA7BFF41BDCEC3C40F3E9C5C22E8BFCFC436FD727CD4B23F8E12FF8323A4377F4197F00
477C625BE61F6BC8FAD121FECAFAD01D3A445FE88AF5DB9C577ADE4B7FDE23C7AD7E5D0F7AC45FD025BA4147ACC374029E8A
2FFA653DE87853B7E8DCD22FCB79CDA598D7CF22F5F5BA4037F49D7EE353E89AEB521778A06B74CE4F740E5EE8ACC891691B
644E409F3667647D9DC7F4857E9007B67FEBD5FD9A9CAF93CFCFB7F77F7DCE529D8043CCFE3BC845740C7FE093E95FEBA13E
EAA57E74BF67FA8C93EEFDDBE4BCEAE6DE8EA923653F65F391FA88BD0EA4F9812ED123BE862EAF6C99F7D29FA7F1E5276F8D
2920AF9BE1FBE86A63D1FED7A53F96FF9A6FF84C89733FD690CF31B14EF22AE6FC29C091F5F1789E6FF35FDEF9AFED7BB8AF
E5A9E699ED5BCA0E7AE63E9FA3F989DF24593ABD91F083BC8DC9E9E8DCBE9F995FFF30E3A4F8287CDADA22C0CFE605F55DF4
88AED0257E435FE8133E0A9EE04B7FE92B3AB1FDC9A171FBE4BAE4D28E3C93F6BAEFF927E7C9F74C7D46CAB987E5BCE6FBFE
C2B35638FF54DF8447C673CD09F86E3A517CC00B1C6CBFAB7EFAFBE1F983E5752678838FBF7FBF5D75E1E30705DB9C74E797
9A4BE800BEC13FF0067F746079B333E00DFE64F359DEF9BBFC1C7B11D5FD9AEFE3F00F9F2BB87BD60AC1857E832738E20B77
6B4C6F24FB2674859E0F3E9F53C1BD0E4E2E29AE96AF1382BE82E3E65DC33F143D5AFEA85F9E3F3C7F973FF79007F01F3EA1
037CDDD6A5EB8C79BD27F0597400FFD1393A30DE2B9EE00BDEC58E4CDA2BAF6FA34FF8036FE0ABE95CD7CBFA6D9E28347FA5
9F87CF360D7AE6FB267AB01CFA604E05D90FD127EE4F3DDCCFE6E64A331EC87AC837708ED96707BE6CFBABFBBD53084FC9A5
47F9C6ED93D7D3F035740B1FE133FC468FF829BA3C726CEE78A9831C8217F00BBEC1EB337B025CE8373967734389F95FFDFF
E56D9E5C73C74BAE8297E96B5EF81A793D079DED6A3FE57DE1A7E5ADFA043E049FE1F76FDBE6FDE8E72AEBB0FB9798F7A39F
CF27B7CDFFCA9F57789EE9F37C9F91BE4EB7150E70B2BC551F250762CE81027E6C2A35E8943F5F9033E81D3DA07BDB2F8D8A
BE2BFD21F79EEF6F7332D0899E73A87E6CBE3FFB6376A9873C44D7E819DEC7E4F4F00F7D5FC4BFE10F7CA27EDBAF6B6EC123
F48EFEE13BFC872FB7D6CD7820FB627C07DF044FDBAFAB2FD377FA8DAF3DAF17BEC6C7015CCCA7745DACD3FC77C494F7E5F5
737C015D833BFADE70A9DF44C90DCB67CD1F748AEFC4D66BCA0F66E5F7E75274847EE887E997794C750C7EE8085DA10FF402
AEAC8375E12BE4D207F13BB8BFF3451EB22ECB53CD577CA8C4EEA0DFE0F82261EF147EDEC0EF3D57A7C571FEA9BE892E2D5F
B70E3AE5FB313E00BEF881F9A0DE973AD02FBA0547F48A7EC1039CD02F7AC60FE01BFCC307A883EB6DAD3D72B3FBBB13CC05
9A6F5CCFE675ADEF41961FB3CBE7F76C3E523F8A79DD28C82D748C4EF159F4FAB46DF45D3F3FD0A1E953F1045F9B07387FD2
3E5D7F3E2BBFF4D3E611F57DD615B3CF0E726547B9E0FEE80A3F33DF4BD76FA2EF7FE6FBA9A7C591F73FA22BDB172DEEFBA5
FBBE06CD6D78023FD0B3CDEB2926D5727FFF4CF31CDDA4BE393B873F7FE243E00CEEE80E1DA23F7C37E6BC27C0777FC4F4E0
7D7EE4BFF2177F8D9B6EF826F11FF2D074A439849E6C9E561CD01F788223BEB67958F4DCE0F507CD73CD01FA0ACEE01E337F
8CDCEC5EFF571E3C7FD8E17FBEEFC303784DDFE1C1D5E4B37388DE2DFFD5AF4DE77A7DEE479DB60F56DF371F8F9C7E407080
07E63B1B837EA067FA96EEECECD36EBFA3F9C6EF793C75B0FE0AF57B7FE3F3EDDE8D49B5E47B4BD0157AB61C541CC0E55089
9915DD1CAB3E4B9FE9BBED3B9597F6FA8DFA2FFCDD786D405FC939CB3DCD01788B6EE12FBC857FF0775BE631A5DCFB13C937
F595E7B5FA7EE9FB297CDD3530B81E3A4197F8F707053B9776E7109A1B97B2CC3EEDF629E4B0E607F55B7EA9EF5B2E768C9E
EBFB0AFDA4BFF01A7FCDB76D664537C72A0F6CFE1E36B984D4491EC6BC8E12F81D7AA2FF19CFCFC9EBE7EDE34BC337C97C6D
7387EA1CFEC25B787C34E5ACEEAFE5AAF6199C63F6B581DF6D08E95BCDE51CB9A77967BA521CC0055FB17308E5DFB3BACDD2
04DFA315CC1DC63FF55D3B47D1BE5EC83427AF3F0F914FE44CEEEBB3BABBB957F1058F2C37E614F573123EC003FCC5F6452B
A21ACBE3C819F81F9BF7C7E3CDFEDECF2DF0B3B9FFE1A8C6725D72C9F23BD580BEBE3F9197391ECDFE5EF8C67ACEA608EAB6
1CD6BC641D9F3E98FDCCCF5BFC38E63C26E001BE817EF08FA7D39AA4F6F1A71FA712CC7EE6BEE79B7CD6F5992F770FFA874F
C15FF80CDF63E68F31A5DCEBEDCA8FBB1D2797907D88CD11EAB33647EAF5F133788A1ED1A7F99D5E9FFB6D4D31FC4BC1C7F2
59F3C1F4AE7983EED917C5F9A26551D12BB9048FE917FDDBF357A04FCB61F515CBF5197DABF9BE8BBED0157CB2F36EF5BBD8
BA4FF6E9941BF26F729875D93981E6CAF3FD9D4BBF36BFA87EF007FA8D4FA00FF442DFC0157DA1B71D79C68D90FE5B0EABEF
DAEB1CEAEFE8067DD35FFA8D5ED009B83C8A3FAAB1CCA7368FA85F582EEBBA58A7E947710677F4CE3AD13D7A425FAC3BE675
88C04FD1EBB36FA21AFBF981AE36E7185CC67D3F0179AE39061F6C3FADFE62BEFF5EB07EF403AEE81A9DDBDC7063F8977E9E
A04774088EF800EB8CED07FC9B75936F1FCCEBD6DB9DBF6B3ED3DF98FD77E0FBB74E4DB921F5DB1CA4F941BDF803F5818FE5
BBE6A5E96B7EB334B27E7466F3FACB31EE3C8DFC37BF501CC0C5E663C5139FD83F3ED089CD23AA4B7493E2EEB4F2FE1C4B9E
D83993F61BBF8C39EF097869FB8D5703C3DCF756680E6F6C1495DB9D6FEB3AC943F48B9ED1F1B3BD2D8BFA3907BF63E6BC71
23640E67FDE683FAFF3C8EFB5AFE6A1E5BEE1E1F5CC6D789F972CBA09FF0F2D73B23D7BBEFB327E7D52F2D97B54EEA46C776
DEA0FCB891645A79B73F8A351FD317DBD79A4E26CC903E580EABAF183FCD97039E1EAA13F4CF7258FB6C3AD1F5B27E9B17F7
4E6FE0CF8DACD7F6A3EAA70FDE197356E6679B53D4BFC88138095B87BBD70F34376C5E573CE1E3F37ADD7ABF36BF28CFE843
B96951B97DDF8A39AF0F7C071D3CF9366472F03DFC413ED30F3B17513EC16BF3F185411FE8B7BDCEA1BE0C7FC13F669F1AE8
FFE9DB03C37CBF430F3697171B9FDCE596E6D5D52FA637907ED16FF2C3F2081C1417785D32E9C8F5EEF56DEDDFD102419FC8
03FABC657DF449A9DF724F7DD7E6984D3332FAF966BC579CC11D9EC24FF88ABE63CE1B021E98BF8C983043E6D5D83AA47ED6
838FC669D473A5D44B1EC25BDB1FAB1E62E6A7E0FA5C97F5DABE56FDD472EB408037FDB6FD40F4900AD23F720C3FB3F30FED
EBA5B23332FAF31339045FE00FBCFEB05113F7B9980F1F868E977D39FE7D3C5B803B3944DF6CCED83A638A9FC3A623F57DF4
B4E7EAC40EEEEF66E20BEA07F4AF42B98878BE5ED035BC47DF9673A783F581D7BD1DE3934BFD368FA87F5D2835638A3C9EF5
9367DB1F8F7E2CF55B1EAA2FDA3CB17FC64E3F6F4F650FD649FE81C7A60503AE8B6FD13FF2C172FA44707FEA7E76BF75B89F
07E8E76CC1E07ED4C77DE02FBA4387E8A97485E893BE2FA22F74052FA9D3F6B59A2B31E7DD811FC5D67D92A53FA6143F27DF
E023FC8447E0098EF8AFE5A0AE97F5C38798FD7AE053E88D75B14E7C001EE00776EEA0F9812E376CE8EBF2C5F25373D37C3C
F5C40E32AFA2FF67D37AAEF47DDCE684CE432A483D362FA88FA0978F16844C96F36374C37DB9DFFA8D11F1041FFC26E63C26
E02BFD44A7E896FC89F345BB05EE5C5073D25ECF507FC70FD029FA840731F368B05E7C08BE94FD76C075DF47D113FAA2CFF8
56F1B8A31FBBF74F6A9FF6474C727F9FD5F25FF58D3FD0277C82FB731DAE7BB7C68F29C56F6CEE501F7F1A3F74BC9F8BF009
9CD0353AFFB5FFF082828FE5BFE609B964FDD57EA36BF40C8EF86BFC7AC3DCE795C86D748CAEC109DF609DF807F9F06158C4
69A98B7CC6A7E1153C83F7F4853EE1275C97FBB36E7CC37C447FE223AC9379033CC1D7E6D74A132FF873ACCD31E5FB5EF37D
1A9CC1973C8B398F0E7C04FFD89567EC7DF73E2FE611F55DCBA1C849FFC99C893EE155CCFE3BF00BD385E60AFA20AF3E8C0C
5BEAFE0E9EE60F3A8B3987087876A8C4E472CE27C96DE52FB8B35ED6BF25EFE0E682B3E5B6FA7DCC9C1DD48FCE582FEBC41F
6FAD9B78C1EDB362CDD3CFF6B65BE0E7017C4567257B0E2F289F53861F0FE70D7B43F463738DFA0B75DAFE4C7D9F7EA6FCE0
A7BAFE3C697AD6EB733FF88FBFA003DBFFABCFC2EBA73F459CF67D137E5B4E6C9B5CCECDABDA6FCB415D2FEB379F3A32F6BE
BC4FC6E6327EEAFF5377CC394AC033F88B5E13868E880E5E2F09F22DEEFC86BB05A7785F34C822AF77E3DF1BFF891AEACEAB
C949CD0D7C37CEFC8EB7DCF987E6CCD194538AFB3E1BB36F097815B34F0D786CFCD49C80A7D79FFF5457F86C39AF79403D96
5F9A2BF8047D8707DBA347D6973C8067F81F7C297D787073DF57F64C1FBF40EAB41C565DA28F98F38080AF3687DD9C1CEDE7
303E6A7D571EE0731F2D8FBC29F5905731FB96800FE8C6F2E6FA94E26EEE55DEAC7F1031CD9D3FEBFAC9D3A711614BFD3C80
4FE0C0FAF129D38FE2002EC7E34DB9EDE65BF253FB7935F9E468E9BBE5B3E68DE591DE973AEE4F1D112DEF4BB19C57DFB179
E5EC94B6EEF51CCDC398F3BE80AFE892BED8BE597DC472FAADA8A1BE3F590E3D9A72DBCF6DF8199B97A7124C2DE3E701EB36
FFFA61FC02D1053ABC94654A5B7F2E216F4C17E6BF813ED01D7E81FEE051CCF941C07B9B0FCE4FFD9F9F9FE8BB78E791F5E5
73C4F42F665F14F0041E592E3E98EA7C0F5CC8810F1B36BF20AF83901B36EFCE1875597C83BCDAFA57B4FBFE2F7883DF5BAE
DF98FA959FC3675304FFB6BC555C6DBE543CF1AD0B99A6FEEFB5DCD6E759BE26889816FC9D859E27829CEE78CBCF03F4663E
A8F88017FC8627F0E6D1DE0659DCEB059C37288F9F548BBC29E76D368FA8BFD83982FA3B3ADE5FF8C794CEDFC867F50F7003
2F7228E63C2EE037FE816FA117F4F3F050DD867E9DD48D4E63CE0F027E5B6E828FE205BF62CE2D023DC7EC1B03BDC0A39DBB
469792EBA3377C141CE90B7D22073E6CD435A3ACDBF2507D2E5ED1FEA7A43EF2F35EC8A8CB9233368FA8CF6EEE37F00DF73D
288A17F9838FD217FA445EC58DFC3EB9BCAE4E9EC49C7F05F8E237E83D661F1FE0481F6D5FAEBE83CE4ABD17BD56E652F864
3AD2DC424FE6E3BB7F4C29F5E023ACCBF254F38C3EC49C4F04FE12B34F9EAAE720417FE80BFDC0DFD16BD204635EBAFDA9E6
3FBAC677D077CCF951800B7E837E582FEB47D731FBFE00C7987381E0F7D4CB3A62CE1B029D3F9DDDFC829F4FF0D5E65DC50D
9F88394F09EAE33E1B73451671E7CF8A237978F0F9C4F2C17E5FE70ED525F78F395778FD3A96A79AAFF849D122A34B059F5B
0970419776BEA23CA65E1ECFF32D9775BDAC1F1E95E93DF00DDF9F1EEDEC7F4ACED5F01D7CC1F6F39A2BE8ECCE83312F85EF
36A768AEA00BE63AF4B1EDEF61BBC5BF6C2E505F47A7F81A7A7D3AAD6B46DF67D1213E91FF8389E5E57545784F9E7C18D92A
B93B17D43CB47974E1B853FE5C4ABF6CDF8C0EF38EB922BCB4FC571F8297E016DBDFCB1D8B2CE2FB13EBB5FD99FA323E94E0
BD4133A50EF2F971C4F7C9FD3CA07F478E4D4A1F9CF7E91CA1FCB2F321D515FAC3770D37C511DFB2B955FB0D7EE0890EF0D1
8F427BD414FE933FE4007552F7AD3AE34E05EF8B0C7880EF832738E223F017FEC027FC264FAE49E9DDE7B1B5DFC657C519DC
D147897787EDF6FD636BF9C1C5DC393679A8BE6E735ED309ADFD79EFFEAB5A6F057FC721C83DD689AFC74DDC28A7E083EFFE
B66DD232DFC7E90F3AA41EEA337EAA5FC353F39DC363AEF87AA67EDB47AA8FDB5CA88FE77EF0047EC0EF077F0D9A29E73AE8
191F31BEFED72BD4E7EDF549135A4BDF2CB735279E46B44AEEFB263C435FF00F9DE1C7394B4C5AE6E73CFC47E7E8803ED377
F8487FD294FE71A3BFEFD9B42C6A863BEF555CF0E37DBBC71595EB9B9E55C727B74DFAC3CF0FFA86FE12C58D2E2575935736
47A8DFA1839D79464C10BFB23C543F83BF76AEA07A7832B5474DDFBFD021F5DBFE52F3207EBEBA0DA53F093EACB343E6567C
D4F24BF1A71FF0BDD486C1C5DCFB2C552757B6FCB8D1ED9BC879CD4572ECB31293FE70F3B9F6095CD3E79DB82578BF5B90CF
F884CD89CAA7337B269DF6F30C5CF1DD0F1FB6FDDECDFD9A1B1B8A443C72FB695D277967FAD29C4067E80B5DC1EFDF0F4FDC
22B859CE6BBE92D3590B0575B14ED69129EFA40EC1EB0AC1E3E823FDC3AFCF1F9ED4C19F1B581FF553F7FA377A8506DFF3A3
BEA3F8D30F3BBF533DE037F0D3F6FDEA8FF86BDC85BD5DBDE421BA4147E8035F375E29CFD0CDE785C7157573A0FA0B7C44D7
E81CBE94AD1E35C3F7D7C77F34CAE9E710FCFBF3717429992F6C1E511FE73EDC37B63F153D3262827B7FB4E205CF13BE5DEB
2DB91E7C47A7E8131ED8F5F43A5CD77253EBA4EE837DC6570A3E37A039AF7E809FD00F7C655BA9E8C2D26778406EA00FF442
1FF089C45F0E5B229F472087F1137C0D5F01D7F2BB221EF9BE0B6FC09F7ED01FF48B9E4D5F9A1FE88CFAA99BFC40A7E8137C
C9B18F1E3471FB3272091F0017FCE0C1D33A3BFC5C04777486EEC0E969FCB6DFFB39015FF10FF8039F764F189556FA8C1EF0
57D68DDED1FFE33ABDDB0A6EF8293E65FB64F56BF4840FB01EFC80DC675DAC131F40FFE0887FF07CD64D5DF806CFB73CEB3B
BE92E08D3EE9377C801FF016BEE2177FFC3C6C89E4007A236FB6B41CB05FF0A6DFE4A1E535F86B3FD05F89ADD185DDE7C394
7FE00FEEF8233E143F5FDFB6EEF517CDED23B5277C10BCAF2CF06574463DB6EFD49C63DDE0804FD8FCBA7AC425D1117304B9
F1D1820E33DDF987E6303E67F3ACF2809C8C5FAFDE72D9CF9127E63BBA5ED68F7ED1273A4617E410FAD8B826B2923B9726B7
35C75847E116A3D2BAF98D394E7FF2FFD4C1756CFFA7B962FEBE68C207C2177893A858F50EA2978F937DB351CEDDF14B7403
1FE0C79DAF465C927A6C1E610E49D0A481EFE3363F68FF4B8F1FB0DF7D7FA7F2727DC65EABDD792FB9A77947DFE9377CDDD1
7EE8133F57F1D1DF4A4FE81EBC9F35F071F8429F93571A7542702727F765191D228FB39C543F7878A86F5BF7FE30D5393E62
39A5EB629DF0815C8117E00ACEF0117E8003B8D8BE577D163EA20BFC057D8003EBC7972D5FAB4456F275FFE4DB0E335FDBA7
29BFD105BC441F96435B2674171EC11BFC3B6EE2663FB9F330CD879BEB469D90F559FE6BAED8FCFD77D457EEBC46F30AFCC1
1D5F38596CC2443F3FE8277D4BF5D69881FE9C633977A1D76AF9F965821E7982D7CBEA2DF7F383FE197F95B7F0187DDB1CAA
3C8057F0095EE27386A7E28BFE8A8D18FAC4BD7F59FB8DDF530FF5C103FA45FFCCBFCE8E0E913AD1FFAF2183D6883F590EAB
1FE35B71E7773EE7FC5473C9F272C78489FE7C71ED9F31034577F489FCC347A89FF5707FEAE171674A4C68EAE7167DB3B9E7
CCD847C257721E1EC12BF4C1BA6C3FAAFE7EEFD6371BDDEB2D9A9FE088AEE12BFAB6BCDC16D4C5FAFF585DB5907F1DAE0BEE
F0073E5DC93AF691FB1E26C585DCB57965E7F8847ECE1F383BA686F0149FC257CE971A7FCDCF69702237E2FDD9708B9C8B90
030FDE8DFA4ABE6F13DFC17F6D9ED83AFE5AF079F7609DBF171D9FD07D8F8ED6C77DD029B9825EE13B3A45B79B5B4485BA73
69C59DFCD97574B83B37B3DC56DF45A7E8139DA0777C0ADD3FFEA3D94F7E1EA01F9B5F5F0EF8D7FDBD13CD49EAB11CD43C83
BFF0191EA1B3923306AD71DFA3ACBAC237E007FE41BE25CA5CA3A3BC4F00FFA68FF40F9FB2FDB9FA323A7B5CB7F339DF4FF1
1170032FFC7AC3F88857EE5C57D743EEA1830259C6D470F39BEA9FFEA053746B39AB78822FBE057FE013F7B59CD5DC45E7E0
068EE81A3DD36F784BCEC15F7486EEE80FF91037B2E555E7539A6FF803FDC027583738E037F4B3CC84A8D0D5AB579B6E0F2F
1CF3DCE53DF382EAEFCF7706FC2B730CFAC1F7C10F5DA373FC32DEDDEE6FBBD771346FD12FBA0517FCDEFAAE3CC0178BE41E
1EE63EE7A7F8A257F40B4EAC8FC7F37C7C065CF01BEA45EFE8FF51A2865BFC9C831FF806CFE7FE3C0FBF609DE88239107D6C
1F3CC4FD1D199B3B347FD06B927C8323E57BBA98236CFE6819F1CAF757F48E2FA07BF8477FE9B7E5B8AE8B75520F75901FF8
41DE3A639EFBF3D7DE7A23863A5DAB9EF133F8036FD03FBA813FC6A773353AFA39447DC673F559F88EAFC7DB1772D77D2E45
FDF8D8C8B19583CF47EB7CA1FCBD7B6870A4F850EC391B5CC0137CC931EAA1BE4749BABF2DAF6BE147F8C2E38896577DFF85
C7E0046EE8EF4E68A58FDDDCA4F9C9FD6CDFAB3E886EB6D61950C8BD2F4FFB419E803FB8E38FE8A378D721757D1F82479F4C
8C3E28BC22E759FF979FF488F479C57D6D5FA8F9810FE56A37B6B29B63951FE88DFB5207BA4FD0B65722F7BA89E698E597E2
068EF01F1EA303F0FB62DE88A1F2BA033A27671276A95D57CE6FF1517C31DE176D72B87305CD197CC5E653ED2BD7E3FAAC67
5399C8778487AC9F9C391135B68AEFE3F0EE76B3E883727F7846DEC04F78894EB8BEE597E68AF15CF300BE27FEAA6A21C9A1
24EDAA0C95F7C1E22BE8099EC37BFA9C22DBB073FE5C072F4A2D1C50C87D3FBAEA617D8D5E1BDDFB6FB50E72EDC03B239B48
BD96E7AAFB9D97A257057F9729E01F7E062FE0093AE0FA5CF7E7DB3D222587F1EFECBDC656F1E791479F87DCF5FD17DE9CEE
33E6A59F1FF4E1C6E961E7043770279F1E8CEA9548CE0F2CCFD58FE02B38800BBCC2AFE117EB659DF8A6CD43E747CCF7F755
E80E1EA33F7C9775B14EFA507673E43BBE3F9133D9FA06EB062FF4855FA0337892746AA58F659DF0853ED23F7C043EC363F8
8A8EF05FF474AEC798067ECED1876B9946CC973EDB5CA0F90A1FD1053A01A7B48B47E5F2F7378FF6B6C9E1FB29BAC53F12E5
EDF386E04AEE592E2EEDB5D1F733723AF380310DFC7981DC88F767E3ECA26B72C6FCFAE5C82632CFE33BF0117EC26BEACD30
6A74B83F175C6C3B3ADC9F27C0E972AD51B982EFA3D5B9437F7F7F48EDBA7EFEC10F7CA868BAE855D23FF84A5DD4896FE167
C673E5FDB62E03F7487EC057F204DEC337F8C7F5B82F759087D449DDF80CFCC66FF075FA485F0F1D1CE9BE9F111FC457F0E9
F8F9DAB773E7349A8777475619EAE72575E003E81F1EB35EDB076B9E91BF5C87EBC2477C031F01177805CFD01F7A84177F1E
EEF386FB7C1CF38EFA32B994F0A3BA8DDD398EFAFA9E90A1BFCAFAD0213E8A9F8023BEB2254BBFBDEEDC58D7432E199F3527
E0353A459FE0C5FAD029BA454F897346EE913E92C3B6FF547F417FE8097D8107FDA5AFF80EF5DB7E54F30F5E971832708FFB
FB65CA4B74962FFFC87B5227FAC437C005FF200F13F7FBBAB0F017FF4EF65E85B372BDE479CABF27F7C1171F256A9CDDCF15
F84A5F6CDFA9BE8FAED133F8821FBA46E7F801FE004EF80675F3FB0DCF2356BBF36AC5851CE67AF80575A04F700347AE6339
AEB90B7EE0897F3C3CD4BE9D3B97579FC547F00FE3B3F2FBE8A25119FC9C40C7E8979C40C73FDFEB31C8BDFEAA7590FBE41E
F8D30FD38BE62EBAF9E3B7C83DCE875457E40A7C297DB6DF5E799F0CBAC5C70ACD18FAABEF67AC9BFFE771E440FC7A4DBBB9
F33FF61DEA4309DA779AEB5E4FD1BC856FE4ABF14E7593F4FBFE21F2F735C87FFC0F9CC17DC78EC171042F78864FE32BB96B
8FCAE0CF53F79FD46DECE70438DADCF041C46A3F0F6D7FA8BA45DFF491FEA15BF8049EE0BBFFCD61114ED7AA677CCB72F7E3
1E83E467A5CEDD1AF8F9409DD47D7CC2A8FFC9FAF06B78040EB61F551FBF33B37F88FBFC89F61BBF060FF081D7F0D4E662E5
CD1FBDBF2EECE707F5DD6853E6BCFBDE68CD879B47CABFE7E7038FE7F9FC9EC7F338FA9C6CF3C0EDD277F2101CE90B7D2207
12DC693042748A7FC3437C1C3EE2AF395A8CFA9F3F8FD017FA817F3D18D169AEFB1EBD587AFEB568FFD2C26BF0254FD047B1
6283E3F8BEF2705ED36EBE3FC2637481CED107EB840FF0E3D4909193FC1CA2FFD463FB4ECD09FA59F0DF61116E1E8BA55BEA
A46E7407BFD1DFAD3203B78B8FC01BF20F1EC23FF4C07DB82FF743F7091F76DDE15E7F210F35673E3E5DFD13793CBE45FF53
56187CD79F8BC8814FBB0478D057CB23C519DCC9AB64712A76141EE337F41F3EA0834D47FA3671E7CF8A2379820E6C1E575D
E143C613E50DFE97E0C316FBDC799EFAFAD95E23CEFBF9473FB9AFE597E6D0C1ADC3EECA7A6D2E50FFB8BE7EF05DF73973E6
11CD517481EFA30FF04BFD28FA377F2EDA5D78F07F521FBCC747C9D52C5123CEFB73CAFA01BDDABABFAFA0F5914BF0BDE4CE
FEA5DDDFE752DDF278CB2FCDB30BCD871F96F592C7E071355EF46F6E1FA4EB21BFA937DDFA619D82BF4313E43F3F33FE38FC
B0F08BBAF113788CAF5CAA30AC93D4CD75B93FFA4257F0FB41D20623FC9C8067F096DC82BFF7E377DDE1CEBBD52FF05370B6
7DB0FA35FC64BE82A7E8133CC117BEC37FF80B1FD129BA853FF006BFB3BCCFD3B789EFAF968FBA5ED68F3EF2971A7657FA87
BFD8FE59FD2BB6BE3F5EDA63879B87353FD105BE8C3E2C777BF46AEBFB1F3E40FFF00374847EE8DF83A72DF6F93E8B3EA997
FA63FB53E1DD83FF13FCC0FDEADC122FFC7C23B7C899841F7D7F54702767EE65ABFE899F8BE0886ED0117DA33FE8173D1F19
3BFC5DC1033F45F7D44B9DD4BD3D62401BB91E7C2297C85BEAA13EFC097EC25778061FE047EC3CBFF5AC6247FF27FFCFFDF0
19FA84DFA077744EFFA8DFF69D9A67F8377DA4AFE81BBD832BFE9AF0EDD64983F72704B90AAEF809FE72AF468F1DEEFD58AA
4FFC9DF5C3437889FFB01E7C083FC1475827BEC4E3A897FBE33F3C7FEFD121B3A5EFE8133FA68E4A53BA85FABA20273F6E57
2BADF019BFC70FF2B41AFEAE3FDFA5185DE6BC9C6BA47AAB7452F119FC06BEC133FC02BD2699D57BAAE88EB960CBD9C87DEE
5C5DFB444EA20BF2067DE003F8027E400E7F92A17232E1313E8AFE8A470E68E3FE9EA4F292FB5A2E6B5EFED679D841E93B7E
8D6EC0131CF15FF0001FFCE0FE93EF8FFAFE0EFF36D68C882779CCFDC83D74489DD48D9EE0153CE33A969B9A8BF42575FD12
2FFCFEDC7FD53AA9BCAF16FF42B73F3FEA51C49D03EBF3C9D3BBDFF5767F8F0A1E93073C9EC7AD6BDCCDBDAF09BDE0BFE8C6
78A57906BF2C07151FF002BF2F720F99EDCF7DFCE4FFC10BFFCB193DECA03F7FC107F2CF78A1BCFBE440DF2AC243F2DCE68C
2C91FB7CDFC28F1376094D1FFC5DF52067E009FC80F7F844A262EDBE75AF076A8E9DFC7E58693F57E0D7CE5B03ABC84FF884
EFE24FF48BFEA123E3BDEAE9DEC85A69FD9C00DF0351439ECA75D1217EC0FACA2D8988277326BE404EE1D3D7FE299DD47D4E
4479093EE005DEE08F0E6E17E85BC57D2E4579433E591EC7EB51C4F721CB3BAD9FF5808FED47D5DFA98B3AC90FFCFBB39901
DEF003DDD8DCADFCA6FF2932F4DBE5CF75B72F564EE6E710F5711FFE9FC79DA938B4A9F81739449FC1DFF685EA23F8EBC755
BFFD52F0C2E7C097BED0A71B17FBED729F0F521CC921E3A7FA233CFD7347BB6FDDB9B3EA1C1F8187F00F1E9327597F19DAD4
9F2FD059D16403ABF8FEB42D7EBF4BEE1C5BFB817F835FAA5903FEE7EFB7D01D7C457FF0A540AF214FFD798AF5C393D8BCA9
52A14B559F3FE80EDF417FF78784A6F7FD175D9DFF34BA98E046BED237FC85F5C6F6197805DFC9BD14D7CB76115EE137AC83
75715D70B3BCD31CBAF6DD00F77D20F4955C27C792F6FCEAA19C87E3A3CC01994E451793BED13FF04E5B6ED0FBFEFE0CBEDB
5CAFFAC49F3E5ED7B1BF3C8E5CBDF8CFE0DDD217E60970E27A19DE1ABCDB9F472E6F1CF4BEFB7CB2D6CDF30E5719E25E9FC3
1FF121FC9EBED0277C0EDEC2E3CDD5FA84B87DB4E2457EE2C7893287B9D7CBC81FF4880EE13DD7E1F9E4A8E951F3095DE20F
F0159FD8D371E049B90F7AC6BFE103B90B2FD013FA8247EB57F5DCE3DECFAC7590C3E8ACC4C37E977CDF5D57B45B52F77AAD
E630F9CB752CC7F5F7F4D9F607EA23F08A5C845FBFA729DCCCCF814B338ABAF79DDEABD4B1BFEC17F111FCF4DEAA6FBFF4F3
86BED277DB2FAA0F829BED5FD5DFD153DE35439EB8BF6FA0BE804ED1277D459FE0035ED44BFDE404BE44FFF02774838EE80F
F9C3BA58277D28B3BC4F887C9E0B9F0257DB1FA89ED1177E145B6749C2BB8C977AC873D60D0EF8077A4487E082DE590FBA67
DDE4D58D9465BBB8BFF79D7EC834DFF7D1253A46D7E081BEC8157466F9BDB35B52DF8F6D9EA8DA738FEFDF77FA7FF5D0CF51
F00527DBC7A88F502FF593B3F80FEBA4DE3FCF857597735CFC1D1FC1770B0D1BE8BECF0FDE701FDB6F29BEDC873AA99B7CFE
B85DFDA12E67356FE90FD7E77EF896E945F5B3A355D41CC111DE939FF0185EC133E6F8B4098B9791F91F5FB93BAECB78F7BD
E4AA43FC1D1F050770E1BEB63FD61CC22F3F4ED6FC883B4FD5FCC49F72FD3E649A3FDFED4F3EE823793C7AC6FFA89BBEC303
F806CFF005722C69D26F9E48DFCD7755079FD4ED1EF5E0C1039B17C8D5146F574824EF0FC2B7B81F3C846FC7FF1DFCDCF77D
F86EFB46D521FA03E7742145F3C8FB44C1DBFAAE7940FFB776EF1BE23E37A638924BE0046EE818DE151B1B3547CE7BE93778
DAFE52FD91EB735D72EBDEC8FA437D9FA5DFF0CBF8A0FCB83DBF7B947B1F92F2833CC0E772BC19E044BFE1317552F7C659BD
DF759F17D63AC82BEB8BE61CFD59D3B8CBDBA26BF26ACDFA2E55E55CCAFAAE7943FFEFDD6A7E44FA886ED10FF7B57CD49C83
17C9BBF5AA203C219FE1333E05AFE943C19B833E92BAD0EDCFF7BB8FF1F3D2F255FF9FFBF17B7200FCE9C7A90B83FE759F17
D47C8277DC87FB822FFC8497E884FCC2BF2FDF2FEEBE370FBE941AD83744DEB789DEF0A18F4F7F77D9F15B7DDAF687CA6FF4
84BE137FD52AB57B1D527389FB705FFCF7E6A05E15DCFB849437E41CF9F069C641FFFA734DFA2B859BC9FBF9324EFF62B0E8
083DC12FF8863E76DDE99F44EA86F7F8E2C1E3031FF87E81DFDCB9F3CD133F3FE8073E077FE013BC48B92BE2B4AC9FFCB73C
D2EB70DDB3AB06FEE2FE5E83E62BFD040FF2E1C6AB0A89DC79A5EACCF4A57AE3F73C9EFBD9DCF05DEF777DBF446FF4977E5B
8E26EC3E46D685CF819BED8FD5DFAF178938EDDE37A37D225FC9C92C5507FEE2CF3BF42B73E28287FDBE81539AA29179FD79
E9428F01594487E42538C17FDB5F2ABFE9A3ED83D5EFC8EF8C030664111E52DF959D9179DDE701B46E729AFBA7FB2DAAAA3F
E7901B499ED54C2975E1D39772465515FDF27CEA253F53BF2AF9A3F40BDFB27DA9EA045FF96375ABD4EE7517F5237C16BEC3
7F786D7AD49C4097F03D7F8E810FDCDF77563FA2DFF0017EE0134592F64FE2FB31B997FC872FAFBBBFB7AD3E0AEFD029BA45
EFF820BADFB63B32AEFBBB4FCA1B72EC5EB6EF2EFBB982AEE05F95295DDEF6E73ED6C1BAF039D6C3F3E031BF675DAC131F4A
B2A56D5DE92FF94C9DD4472EA253DBE729FF8E4C1918E1FB29FE71EE8F8287FDFDCC851FBE18ECCEE1950FE43DBCC037E001
FE816ED011FDB75C268F15477C08FEC027784BFEC15FFA66FB09D533BE9EB85FA392EE7C51F38D75B02E7C7CF3C6DE0FFD39
85DCDEDBB77F3EC1051DE2BBF08AFC865FE8DDF6618A3BD7B77940731ADF002FFC037DA377F05BF77ED7367EEEDADCA0FF4F
0EF37B9B9F0A47C6F5FD129CF00D7C64FDE11E9BDCFBBAB53EDB27EBFF5BBEEB7DEE966E5BD79DB3ABCEF1777865FB5AF535
749FA7F1C0087FEEC37F6CDFA638DE8D5333A59F5BF0067F03AFD83EC7BFC18FF5D93E06BFD3FE935757DF2EF9A3FCA4EFF4
9BFCB03CD57AA88FE7711DF20C5D7E12DEA1877B3F95CE17E899BED027748DEFA0EFDFC2076CF07302BD821F78E22B37A77F
79DDCF51EAA62EFE9FC7317764F8AEF070E927FAB6F9A95CEF87BEEFDADCA47E8DCEFEE8DDA8A49CBFE383E8D9CE8792D6F9
C1FD3D66F57BCBE9175DDBF87E89EF7DD1A77F3EDFFF581FFF1FDB3F2BE4EDB1C9F751FCC3F8A97CBD3DAE430FF77DF1B1E6
FE9D0522DDFB4DE9377E8F6FE51C17F4011E831B7300F89DEE9BAF8EAC031FC577E83778822FB99A3A4399ECE267F80ABE08
9EE08B9F2569D764BA3B8FD4BC22DF92FD5BE5BAF01CBF843FB6BF547DC2C7E4DF74EEE9F6059AE7077244FDEDCF95F8DFC9
22FDE7F9F9011FC1CDF2517388BE561DD7E91D77EEA1FD65FDF48BFED93E53798F6EC01BFCE13BCFE73EF41DBE9113F0EEE6
8ACE3DDD79AEF69BFC30FC35CFE803BEFBD9AEFEF3FCF965EB5F1183DCFBD0749DE4067DC6572ECE2A3CDCBDAEA478DAFE55
FD053E1A0F9597F0BDE881C8A6BE4FDCB953E707DF07E1077D4BD9A45BCD870F1F5A0E33A7664D56E02DD1013C8527E4167C
39F3765459F7F71234B7E833EBB57CD43C637DAC171F5CFD5FE7103F9F5685777A47EE73776493E9724E8D2FA04FF4CC7A59
3F3CC75FE03BFD2F703CEA6F7F9EB2BC507CC00BDEDA3E52F5409DE427F56E3CDD73ACE434EB240FAF4FEE56D3BDDF5E7943
8E9233595F4595F5E71DEAA24ED6C9F52D67355FE957B63EF9EAF87DFBF974B7417EEE5ADEEBFFF37C7E6FFB58D50FFA8617
A9E7F778E9CF69F01A9EA3B3739B2327B8EFBDD23C861796779A27572F9671AF4BE3EB497BD65B2FFC261FCC9FA636FBCBBD
CF51F3EDD69B55AEFBF96479AAD7E3FAE494E59A3E9EE793F799CB444E90EBD28F4337A2DAF87E879F5DADDBE3A57BDF8FF6
8FF9C0E6A4F72206F93E8DEFC25B788C5ECD4FD45FC0355D82DEADFD7DE1C5267D66B9F79DE8FC01AEBB3746B690FAF005F2
C7E6C1C97D66F9FBB94B0F7AB7769F53D2BAB91EB99DBE7AD133CEC7D4CFF0C7AFDEE81CE2FB317AB7FDA8EAC4F23E5BCFB1
BE4FA34FF0012FCBDD6CDD06095FB99FED0F957FF813EBB57DB9E607B9898F9EB955C07D4F037A425FF0947C4B76F19B3882
237E4FDF6D3FADFE4BAEA6EA5FEEB9F41FFFA36FA66BD5F99DD0667FC9B9323E855F1BCF3557E03BBAC99722AA8D3F97D117
E603FA433FBE59D2C19D7BD0177467FB36ED1378800FFE8BBEF0477476A77FBDF57E4EC063DBA7693FF00374838EC09B75B0
2EFC0C1EC0137883CF152E17D9C2F775749FECBDF0326E1ED69CDF5E3B62AFE4283C2057E9033A45B747E7F61BE9FB38FAA6
2EEE4B1DF012BCE12775F37B788C4ED127FDB0FC553CC117BD8323BA47AFE8175C8DB79ABBF0177E327FC053FC061EC24B70
B33944731A9CCC1FD42F6CBFA4F55117BFB77D98AE1FFFFE24C37769DC39A8E630BC23077E5F56F48CFCFCAD41CEB5EEF36D
AAE3E39772BF293F579EEED844EA206F9907F6D589FC44AE8F3EF13FE694CC570BBE94F3797473EBAFF032EEBC3ED61CBFE5
ED5E894587AC9FDCA62EEA2457F08FDCF5031EA1DB5B19BE89E3E710FD866FF00C7FE1BE362768AEAFDDD225CCCF69D64D5E
717DEEC7E3C9739E876E8A2F8AD82B3A406FEB0F77AFECCEF9F57EE43F7A4ABEB84D2737CF6BFEB36EF2F35ACF72CFDDF76E
E9756C7ED0FBD337FA88EFF03CAEC37A8E87457EEFE707BA44EFF811BAA7BFCC37F4193EC0037282FED02F7C089D830F78A1
F3AFB3756CE2FB2B7CC6AFE135B991ECDF1AD5A47E7CFDF6C5EFD2B87DA4FA263A447FF0105E92B7E9224A3416FDA263FA56
FA55AFC4BE5FDEACD5A6933B07571E9343F86E8E9F22BFF7E732D6FDF9C2C84FFCB9929FFC3FB8E06B954B7709F3FD8FDC4B
95B9E263E73FEA433179DFBDB2AC87E7C3A39403DA0F76FB0BCDF35369FA5CF1F389FEE38BE0068EF847B2383F0C72E7979A
57CCBB9FEDCD335DEE072FE0157C4287F0111F3A97BAE04BF7FEB7FD11ABE5F7F0863CA12FE42EFD8137F4853E1DF8B76F3B
F77950D539FE84BEE0313AA33FF40B9E1AAFD4DFE117EBCD912EF79B820BEBBEDEA3FD60771EAD3C20FFC88D4FAFF4B9E2CF
41AC97759207E0013EF8203C4017C607ED67EA919D3A497F6D1F094F3487E00B7C37DE2A8FCFEC8B48E63E4FA279092FC017
5EC5E61975C1B75B6FD6A8E6FB2F7CA22FF403BF583EBA7D1FA98F9C5951B34332390FBEDAAE5327F7FE3FC5915C24F7B27E
1E91CC9F776E3DFB6190F8526C9DFF5AB55709E91BF89233E0947657D77CFE1C852E6CDFA93CB63C523FBE1459C2BD7F33E7
9C2037733FCB7E4CFA012FC807700097F38D7A7DE83E3FA3390DDE366F158C58EDFB22BAC3EFD01F7340A669BD3EF4E714F8
5BE0CDBEED245FF09DCB45BAE673EFB7531C991B6C9FA97AC027C02343FAEE9BFCB9083CC1911C050FF0216F2EFEDE7D937B
9F9CDE8F75A3037481FE564775AAE7E72EBC20D7F08F935FE4992E3FF14D700347E68E4C970B1595FD1DFE71ED5CC5C77EFE
D1279B8F1537AECB3A783CCF276792FF507BA9EC53C825EA672E601DF865F23C8DEFB9CFE96AAEA23FF30DF5918DB7BAAF95
7AC09779E1E7595D93F879CF7D787CCC5C11CC07F845B536EDFBC8CFEA57DA95095EAF0D7C1DBDA01FAECFF3B9CFE1A57DBA
F83E8EBFE21BB61F453FDA1FFA85AFD117E606FA63FB4CF8A77E66FB19D53B7CD95333A2A2D4818F901BE8A6E4AA5E25FC7C
20CFD3DE2C9D4D74802F929FA9FE57E513F777D7D4EFE983F98FFA91E571AF4EF57C3F3E3CF1B371FEDC7C344EF663C1FBB6
1BDF93D7DDF023FC171F30BF551CE1157903BF6CBFA7B8E05BF4C5F24E7DD6E69564DDD70ABEF835B859CE6A0EE11FF80978
A32BF0045F9B33BEEB9AC4CF01749CB74650177E047EF80F7E842EF035F471737AEDA5327F901FF0153FC147C011BF0297D8
BEC5BF591FF7375FD2EBA2CB14A39B9573EFC7D27901FE307FC023784D3E9C4F5BC8FD9D31F881EEE009F3D1A7D3F30D13DC
D00138811B3E74AC41C4677E0EA13FE3ADE63AFC8577F0101D5A2E2B9EE08BCF155A1251D1F751FC00DEC263D6C77AF17DF8
01DFE01FEBE0FADC6FFBF25E09DCB981F28F7CC64753BC5D2F99E41DB9076ECC5BE007FFC01B1E5207BF87AF37DA342B27E7
F5E810FFA51EEA20E7AEFD5F954FFC9C0047D64D3E5C4E5E3A9BDB2F7FDAA189F088DC62BEC1B772CD89F8CC9F372D8FF4FA
DC8FEB717D1EB7EF76EF816E8E529DE37FFBE76675DF67810FC103F8036FF00B789D2A7BCB816E5ED4FCA77E728C75A00BF4
893E4EBCD36B8CFB1C88FA3E3C650ECAB52467A8DFB72DA9BBEF93DF832FF9463E674A52648BAFB3B5353BBDF0738EBAB88E
E5A3E61FB882337A4207E8F5D40FF986B9FDBB5E9FDCE33E368755EB95C0F73FF22155E6AFDD3A6DDFA9BCC5AFE12F38301F
82C7FAF95D5BB9D711B46EF296FED257FC9AC75B2E6BBD375ED54B269F9BC11FD12D7A85E7F0FEDA899603DD79B7F29E9CC3
EFB3BFEC3546FC18BE929369A3CBA5175CF01BCBBF531D9AF83E677DD45CA79FE81F3CC117DEA5FEABF592478F1E59FED31F
FA85EF9C9AD663B97BBF9AE625BCCBDB2CE07DFE97D9F6FA78D3FFCF3FE93D509E871FE17FE0068EF84D8AEB0D87BABFCBA1
B96473A15E87F5C07F788FDE6CFEB8DA7D9FEF8BD459B07EA057EAB51C5DD2E985EF67E814BDA1D7ABEFB55EE2CEBBB57FE4
2BFDA871BCED4A7FEE23F73E6DD463B93F4F199F353FE0B5CD0D75BBB6F27D94FEA47DDECEF58B79C478A2B9085FE011BC42
97F813BA3C5633A7FB7CD3AE4C3D8749BFE03DF970364BB7C1EE7D789AF7F4FF60FF5EED7D5FC317E13B3C87F7F008FF3E7F
B788FB9E98CB1FB45BE2EF3B990F981BB29CED36D8DF6FC163FA4BBFC123439A8E7DFD39EDDAB9AFBFF7F303FE9D1FDB25BD
7BDF9ECE1FAC8FB920DB82023B9C8E544F5C2F53AB2EE9FDFDD9C52B1DFBBAF7C968DD5C0F3C2D6735E75837389037F00F1E
C143CB41AD9FF52C5BD16E909F734B73B45D29E7FDE8147DA2931B291B0E9573FAD83E45FFE103FAC71FE03D3E613A559F45
AFE420FC3AF84EB6BDEEF32CDA27F80CBFE18FE5A9FAAFE5AFE6EEE5CEE5D2FB3FF97F70E1BAE0875EE823F30AFDB4FD8CEA
117D8027390CAEF0D7F4AEFA277F52F5AFBED57D2E4BF3079DE5EF19AC0F1FC4CF52BD15F2B73B6FD63CFF35BC7B7F37FF28
3FC8677CAEC8F99EC37CFFC62FE9237DE53EDC17BF44BFE6DBCA0F7ECFF5B9DFAAEB1D1EFA790F1ED469F383CE01B66F547E
E04FCC1119A38A47385F527F42C7E89A3E1C49D7F323B93E3E8EAF500FF30275D107F30DF51178677398F29BFC4C1B5EB1BC
AC1BDF456FD5BF6937C8F7EF8D7BBBB6749F8FD075322F304FE51C9DBB9CCFAF9FCB762EE2CF07D48BCFE02FF483EBDB1CA2
7386ED331547FC0C5DE2A3A7BF2DE0BEB798FBF27CEE6FFB2EF52BF00627F32BF52FDBBFE9FD781EBF379FD4F55CFB27E46F
394FC047F07DD385E619FAD833B7C74A5FE7E40C3C615E812FE0CCDC06DEF4A976CDD6AFA44FF4CBE6C571DDFBCBF3F19DDD
5199F6084EF813FC003FF0C497D129BE865EAFF5ACBE55EE437EA043FC23CFA59E1FF9F32F7C86C7F81D3E5535658787BEDF
E307A9EB37FEC17D8E80B903DE6AAEC35FF00447F2063D8303B83097E5DBF159551F579B93BEE8DAD2F7697286F5B27E9B4B
36752EE2FBEEB175DDEB099FC91B748C3FC11378034EE0866F2D49D96687D44D9E2D5AD2FA95E89BFA99DB5807BE9B3A43ED
BC32EF904B365F689E5CE8553CC2BD5F40718747E08F4FA0B3DFDAE42E177CCEB3F10FF23EB7D8733CD7E3FA9653E32A96F7
F3091C6D7EADDF63A5EFBBAC97FC61DDE0C1E3793EF348D61B5FB8B907DD583EB28FD53A98770BE5C97CC7C70FDFCD55A97B
3DE12F3CE6F1857B05BAE179F0099EC09BEDFBBBDDF2F9478EC1DF34B79BAE7273AACE11F0863EC01FF88FDED001F7E379F0
8FFED02FF40C0FE11F3A3C51A2EB776E2ED2FC83A7E069F381E2CABA2C9735FF56D46EFF969F8BF46D7FCEEE51BECEF13972
266D74E5E1F27C7C11FED5BCDE6687EF67F8137C3C54EC33F7B957788E9FC277723563F6521D9D3FA84F5CF9A4E92A39D786
AFE425F593ABAC035DC273786FFB07E5097ABD7AB1765EF9BC013E88CEC999ECDBBA7E27FD8467F43FDDC21677DCFE45E705
D6FFED9FAD7E90F583C3E6FFBAC677E7158A3B39BC6673C7F47EBEB28E536B3AA776EFFFD17CA5CF5CC7725C73975C00EF3D
4732DF7139A17DB49C529FB539AC60B75BBEBF520F394C5DF825F880D7A53A2DEEB8736DED07798C6FA57E5577B1E045FEAC
AFDA39A37B7D41EBB6F308FD7FCB7DBDBFF14173145E30A7E4989BB7BDCF77CB41F5CB3329BE28E1BE875CF3FFD32A9D534B
3FE92B7C849FE8907E6698D8EA9FC78F1FDBFC820FA067F30374A179863ED0A9E588EA191F047FFA818F7CB3A8FD5BBEAF9F
4DD0619AAC9B79055EC0DF82BF057DC66717B66C9D5CE67DF26641A2563F489F2E366BF58F3B37D73E310F31276579D0619A
3F2FB1EE4C2BDBDEF4E737EEC37D59CFF9AFDBDE74E7E07A5DEAA50FF485BCB17D8BEA0B5DA257F48EFEE93BF313FD87BF65
DEE81ADFF7F5CB9D2B0F97F7CD935BF0187CD12B38831B39067EF87595321DD3FB7E7F704DB777A4CFF8327E89FED0233CB3
3CD57AA8CFF25773EFC28952EE73FCBB3EEFFAD4D72DB964F3CDAACE19FD3CC197E025FEC4F562E73A3E6339A27CBDFA76DD
C5F279097C103FDDFE7906F7F94FFC86759ADED5C7D13D7359DE753932F97DB3FD92EA91FED367FC047F4127F4977E933F69
6F567D243A204F9823B2D62BFC97EF1FE069F9AEB9080FEA8C6F9DDCF76FEA673E631DF859DA840DB3B8EF93D09CC73FF0A1
E3F5F3B6779F7B50DEC03BF8B37441DB5C7EEEC23B78059FF0717C89BEE24FF883E585F6D5724DFB485FD159FE2ADDDE91F7
05E023B6BF521FA03FF4C17C437D847EC04378496E678C5B7682F343F545D64B6EB36E9E0F4EF0D9F68DBA4E7C88B9EF8B0F
B236F2FB83FF984F2A2EB60FD4E79BCFC1139D93E00BEB331FD3EBD24FFA8BAF510FF32575D97E40718EED9745F6757DEAE7
C0AF5F7669E0CE4F9497CC0B470A7629E4CE85353FF09B92E9D387B9CFD3944C3BCB7DEFD8BE409FDC87EB733FE6E062A333
16F0EBE5FF793EBF377D699EA1B3CBF71B6691F7EFE00BF8F5AAA4EDD3FB730C7D665DAC8739065F4437872BE5C8243FB90E
F308D7433F35BE6D9BCBF76FE39BCE07F0EE72F2AA8F645EC1F7D103FAC5A762EB385D48C824770EA439BFB175A755EE7DA5
5A37F3CBBAE51D9EFBF30AF5CE9B113E5474429E90BF5CC7E61F9D7BB80EFFE67AF08B7C3833AFF05FCE77E77759EAFB02F9
808FE6D9DFA5903FE78217FE066E968F8A0F78314F65FF23FF0DA70BD5073CB4BC539F2207E8E7DEE7591B05EFBB0DFC069E
C01BD6C1BAC8279B477FEED240FA88BFC00B78820E37F549F3CCFD1D30CDA5CDBF063A38B6A893C3913C4337E443BA88EA2D
5C4E695EF17CD6B1A34DC602AEAEC76527F879095E3607DC69FF9AAF5BCEEAE379FEA5192193DCBE457542CED147CB77CD27
DB9F681FF10FFA572F247CA8FCAC3FA3E525E91BFD621EA26FE81F3CC1175EA7BF12FAB7F821F382CD4F633AADF27D7DD1A5
D65DC4FFC82778CDEF2B55EBF0DCF76F722357ED4E19858FF012BCE023B8A13B7488FEC8EDACEF14BBEBEB1B3C0AD50DF0C3
EFE01D3C8437E6BBFA7870615DE412EB3BF151FB03EE7C5973141E916319B357709FAFC7177F4F13FAB7BC6F9FBE9297CC3B
794EE5AAE4F7077E9A9FABAEE0273E0B4FD10DFC814FDBA774AEE6F386DCA09F19C7352D2BB9C9BC7029B27A0BE10D3E883E
E1077D8707E80C1F3A9138FF0D5F1FE805DCD11F7AE271E46DF627ED0FF8731FBC8117F007DC98B7C0CFF609CA2BF4B13F79
E724BE1FE18BC61FCD757874AA4EDB82EEFC5AE703FACCFAA90B7E83B3E5B2E6D985F0A665C5E7E93773C3EC2D2D93CBBE81
1C981BD2F292F06379B5B6C5E4BEE422BC430FB5D3B5EEE2FB347DCBDCA9454BE963EC39E8D3856D0BFAFB36FCDBF6E7DA6F
FA461FF15FFC2F7DF55AF7DCF72B9387F05C7314BE9FE9D7AA993B37D77908FCF013748BAF9C1BDAA2A53F4FF13CD691B577
AB66FE7AC087FC072770C1F7C087B9EFF30F3FCDE5F305FD9AFFABDEE089E5A6FA6FD9BE418E54389FBA863F4F6D9EDAF1B6
9B5BB4DFE43FF765CEE3FEE8121FC017A8BB61E996AEEEEF5AB6E826BEBDFA41BB6BFEDCC0BA6D2E6CDCB99AEFDFD4637388
CE13B68F527DA10FCB59CDA5332F8BDD0DFE5DE192BCFF9C3C8337CC7FACBF74AA748BE571CC89457B659AE9AF9BFA995F58
07BA2C78B37312F169FC1A3FC6578E7C9AAB923F7772BFD8F329FFE6FEBFACEBD0CEBDDEA4EB676E8A595FB05EF2D7729CFD
B6AE9FEBC49EC3C809F0A71FF82AFCC65FE7BF0A3F2BFB047211DDDBDCAF38705DFCAEDAF2B6C5FC3C64CECADEF4F39F7CDF
FA7D59AD7B6ECE56BF260FF031CB29E53DBCC397F029F44EAEA0FB83FD3BE675739DE6103E6AFB40F519F845EE658CAA5C51
728ABC020772153CE8337DC78F7E7E996AA59F1BEB3305FAC3F7C017FF032FFC00DC6C5E0CED78DBCF31C353E7217025AFE0
E3BEA79FE6727EA8FE97E1BBBA1DDC3A751EC1F72C2F147770C5DFF0BB9DCF3A3EF2FD857C2657BE4AD0EE9A9F7BF81F7E08
DECC3B596E97FAC8F763FC135C781EF7C767A9DB725FF9039FC0A962A50EEDFCDCA31FE4F7966BE916FBFF06BF9D5199DCDF
C74327E04F3FE095E5B2FA2FFDABFB76F859DFD7F1ADFC3D3BE6F5E75FF0B43947F37875ED14937D1EC19F25AD5B3BDCC857
F8C8FDF061EE6B7CD69C83D7F409DED1AF8BB3EA7690CFBDA0437C1F1EC12B784F1F993FE827F3609E2C79D6F93CB5BC56FF
3E31E9F3E07D99D7DABFE3E646CD33F46AFB04C50DBF010772153CE8137D8707177A55AE28730379006FD0057E813E783E38
C11F7493717A8347EE1C52E702FCAF689CA09FF086BEC13BF866F8E8FC074EE0013EF83EEBE6FADC6FEB990EC1DF9154DE90
E773DF68595FEAC7AF99DBE0353CA7AF969B8A0F78C5D607CFFFF29D80975F2D0A785AE5A79479DDFB7BE3B7F9537C9ABCA7
3FF86BDE54EDDFF1E75CEAB73944E700740C0FE1A5E5A0E6DED94F4A7D24B9C74FCB455D0FFE9D2949D55DC21FF201DED71C
D3BAA3EFBBAC8BF53037B12EE602D677E187068FE47DDCE8847C625E2BF845F6CA7E9FE90B79427FE0BFE594FACB9E6F3A94
F1FD05DF3D96B7CD3877AEADB9886E361C6A3FDA9DE7288ECC11F03473E29046E2A7E4FCDA0BED2AFBF301EBE33A3C9F7982
753117B23E70C487C0939CCF9EA4F060DF27B82FD7E5FEF80A3E71E46C9E75FEDC463FCA1F4E534DFA0B8FE0153E81FE4A65
EDB0C7F76FF227D7E136E3FCB9125DA37FF4CD7DB91F731C7995F59DB275DC7E4CFD95F9B1C8882C4D7C9CC0A5FE7F2DEBFB
3E4D7FE92B79C0E3420A37DF29E7DB3CFEDC1F218DE45C179E91D3DC973992FB335F5277A91419E2B9CF9F360BBFEDCED335
EFE1117EF7F5C3367FFAF9C075783EF7415FE40A3A23BFF0BF3561819E8C0F3AAFC00BF898E597D00D4F9E3CB1F908DDA32F
F4860FD25FFA0D5FE00F7EBAB06DF85ED96F9197E8EFFCDDAABBE4F31BE401BEC3DC916362F86D7F3F0A8FCAE76B3FDACF19
F4880ED1F3D98AA11BE47D22F48939069E165AD1A18C3F479FECDBE25DF7BAA2CE37F487FB54CED8AEB29F4BE0946D70D853
C18DB9C9E672BD3E78F1FBCFFAB478D79FAF66FED9FCFFFCFDFD8CDDCD77CAF5C90774B97F6FF6CAEE7582AE614F053FD643
BDF0063F803FE885F900DD2C1B996C8DFB5E25E50B3CA15FE43C7DB37DA0EA17DD9357F491BEE2F7992E57AF268F23E7ED1C
4E758E2FC257726F43DE34EEFDBCDB93B7AFEEEB93BCDD37BB5D79DFAFF169CB71CDAB13770B0FF6F3135C77B5CFD224F8DE
9380D7E81A9DC36FF806FF98A772CFCB1FE6EB89BAB9CFD61B19E2C94FFA4BFED067EAA24E72147FB07C517EC33FFA0C0FCF
BC2C5B47BE3F88DCA51FF005FEA04FF00447E600CB71BD0ED7854F754685EFF5736659C2D61FCBFE903903BE8017F5A21F74
8CAEE111FE030FF021E69BEAED029E56EB9EBCA4FC047FE636FA007FE111FCA16E7C9DFABF4FD4FCFF64BE09ADD06CB8F83C
7E0FAF2CF73527E0057E82BF8003730A78703FEEC3F5CFA7AD5E4DFE8DDFE1B3F80FBE030F4C47EAEFE809FD7DDEB05D797F
4EC757F019F0064F7064CEC2C7C01D3FDB74BF5D7DF7B974E507730D799BB55E856DE26FE40F78E04BF8143A4377F0D2E6F5
9B011FD11B3E491DAC83EBE177E0C55CF3D9DFC516FA3EB7EA4A9BD29203CC55F4075EA05774CA7A59277318BE95F96AAD4B
EEF39A3AA78037F8E3F7D37635FB5970C0C7A7AE6F365C7C85DFB35ED64F9DCC65D43BEFDB9645659DE425FE41BDF80CF81E
B8D7E62D99E7C83F7C91F9B7C0E89C2B7DDDE0D3F0105EFE7CA66D1D791EEB674EE4FFF937F5E247D5EFB7FED8CF1F9B07D4
8F8FD6CB1FE6CF91F0AFF2EC5499C56FA893B9967AE92FF30A7D26BFD0FDF28101AFC18D1C033FEECBFD982B8DCF9A8BF09A
3EE2DFF49375E22BAC7767FDB6437DBF20AFE075D984EDEAFB79722E75AD4BF2F91FF44F0E8113B89113F865818FDBBCE5EB
7FE192A4FFB8CFBB281EE0C0FC5B784FB6EE3EAEF872D5346D4AFB79C29C0A2EE52EA55D249F9F8057F0097F6C54A4D9CFC2
97B0264D3B89AFE1A367E655D82673003904CFA887F998BAD017BE165B675993D50D77E7763A2F80DF9759DBD6F17386FAA9
3BF6BC5CB268C6BDFEFDF165D641FD3C9EE7F33CFA5D6F41CBA27E9E1CCE18DED3CD9F9AD3E80F3ED8DCA2B96B79ADF88017
BE004FE00DF7A35EEE0B2FC8ED93EF165B28EB8717F0047FE2F73C9E3A6C4E50BD1E6893D3EDC317BD08AF2BF943AEA3A7C9
D14DC7CBFC82AFFD34B96927F9DC0F3AC45FD7364C9539F8FC40E02BF01CDEC33FF4040FC993AC372ABBBAC93DE6803A3503
BED7EEF2491F370FE8FC95BBFEE7D5DC7987FAC4995B75C3A52E74457E9387792F84BBEF754057AC973C67DD36C7CF6D3BD4
F76F70E5FF791CFD839FF07576B6E687648EC6C7F115789DAD4F833BEEFC55E7209E0F4EF0FC68B9E6EBA53EE60378B73571
9BC33ECFC86DFA497FF1177283FEEF2E94ADBBEF1FF4933E823FB83347D06F726963BAB48B5EEB8BCE67F4877CCE9EA4E433
E7B7EABBE808FDE03BF00D3F8677F8037E814F5007EBFB7567C6E0F32C7D1BDC91EBC003E61570032FE680269D9B8E17DE36
1F1BD653FC0B9F631EC9BDB1F97AB91F7C5ABEBB551DA99F39031EED9ED0DAFDFD267C937CA0CF9F250EA9F1F4E9539BBFE0
0939065F6CDFA5BC45DFA617CD2D7483AE6ABF1F5ED7CFABDF2A856D973E303FC123EA671E611DE8151DA12B7C139CC0E74C
8ACAE5DDF731A8EFA3FF937F84D410BE833BF31C735ECE7561DBFD790F3CB2FF163A55F0E1712772864E155C783EEB58FF7B
9B387EFF985F26960C1B2BEBC7A77E6C15E6BECF065CF057F0214FD1C1A22101DFE12DF307FC5DB3AD7511B90FF310F8510F
75303FA12774835ED01F7A44B7CC49E8ED9B9F527C269F37E2BE765DBDFF9CE549BABA7DBCAE83FA99FB0A24CFD3C0D7A5CD
1B9A2747E77EEEEA6A78BAF921C1833E933FF0D9E61CCD3970E4793CDEE6FB3F029DE18FE41E3C8137F89BE5A9FA023E89CE
F14BF221DB82AAA7450FCC11E0065ECC7DAC8FF5928BF401FFC037F8899FD027748F0FA04BC359E757F0C6C7AB156E55C7CF
3D745CB845EB6EFE7E825C013770447FF818BE367F7B8B4DA26B721AFF3871B7E433F91E25F21EDC9BFD1A36563EC7127EAC
C92CC11F9F030F1ECFF3A98B3AA90F7DE1D3E86CFAB5667D449FE41CFE4D7E37A816F0B4FEF0248FE427BCA890BE4D1C3F3F
F12B7C0A7EE09FF0001FE579554AB42EE2E727732D3CA8B4377517F77791D46FF01FFAB97777786AD11539847FD33FFC833E
325F17CEFED959BFFFF485F982FE8CBBDE6495F40F3F1A97ABC92CF157FA894FE29BF82C7DA35FAC8BF530FFE25FF8CD8AB0
149F05EF33AE7A5AE6307C197FDF56B4D56377CEA67A668EC087C1173FC687F15FFA6438E87C0F1ECC65B94B156EE7FB28BE
CE75E903BA26CFD03778E0FFD4C11C91BD69D9C62E07350FE11B3C637E819FE42D3CC577F1C50337F334909FEC07C0B75CC2
F40DDDDFB3507DC25B78CC3EA344BFCC277D3CEC3AFAFCD8FB12FECDF3F06FFA491FF1B94FA7576FEFBE975FE720FA6873A0
CE37E801FEA20B785CB7788B4D7E1E86A76CB24A78D63665E3BBF2FA3E3EF143AA667D84A7F0917C402FE8849CB07940F906
FF968C6819EDF6F13A27A10FE389CE8BF0057FFDA270786A7F7F03EEAC0FFDB31EF228763EB10EEAA7CFF49D9C40A7F81C7A
9D3E257146F93D793BB743C02FFE4DDFF8FFD1BD1A3F91B9159D8FBADEF8AED4010ECC5FE0617CD639035EA367F26ADD17A9
DDDF87D97FB3793FF7FA87E625BE054EE0867F9363F07DF789CFCECACF29D14D6B0A2FF06B7C7DE39D562D7D3E3167598E28
5FE009F9506267ABC77E8E31C7A0C35ACB92854A3E9FFAA17A7BB74F533F22E7A897EB705DEAA10EE6AD399D9A9796FE937F
E87BE590F0628207F319789F9854B6B19B3F3417E12BF9F443E3A0CF211B134FF071672ED8783F7D43F9D926AAF113E17587
B68D6BC8F3D10D3C22578F6E2DDC4E7295E7D3876DBD339F949FE8F8B3BDB5DF77EFDFD17987FA99FB58073E803FE207F804
7C8127EB16859FF3F1627EB43941710017EAA57E1EC775783E753157E63F9ABFACD3AFEA98BED36FE632F2B660F2E6FDA46E
748E4FD56CDF32DAF7A911711AD7115EA09FE1A31AD710BEC05B7C1DFEC26B788ECE4D8FAA6F74097FC86D78B4A361CB7DB2
5E7C8D9C4477E80D1FC4E7E033FC3EF8306CB1FB3E3C9D63D02B7995FD8F0AF7FCDC6BDCB9694D791E3AC377C1059F031FF0
67EEA40F27BFA8FDBEAC135D3167C073F84D4EA2D77249039DE1B3E81CBEC373FC947AA913DEA11BF4429FE93BFE04DF73A4
AB3BCEBD1EA1F31DFEFE7597F0627E8E3147E4AE5F6C93EF97CC11F9E3872DF6E7F4F0A28912C8FDDBBD9BD0F1A0E5B64489
7DBF46378BAB270BF575894EE11373243CFEFA78CA8E729E074FE007793C295DD897522F7E6AFBD2CE1FBBEF21254F66940B
FA405F1B0C6D5EDACF31FA50A976F8393FF7781C7DE0F1E08DDFB09EC397426748DE32B7C1477C96BEC303720C9CC087BAF9
3DF51FBF54779CDBE72BCF98DBF0377C0A9F1B392DE8C7F812413F46FE1DFC1BBDA357740FFEE0CEFC0B5F72CE6918F7D9B3
6736B7328F17DA90B38FAF8F85459B5F123D327FA07BE6D1BCE94267F8E713A653CD63F47A7479480B99839947C1F5B7060D
E3BA730A5D3FBFA7AEDCD5425A489DDCCFE610CD8D83B9F397157EA18F62B35BEEF3F761E4535874C0A326FF7EDC567EA213
7C02BDE03FF80EFE859EDB3F6B5CC7CF317C0BDFC0474C773A9FA13FEAA24EF804FEEC07E8C389C415EEC9BC4BEEE29B3313
34BDE59F9BE2EBF48BB989BEA12F72009DA17B7C00BF69D728E057E79A094EFB7EC07C0F9F2A7D98369DCCE94D2F857D2938
E07FF8313A418FE893F9063FAD7731E9213977DFFCA46577F7F748D41798A7E8337D2777B92EFD44EF56A7D6C77EC2E619C5
137CD97F147F982DB7F48F3E448F0CBD207D4567E46974DB04E77C1D462F0970DAD5B1793DC937728E3CC4B7D133FE4D8EE5
985B79A03FA7500FFB1CEAC27FC98D953952BAF7DDB32F629D65C332864A9DE3E73509115DE1A7E43938832FF322FCEEDC2E
F4827CEEA9F3C1D00B8213D7E7BAD4459ED17FFA0E3F9917E0297A4487E42A7E8BCFA243F4844FE29BF49F79EBE8DC629BA4
DFFCE4FFE9F7B2FC2DDCF73333D7A28FCEA3823EF658906090CF777C96FC99FC66502F7E516767F34BAFCD17CA6778CC1C87
2FE293E86C40E104C37CDE0CF836B8BFCD63EA737BCAE774D7A57EE65CD631286E687FD96FC1D3E8FC41DFC81FF8024FF07F
F8470EB01E729275FD32B165169F1FCCC53D7607F5F7F92F41ABD7F053DE747B1CDA5F78D23D77E820E123BA400FCC01CCBF
F97FF93CD4F7CBFF77DD0E3E2EDCE7FB074D6F495DF08E3C275FC815F8045FC817F286BE92C7F497F5B24EE678F28379A865
B544BDDDEB389A4BE4147D0647EAA5CED5B55A1C90EBB11FA09FE41DFDE07A3C9E3A781EB951E6A396DDFD79A465BD2621C2
07F20A3F66DDE42B78E1F7458635AFE7EF836D2ED5F968DDD3B4E9E467D4E0D06F641DF06EE0D1D041721FE6CADCF34AE5F2
F3FF78FDCA038507E41679093F7A750DFDC6FD9C105A2BF83C4EE087CC0DF8223E468EE167CC55F07D5E86A487FCBEF27BFE
DFE655D5D9F6F8D97207FC0A2D217EC6FAA25A84D612DCD01F7A2487A89BE7FDBF9F65659E27D7E90B7D025770E671C67FDD
87A08331691B4F12DEE2D7CC775327862D97BE91F7E4397D627D9B7ECA18EACF39D50FB6F8C19F77E61E6E56477C96790E9F
DED3236CA1F091B9871C638EC7D76ABE4A3E5BE64AF22AE7E8AF3F75DFF7A3732B79001FE1217989BF8E5F1EE8881C8627F0
065E576C16E892BCC7CF6CDFAEFC349DAA3EC919E3A7CE61F094F9831C0AC998245CFC9BFB7CB5B8C5017FBE0037E62CF003
1F7061EE6E7DB9F124593FB98A2EE907BE4B5FB6FE5FF3A2EEEFBAAA3F3207E13BCC616D5724FC557EC207F2165EC06B780E
8F6D8ED2FC3C58F1F350C94972A9D080B085FEBE937D40A12E799AB8CF812B4FE91FBE491F470C0B7D2EEFFF83AFCC4DBFB5
F9FA53F989FEC975F806CFC809F88FEED1017D86EFF01CBDC1137813DA2C6CB9CC0DF091DC385AAF542E37FFE97C00FFC003
7C9827E80BFD607EE2FEE80EDC8D57CA27E6827DAD42B7089F986FD099CD2D9C2B697DCC07E866D43741DFD159AE25D50FB8
D77774AEB47EE9FE81BEA153FC14BDC2BFFA799BD511DEC4CE2DF00337F40CAFE133FB0178FF75B6D4EEEFD6A31FF88E8EE0
177CC317A983FB534FFB8EA1CF853FCC17E869D1AFCDE60A7ECC73E8923CC02796BC9D7CB6AF1FE6BCCE291216929FF0867A
D1FF4F719BFC247989CF923FB65E5D27FB8CE847A15D052F78CF5C461EE6C95276939F6FE887BC9A712149B83F67A073FC0C
5F4447E88A7C63DDAC9779E4F3B1A15B8427E81CDC982BC18F5C2AF5BFE645FD1C3B56B3FA01774EA1BA65AE629F547CD067
C77CDEB11ECB635D17F31CBC8EBE11E04FFFB92F756C18D83CBDAF37E650E6C4FC470BAFF07D1DFFC13FF0A103E54366BAF3
449D0BD11F3C607EED3134C12DF7F70A554FB99FD56E227319F3B2E597FA1D7A341C749F061E9DE3857695FE30EFC0DF5995
C31E483DE42EFECDFE0D3E7DD9267D5A791E38B07EE654E662701DD029A87FD037A1A9FCFD2BF31AD7E7BAECE762E7C39E21
79DCDF3DC1AFF029FC9ABC245F9B44242EE9E3485DFFEFE741F7FADCA56617FD7917DEB18F2CF37BE6750E37C5CF7253758F
DE793C75F33CE6E5021B4266FAE75247E3D46E22F9015F6DFE575D939FE89B3A99BBA917BFC537F1DDB5999ABFE3F783391D
BEE46D56F78FE7CF9FDB5C7FE88D8675A5CFCCDFF071DBB3B05FE479E41FF9DF6D45682AE12FF3263C0A7BDCE427D1237E4C
BE319F9137F54A277B4FF28C3AA98FF5998E74EE444FF875ED92CDE64ADE9313F8077E827FE36FCC55F81CEBCEF75FC3BA82
037D8A5AD0E8A6F0135E32771F9E58F70F3927A04FE064F39EE6FFCAD3A993BEE6E7EAE3F825FA625E4367BDBE6D7453F064
6E8F99E31B9D101E5107BF4707E882F900DFC50FF147F37FF52172E0C8D9B29B64FE65CEC0C7C957FA491FF147FC059F9CD0
B1F142E12D39C47C449FC84FFA857F97EFDE3CBD3F1F194F941FCCADA66BD533F38BE5AFFA193EC6BC837F6CEFF6D931B92E
7C86C7CC4FE412FDA48FF833BE8CEFD89CA6B88123FE4B7E4E8E0CEA851FCC653F8F4E9FD6CF7BFC1DBE7CB736EC81E0451F
983B9803AAA56B76D19F5FC097E7F178F2025F20376C9ED4F9EE60EEC22B64BE23FFF39CAAF89D7CBE93B98CFCA35FF40F3F
C3DFD09FE5A9F203FFA7CF95CF07BA673DE042DF36A7CFBC2EF0FDA02E3B4753DC2D4794D7E409BA2F1127EC177F5F8B0EEC
9C4EFDD7E61C9D03E897E5A3FA1CFEE6E5C841BF6FEC3F0A85E4EFE7FB96E5A0AE9FE7511775822FB9C55CD5E2EEC75FCACF
F9F99BC6959C638EC43F5A0C6BBC50D6452EE1BBF8393E0E8F997BE0FBBC2DC9DEF3F54BFEA063F202DE901BE8097DE147D3
5636592CF7636E20C7F057FC16BE803B75518FF9B0F6193F86F7F80AFCA79FF8337E6D39A5FE8AAFB21FC05F6A364E594A78
B5E36AE813F77E5F9D27C83DF282F5B30EEE431E80FB98FDA1E7E5FFF141E63EE389EE7FE00B7E834FE00FE6FF8A0739001F
9897E1C5A6B66147DDF9B9FA1AF3694CBE04D72347984F0B242F9EC9CF49748E4ED12D7900EFE0DB914F2B7E2775911FE490
F989E62EBEC275E9A79D67D017ED07FB18FC007F807FE410BE3B2149A0177060FDCCA9E08ADF705DE61272E5FBA94967883F
C267721B5EA353F4492EE19779D755F940DEF7CA5CD6BA60E879C9DDFFAFAC3B8FD7B22EF3073EF34AE737D6D0E206822808
1AE68682EC20ABECFB72D8E1801A96A21162A1CE29954A1110A7500C3732C5845CC748C5ACDC2A6BD2DF94F6CBA5576EA963
4DB6BAE4F9CD753FD7FB7BEECEFCC5EB3CCF737FEFEBFA6CD7F53C9CC5BCE2577AA333F369E7A4D3F68DAFEDA378960BF689
B396EF5BE591F73743C71E7377A5ABD4179FCC7DFCB40EF5392CA787745DF687FAFB5D730C8F78A5137B18BD344F39E5D638
8F2ECCABB29FE49C79F4FA3E17467EEBCB3EABBF3237337FF1B8E9C0657DE27CFEB1BF781F421F93B71EFAE198FB7C27CFF8
4FFFE6231CEEF9E3A9FF55D7877DB9E099EF97E0EA7A7DEB977EE8C6FE22D7475D7EEA4FEBF3F3E1FB96FC39F8B44FC81BF5
A8C37E8C1F7C991B656FC9FCD8B9ADCBF0F8173E70F1BEEA8ECB4FDBBBFA7D39B967E3013FF83277BD5E1DAE33D7F9EFF26D
0D1D9679947AA773F5C2012E3F193F719FC82379614E7F75EFE59785EFCC0F73E9472F0FAD7E9F863D805EBD3F1B7941AF59
75DECE7E75599FA8C33CE7BFB257E439CEF5FE0D5E63DFEC71569C6B2EC11BCEE5BEF97ECFFDF9D33EB7FA958F54DFCFEA7C
E7BA8EEECDC31B961D787DFCFBB50DA76C0D9D98F77248AE9CF0FD29F7C47CB6DFC91F392B87D6CF6DBE2DF4CA3FF6A69227
9923F24C2ECFBCEBB47DEBF3D6DC1E347AC99FE37E72A4CCF1C449FDF61E7A59DFA981C30353970FAC3EB7CCBC33AFF863FC
074EFDAFFAFC29FB46E6CAF7761F7377FCDE1DF3AAEFE6E1DFAACF99921B9907F2E3D1F3173E1273D23EC47FF6B5016BFA3D
51CF393EE54F794DF7728AFEE5253DD3371DD817D7AEFA48F573CCAB6F6EBE2DF8B0A7D029FD9B3F7CC0E7E6D89DA71DFAE1
BA4EEC05F402BF2967357C8CD71F0F99724FDC97DFEC4FE69F79764A8703CE0A9CE9408ED283FD49DECDFDF641D7C5F7957D
F1CBCDD3C30774668F930B7CC98FF6627A59B7A181CF296F2FBF2C78942BE60C1C5DD7D2FA914DD5CFDBA60FFABC336361F5
FF95B9B7D229DDCA517BDC808B163E52FF9C863EE9D23E21BFE8935E3FF3C9E6E9A12F7B163ECBDE92FADD7361AF59D5DFAB
A1B7DCFFE8EE07DDE79D197CDB23E994DEE95F5ED381FBECEED0E3ACE0CF7C821FDC1EDF6BC6C2EAF38BE4DBFE7AD18F9A3B
C5D778B26FCA31F92527E07AE2FCD9BFFFCB5FFE52F662FFF67B6EDE99F1387CE525FDD2ED95E3965D113A92EFE6365DD193
7DF0875F9FFDFBEAF3A9AC1B4EE7F76DEE14F8DB6FF150F221F71539213FE5A31CA5EF85AB4ED95ADF17FEE7DCBD031FB8B8
8F5C9537F2D5F578704E999B9917F2B0F82EFD668FE05BFEE3DF1DD79E7246E490BD90EFE9D275FC2FF7E49D5CE38311B72F
1F587F3FCA9FF2492E79DEEBF941CEC913F9F2A395C3BF15FB87392D1FE5A6B9F7D5370F687C9FC377965F52F79BBD8C2FF9
948F4BAECAD3D4297CEC9B702A7B4ECEFF472FE9F7447DFEAB53DDFCCABF7C68BF1E7A42EF27EA3E965F789563F2905EE4E2
754F2E9F10D799EBF21E2EE6287CF886EFF95D1FFAD24FD95BD27F379F7CD075F59CA47B79B9E25BCBAE88FAE4B4DCC50B3E
ECA90F765BBAACFA3C2DE7A279CED7E6187F9BCF27F61FBD24F6137B8A3D9ECF67BCD0F5C4C0DF1CB2979CBE69FFD5F1FF9F
656E66EECA5B78E95BBF6DF3A23127E8057EF213DF651E257FAE93CF7299AFE85DAECAD95DCF9EF2B1C0C99E4D8F453FF9BE
828EE4A5FCA4F792CFA92339ED7C7B7AF97F82E1CB9FA9F3647F87337CBD3FF8B7F39BF7045FF2D7DE54E65AE6AB5C75BE73
DDDF5CA027F3C11CA22FBA926B4DCDA79C51DF2FF8C25EC31FF25FEEE303DEE682D79947EA806BC9F3C455AED7E648F5BC9C
B70F0EB877D0B8FA3CA17FBAB737C8C5B1C3975F52DF17EC49F27D71AF4E7D62BE79BF3272E45127D6750117F3043EE6871C
DD72F6FEABABDF3399EF5BF03AF9C8EE0B822F7979D2F34B97D5DF47167FA5AFCC1F3E352FF9F5878F8D5E52FDBF41E6B8DC
2FF5E7FB1F7D2C3D66F984A817FFE661A933EBF37EE88C8B9AF7047FE68F1C326FE0095FB8B88FF3CBBE9479FFBD1FF77E22
7ECE5DDF7070AEF76163B71EB14FF5F76CB20FB963EF39EBD0FDAABFCB69BECE38ACE13F3A2EE7643FCE93BFF27BD7C15D4F
6C7C6E9CBACAFD8FBEE865D277963F539FE77C2A5FF9F5BBCF2F9A597D5E99735D1ECBD77E1F18B75FEC39F69DA2C3D49F3D
ABCCE59C3774B5E9EBCD1FAC7FAE69DF29F33DF5401F373EB5ECB5B8BFBD848FB7BD6F59B7C0A9CCB99C0FF8C1977926BFE1
073739422F7472EF854B7F1538F1A13DC23CE69BCB7FB55FF57345F8E3333C9E3DBFF98371AE39438F8F8E1934AEDA0F72BE
EA1F5EE6CF0D4F75EA53FD5C61E2601F8147D92B52A77BF61C7562ECBBEDF711F7A37B7AA00FF36158B74533EBEF6FE8CC7C
BBF3E7DD17C4CF5DA80FAEEAFCC11FC7553F0FCDCFE619DDDBAB56DFBB5FF5F31BF446677270FD434BD7854EE8C3BE537C94
FE9137F8C2933D439FF8A17FFA311FE8A8CCF1C475F76947EC13B87EFF8BF32B3ECD69FEE083FEFB4CFC61CC2B7B10BFC83B
BEB1DF0D1DDBB75BF4C9AF7CC197FC618FE1ABB9E3BB6C8CEF13B04FD0E3FA31FB553FE7671EC8DFE53B0FFC5BEC71F8C4AF
DC2AF99FFEA4E7E5EF2EEB1675F087DC35AF061E35F4F97AAEAF1EBC745D7DCED33BBEED496B17ECFBD3D87FBFD161D9EAEA
FF4B738FA16BFE59D06BD96BF57945D774CEE725E7D367FCC5DFF44BB772853EE9D53C1FF299F9FBD6DFC77FF199A53D63FE
D39FBDC99C80137C1EFBF3C41F464ED3BD7DE2EE379A3FFEC20B2F145DDAA7EC83F85C77E3BE3FAD3E87A793DCE7E8C59C19
73C1D25FD5E718BCD4053738DBEFE08D0FF76F69DDF7EED07DD15BEACC7E472FEEAB0E7A1FD07BEAF9819F7DEDE103E6DC52
7D2E967B15DFC041FFF6BFAB7A34778DE7E5B8B9F4991E4B7BC65CB24FD1D76D6736FF779C635FD45FC999CC173971D1EF96
BC1A3CE2CFBEE61C7538CF9E37E8B539B7FCF5AF7F2DFA78F43FA69E1F3AC3ABFD504EC80D7AA53BB9E2757CCA9FE627FC06
9E3FF390B8AF3A1EB968E621910FEE07D79293991FF2525EC92939CBBFF6033E3EFFC34B5E0D9CEDAB6DFBEB929FC5E3F0F2
7C99839917F2471ECA3BB978FD5BCD1B4377E6B7DCBDBFFFA2BF54BF873C73DFBC2D733373459EC82973725BFE7E077EE55F
3967FF1D79DCB1D3EB7EB2D7D2F78C9643CF8B5C28FB43CECFEFEDEEDB2DE6A73927576EBEA7CBC6E0D35C99FD66C367F2EB
91FF1CFA7CFCBFA8B92E77E48D7C90177C3C71DFE68FD7E7335F9857FC51F6849CBB7028799B39CB070587DCA7E1E17AFDE9
ABE093B8D8AF3FFE4C73D7F0A5DC351FE00437FB46D9731207B8945C4DFDCAD732D7524FF2907EEC417464CED9034E1F7140
6BF007BF695734FF777DFECB2579C42774487FF68932D73287E40FBEE52EDEE521DDCB457372E0A32336D4F782B61C6EE42F
DEBC6F18BBBCD76F0227BC957991B8B8CEFB043C4D7ABAC79BA12BB925C7F8438ED31BFD996778C6EF97FB2E1D1E78C9657B
4AA933DF9FA8573DEAF03CBDCB4F795AF2397521A7F978D4638BFE527FBF55E642E696BC9293F211DFF62739BBE8B1CE8784
2EE52D9EDAE7AEAFE1675F1BFAC5FE57C7EF5D93A3FA90CFEAA08F25FFD4BCB13E9FF5A12F7C781D1EBCBEE470E2248FDB72
BBF1BC7C2EFB4FCE9F5D9F3BF4BCFA5C92737B7E7AECF4C0834FE89A5FE4BBDCDDF2C001ADF5CF97F9CD1CBBA9E7D297E23A
FB1CBFD299FBF0E79E5FCE7FAAFABC2FE7933970CFE6C5BB622EF3A1BDE9933F5A3A3CF0301FE450F157CE333E7B64C088EA
F705B6CF5F7AA003FB4EE137F72D3CD3BF3CE603FDDAB3F42D8FEC1F2B5F3A6060FC5E8D6BFECFD2F9D1A7F927D7E52EDE0A
5FC90FBEE42C7CE062AF2AFB46EE19777DB4C79BC177D94F92FFDDDB7AFD26EAC1ABF981DFCB072EB934FAE473FB85DC1AB4
71CCF8FA9C97EB7086BB1CC11BBE6E3D7FE9D151AFBD0D7EF4400772596E8DEC39FFA9C857B970FB5B4BCE8873F46F8F2B73
2D71DDDEBFF321818B3E3DEF71F5D807D5E57CE77ADE7CE28FCD9D0F1818B98D5FF30ACF653EE63CFCDE67FA5F5D7D7F69FE
EB71389EF5C8924B23BFCC1F7AE7CF79BF6CE84EBEC88FF12B17EFAAE7BAFD86BE9B7EDFF5C9E0E73BA7377DB1FA7C27E71F
FFD98F468E3CFE864AFFE9033AB64FACBE62FF5BE3F3153CE0459ED12F3D161DA6EEE890BE1FFED498F1B10FF1A739C147FC
236F2EEBB3F86F91EB74698F808B1C858FBAE003AF657F5D3A3F70A053B94BFF7CCF07E687BC5BBEA0D3CBD5EF77CFFB39CF
39FC228FF8A6BD2FE1C24F83CF1FFFFE984BF6083CE2CF5E56F8CA7D0E6FF6033A5D7FE6FEB7567FFF2EE7D5D0FF1E7C623D
2F4BCEA7FFE9D9DC1BBEA5E98BE13F3E94E3B32E5A7A747D6E7CFAF1C57F0BDD9ACFF4CBFF744DCF743AF59F969C519F0F74
661F5A3B76FF3981D7D6F72D991C3A959BF2F8C1E7667DB0FA5C26E73D5D7FE988C557C5D7F461BF904FF4467FBB7B2DB8BB
FA1C2973C7FC79E8A2F1EFAFBFBF37EFED5FF05EB77BFFEAF7266EEFB1E482B8CEFC93077CC76FF24E9FCE6B69DD7F60E802
FF43364F7AB4DA3F724FB1E7D1C18C7E877D28F61F7BE2D8E947BF2FF0A08BA2DBD4ABB967BF39A9FBAC0FBEF5D65B85EF73
7FB1F8AAE0C97E81CF9203E97FBEFC5EBF696BE2F5F61E3C14BDE51CA5BBEFAF9CF468F00047D75DBCFFE205711D9EEC3BA5
DFECD35E59E666E6CD8E0F767DB2FA7FD5C469E80FA6AD79FBEDB74B9FF42A27E8B6E0667F4DFCE48C3C90377CC017F24FBE
D0337D97399E7360CF9EE36F08FCE98DCECC75FAB15FD0D179AF2F5E1038DBF7DAF6BFC563A34E7879DEDE0DAF49371FF193
AA3F8FFB3AFB9617FC2F37E482F9B9EDC64E2F87CEE5955C935B7C6A1EF22BFFCB09F9600EE3138F70C49B7C801F3CCD51F5
95799F757EEF43834F8CEF1F345FF9FFB477974C0E3FCA2773401EF03BFFBB8F739C2B7FE40E9F971CCBFCE21B7A19FBD482
BBEBEF03FCEB717A326FCCBFD3FFA1E38BD1B77CA35339670EE94F5F74BFE89986BECCAD9B3FB57853E8CE1E4677658EA40E
E49B9C901B7C70E7BF2EDC5E7DBE973AB6977DF9778B57840EE4B5B9CAB7E6E1AE1F1CF6A1C0873EE490FC51B7FED46F8E0D
BB70D8ABF5795BF68DF4DFEE6F1EFDBEF87E56B94AAFF2D59E241F179D77F0E2E0B72D9F1BB98C4FB8D91FE10707FDDB23CB
BC481C9D475FF2503EDABF465E7062D72A3F3347E5AD9CA5BF6FDF3EF7E2EAEF98E45C31C7DAEF5F77CD3DE227A127794B47
72575ECA4F3CEFFCCDA2ABE2EF89DB83F5ED717D7A5E0EC3D779FA93FBFAC1B79C6B9F6FEA9573C517E907F34D7ED0175D15
7E9357FBD4B51D160F0BDECC69B9D596B78DFBCB55F92817B7B4767CB1DAE3525FF426A70A8FB99FE1F3931F5EBC22742C17
E510BFC851BEA1737394DEE5C4E4CF2FDC5E9F3F7273EEC6C59BE27EFC4F27F4610E7DF78261D5DFF7958F7250EE98F72BEF
EE7870F459E642E62C1C2EFF8F453BA36E3E373FE5C4C953E75E1CF92707E16FEEE1A1E473E2041F7CD2353DEB7B66C74557
D5FBF73CFDE349CEC20D8EF2E3A4C12377C5CFBD996FE60A7D6DBFF8E0C5D5E729E64BE20647B8CB2DF8EB1B0EEA80A3FB38
DF9E47AF4D77745B1EBC9439987369CF852776AD7EDE38F7BEB1CB8F9B55D7DD7D3D661F5EDF43E584EBCAFBD2EC47EED3FB
E6891D0F8E7ACEEABD68675C2FA7E9B7F9CD068F7893D3A5FEACDB1E0717D7797DA93FF7307DDCF5977957579F7BA45ECD99
A2CFD4A53CA0577377F5AA032F8BCFA71F7C68E4AE98EBF42D7F2FBB7151E7D02DBD9AABF441577456749E7EA577AFC33F1D
9A2B726A79972E3DE2FB1DBE36795197E8C3DCE503B93BFA99D987478EF28F7D0E8E333ED0F3B56ADFCB3D6ED29547BE5AE1
96F8991B239F1B707F3D3FE85C8EB4D7FBF0CE279F1BF3B1EC752F2DFA68D4251FE553A927EB68BF57FA5A5D74487FF208BF
7835474BDE66BED061C9CFD4231DDEDF32FDDAF8DA5CA477389BDB65CFC9794CEFEB371C7859B5872C58D439AE3397E8BDCC
97CCA11D53BA2D0FFDD199FD606DFF03FBC51C9007F4467F450FA90373A9CC11EF13B61D372B7CF99D974E3E37F6037A3547
E9824EE4CA97562DBC23F8A16373DE1C1EF5B9E9D7D6E73B9E471C3C61583C6E6EDB53E0B4EEB103FB45BF455739CFE8EB9B
93E63F18FB00FEEC053BEE5BF0DB78DCDE800778B94F4BEB8107455EEF7970F25BD5F7ADE77E804F3A3437E8D17DDDCF5E22
5F27FDF3BCABEB79AFBF91274D7EEB9D77DE29B83CF0C284613167E0EBFEEA77AE3EE4345DD0C3B91B16DE11FAB147D051C9
A5CC23B955F22773450EE153FDFA912F0BEF5CD4A53E67E8857EE4FBC58316AE8A7AE8C03EC557E6CFB617BBF488FACA9CCD
F9B9EB8F3D5F8B7A4FEBDCC09FFFE47D99B3A9DFBB561CF96AE45C996B8993F3E58CDC9137E73DDCA8D3DE05BF328FD2FF7C
5AE651FA954F5B5A172E8AAFF5E93CF920174A1E751F707FCC6FF38F6FE16B4EC1B9CCCD7CBDEBCB1C4C9CE0C3F7FCCBCFFC
C78F7CFCEFCBE62CA8EF6BE69EF961CEADF871E76E51C757762D1C54FDBEE1CC6573C9DE205F168D3EE48DF8BE593E987ED7
FC07EBF3ADE458E6179D7ABE69F482DFD6E75BC999D4AFBCB1778CDD7AC293F59C307F469D3A644D7DFE941C4E5EE58C5CA2
0FF9C467744A9F72435EE0BBE479FA921FE9946EF9E7BAFD173C51FF3C46DEE04D3EC80BF9214FF059F04FDCED416DB9D7C8
3B78E31DDFE662C9E1C4C575720C2EEEAF2EF90597C24BEE55F8912FF4401F70E56F78953C49FEE50A7DF00D9DD02DBD9A7B
2527D3F7FC5EF49C794FD77831CFF1F389198DD7CB03F9242FE4DC952774EE16F8F1D9846BE62CA8EFE3F4699ED169C993EC
5BAEB4E555E37939E33CE7B7C7174E7C75FFD543D6842EE48EBCB0E7D14DD3A53D8684FEEC77936E3E6670E4173ECB3C4A7F
6CBFEF9037AAEFCB95F399BF70D8FCF1059F09DCE485F923BFE00ACFB25FE7DC5AB9E1A0D9F17DD8B7DDD4F48D988BF45A3E
0FFAF68CDFD6DF7F943C33A73247779F76C293914FA5DFECD35EA96EFDE9ABE090FB253CE4313CE14B974B5F5FF0443DD7E5
D0E88F0C1F1EF3C85CBA71F5FC95D1977D85EE5C0F27E7E8435FEAB007EB67FABBBDF6AAFE8E72E2062F7B8CD77B9DFED463
1F51D72DEFCDDB27BE76BD3D65E5550D7EE51F7F78BDD7394FFEF2CDE65507CD0E9CE5D9F89367FC36F2867FF08B57F965DE
F0F9B237BADE1E3CD01B7FD3DDD58FCEBF23EE2F67F9FEBEDF0E1F1EF3921FE46FD161EA4FAE983363DFEC774EFCFD643E2E
F322F36CC7393D86C4F51B3ACEAF3E7FA77B738CBECDBFD5733B3D147A2DF322F57BD7DC63AABF2FC3DFD3E6357DA3EEF36F
BD33B5DA4FCD213EA037F9D25E77A39F1FBD5FF5FF0C39F7E00D7FB9224716AC9FBFB29EBB7C8B7FBCE3B5F83FF9759EFAF5
53F23FFBFEE6FB7AED55FF1C0DAFE61B1EE5343ECDA5717B4D9D5F9F8FE619DDADBFB9D34371BF55BF69F0623ED06FE137E7
169EE5019EF1CBF7F4434FF899F38FF3F6A9E7F17DDD46EF17BEA3337314CF78E7F36FF79C58FD3E7CF315CF7463AEAFFD58
A73303373CC90F7CE179CC4FC6BE10F898E75F7A60FEA70307BA346F3D7FF22F27DEF2B7BFFDADE07BEFF1635F88F9A96E75
D917E0BDEE679DCE0C7C4E1D30FF8EBA5EE4D9ED6FCE5A1CFCF1B35CA7673AE67F7DBA4F4B6BA7EAEF3EF117DFF319FDC9F5
6BF6ED7A7BF8ADE47CF2E97178A9579DC5A7E94F7970EE88066EE6331DEDEED0EF9CD0AB5CA68F92DBF9B8D7153FA63FF9F2
E2A5F34F88F3F1612FA03BBEE1177AA77F3EE6373AA44BD7F11DFDF1337FF085FD800F17F5E8FE70E8D7B97447E725E7D3FF
7C777DC7790705BEE61C1F9F775DA35F7B0E5CED1D939EEE7D5FF57E3B7D493F7463CEB5B4CE3F38CE8197F3F892EFF870CB
27E6DD1FD79907E686B922F757EC38F8D6F8BEC38243F66F8F3007C6AD1CD4B79EBBFCCBB7F45DFACA7D427FF03187E1C4A7
F4CEAF2507D2FF745A723DF5C4F77CCCD7F4020FBEE673B93DB5C3ACC5F5FCA62F8F7B5DC9DBE45B0EDC3A60F6CBD1179EEC
112587D34F7CE4FEF2058E72493FFAB0F7E0A3E98223E6C4F5F83537F06CDF99F181A3E6549F33260F5FDFD27441D4695FC1
83FAD56DBF293A497D9857A59EACC39ED596638DFCC29BFB3AD7FDD5E97AF596BC4D7C9D77FA57E6DD1FE7F0B1DCA54FF38C
4E4BCEA4BEE40D3F2FF9CDBC83EAF3A4CCC1CCA1EDCF747F383E07BAF3D8692F569F9B675ECB75FE97C757361D7C6B5C57E6
5DE6E25D1FED7D5FE462C9D5CC0339F0ADCD83FAC67C95A7FC2CAFF5AD5F757B5EFD251FB26F39610EE84F5FEAF6BCFADBF2
AA719E9CD9DC3AF70FC19B1C3197F4A75FF9EF5F7DE953EE98372B3FD5E5F9B89EFF663D3AFBE57AAECB9BF11F18B6777D9E
787EFEE94D17D4731A5F7892A778351F765C78C49CB87F9917C9DBAE3F1E3527FA901F939F98F662E4815CF8EACD4DBF0BBC
E53BFFC90BF5EBC75CE19F654F74FB42F02BEF271FD9B739744BBF2BFF61EE1FE271F941EFEEAF1EF7BB6BF6E4D1D5DF35CE
FCA5537CC909BCC92FFADDBCB1CBF3A1877BFE386CEFF035FDC94DF8D21D9C0B8FC99F9C808FBC2F9F07A5EEF0007FFCF22D
9E375C32774ED44D67E6037EF0457F05E7C4D71CA397F18F8D783EDE07981B78D21F3DD1A579B37A4C9729F13954F9BC27E7
19DE4A2E256FF8DA79EA8C0DD5E7EAA93FF341AE4FFAC6E4D1EFBEFB6ED159C993F4215EDAEB54FFFC842FFCDD7DEAF82F55
7F4F3BE74AD145FF11CFC7FCC3AB39443F74C33FEA57B7394107A7CC6DFA5DF42F8756AD993B27FE9573F40BEF095F1BB33B
FA35DF5C0F3FE7DC70E09CAEF1FF7EE6131D9907F4BBFEDE2E53821FE74DBC7AFC97DE7BEFBD82EFBF2F1CB33BE68B3EE151
F230F578CDB1DDBE507D9F68CE57BA59D4A9E773D5CF7B271FAEF3FA4BDEDFF452D44797E60CDD987F6BBB77FE737D0EC8C5
3B7FDEB7397251BFE6A4BECDED19FD8E7D3AF44377F4423FF2C079CEA78FE2A3CC037EBAE9B65953E27E78359F3DEE6B7515
7C1217F3BFF82B7DCB67E636DED63DD7F9CFA1E3D257CE7BFD15DFA5DFE40A1C9DD7D2DAF917B1975C79F89CCA6FF257CEC9
C59957CFD850CFE3357F6A7A295E6FBED169C9DBD42FDDDA4BF4D974E9913BA2DEF2787EAD8F929FE92F7EC38F7AD477F117
9AB647BD74636E9B13F274C515877E2E7EAE87AFF9996F4B9E675EF085793279EB809E916BF2F5BCCF366DAFEF0B6DF3BE69
7D3CAE0ECFCBC7C5AFCEE95A9F037CC677F4C8EF74C0F7F8321FF0869F79D3664DA9E731FCF8979FFDCBCF70297E4CBD9779
71D1D47FAE3E57CD5C93C7651E65BE6E7FA5E773A12B79A26EFD95B9933EDEF583639F0E9E8ADFD3E7F45B784CFECC2538CA
3F7816BE72AEE08D2EE9976E4B1E26DFFCB5E2FFCDF9437D8F90AFFAD4C78E738EDC1175AAC3755E2FE7E0041F7DC341EEB4
F9B4E14F78C01FEEE6CA9DA70DE81978C829FEA163B97BE599877E2EEAE61F3AA19B9233D9BFFBF3175DF059C9F3ECCFF957
7C654E8FC86139222FF5A74E75973C493FF191FAF04CCFDB2E9CF55A9C6F9EF101DFCC387FEA3FD7F73B7E9DD26BF03BF57C
2D3ECDFEF8555EC8CD95CD87DC1DFF1FD0E6F7C6757CED3EEE0B5FF38F5E973DD6E3BAF87E43F59B57FA309766B49C70E12B
AFBC52F44B2F74476FFEA5437AFAE6A6898BAAF713996BFC7AF3FDD35E8FFBD08D79852FB98237AFF73AF3ECCC4FCCE91175
F1137F947EB34FF3ADF495F3497FDB3BCC78B13EAFE171C75383DFA9723F752CBFBC5E1DAEC32F9FE3B9F098FCC90F3EE78F
CDD71E7277E02B9FA69F3D71517DAFA4BFE517CC7AEDF5D75F2F3EA7A329FFF7A467AACFC533A737BC3DFBF35127DDCB51B9
003FB8E10B4F7C7EDB17C66EA9FEBE75E6323ECD51F82EEAD1EBF8EAE77EB23EF751A7B9DA74C731D5F78B943D21F52D6757
0F3AE4B0F0275FC21F1F743D77D4B4D7EBF951EAC93ADACF735FBB7FC9A5CC0BF8DF71F449CFC43E806FF3A9E467FAE09AFE
3DAEABEB5DAEECFADC091746AE9837D33E3B764B6B6BEBFF9AABF8E557BC4DBD7AD403F1B8EB6F3F7554F57B06D581077D2F
7A73C68BF57C727FF5C0519FFA50BFFBC313BEABF69EFDF9E04D2EF2995CA6DFF50F1F7258F076E5B3335F0E3FC9413EA643
FEA6C74BE6CCEE105FD3A55C2F7A4EFFD0357D98436BBB74FD61BC3FA44FBAA463B935E3857EE3ABFFC74F5F16BDA5CEF8B2
F83D79A783927B9973DB9FE9757CF8B9E45EF2BF63CA3147C7797296CF575C7AD8D63877C7FD93E7477F7C25FFF0C7377834
F7E0BDEEC5AE3F8CBA4BEE253F7881BF5CC703DDD121DDC1D17D5A5ABBEE8A3ED6DC32BB439C6BCEE09F3ED44B7F17EF9CF5
68DC177FE658C9B7F4159DD0395EE99DCEE91BAF70839779801FBE703DDF3857DF450F99FB7471DECC46FDE6A9BE0B8FC99F
79D2D23AEBF678BD3E5D877F7AE0433EA323FA5971D8CC97836F7CCACB929FF9BCC78B5FD24F7C734DD3F45BA21E730BFFC5
2FE9137ADF7570BFF191A77281EEB7BF33F5DCB80EBEE65BC981D4053DB8AF79E8FEF4D5346AF2FCFA3EE25F8FD3537B3FC1
B1E44CEA86FED5A93EF72FB9977DE993FEE4DF95E71CB635CE297998FC78FCDF8E9CF9A1E84B7EC92DF3AAE9F7BD9FAEF493
3A6AF351C33FF03617E966D91347FCA1F17B571A393473F9A0C3EB39577C943CF353C981ECDF7D6ED939FEA6EAFB2432EFE4
4AF16FEA9A9EE591BC5BB9A0FBE3814FE92BE7A2FE4AFD59B7B9CAFFCB764CBF25FE8E3DBD165F641FFC51F0CC7C826BC993
BDA69E5BCF95363F36CEE32FF37BD1E8A3EEA97ECF8C7A3DAEBE7C1E3E70912367FC7CE687A25EBAE6BF9DDB061D5EE57EFA
41BEE04B5EE04DAECC9939FEA6D03BDDC3857F8B9FE19FB8CB15AF53977AE862E6BBC3FE25FC222FF9076E702CF9933ED8F1
C1DE4F870FBE71C9989F577F9F20F3945E4A2EA5FFAE39F6883F846EF98D6F36DFD8FDF1D09BFEF4AB8E8DA7CF5811F5F289
7CE23F78C2575DEAF43A793D7BCD989F07AEF0DDF9BE61FF12F98E0F39AF7F78D30F9C667D73C49FE23CE7F2B7BEF57BD5A6
69DBA27EB949A7B74E1FF1A7D095FBC1AFE4499E071FBE92B3ABFB751F54FD5CA77C4B7EB6DF77D43DF1BAAF3D33F1D9B81F
5DCA5DF5C857757D6ACB8C15AFBDF65AF137FDE217AF745FF4937AA5237CF165C9B7CCB3B9E3FB7E2D7448EFEA549FFC9547
F4B2FE07DD07559FD767DE157D5F71F81BC11F3EF14B47975C3FFDF7510FDDC8BFE2C7E49D0EAEED39A54F7DEEC0A93C9EF5
791EDEF0A73F3CC9EBB50775FB753C4E1778733D7FE201FE85F7E4BBE46AE2294FE14AAF1F3F7BDAB67A4ED01B9DE1959EF1
47D705E7C4574EF3E5C21E139FADE784B981B7752F77FB75D4CD1FCEC5CF9A25D37F1F7DCA63BAC3035FA85B3F70D2177CDD
BFA5B5DB43D5EF4F9333E91FFAA10BF5EBE7E227A6DF129FABD1873973F33D7DBF163A9323FF2B67F271AF8347F32FA7F489
DFDB818FF38E9D7E4BF46BEEE8DB9C5974DE71A7C67DF149B7FC44BF2DADD337479DEA739EF984A7656F7CECB3C16F393FE7
96FB94BC4A3D5E79E6E16F045FFFB67BEA83713F39C5C737DE36AE5BF5397FFA5FAE94FBE6FDDC07EE6DFB50E3FCE28BF403
5D97FC495DD0839C98FBED0117D5F3A6E443F24A6F70A717F8EB4BFEE98F0E5C873FFDEA534EC91D3E5CD97C78F7EAFB198A
0E1BFAD35FD155F2495FC52F59B7EBF863C1B471DDEAF3B0E44CEA6EFBC5C79D1AB8151FA5BEE8AAE44C5E77CDBE1FFB6CEC
0705CFCC0BB8DE7CF2808BAA1C4C7DF3E11963A73E18F8E187FFD4A54E38169FE6755E5F749B78D0EFD7578D1E57FD5C77E6
117DB5F9A5711D9D6F1D3CF9D0FAFCC267E125F9905BF4326FC490B3EB3977EDCB1336C7F5782DB99EB8C299AFFD0B1F78A9
C77C28F322CF77AEE7F9003F78D9F8B3A9FB479FFCC6F77273FE86D1E3EAF38DEEF96CF3B587770FFDDCF4C090B3ABBF679D
FCC9CF9B0E1BBE34F0F535BCE131EFD9E14BEB395DFC9BF53A5FBD74052F3CE28F7FCABEFF589F3EC127DDF10FFFAF1ED473
63F5776A922FBAC45BF1E9411336C7D7FCEADFD31E9A7CE8DFF93E73900E56EC387275CCBD8273E601BCE10A67BAFED4C7A6
EE1F9FBBF2191F94BE324FF5577C94F8E3A3D49375C84BF7D1973EF5C347FABAE43753B6C67DE843AEC811BA58FF70CF8DD5
EFDDE6DFEC43FDF44007F893DFCBBB1C7B64F579B7FEF271757BBEE827796AAF235FC3B5F09E7CF32B9CDAEB953EE4E9DA2E
3DFBC55E42777CB9BD7F9F3ED5DF23D876F2E3D5F703A4FE8ABFF3F9E2E7E409EEF0A6A7351DA76CAD3EBFC91CC27FF15DFA
01CEC5D7E98F2B9B8E5CDDF8FD458D5CC7CFBA177BF68BC7BF3C62D2D0D0A1FCA2477DBAAEA5B5E74151EFC57F9A32357C83
6F39AA5F7DCA9DA293D43BBDA85BBDEAE4CBC5BDFA778C9F17E43F75CAB9F2FE410E240FDB5E3CF6C8A8F3BCF74F99DAF87C
B191DB6D393EA55F9CAF6ECF177F25DE702EBC27DFF4B2F5D9F16BAABCCD3EE56A79DCD7592FDD2D5D7EF2E3F1B91DBFF9D7
E374202FE872E5A77A55EFFF8BBF926F3CC28DFFE1474F9F7C60D2D0C0910FE10437BABC7EF4C837AAEF9F481FD2EF0D4FF5
EF18F5D02B9F38177ECE2F3A49FEE905CF7853479BDE1A3A8377C121FB975B45CFD9BFEB8A5F927FBCD3D1E26B06EF13FFCF
2F278A1FB32F7D163F1E367E4DDD977264C97D23DF88CF9BF1BA69C8C4CBE23EFC43EF456F8907DDC9B3E50B8E7FBBFA3A79
833FDCF9AA4DCF8DF3E8F6FA6786BD577D5F4CE614DE6E5836789F2ABF9227CFEB7B498F61EF85DEF4434F7CB97963AF1D81
67F169F2E671F98B8F153F3EE6886A0FD757E6B2FEE81B9EF02DE7E4F5AE832B9CDB5F4F777829FCA6EEF18C17F9811FBEE2
E7D5637A750C9CF00F37389EFDFD8997C5B9F445EFF0873BFEAE3E67CCAD711FFCC989E2BBD4DDB61B8F7F3BEAD08FFB38BF
E82DF982071EF5E9793CB9CF95271C7344F57B661E1AB731709047F0283A497DE847FDEA966B971E33B143E88FBEF95B1FFA
8263C139F1E5E3524F9EAB2E7EA6C7F5F7F6EA58CF5DB9B0B6FB47BF13E7D3A3BEF5CB4FCB779EF8EBEAFF0D531745E7799D
F3CF79726287EA73D9F4199EE9EED44BC7DC1ADF5F437F70A0AF36BF3772481FEB9EFBE877E2FF6BF1E61CE7AA175F702C3A
1F3C6E63F04B1FFA779F96D68F6E8FD7AFEB34E1BB512F3EE44DC1217D080FBEC5C7CA0D472F08FDD21F3CE0B0F69509DF8D
5C973F70FAEA3F8E78A4FABE8AD43DBDABDB79EA6F699DB0A3FA3EF7ACCF7985DFE4955EB6CD3CF1D7755FD26BA93FEB961F
F0F37AD7CB95D3FFA1F75BD5FF3F251E456FA9037CE1491FF8708EFB398FBE4EE930F0E9F8390FBE840F5CF88A5F4F796FC4
23F5798247FAC2E7A66BC7FE22BEE65BBA6BD3434307F0FBEAE7873E597D1F6AFA193F78C567E1F1CD814F87FFE1EBBAA2E7
C4C57D8AEF52BF9B571D5DFD1CB7BE4FF9D7A14F067EED734BBFF02D7D651EE80F1F650EE5F95F19367A4F3C0F47F94157F0
A32FBCD03D7E8AEF522F5B5A7BBF15B839DFB9EA6AD37FE33E745BF499B8C191AFE878F5DCA37E54FD5C355F348FFD45FC5E
3AFE287C254F74E37578F37AF5AA5F1D4557893FBC0B2FD9BFE7D58B6F75E0E713DF1DBD277EBE47DD78809F3A4ABFD9275F
5E7AD6D8A6C08B5EE99AFF4E1FD1F7C578BEF093F5BAAF3AE88F2ED6DF7CD48FE2BE45275997C7CBF9E94FF7C1079FACFDD8
519F8E1C2975A61ED5AB0E3CC081FFF5B9F2EEDEBD634E7EE5B6E1FF19F7C31F9D9E73F9D8A6F83C8F8EF059CEC9EBD55B72
3DFB5EF7B3A33E5DFDBDC88163F789EBE148EF4527791FE76F79A0EF8BD11F7DE0BBE8271FF73AB8B87F4B6BE3FB8AF4A51F
3ED9F4F6E8BDEAFE855FC139F1A50BAF777DF1D72363F7A9CF0738D1C5E99B063C1B7DB5D7EF27A60DFFCFBACE5B5A4FFE5D
F8032ECEA30F7D6F9ED8BB77F57725E82DCFC1F7960386FCB2FAFE92F4033E8B0E536778D972F680674307F0755DE12BF585
377D9CFEDA905FC6EF83D18F73D5ABCE8243EA181E459F7B8FDE2B3EFF68AF03B8C2B38DDF06AFF8711FF8D043D161F6E93A
CF173D65FF05B7C48B6ED4A72E75F20F7DAC5E75DC0D8DCFA593C7C40F9F70821BDEDAF4D3B80EEF9E777F7815FC1347E715
1E933FFD943AB3BEF67E5DF952DFEAFB99F073E9CE512DF53CC17F797D5EEFBACBA70EEF183E80177F9473D237CE2B3A49DE
D66F38EE86E0034F251FF3F1727E3BDF16DE93CFCD9DFB76AD3EEF4A3EE8636DFFE32654FBC3CC512D917FF0C6A7EB9C03C7
C27B5EE7F57839EBF6E11D0377AFC383C7BD6E5DD3A8DEA17F38D205BCCF3A7440E7781F884F3A2DB9F9D87113FEAEAFEC07
9FFA745D4BEB71C746FD976F1BD225F8A003FCC0157EFABFFC57033A57DFD7DBCE0F6B778CEA1D7AA74FF8A9F7ACE543BAC4
7EA97F7DD187FE5A5A471D50FDBF53F65F7C82773AC87E0B6EA923F8C1172FF0D7AF3EE9865EE1B5FA8A3E9BE2F1C27BF6AB
FFF2FA7C9D3E8A7E128F8203DE1307F5B5F1D2E0036F85DFECD3F397BE74D2A6B80F7DEBBBE09C38BA8FD7EBD37585C7E40F
2E05FFC4150F6DFC369EC74BD149E2B4FECC3ED5EF6F29F8773E6953D40117FDB8CEEBFDDB5E67AE777FE7E8CFF3EE57F49E
38AE1DDB67503DBF56DFDB7F43F405A7F59D066F0E5CF40FA775AB4F1A56D7153CCAF9796E7BFE7DED3EEBC7F4DF103AD14F
1BCE0D3E8B8E77F71914CFAB6BF52B833747BFFA2C75667D782D75E6D7EA6DD351E33E2DAD7D0E893A0AFE89B7FBFA5A9DEA
5EBBFEA46175DEF1A02E75E2A1A5F5A41E75FDB8AEF098E79573F49B7839AFE09678E9BFADAF463F782878663FA5DF76FD39
B7D499F8157EF3FEEA6C5FB7EBF485EFB50BFA0FACF66F7D783CEBF3FCBA0D8387054FFA293CE6E3ED792DF567DDEB6EEC3F
B07A5F9FFF7ADC79EEBF76D5E061557DD94FE92B1F6FEBE7EF71726E1BAE8DE75B5AFB77ADBE0F24AF878373CA7DB34FAF2B
75667DFA6C691D7C58D55FF659CE2BF7CDFBE579E5F5797DDB75597FDEE77FD7D9B8AE8DB7C6F305C73CAFFD7D0A4EE5F1C6
756DE7E5397FDFCFFF075EA7489D
}
\declaresourcestorus{band}{
78DA6C9D777C5555D6F76319C75EB0F7C732961967B08E05B0204590A680481150447A912288E0205D7A6809298490463AE9
85F40E8190404880145AE804B08F0D5FD7B9EBBB721E9EF7AFF3B937F7ECBDD7FA95B5F6BE255E73AF0E8C4B4C4CF49AED15
98A5D77CF7E3993F79FE3EB33E3049AE0B6E0A5C2057EE9BFB88E7795EF7F905CFDF3FCF0C5C2CD7AF6F0D7CD37DDF828E9E
E7677B9D5F5D5858E8F5C5C0F39F6FDBB64D1E2F2B2B2BF39A977BFECB82828296C7A1E73F2F2D2DB5713F3BEB19EF33FFC0
0E72FDE2F5F35F161515797D967EBE8D8CC3785FCF3AFF5A7E7EBEDDFF75E7F36DE4BAB455C031F7BABEFED8330EE34E6EF6
FC7DF2EC80A3CEF373CEBF26E34FBAF7DCDEEDDBB7DB3CCBFF7DAE46C667DC654DE7F696949478ADB83E609CDCC73C4BFFE3
198771277EE7F9FBC4B10123E4FAE9B67335C5C5C55EE3E69EEB585E5E6EF3AC3A79EECDBCBC3C1B77D5CC731DE51AD0FAEC
A6DCDC5C2F1FEFE685CEF5DAE6E5CEFD779E7B53FE3EF287E6E53B76ECB0F1468E6F5E28EB1B56797693E4F3E34967E376EE
DCF9677E37DCEBC60DBCE79CDF70A3FB3138331FEB60DE59377B5E0F3E8C37BF3EF084FB7EF09FF188E779F2CD7DE004EE7F
5E1F91752E2A0D0C77DF0F6EBC0EBCBF987EFEA68A8A0AC383F1E62D3C7F938CC3EBC0EFB3EFCEED91D76FCC39B3362727C7
F2EABFEC6C9CE03DED25CFBCE0C67AC9277924DF4BE202BBB9D707FEE0C7BABEBEFEDC1E793DEB00EF4993CE4DD8B56B570B
BE9A5FD6037E43DE38B356EE1BF2FB997059FF94773CF3C22BD6DB82A307E7167D7B1E83DBA578711FB82D5F11B0CF1D0FFC
E67EF00447F841DCCB969D9B203C0447E6014FC351F10377C60507F2EFE289F33AF0258FF07EDC4FCD472B2B2BBD3E9DE859
FFF8CBFCEBDC79323E280F98DFF2ACF90567E395CE078EE04E7CE0EF3D2AA0AB3B5FDE173DF387BD7A7A88F00DDE6DBCFC4C
B89B6FE069BC8327CA0FD60B3FC8F7AAAB9B8F4ADCC6375D1FF9673EF8C4BCE04B1CE06CFEA1FC1EB5A0F969C927FC1B907F
7A88E47BA0CFE9E9B2AEF1EB3CF18EB9D3BFB53BDFC4030EE4DFF0529C886FCD2BFEA1EE7CAD39E919CF78A5F1102F79260E
F24D9CE469DDE7CD4FCBEB593FEB0667E389C6C1FAC7147BD633F225BF51EEF8C80F79210EF2056F863F78B66F5555955764
CCA9CBB3B3B30DFFD091A7A7CBFAFED463A35B4FE873CEBE0DBBFFD763D5A5CF657E4DEEBCF8947AD665BC525C19DFF2A979
22AFC61FCD23F999F5B8675EF4C7BAE6C76F8871AF031DCFE8E5791EFE729FF141F347DED00FFA9DED75AEEFEEDDBBBD16CD
DD30CE3D2E7AE275E8F30B9F736FC8EBE1DF7BBD4F5D2EEBEEF7E6A987054F78CF3CF3469E7B43F24F5EE085DFA1B37DE57E
C64557D39E3EF7D73D7BF6784D9BE9590FFC25BE91173DF91FFEF1FA256E5E803B78C38B25BD37DCEE8E079D198E8A0F78C2
6FE258B4EBDC5F653CF8045F87969CE92379589FE377839B077EFE9E751127B8122F71A2B74915CD61122FF9248FE8634A8C
67FDE8837CB8F8D9E8E6A9E945F3447E8CE7CAAFFFC32B1D0F7E2DBF3F30CC9D2FF46DFC54DEC137F33FCDDFB2679AC3246F
C62B9D077E19AF944FF0D0F8ACF8837BDC5F4E7C29FA45C79BB34E3D2CFD2B8FD11BCFB7F0DC332E3C24EFF8C5F849CDF757
57577B7D7AC4132FBA21DFF0135EB2DE4FDEF0E0FED106DFD56E3EC213F8011F4D17BA0E7865FC513EC223E3ADE60FFE1227
7A0A7AF94C1FE1CDCAC381F7B971C26F4C47CA23F8633C54FE115FC0B8F5F7BB791CF8A1273E786CBEBFACF97EC76F9567E4
833CB03E78CF3AF18B777E39F1A5E4A3F7152722641CEB53342FF0C778A8F9808FE8109D0C9A79FA9FB2FE090F78F2800EC0
13BCF185D10F9E9DB077EFDE167D29FEE00E4FE007795A1D1830232929C9F28B7F1017F1A063D385C643DE0C2F8D0FDCC83F
BC008761AB3CCF0F39E633D7CD3BF20B3E6B0F9D9DE0F4099A1FF2021F89933C919FB11F05CC70FB06F90357E6675EF2491E
C9133A458731BF9D88906BD04DBE77B8F9B4F11E4F1CE0814E86D79CF193F1B81FBE87CC3DFD4FE91B12DF3A562DE3330FE3
AFEBE19F2FB880077A9AEDB5F14A791EDFC20719175E32FE9CA2A09F929393ED75F822B8831BF893777025FFB3DA04FDE4F6
49E69FBF20684F4A4A8A8D8B2F82073890FF199F07ED71EB98F1E03D3C7BFFCD93870467FC05DF9BEDD5FCB79A9A1AAF459D
82D6A5A6A6DA7CF80AAFC3D7BE486B6E555B5B6B78C14770439FCC33EFADE656323E79848F7E4F9EF19375A1F3EE69C7AAC5
177BEC3BEEDC3F2AC13FDFCD2F78C73AAC0F987876BFAC675A46D03AB78F910FE38FF2063E0EBDE059FF076DD70D74EB66C9
E541AFA7A5A559FCF80AFC01677864754CE35EB4E2EC7E89D7F707BF0BC21378865ED109F3332FFC46871FFA9DDE28E38023
F8A157F2089FC82779C18F26FFEBEC8C7DFBF6599EC92F3E34E5F7A0D7DDFE46FEE0193890FF16FE7B74830E8237ADFBAF5B
B79B0A3D79357D29DFFF0FCF753CF8BE2C63C31EC93F79C7FF4C2FCA77788EAEC06169D5D91982AFF159E781D7F801FE1191
75F290E89AC7E89CE78DFFCA7BF462BA53DEC1B7167D795E07FFC1031F1CFFCD991F059711D7FA5D70F302BE4FEAB4618F5B
F7E063FC57DEB32EF8036FE03B7E48DC098F1FFFAB5CE127BC86A7A65F5D377C365D685ED107B8E307E04F1ED1FD86E1A737
0A3F578EDAD0233D3DDD70C52F4DBFCA4778983AA9E92F5959591607EBF75FE9F757593FBAC287E03FBC274FF820BCEF73FF
89F18203FA8147DE379EF951F206EFC1075CCCC7546FC435F8034F1E0684AEBDCFED2BE60F9A47F8693CD7BCC277F48FCE3F
1877EA2BB90FBF601EC69FB06E430FB75EAD6E2BCFF0BBD1D967BA4B3EC82FFC8257C64FE525F95B7D4B6042464686E184FF
914FF2883FE1E35D9635FD45FEDEB5A4A9BDC4437E899F3C83237C044F78429EE0CBD6F4C3A3376E4C7578111E9EE59539E7
C82D72FD7882DF5FDDBE8A9EC0099CD7B63FD3DDF17FCD33F94537A63BCD2BF9242EE2C147C79E0B4C70FB2438C01BD6C57A
363CB0FE6AF123F88A6E0C2FC5091CD02D7E147DE4C478B9863EBAB63030D0A33FB98ADFCA95D7A1BF4DAB4E7D25F1C3037C
E19335A7BFDCBF7FBFE98979187FED9680DF0477F0C6271817DE33FE6CAFE04582173E8DEFCF490C9E29E3D863F57BE395E2
0CBF4C8FCA0FF09AD53D78A6F000DF659EF9C382FB676666DAB8F8BDE942F10117F433B06E6DA15BAF330282FB8BAFE143CC
830EE177FF5E27DB899FE08FF8FC6CAF33EBEAEAEACCAFC86BCACB4DEDE5FD9245AD82EFDFBA75ABAD0F1FE57E7CFC8B7D67
E6CA38F0047DC017FC0BFEF5BAE1D83EA7DEA9BFB0AE798F9F997BE0C001C3091DAC1F73FA4BC70F7A06FCE6E629FC651DD6
57259EE926EB99D61C7CBFAC07FF204F1F1E5E7FB53BAFE80FBEC13378BD387BE376C90FF1E37FE8185D4BDEA2A30B5A78A8
7C818F56E735EE45DDCF749378D11BEB623DBE1D029E119EC063FC063DA0FF8FA69D5A22EB0577F0C627CC4F9497E4B343E7
C3A3A5BF7A63F1C15F04978EB38EDC22389147FC74F2CAD3DF483EC111FCF0D1A9ED376E9775E143E47B6380CF77927FFC82
BA0BBFC1131C5B74E7D1654BBFA63A553DFD1F1DE97DE869D9471B8748DD0527FCD8F4A8BA412FE16BD74C171FC29FC206AE
BD4FAEE819DC964E38FD8DE065BAD0F9D1073E849F85C79F6C2779E731BEC3F3A623D50F7A345D2B7FE16D8B7E3DAF432FD4
2DF0EBB4F3487F773F812F4F987A7A627D7DBDD788AD01CFB87D1E3D4D0ADC38C45D47C0D3F4A5BA62BDF80571C77F7B6C9F
5CE127BC444FF0197DC06BF3018D075DC01F7C051E991E150774497EF18FC045A796083F575E16D420E7D8F000FFC64FFB8F
5E335DAEEF3FB8BA44AEB9F907DF927A8DBE73A61EFC252C2CCBF24C7EA9ABE0840EE03F3E4A9EC88FDF0EFF9E9207F48C5F
9AEE546FE47BC8309FEFDCFE8F6FE1EFE8B06FE7E3D582333A87BF2B179F9E28F934FF54FCC19D7CE10BE4CDFC4AF180FFE8
0B1CD019FEC57A59273E851F0D1E79728A8C33E1625083DB0FE01D7CC57FC77C786A87F0D67C46790A3FE139FC267F21D7AD
2BDDB831C5FA0DF4BDCA37A883BC1F0A0FF069F2491EF151EA167877BBADE90189037EE0FB19CF1EE92F57F30DCD17F937DC
5507E6E7CFF9F774FB3FBA34DE6ABEE12F3881F39A0DA77688BEC833F9457F9BBF5B7D97C48DCF451C5A5D2257D3BBE240FE
899738F1FB7123823AB87D1BDCE01371B07EF22F3A0E0A6AC121B0B35F3FD12DBC47DFE00EDEE089BFE09B51E9C7ABE51C10
5D8203F98737F8CE88F74F5D29E3703F3E11EC73728AF01C5D320FE3CFF60A1925F9A76E5087D63EB6618DF007DEE0677342
433E90E7791D7589F9D015F31A0F1567F808AEF0067C670D0CF94078407D605DF33B85BC229FEB603EEA12BED6EFFAD577C9
B56FCC2A67FF35E8FB75A56EDFC357D02D7FE7F9D75F3DF896EC23DAA534E409FFE917666484BC2238E27FAC0B3EC1177882
BF53BF667B9D8A696C6C345F406F03469D78595E877F824FD299A607E47DF4852742BCE47332C4491D605CEAD3173F9FF293
F1E11BFA8377F81AEB9A7795E7F5F829BC7FE7DD635D844FE087FE7C234E5D29B8309FF59D3F9C1A7EF0E041AFD1FB37AC71
EB03DD4CBF2BC44BD6819F91BF8FD2FDFAB971C00FE0333C461F61A3D62E147DE16BE89B3E009CDA7F73A887AC7BF1DC4D71
EEBE1EFF4637F898E479CB9632E3393C85EFD60769FE165D7B6AB8E40DFD1307EBF7890AAC107EA213FC8F75B23EEA2ABAC5
BF860D387957434383E5177F9FFCC3A9EB25CFF007DEE07FE0847EC00B7CC115BF9F3A73539CAC0F5F049FA0E4F55D245FF8
147D498BDE3D3ED1D2B7AA6FA85ED1A9F9A7F2023E987E753C745C18D570B9F421F85C41D7863CB92EBB675327E927C08FFA
60FEA0BA448FF0043F862FF80B782EBDF6D4F582A3E94BD785CE4C5FAA2BF48B5FE2AF61EB4EBC2CF96DD1BBE775E88BF5A1
03F8CFFDF823E35077E1C75B2F1C69EFEE93A813137E3DB952F8607A541DB2DE49C73675121CF027701ED937B0C25DBFD02B
BE437C71B1C7BAC8157EC24B74697EA271A20B74475ED19FE94275863EE01BBE05EF06AC5BBBD0EDE7F8377ECDDF791EDEA0
EFEC1B0FF5902BF9C5CF02C24EDE25F1AC880E2E91F78FE01375C9FC4AF5812EC0053CE817F069F247DE4CBFAA5BF0F17B24
E05BC903BE804F0F7D7B7D17777CF8227A86BF2BAF3CB952F25C3AE980BFE805DDA017EA0A3A7F6FD6B13A59073A047F7037
7F53FCD005F9C247C89BE951714697F82071B07E7C0D9F1A32FAC435C223788CBF8FBDFAC437929F897D824B44C7F807FC8D
485EF3ACC489EFD3E7998F29DFE1B9E94575020EE409BEC0934DAD7C1F0B0BDB6AFD1BFEB1EAF9E0A7E5733EF083FA409EC9
2F7E4C1D851FDDA71EF577E3477EC111BE5287D2B61F692F57E3A7E61B9EC21FF4048F86D7077CEBAE47E81BFCE0D1EA9F4E
7C23F3924FFC92BC821738A17BF20B3EE0D2B66FC3E5C2F797961DF097FAF68A57DD70F12DF240FCD49F71E5C14F8BBEE1A3
E1AC7C8007E0C6F88C4B5F663EA671132F388AAFB8F10C581EE0F8117AC36FF03BFC3772CEB13AB9C20BF00117F88A0F8E7C
F0446B5927F7E3371BD79EB846E29EED15B651F0A35E512FE7F886ADF17C7E4D1F6B3DC43F5817EB5993BC719EF37933E51F
3E3A6B44D81AD121758679CC379467E0CEFAD033EB84E7F00EBED35791F7D79E6D7C48709DDF3A6C927C9E9075533FF1E5F7
DF5EF3ACDBBF3F68F67DCCEDDBF81CAFEB53E01DEA7EFD8CCAB04982177E4C7CAC8775D04F5187A8BFB3BD4E8C3D72E488F9
11F333AFF98CF2111E723F7574E6BD2706C8380B8BC25E92CFAF1227F50D5F43FF03971CDB277AA12E8067E2E2A3FE72C567
59EFDC26CFF8CC479D98FEF289C7E579E3ADEA1EFE4E6F13F692E41DFF234FE08AEE7D0E9D682DE3520FD0E1BB614D03C5B7
E0333C467F63DEDE384FCEABD0197A1DF6698097DBE7F19BD0C1EBCE8B9FE29BF8CAE2BEA127DDFB18EA06FD0DF8759878E8
4D5917FE489E169678F21079765599F815BE1FDDCE3B54AED4417CBDF88FBAE17245AFF8AAE43F25A5A24547CA77F444DEA9
1353061EAF9279C923FA219FE81F7FFCB8FDF1AE729FCF3F820A4417E8139FC6A7C81F79230FC44F7FD1E20F1EFFC027A646
859E942B7E0CFEF010FEE1DFE63FEA13F8037C8007D431D3B5CE8FBE379CF5FBAFE0876FD2B7A11FF8024FCC9754C7E877E9
E1D039D227C103EA1BF8E1BB793B1B1F922BFE052F96847A70810FD43178010FA903F0D1F4ABF1A163F48ECEF18D163FF13C
46A7E64BAA3F74076FA85B13E71CEFEFD6357A661E7C9DBA113AE5D83EE1198FF1659E9FFC40E81CC10DDF027FFA0C78D325
FE70B9BB9F63DDE8C8F4A8F9409723AB830ADCF5D4FA44E539FCC627F035E2881DD03450AEF01F9DA103F32BE52FBC5D3122
F44AE1033CA0FE94DFB9EF56C1157DA3EB81C1EBCEBBEB06F5847CE17FFED9C7BB4AFCF0161F85BFE81D9D830F758279181F
5EE22B5B571C7A53AEACE3BD5B5795B9EB18B880077D103E02AF57CCF2F0C4F6A38A133CA76E906FF24C3ED00D7A59BF6EC3
69E7FB45EA3BD40DF2827F919F0F6FF5FBAFBB3EE2EBC4F3EEB295F7B9E332FD2A8FD0B1D50BE511FCA12EE22FFDB29BE639
FDBFEA847A3376CBB15C8977A26FE8956E5F81EFF82CBE3674D2B10E721F7C421FE8029F253EE242A7E813FCC3FDD686481E
A95BF88AF7B1904FE5F324F092BA41BEE1057C0057F003DFE057D637444515585F8A3F911FF2827FA30BF28D3EA8DFF0A947
D3910CD167CCCDDEDD64BDD4C5D8492B9DF7F5E02BF52FB5D7E172B91A4F54AFF0057EC0D3D53D3DB87C326AC36977FDC227
CC971437F0623FF0FCC97DB74ABFF6EFF1FB4E0B6F799EC7F4F5E00BAEF80DB8E0EBE063FC81378A277D2AE3B6595D572D3A
1C7F4FC8A7E2D3F0147ECDF6DAFC84E4977A45BD24CFE4973A095FF1C191A78E5D21F3CE59B9F97EE109F7533F590FEBA06F
4527E4957C8223F50D3C674DD87CBFF835F587F5C22BF1B9A8A87CE397FFF78117E5F5A67FF53FFC9ABAB2B97DD33C779DC1
4783961DEB20F8CF7F30E257F97C2CF1503FE13BFC8257F006DEC31FD391F21A3DAD19B8A9CCF95CBAEA8ABA81CFB15ED6C9
FAF02FD639E350C4AFB22EEA0C79A2BF24EFAF6F6B7C5470A52E52DF677B35D51D3B76CC1E53A767B66BDADED4D4E4B53026
A2403E3F40FCD44FEA44FFE16B43DC7589C7F828CF0F2E5EDFE0AE2FF82EFA81A7F0937AC07AE716346D97751217F1D03FB2
6EEAD9F4594DBEB27E7C90F999173FC577062536DD77E8D021AFE9BD230A0407FC95BC528FC027E1DE2319A233F0C03F7CEE
F0E882F87B9FF7EEE6AE27A623D50F3A5EFC7CC46479FF85FC5267F053F2B1704E93AFE4019EE367F07D4CE8A632399F4667
E89E3A849EFB5C71F426592FFD0675759BF7BED372258FD48929794D1D259F1F5F1778D15D5FF045E2EC357EF96C77BC219F
FBC6BACF99F033FA36F0EBF4978307A41EB4F889C79FF015F4829ED08DF998FA04FEB0B9604D0799877A86AF102775AC686C
5DB55CCDFF55B7E038B53C62B2CC8FDFC30B744F7D113E6466EE317FC4C7873FD5F4988C637EA2F1E12BE66FAA63F4BBAE64
D335F27A7C817A056FE12B75033F051F70599A1E71A3F0163E5157C93F79A75FC3BFE0DD92D79A3A0ADFAC8E28BFE095F983
C6874FC03778469D379F517FC17F5A7C491FABDE37740CF84DF2838FD3D79AEFA95ED129B8E2D3B9FF6E7C54AEE60FEA0BCC
0FBFA9A39F5E73B44A9E873FD4617864BEA4EB438FF09CBA04DF27778EB851F888CFC123D3BBE60FDDC73DB0225BE6A18F89
F75E3EDBD9DF699DC2F743BA37DD27F9BCB48EF13C3C847FF88DF99BF2DA7C0EBDABCEC9D38A7F87AF92F76BE00D758C3E0C
DEBCFDCAA11FA55F18F5F2A66BA40F828FF890F159750CAF89071FF5DBE359AFF99EE6171EA24F78BAFCC7A355F23CBE4C3E
A27F3B7A935C2B7CF78E917CE273F81B57FC001F20FFF423E0000FF15DF83868816FACBBEE51DF4C778A33FAA30E711FAF87
AFF84AC62F070FC8D5EA8BFA11F961BE7EEDD67470D75B70010FFA38F84DBD19D7FF685BC1973A46FEC8DBFAC7835F105EE2
3BD419F40FCEE06B3EA97C87E7E4155DA247EA2EBED0FF8A237BC4CFD0233A848FC4499D21DE89DBC25749DDC39FD0D98799
9EF75DC937F50A5FC6373F9C79F449379FE019BCB27AA7E3318E774CF865D2CFA307EA5258B8CF2889933A877F988E9417E8
8978DE39BC22DB5D3FC90F79C1EFC91F7C8127F0037DAC0A3FDA56D6BF71A1FF71791D7D32FE64BEA4F9236FF407F0A6D7E0
C333847FF0011DC38BA89A55BB253EEA2BBA816FF00C3C677B45FD26F153F7A8DFE0021EF8D69C2551DF783E47E9791DF513
9D5077938B0FFD28D74FF605BF20FD23B8E143B3A6447D23EBA76E303FFB25DBA7DDBE7F84D445F66FCFFE51D34BF62FFC7D
7CEFF0CBA4CEA0037887AEF0D9518947164ABC8CCFB8EC87E6DF1EB54BF8433CD44FE67D66C4DE31B27F627E78459D815FE0
0D5FC19DBE9FF9DABD51774C743DE374D42E99873A403EA8672DFB81C37F9C3871A2C5CF35AFE4139EC00FEA3CF7535F670E
397C5EC6A12EE1F71BE67A70361DA93ED053EB1FAAFF258FFF75FDEE50C93BF9C0DF192FFCB7237BE4BA3030CA47F2461EA9
ABE481F8D937ACB927ACBFF3797ED52D75C9FF8B8DDDC4AFF13BFC1E9ED317C077E2C197898BBA455EE76EF4E4C1F4F8673D
898F2F355DE2BBC4475CE4917A3C3DEF70AE8C83EED101FC373F549EC2CFE91F45F9489CF834F8D37F9397F69F356C135CF1
41E258F8DAE1DC93274F7A2D7E28EA2DC181FC52DFA86703FAFB8C72D74374811EF003F2856FADEB7E64A1CCC3FDF83DE3E0
D7F8E007F547AE93FB864CF73FEEAEA7F83DF9A27E4CFDDBE109C78F1FB77889D3F61DE6331E7FC26FF07DE6617CF331F515
FCC474ADE3A1EF31C7C2FA4B3D444FF807F517FCE3830FCF90EBD4C6A8B7A45FC0A7C1CD7C4C758C7EC957DF2757ED76D749
E3B3FAB79D7F69DEF10BF26FBA505DA20FEA343AEF3BE2F018C161CB1DCBEF131DD037D1AFE283F8F827771C4E141CE8BBE8
0B4A4FEF1F2157D3AFE60F1D9B9FA88FE04BF81D3C5D7CE0F004E1E7D2A0C826C927BCA4DEB6F89BE77EF4F8F1BC8DDD641F
03CED40DFA46D6BBD3ABA6975CC92FF59A3C9BAFAACED0D7BA2F4303252EFC85FAB929C32F5B7C807A8CAFD34FC3CBCE5107
6F759F0FC27BEAA4F983FA02798AC85B77AB8C47FF816F9177FA88829CBA6372B5BAA63C859FE887FAFDE98A437F888EF01B
EAAFF03637777F8BEF697EF10DABCB9A6FF26C7EA2F8E32BF00D9ED107911FF2429FBB676A95F33D3370ABFA6E77A85C270F
8D6C125DE28BE8097FC017C89FF9A4EA0FDD051604F5143E510FD80FC073F88DCF997FE29BCACF1577458E723E77AFBCA4EE
8107BE9EFD75C336B952E7A82B9B1EF1E08F6ED1C1F28987FE90F1A91BDCC7EBE1037D04BC80C73D4F2D77CE6B7A8C58EAF0
1B7E52AFE0694C8D770FB98F3E0FDDB33E7C7DFDA9C389B22ED3A3E28C2E4D17EA37E863D457A181D287C217FC12FE513FC6
7B1D7270337E2A2FC19375E06BACC77C5579011FE877E153B7CA83C5A2377C99BC46F91E1E23D78927224789BED10DFC329E
A8EFC217FC033FC1F7C00F7E80A3D50BCD0779E03E748C7EC1853E057C7C6B437648DEF11DEA00BE8C8F7CB4F050F6D1A347
BD3EE8E497EDAEE3D455EA2E3A1CD0EDD04E77FDE53E5E0F8FF189F4BE076FC527DCFB73F249FEC93B7D25EB78FFB575B7BA
FB02EFC59BF74B5F877EA857E81E1F20EFA63BE523FAFBA85D504F777DA70E905774891E59077580F5C03374E9FDC72107CF
D95E71B749FEA9DBF40D73BE8EBB46D6638FB55FC01F5917EB311F53DEC13774811ED0C1ACCFE2AE713EBFAF7589F9E7B78A
3D2179633EEA2DF8822B75069EC00F7C2EB46EFD6CC90B750B5F419FF8E0E81907F3DDF580FA3BDBABF1D0A953A7BC6634C7
9E90F9F15DD6CBEBA86733673756C9EB8DCFEA1FF0DAFC50F543FEAD1F7AD2BB87BB6E8F9FBA79BF9CBFA233741CF4CB8657
5352765A7F8EFF2DF48D8D77EFDBA857F002BEC3737C9FB8E7FEA7B1EAF4E9D3D63FC1EF770E1E3C22791FF144C80EE9D3C1
1FFF237EEAD3E797376E72E283FFCA6B7480DEA8C749AD0F16CB35B26C4D94E0431F804F4C1F111B2FFCC317C181BA81AF07
4E3F942DF9258FD40BF2B9FABBCDCF39DF47531DE2EB093F2E1D2DEF5BD22726F82EDDED9CA7E9B8F832E3B35F65BFFCF2FE
FD8BA4DEB3FF7E3EB3A65EFA77FE8E3F923FF2869F82C382DF1B3749FE199F71D917323EFB6DE6C1F7C9475892471F8BEF89
FD50E68707D437AB0BCA2F7805EEE04D5F838EA8B7E8C9FF914DC79DDF6953BFA69EB4F886C757F00FF841FD9BBAA8B1EFFF
F21FF5037C009DA24FFC837D207979EDAA3A3FF165F6DFAD53F7F84ADCE409FF264FE4C77C46D78BDF80377EB976FEC17C59
27BA93FA999CDCA23FD6C33AD8B799BFA97FE01BF09C7E07BEB37ECE138803BDA21BF432F558EC87926F7C177D502F6DFF79
D5C148F145F32BF529FCADC5F7F4B1FA80F989E6095FC1BFE0EFE2698D7DE5F5E40B5F236F4B7D627E73BEBFAEFCA3BE99AF
AABED1357E828FB0DEB1D76F7E4EFA20F4842FB20F21EF1D5E686825389BDED587D03DBE4F9DF8E462E377EE3E98F8D1635C
DE8A0CC1873E129FA09F18F8E8FAD9EE7E64E85F36BCEAEE37A8DBF09EFA37A955C32989F3A9C55503655D7F3FB02B5FD609
FE8C4B3D677CEA2EF3303E7AA4FEA04BF2435ED897983F28CEF804FD0FBA89FB9F8347DCE7BCF8C27B1507273AEFD7A06BD5
33B8C143F8477D5EF76EC474990FFFA54EE247E00CBEE6FFCA77784EDE9EEC52315AE2257F9347C6FC263AC3E7D007797CEF
C53551EEFEE8E3FA4DC7651F495EA96FF01B1EC073F8403F0B2FE87BE9CB4A1EDBBF48AEE801DD2C6B6E3875E6CC19EBFBB8
AFBC634DBD5C8D9FCA4BFC6FC58D31730417F4435DB5BAAC71B07EE2A41E12EFA64702EF963E80BE015F675F013FBA0C6C70
3E0780DF50AF850F85852D7C8007F465A677E517BA0F3FEDFB8CCC431F86FF912FFA94BC9FEBFCE48A5ED009FC423FD483F1
A10D43E4EFD43FEAF4C69F0F463AFB297C5279011F781D7EBFFEB2C6EF647CF240FCF4FB81C382CF4A1CD439F651F4F3E059
D9658FAFF37928AD33AC8BF5C07FFC091D98AF2ACFE009EBC34F59E7A8D888E9D25F833FBE38F19B9839E2A3E8091D905FEA
C9D6ED0DADE40A1FA86FF0C2FAD1D75664B8FBA9D95E09230557EA1BF57BCE92840F04777BAC7515DCE923C0DFF8AC7C81D7
F8327E316C79E355127FF4D6D54FCB7DB67F54FD595DD07C9007FC0B1C229FF1F074D694840F9CDF21D3FA431CDE5F467FEB
7CBE407545DD987F7B421B5927F1500FC1037D7B0F6C182278D137C39B1EE31BDB884E8CB7CA5778ED1B1EBE59F2890FE2DF
C637F57578475D476F03D7367E2679AA49AD182D7D2B7CDCFBB75DF972A59E51DF677BD52D90759A4F2A5FE009AFA3EECEDC
58374D5E3FE374421B891F3F267FF0A3C7354B47BB7942DD0007F28F9FB2AEB943EAA69D3D7BF6CFFA937085E49FFC5297C8
2F7A20CFF825FE89DFE20BF8E69876F5A9926770010FFA5CE2A46E7DDEBEAEABAC8BF1F0277CC9EA7B7DE0DDEE3A0F3FE9EF
E029AFA39EF37A7488BFA686367C9393B3CFF80FEFA9031F05049F95F315E6A75E316EFFDB7D9F71F70BF0047E5087A77F94
7085E80C5F044F7C0B1F23BFE3BF8AFE56EA2D3AC357AC6EEABA580FFE05CE0BB2EBBA0ABEA62FD53D3A23AFF813BED4E227
1EBFC1578893FA46BCE64BAA6F746DFA52DF45678B1FDA9222FD223CA35EC103EACDD41307FE107E980FE8BAF003F313D507
BA587D3CDAF91E2E3E82DF8FE81FBE59F65BE40FFFB67AA4FE1EF069E35532AFF9B9E260FB36E539FCC657A8338CC738A66B
D5333ED0E20FFA58F518DAC6FF88E0417F807FD3FFA19FDEC31BBE179E05ADD9D89891B1C7F647D413D3B5E60F7D4F6DDC92
2238E043F011DFC2C7E0E596D5CB6F75AEECD375DF4EBEF0E33505F5A91287E95AF54CDCF8083C5D7CD7813F9A9B9B5BEAA9
E2032EE67BAA1BF4024FD131FA85977D3AAC7EDADD7F5187C0277474E36772353FD4BCA3177C817E2BC1BBB18D5C97066DE9
25FD097CA57E5A3D527F240FF0987A3369D5011FC90B38821F7D9CDFF9B02219973A44FDD97C7C5DA5E04FDF86CF998E1467
F484DF1337F19ABE5457E003BEF4D7E04C1F807F0DFD67FD6FE28FE497FA439EF17D7C73C4E1BA0BA213CE0B38576993B2BF
9BB34F44A7EA4FE895F39B7FAFAEBD5EE2E33EF34FE52FBC35BD683F856E9897F93827309F54BCE1A17DBE43CF77987FECDD
D1DDA46F85A7F8B1E9EECF7E212363B7E96FF2D02DBDC487F00F74C63E9975BCDE70E05BA98B9C973C33AFFA398993F9E113
BA5936EE808FAC13BEC353F8C9F88CCBBED3EA88C6C1FAE1157D22FC623D9C83B0AEFD15E51F489F43BF835FC013F841FD34
BD284FD18DE942F5001F57DC15BF57D68BAEA857EBDE8C5A25F7E3E3D413D38BD61974030FE11F3E8D0EF1EB09ED0FB49275
D0CFA1C3F7FFDDF03789C77C52F10447F669E4B7D36575179DEF7F69DF85FE62FC3CFA1E3A6663A39C3BD28FD047981F2AFF
C0D7FA0DBDEFD2BE676091FF1177FF11EBEB1D26B8D9BE58753CDB2BD9F9BE3F75923A3D6765F271C1C51E6B5DE5FCE3A935
9545F218FC998FFACFBCE481F8D9B7D9BE65ECF25BDDFDE9F09BC38AE45C81F8A8C3D449EA89EF037517242FB326241F17DE
51078883D7E1B3BC9EFAC13C8C3FFFC1E422592771520FC11DBCE9D7A813AC2B6877FD6F4E3F70227EAFD41FF484BEA91FD4
F1D95EB53F098EE493F311F2CAEBA88B33736B4FC9FDC67FD50D3A30FEAB5FA203F0E877F7BA4A775F36E350B2330F7E49FE
A803C4475CE8023DA04B7C853A287C2E2D3D62BE49DC735FAF3D75EEDC39C303DF5A997DA0953C0F3FE9D3E169706A50B6D4
6FFA15FC9BFD1B7C7AFB8EFAA9522F17C6242F71BE7FA3F851F7C823F5ECF3FFD4E6C8BA466545AD927D15BCC0BFBD17C777
74BE8FA3FE42FD61DF415F5CD4757F37B99A9F2BFEE01EF8F7D072C933759A7D23F50F1FFCD8BBEE3AF121FA7AE6D936B6F6
7AB9A2936E6F2CDEECD60BB8526FC1D7EA82FA0579B63AA571132FF9248FECC3C2DFF0F395F8E82FA903C44F9F95FBF0816F
E58A4F81FF82D9B539E7CF9F77F989C76FF095E9BD9397387EAF3E0A1FF15F7C79CCAFFB2778CEC1D5AF54F7E89D73C6C7AB
7614C863F464BAD6F9D1373CA0AE7CF666ED67922F7C191D446C6BF89BE0429F0B3E155F543FE7F457CF27B795BCC133EA98
EFBAC842998F3A44FD31BF52FDA13BFA27FC65505DDD1ECFE75DD50754FFF8468B9FE863D517BC824FD437D391F623E809FC
E817322ED65D94ABF15FFD061D98BE34AFE8CCF4A5BA62BDE3A7C67794FE0BDEE18BE882FA8C3EC0F1B167CA3F101E8027FB
1978D9EB96FA81CE79053EA3F1C337F337E517BC323E6B3F08AFD109FCFD3AABF6B30B172EB4F050F987FFA1B777477887B9
FBDCA9E5C96D65BDF801BCA6CEC32F7805EEE04D7F14B56BED07B22EFA787CD1EA8BE6957CC2637C79F2F5B5AD1C9F539EC3
6FF059FD5D5C90F0165FA36E24F6583A58BE77C2FE2D2967F166F91D18F24FDED9972C4D4FAA95E7E13D75C9F8AF75091D50
2FA9C7FEE3EBAE73FA6B74A4FC424FBC8EBAB1E6CAFD130467FA4A74B9674C65915CA957CCC3F8DC8F0F328EF995F206DC3F
DA1B5A2EE794F45797F65B83BB0465CB393EFD16FAA11E269FAA9F2A57EB67743CC631DF53BEC013F884DFC32BE6A3CF625E
EA03F5827AC4159FC69FC197FD1C388F181559289F5363BDD44378C83E003EE2CBF8CBC89C7D5B9CCF7B534F753CC6B1FE30
C7CFD7DD278207FA5BFA5D6DAB6FBEF9C66B72E7A45AF145F86E7A529EC36F788D5F837FC8A3757BE48A4FE3DBE48B7CE07F
F89EE942FD0F7D18CF95A7F09DFD007AE873BA2ED5F37E6E7AB8BC0EBFA76E113F759B3CD007A0FF0F7B1CB8561ECFF14DF7
9379B99F3A16D42F244FD6CBBE125F278FF0987CAEF877D2C0C3870F9B3EF1F5B1D7C705399F9B53DCF0697C04BF9B5057E3
2D8F678D48F7FC3E9FFA2BF199EF296EE0E5D7B0F9A0E3075AFFF07BFC0FDCC06B7EEBF4E9B24EE2A4CEE0EFD4BBD95ED5AF
497CF010FEE1D321651B16BBDF3FC1EFB99FFA34B3BEFA191917DDD2BFC69FAB1F2857F00557F60156D7347FE40D1EE107F8
00FB627C6BDFBF7614C815DF21BEB98F543F2378CDA84C9F2EF1E313E41B9EC00FEA7CFCD095D7C8E72BED7C46F58D2FE293
E881F8A9079F6756FF55F2B2EEDDD826E12DF500FF5E5894FE9AE77BA81E3CF06978DCF7E9B51FB8FB50E38FF206BD1A0F55
4FF0D1F8FC675FE0E6750BFF3D7A4207F863F784A583DD7DB1D50BF5779F37F66DF1F046F5A8FC85B713B7250D74FA4DC513
BFA16E301EE3A02B705BD0B1FAAFDF7DF79DCD879FAD7CB4C65B7034FE6B3D44079B5F58DF4EF8483F8DBF1AFF356E74009F
E9E3E0F5F436E9AF493ED01F7C319D2A8FE00F3C61BFD5720EA0BA503DA02FE2C147880BDEE0BF9FF9EF2990C7F0803A061F
5AF4EB19179ED3EFA2ABFE7E07968BFF99FF2B3EE0C2B918E775EDC2F6A5489C3CCF63CEB738177C697E8D73EEC27DE812DD
A3F7C57DD3CEC93AE0377EEF1D93B8DDF99E9FFA26F5817A035F36241CB856AEF6B9333D1F647ED39DEA8D7C9BEEB4FEA03F
D39DF202FD71AE449CED0BF79F17BED067E2AFD1B7D7A5CA151FC71FC726EEFDDEF9FC04FEA378C067C6675CCE6FCC979477
F08D3E003F1ABEB4B65AAE9C533EF769F554B98F7C0C7F78F341395FA42FA2EF800FF080BED8FAB1CD2179F2FB47F461D40D
7020FFD6CFE83C8C0F8FD1F1D71FEF29F8FEFBEF6D5CFA2CC6379DEABE01BDE2DFF8393E8B1F755DB5A8D6ED4BA617D509F8
931FCE4FC9D3A8D8D826393F200EEAF3D4A8B473CEEF5AA93ED03DE720E0D4B9F6C06AE18BD5771D8F71D02DBE3C79F69E0F
9DF75BE0BFF20E1D9097412F6E58EC3E5F345DA8AFA30FF806CFA803E69FCA33F20F7FE00DF59CB88887F389A587D316092F
D027F584F5D167B14EF349CD1F7943F7D47DD1CBF6EDC7BDC6F74EDCEEFC1F24CD1FF5212664F501C19FFD183ECA39E8BFA6
544D9479C033F0CE881382277D0FFB5EFCB157D0CA6BDCFDFE6CAFAC59927FFC9BBA65754AD7C57A7CC3636F93B8A9D3D499
39A1599F7ACE0F3CF7534FAC4EA98FAFEEBEF7FB6FBFFDD6EA1B7FE779FA517438F89E7DCDD277187F9437F03A78CAA660A9
DFF4A3D42BF6EDE0D7ADFA4098F3798F81599F0AAEF80A71E35FE83B7CF881E57265DDE878E97FF67CF8C30F3F987F515F66
7B557E2C8FC1857355F099DF29AB87BBBFC65FB91F5F9EF953653F79BCE69EC44704177C1F9F361EAA6EE0A3F5F9DBD7B773
F7FBD409F2ED37B9B65AFCD77C49F5405C931F485B24EB067F744F3EC9A39D87A8AEC8C7DCAB2BFBC96370043FFA47F8CC3E
0F5ECFC8C8EAE1FC9E8DEA067CC80BFEF5F985CAD6F298FD26FBAA8201FB52E48ABFA0FF8963771F73FB34F9200FD427DB4F
AEDADBD6F97C19BEA7F9200FA62FEDB3D0590BFF3D7A4207A617ADB7E826E99BC5EB655FCC7947CAB845B572459FF8FAFFD1
A9F217DE2E3C91759BF0163EE1A3EC3BC953E98C1A673F023FE89BE04958696007C91FFB16EA0FF9659F91DD76FF79B9AE18
91769DF37E9AFA207E6C3CD77CC077FA027486BE4CBFDA8FA363D385EA017DC15778B7E0A6CAD6C25BD38BD60774D3A253CF
FDF0D37C407907DFEA33CA5E907D13FB27EA2EFB6D78FF4ED281BFCB7CF00D9EE14FF00A3E5197C09D7311F0E7FCF849FF9D
7B8577F807BCC7EF3E3BBBEBAC7B1F891EC8D3F4BBB26EF3FCDE9E87BFE809FCE843D39F38B05AAEECCFE0CB8EE5D553E56A
BAD67C9167E39BE20CEFC80F7961FF6A3C51DF822FD66FF48FBD4DCE135BCE6F3C8FE9A3789E3E0F1FF968444DBDF379AE63
898FC8E7F4198FFAB978EED602E7FB14AA07FCD2FAC0931127E47D08FA3C7C1F9F083EB6AF59AEE64B9A57F2495D647EE6B5
7AA17EBA62F5EE633FFEF8A3D537FECEF3E88EFD16FA637DF489ACD3F8A03C007FF8C3F9353C9AE89B769DF813EBC5D7591F
7AFAFAD65D677FFAE927AB13AC77DDB8BD6DE56A7541C7631C7CA4F7A0D507DCFB30EB6F976C0A96DF79A49F33FD2A9F9817
7CC195BE0CBED23725FDE340985CA36EF7A9923CB18FA59ECCF6CA1B277A43DFF8A6F141F90B2FA803C44DBC7312F33EF47C
1FC5733F7E860ED1F9E4E65DCBDDF5069D8F1BB2BBBB5CD967A287F7C6EC7FCA791F7FE6D602E95BC92B7A255FF495E46DDD
97099FC87DF401D4137C14FEC2DB849B970F91DF69E05C06BF41CFF8E16CAF1D93E43AAB7BDE879277F84EFEFCCEC73C2D8F
E9F3A87BDEC752A7491EF17D7C9771F1A35937EF182EF9355FD2B889977D0C3E54B5A46AA25CE70FCB6B277925EFF80DBC21
8E39E7770C173E192E8A07FA5876CFD6019217F48F0F810B78B0FF0347F0C3D783DA87D7171434D87911F5AA85271EDEC117
E2479F331ED9F1AAE8169D50E7D18BF1567120FFF435D461EAEF8C80BC76CEE73735AFE0469DA7EE733FB88007B8C3B3012F
057670EF632CFF1A1F38B4F0C7733FF9876F9CF7C03BD3A9EE6BD1AB9DA3DDB878BD7BFF4A1F86FE3F99553D42AEF09F3A6F
E7E6F050EB8C9DF7E16FEA2F4B5BED5AFEF3CF3F9BBFA3FB40DF9A7AB9E22FECC7E2BA1DF8BB5CC937F593BA897F330FE32F
6A9577B9EC53E02B7E63FA52DCC0CBF8A9FE044FE10FBC0137D3A3F695E89238D1C7FCFA1DAFCAFA785FF3D1EFB6458ACF50
1FF077E258B5717777F7FDE0CA38F083FE119ED02FD03FC03F788F0EA7BD54FEBBE4057CC195FE28644EF044C92FFB37EAA1
E553F9485ED91FA2DB013FD52639E778EACBE876E277BBEE902BBE4EDCC43BE9D8D601CEEF8D695D333F569FE2EF97D6BFF1
F7A44E937A4B5D357FF8737F5050506F3E015FA863D42FEEA77F601CFC1B9EC2CF69CD79973BBFEFA8EB42F7B3BD8AC63A9F
4B527EE123738A8A867ADEB7D0C7AA63F8C7788C035EE004EFACBFFA2AE113390FA5CF821FF8C06CAFB24FDDBE025F673D5E
36DCF93CADDE4F7FC7388697E2847FCF6A5334D4F9BEB5AE97F8C01DBCD95FB6C4EBC907712FCECE8D93E7D1213A432FF07B
5169F9EFFFFDEF7F5BF2A9F3332FAF23EE39FBCA864B9CF31714B5753E6FA2F985A796078D9FFCE013CCCFBC2BA233EBE43E
EA137E66FDF0CD314FCBEFE0D0BFB6E4DF332EEB41EF8F742A7B41F6ABE89EFD17FD42E49AFD4F39F1294EE03DA357D96BB2
0EF63DD4E79A8F77EE95FF6768F16ABE899BF5D1B7B24EFA21743E7A6095BF1B17F40D3E6B06A654C8DFE927F07BE20417C3
89F5EBBAC937F968E1A5E7F5F83E7EB07E4EF508A967E63FBA1FC5875AC6F1DC4FDCC6CFCF8B9CB8E069DCD0556F890F724E
473D44F7D419EA8BD511E5C58AEB77DDF1CB2FBFB4C4ABF1301FFD11FD12FECB3AD035EBC17FD1D59477CA43641C3B0738ED
53E53E0F309EC3EFF8B2D78467F0937A054F591FF3B34EEE631CF206DF9987D799BFB5CF8D93BA86CFC11F78435F867F933F
F246BCC4899FD97E2A3BBC5E7E3F94FD896FED96B51207FD34757551A7C28B4EFCAA67784ADFD6E3FCF221EE731A8B57E344
271129FECDC223CE53E8838893FA42BC964FE5337965DF8BDF0DE9B5B74CAEB64E9D8FF5A267E29F36B3CCCBF97C2AFAD2BE
099D5167D007BAB071F47ED6451ED9AF904FDE1F7A2DA4B650E2E3FD1C9EB7F7CFF4EFC6DB8CC28B120F3CB0BADC27B3CE79
DF54EBAEF559CAA32571E521823BEF43305F87A47D6B9CF31D5D2FE3331F3AE53C12BD9267F2CB7E149DD0BFA3974573CBBC
DC3C0607E3B93ECFEB783FF095A5358364BDE463D947B91D9CCF07697D4267C44FBF401E88973879BF82BE079D8CBFAC629E
C41BF8F7E8979CFD9CEEDF38C7B47338ED3BC2AEAE4D72FA35CD173C0517D6CFFB8BC4819FD13F4B7DD9B9F394E912DCC0CB
74A77D3DFAB37C9247E5377E8DEED78656F94B7E2D6F5A4FC89FE1AEBC251EF2435EE83B789F803C7649D8B7DBF97D555D2F
78B23EEBDF42532AA4DED2272DB9BC305274C3BAF10D7B5E1FB35EDBEFCF0A9E289FD3623F866F703E847FB04ED6C7793CEF
73BE30BD7A90F0097C581FFD17EBC417A6C494463AFAD33C5127C8CBA713CB5BCBDF573D9FD9515E4F5F827F988F05E63ABF
F7866ED9BFA19F6153F74C77F4CDBCEACFCCCFBCCCC7DFF121F00777FC97FD233E0CCE537E2F8C747E2F4DD747FF423F431D
256E5ECFFDE493F747C9EB92DEA5916E1E80B3E94BF767E80C5CD93F81EF43DB4A264BBFFF3F0B8B7E90FD18FD9AD535D58D
F7C58A79BFFDF69B57F090B01CD93F70EE409FC279343CE8B17ADF5CC73F150FFC035C8CCFBA7EE2317CB5CE833338829FED
DBE8E79FD8B2567E57823E95BE125FC14FF017EACCC6F8BD6572B57C6A1ED10971110FE7CAE4853896AF286F2DF9311C954F
E0C9957CC007D64DDFCAFAD795A4DC24F7D377E24FBC2F405F5C77FDB648B9C2477003AF9597E514487CF838798D3EBCAEAF
E898F34AEA0F711017794999B1A8343474ABBD6FC0FB44F4D3E031E2A3CA5D82CFB28CC2B6526FC8233E030FC69567769473
37FC8BF78D9F9E52E5F41FF0DDD6AFFE4A1CD4BF778256BDE53EC7625EE6E33ECE3DC1AFF73BB5F1CEFB349A27F0247EE3B9
EEA7E13BF11227FD203EF1E991D2971CDD2BFE86A3F63BE069FB94BDD12FC9F766E90BD1E7FB5DFD9BDDE71CE40B1F256FCC
33A953615BE777E89567C44B9C9C2FDA3AD5DF6CBDF053C781A7CBEF2F7DC9CD03F0443FF862C0E23DD3E57AA9CED017F8F2
3E3438D31F90C73177EEDCE7FCBE86CE3FE1624E81F81DBCE11C91F3B0BC41B585720DFB72E33FE5FB7E9CD750AFC82F7D22
79E67ECE45B676DBB746AEAC8775509FFD1F898F973CD197E3AFD4A7F1EBB6073BBFA3A379351D69DEC8A3BD5EC7E53ED3A3
9E53A24BCBBFE61D7E833BF8811BBCA14F843FC4C9FE2FB5C7BEDD72E5DC0A7EF43B5633581E737E46BE8B27D70C723E471D
B9344CCE57ED7D43F5037082A7E06575FFE5949BE4FB0DF835BAA57EA2DF35C919EB9DEF49695D4737C4C73CF003FEE17FBE
8195BBE44A5CC4C3390D7C84BFF0D6F0D57E109CF107F0F41EB53DF8E2C58BB64EEA15EBC58FF167EABDED533686E5C8FF3F
A6DF62DDCCB3E6E4CE7DBFFFFEBBE1443F9B3076DF5CB99247EA06F95C39AA2055EA28FCC0F74C177AAE803E8883B8E019F8
728E0ECEBC0FFC8FA515DEE2CBE8D8E6D5F9F0EB55BE394FCAB931F5863C714E00BFB62FAC1E2457F693E034A86BB59FF3FB
A4AAA7090F94383E06DF83FA453D555272D8CE39F169E6651DCC4FDEE983C83FF8C323F8C33C13D615A44A3F05FF561E2E89
FFF5D75F2D5FE4DBDE3F98BFC8F97FB16FFDBEE027E77D338D8B78D8CF199F757C780DDEEC5B62E26AE3E50A2FE009F8595D
783B63BDECC7D01B79C35FC81FF1C13BC6854F1FD7C7C7CB7928F509DC791F17FCD73F9EDC467C92BE04DD18AF343FE4CBFC
6144CE93723FEBC6A7FA3CB0AEAF7B5F4BBF4EBE86CED83D41C68F5CE1B740F8CBF9057EC93CE0012F8D277FF6C56EBE18BE
9A1770605EE6A3DF24BFE8CFEAB0FAF898E26DBF3ABF27AAF9268FD41FF2495F4A9DDC7C4FCD60E7FB628CA37587F12CCFDA
4790EF2D57AF0C97DFF9E35C0FBD12077926EFD6FF7FB5F19FF2FFE9E9330EFD51F837E927E92B0F4E2FFA41AEF473C4FFF1
D45D27653DD45BEA2775932B7586FAC238DC7F69DF34F2A59DADE479E2663F844FEC5A525520575ECFFDDC479CE003AF42C6
842E9775B03FC1FFA863D435E261BDF82E7EBBFA9682BBA59F864FF80A3C647F0C1FE131F50D3E8337F52724A5DA4FAEA66B
EACA2BDB7EFDE38F3FECCAF3F0DBE71FE935CEF7BAB42EC20B70A5FE83AF9D0FBEB734CCBD4FB5B82EF16B78C5B994F14BF3
813FE14BD65F6D8E7A4A3E1F439D4307633F2AB9C1F139D599AD5FEB2771A0AB4FF625B791EFCFE047DC37F65CC1DD4E3FA8
FC0E1C1697223CA00F63DF4B7EF1F5A0F9BB27C8757560C90DCE7E48F3089E6B1FCBFE5CF67BC44F9DE2FDC9C7766CDF217D
0CFE67AFD7D71107758CBA06CFF05B70F629DDD94AAEC64FDD67C053E2235EF84D3C8CE7BF78D749B98203F542FCA4B2F2AC
E51F1F010770A21E8217BAC707E863C0033F0317F24FDF040EF0963C9377E2260FBCCE70513CE0257A1E599D5E2338C00FE3
83C6439EF0517486BEC8DFE8FDD99F3BFFBF5CC7055FDE970267FA357CA87A7285B75CD1D57B13FD16B8FB6B7C6BE4C56D6D
655DE41B7CC0CB70D775C153F8893EE91FD069C4C80DFBA57F633F83DF1037711167ACCF9A6DF2FD66F6F3E8181FC557E12B
E78C9CDF76DA5CD32479E5791E73CEC77AAD6FD63839EF7E23A666867B3CCE83B8FFED253579CEFFFB515CE117F8F27A5EC7
B90BE3735ECF3C967FAD07E060FDDE9AD0E5F27B49D43FFF2F92DE757E4743FB15F8CBB900F3F5EA5E73873C0E9E127983F4
0DF497E875ED96FCE1C253F88A7FE0FB3D7F5A19EEF65FEB0702E252E4FFA3E09BE0473EE027EB611DECE339A76F1BB27796
3B1FC445BF407CF095BA046F7D2EDBD6D65D67D021FB40E6EB73B2BA59AED4DDE11FEF784E1E936FD6C3F93BEB32DEEAFE03
FE9237F2859F9A3FE87A581F7C667F0CAF5927EB639F07DED443704797A3128AFBC8FAF00FE2A23F213EFE3EBA67BEF37FD1
D11BFB13E6EBDF67CF3BCEF95C8FE23EB27EF2082FC01BDF0577E395F6A7F08BF119977D037DD8876B2B4EC8F3AC77FDBAB4
3FC4EFE803C0DFFC41D70FCF1887FE96F1E8B399EF8315556F3BF53A2A2B51F00127EA1DAFE7758CC779FF8B6BF79C943CC3
0BD649FD67BD8CCFB8FC1DBCD033B8313EE7E0CC83DEF0517467F56B5ED2BBCEFF09501DE2EBF83CFE68382A7EF0159F83A7
7EFE3B9E932BE7E50FBF51325DFA03EA09FC850F8C0BAF88077C929F58DC28FDBB9DAFEA792BF3302FE39A6FF4CD4A9479F1
8FA035B111DBB71FB3FD0EFC32BF1D95F687AC97BCD0F7D2DF254FA9C993ABF5933E1BF6BBCF23F17FF2B16174C509B9860D
0FE926BE489F8DEED917304F46BF9A26B95A1FB224F206F97F90F8357D077578FF73DB77C835BAEFFA1BC53FD8EFE107AC9B
3E283EB1E60EB98237FA0377F24DDEC823F9E7FC031C38877E76655596F018DE311FFD56F49DD5CD72A5BE93A74D13ABDE96
2BFEFBEEC835DBDCFB037CF49337CA3E967A06DFD99790BF9CDE3533E44A9FC33CE1D17BDE71E603F761126F0BFE47FD7273
A59F639F479F076F79DF07FEFA76C8F373CE0794F7E89E3C706E4B3EC00FDF0547E383F6B9F08278F159E28647F00ADEC33F
EA183C840FECC3E005FC067F780EDEFC1D3DD0A793EFC2417B67C9757D4ED9C7EE3A847E2D3FD42DF52FF43C626B9EF37D5C
E399E239E2DAA2B9922F5E175890D8D7F9FEACF61FF4D309BECB7F94FD979DBF291FCD07743DACCFF787A2B9F23CEB033778
862FC237F8083FC93BEB65FDE4CBD6A9FD05EBA5FE7CB4A1DC39C7024FAB77636223E47701D0B9DF235B6BDCF9C2E7D1377A
C797E1353C67DD368EDE4F5DA19F05CFB2D17B4ECA957552E7582F7A45BF32DF9E3DE7BDFCBF4F6BE58CAB7594B8C185FD31
F8700EFACF0D153F4AFED0073E844F309F8DAFEB661E78821EE00B7E8CCF316EC89C887AF151FA5174061FE0017CC577F161
C68547C66FBDDFFA28BF906EF23B71D403E3B9AE0FBE5077DFAE5DDCE8DEC744BE1A382537779FED5BD02BFCC03FECFC40FB
88C1453B2B9CFFB3A4B8E17FF821F1587F1BB5FE46779F0B2E9C57820FBC1E5EBFB546EA347EC0BCF433CC6F75A45D625FF9
DE29FCB6FC6B1F070EE809BD067E58DED67DE579F24EFE580FEB23EFF81FF9C79FF02BF262BE775D5A2BD947B16EFAB90129
954EDF417C0F7A15FE4DCE2FEFFF293F50F61BD401C34BFB7770631CFA52C6A39E50B71B724AA60707A779C5BFB8BAB5F4ED
EC3FD10D7D167ADD39A12A4BAEF8C7B055A553845FE8803EFBBD4955CEE7D199DF6F476EBAB3BF567EE32BBC9E7D01F7C10F
7C119E802B7C045FF63FEFFCB16798F33B283A1ECFB32FE3EFC1A9310BA4FED26FC10B7C0B5E04B7D959E17C6E44F799DD87
560F73EF37E1013882ABED27472CFFD1DDE7A2B72EB72F38E1D61DF967FF040EF8223E099FC82779C467CDD7D537C187F5B3
4F278E8071A553FE7FF509BF18FE5C6EBAE812DE5BFF302BA25E7E2F09DF217EF081CF9C23BC15507DD43DBFE55FFB207080
471F4F2874BEA7CA3AFC5716FAB9FD0E1E71CEF7C4CE723FF171FC8379390F617EF0045F748EBEEE1B9E9B2BF8A133D3BBE6
01BCB98FBC9147EA203C0AEB5A79835C237A04FF2CBAA7CF46AF9CBBBC99595DEF5E2FBA44A7D4DD80E5994D6EBCA94BD677
E5074E71F75F8CCFF90BF304FD92D041EA12FD2738802B7E06BECCCB7CD4451BA740706F19CFF6A9C7B63FE5BC4F73491DA7
7E472EAB3A27577087AFE06FB82BDEF01EDC3907047FCBB3F6B3E41B1E312E78D2BFE073BB3FACF8D1F37F8B536F96CFA310
277D0538C343D6872FE293E8D6C6D1FBE9CF38077AB5A8BA931B9FD8AB7C16CA7DEC97F0037C04FFC04FE0177C83DF56AFBB
C42C90DFD140B7F423E010E7B567985CC18B3A0A6EAC937322D60BCFF069F8863E867D9AD924F5095FC0C77B95AD6EEDEE73
E96B58576250F530B95A1E74DF403E52567C7D5D64649EEDB7534F2F38111A9A6978322FEB805FF04DF25553F35D4B1DF94B
4207F9FE2DF9E2FCE3E5C23DF7499CC4BDF19EED4FB9FB027840BE799ED7857D19BE5EFA0DFA3CF8481F44DC69C3AA8FCA95
3C7D98997AB37CBE8FF5811FBA0447FA824157EC089075821BF9617F4C9E023BE7E4497F897EF00FE2E55C86B8AD9F4908FE
59DE77C75FD0F1D00B25E3E43EF4CA7AE8335897E541FB2FF201EFE02138B1FFA54F3972757EA05CE99BC85F56C76A673CF4
8A7ED1213C4337F08DF39E479F288D16BEE01F961FCD0B7E079FC0C5FA35F5918FD273F29CFE54F9079FC00D5E25DDB4EC65
E9A7D81F824FD04D25E3E4FD42F8849F445FE6FF87D479FA6F74B0E181027F399FE275E0090E1F1E2EF0977AC5FAE8BFFA35
57A439FFEF45F10147740E9EE80F3D9A3F68FCC4073F899B3C902FCE5D9E2FA81A28EB825FF009FF8057D437FCB8F6D972BF
E8E842AF90B2E8DEA20FFA337864FDE1CF3E0BDD7D22F1D24712373E895E437EDB1120577061BF083EF459F02EBF4D7527B9
A26FF0023FC34BFB7170DB7036A3C1CD27EA09BAC707C8BBC5ABFD1171A30374417F401CC4C5FAC83FE73EE000FFEDFC5375
6075FFABF0F5F2FB53E81FBEC1337462FD7568E524E7F7E635DFC44B9CD43FF809CF2EE5298FE191ED4B277E7D9DBB5F26DF
F09FF5B11EFA7BD6B5E9912D0F4A5EE96FC8AB3DAFFD0B7FDFB830F997CACA6667DD72A5AFC16FF01FF0C3773F68BBED4A89
0B3D310EF7330F7CC59FE0ED963CEF05F2FF16ECBC42F9882F52B736B7AA48932BFD053C2D69BBE73EE773DDBA2FE9F95495
D367930FF3B15B331AA4BF45CF918F6C5C2BBCA61F85BFF000BDC20774894EF15FE6655FC1FCE0C9BCAC837397D63B775D2E
79B3F349EADD8BD1BDE57CC1EA9EE699BA42BEC91BFB27F267FDC945FF3FDC7D0A3A46D7E04A5EC9B3E8EAC081FFB6E447FB
4AF2445DE97661D9CBEEBE8CFD4DD7A5BB7DDCFBACA0E4EC44D10DFAC10FC061C8F4E45F64BF84BF6E2ADC76A5F01A3DC37B
FEFE41FD960765DF4BBEF069F81233B072925CC933E729E49B75B22F62BDF809F3B20E74D8397FDE6F6E3DC2F7C11F14F794
79D037F8C12370240FC48F5F709EF16076616BA9C7F813F51CBE97B7AB1A2857EA58FF3EE5EF38BFA7A6BCC04786BE9D9D28
3E07FFE282D77E2BBF0F443F0D5FD86F752CDBDDECCE077C051770824FE8185EE1C7D4B1BADAD2E8F0F0EC3FF912D651F845
FF078FD02BFAB5FA1290BF403E8781AFC197E04DC53DDD3CE0EFE47BC8B0FC05E2E3C40FFEF0017E132FFB40E2A6FEC0A32D
7BAAF2E44AFEC8277921FFD45770800FE8135E581EB41F241FF018BE189FC14BF701E0061EDC871EACAFA8DFB8567E0F11FD
A34FF88B4E6D1FF29AF702779F18FC4A7A7E75B58797CE55EB003A20FFE0C179C3DF2FECB845E24567ECE75E3FB33BD39D6F
F00557F400FF3977400761A55137489DA1EF217FD41F704B99BCDB47AEB67E5D3775085EC013F0649DECEB582FBC023FEE83
47F8637874F93B070EFC647C87FFD409FA82BE2FEDECE29C0F2A5FE0177C83AFE4157E93AFE35766BD297D0FFB03FA21CB8F
F64DE48975B14EC1AFB1F137F30FFC041DA6D62DFC87EC9BD84FA6BF3AEFB790904CE313BA875731F17E5FC87DF4E5F07453
4652953BEFD47DEA247E56F9ECAECBE50A7F38BF8047A16DE24E3A786ABF427CE493FE86BC322FF3D1DF300E8F198F78D92F
12B7D57D9FB08EF2BB51E80A9DA027F4821F0C2E4ECF97F333F4891F0F082D2D163DE007ECA7DBFC56F55F37EFE023E3302E
F519BE67BEB8BBD979FF98FC6BDF0D0EF8F23B83D77EEBD639FA402FF02F79FC92097171252DFB65CD07EB643FCC7AE9A77A
9DA968707EBF8C3EB355D68A7DFB7E74F4BC7FFF4FE61F56075F8ABA417ED7C0EAA1E607DF244FE66F9D92AA3CFFA7D3E373
AC93FD0AEB4537E832AA746717B9821FFC05C7B081A5C50D0DBF9A1E446F070FFE6E71D1B7111FF30F2C8A3B29DF73239EA8
63410F897FD32F8203750CDC726FDB9D2957F80EFFD12179235FF8177C195857748F73EEAC3E447CC48BCEF1A70F9AB35688
8FC2277C82FA7CA87D616BB9721F79202FF006BEC31F78070FE1137C052FF063BFFBE25F2A37491EE111E70D8FBD54E6FC5E
2DBAB77C6ABF495E43AECB1B505FFF8BE383B23EF8C539F5BD3E39ADA53E5337129F5D51999CBCC3F60FF025F4D1A27BDC78
331EBED6BB97DF17EEBE867C0FFA3E6F808C7F69DD20CFC44F9FDB6D52E59512273CB2FDDBA30BFFE1EE1FD1033C4517D40D
FC72EF4D3B6E912BF9B4F304CD2BF5E4BDF46D179CDF33553D9107F607E48375D2EFB25EF083E7F01E1EA24FF8887FA0BFF8
DB2A1AE4CAFDC4875E423E4FF376F39BFA10F96AC821D1077D21F8C3537B1F40EB2BFCB1F32BE5113C811FE88AFA862E8BAE
A8FAAF731EA97C870FF003BDA25F786FEBD43E8EF5927FEA1038102F71529FE011BC8247B68FF25E32C1DDAFC14B7479293F
C183F1C0015CF00BF472F7AF596FCA7CE886B8D10379641FD0B9AA32C1CD17FCB4EF3D410FC9EFEBA19F2D7F5BD343FE0F80
ED43957711299B9DDF35A06F237FEC779FBDADC2F9FC38BC0EAD4B5CE3CE17751B1DA00BFC9475D2C7B35E9B57FB23E68FEC
BCED82F003FF830736EF257D03FA1EB420CD5BF63FE88478E9D3899B7C333EF38133E3306EF81BB141EE732FD68BBFA2B3A4
6595573AEF1BE9EB591FF75187DFF52D77DEAF808FF878F79D2B2ADDFD0B3A4017F01E5CEC5C42F1C157DE7FB06482F00E3F
B03A9B1F72487EC7047E50DFD05FD92F959BC4CFE10DBC863FA183B7F6DBB7EF07D32D3EC1BEE1CD5B2A1BDC38A3137C15BD
C475F23D25F3D18F820FF9239FE884BC51EFC99FF9E1A3896BE4FB03F81AFBF287128A264B1EF021D669FB1D5D2FBE85AFEF
2F2D5B2257F3FFAE9B53E5FBBDE4913CE323E49BFC9017FC25E250C904D1ABF158F54C1EE813C9073E3130786B3FF93E36BC
24BEFE39B141F2BD2FD687FFC1C7F47F5526C815FCFB8F2ED825F1E343EC53E83B9A46E6B4962BBC279FE4D7FA8D036B7AB8
FB8EB051B9F7D4D5FD6C7E042FE0193C826FF01B5EC0738B4FE3651DE892BC9147F6AFFF7C76C73382033A20BE01EB72EF91
FA4F1EC3D716EC72E79DF5C277F88FDEE01D3E113BA2BC50AEF40B3DBAEE74E6858FEC0F5E7DA572B69B5F31DE811BE4F3D7
F49DE00C3FD93FC3D3F4EBE7BFBE7973AEED7BD8DF313FEB212F77F9663AFF27FCCEC2F424E9ABA81BF0173EE3CBF0029EA0
37F80D6EF03C2C3CB54775F505E3317E0F5FE00F7D00BC819FF0953A80DFEC3C53912A57F2435F499E5292BE5E2FFB08DBF7
693EC83FFD103898CFA90EE13BB8D83E4EF181B7F015DDC22B3B27517E110F794277F831FACB3A57D9E0E850D7CF7E8838F0
FB3EADCBA6C93CF8073C85B7E880FC9377EA1F3A80CFF0DBFA900CDF53EE7E04DDA23FF48B5EC817F3475FB669AC5387B5CF
031FE6413FE4D15EAFFD17F7252D597E6352D28E96FDA0E623FC744247773CD449F283BF9327FC6040FFD41EF2B95B7846BF
D3A5B2E269E77B41CA0BF6312FF7D9B5CAF95D27C59FBC320EE3C2473B77515FB075EAFAA8EBE8069F4848D9F98C738EB922
A2A7E04FFF443EC803FD3EF9802FF0077DDA38DA1F319ED5E5F1811BDC7DE7E61762DE74AF8FFE84E779CC78F0035EC013DB
FF7C37FF75B9765A38F72FCEE7A22BCBA6EDDFFFA3F90EFC236EF2804F810B7D1EF8E013D4DBC61E4593E50A5EB60F55DCC2
FD32AFAFADFDDEF4890FB11F7D62ECB6DB453FE895F95927EBC6D7D16B7E71E56CB9B21FBD7F646E95F01F3F4E58B56A83F4
8DF4C7E0865EFBC6145F21BAC16F6C3FD0EDEBF5EE3CE2BBF8DFEE9D3B9E91ABF9E4ED091DE5FB03F80BFEDBE7E2A6B1F23E
1678DB7AB49F655DE487BCE017D4A55E376D7F4CD68B9ED077FFE199D74B5F8B1EE013BC8057F4711D9FAD48717E874DF104
373B4FD57A032E766EA0F8587D9918D153BE174A9C51BD8BAF10FF804FF80BF9E8B73DE64DF95E93F9AEEA96BC91C788E4EC
62393FC58FE00B38DB7E51F5068EF8CBA578F6BB3E7FA2F831BE882ED17D6AEB8AA7E54AFEC827EB639FF4FC988A5BDCBCB6
FE7ECAF21BDD7D507C994F979C9C5AEB83C199F9DF7F3BBB58EA3AF86DFE2E7FA2C449BE881B3EDB798BF21A3E58DFA2BC38
15913A4DEA37FDE0C9B6E949720577EBE7157FF86BE764EA9BE80C1CC1D5E2D23E97F8E037BC83E7E4953CE32FF8313A2E89
DEB5CA799F2A2FF9A73D7B2E18EFA927E4DFF6AD8A03FAC4DFD0297C826FF00F5EC127F4469EC81BFD10FC263E788ECFF618
B76A83BB0EC337F807EFC019DCE157FC85ED8FEDDDFBADF99CE55DF340FCD42DF2C7EBB93FEE89801BA4EFA21F051FE2270E
74079EF00DFED1E7BE3EBD22CDF9DCA6F2C5C6D7FE8E79B80FBE925FF4003EE802FE9317F2419D7F3BAF7CBCF3B961E535BC
B07DB3EA9FFDDF23BB8A3F11DEE34FF801FACEDC5991E27CFE41FDF29DAC9253820F3E81CFE1EBB5ABB7DD9E90B0CD6BF3F1
F82677DEA9C3F8C1FBAF25FF24BFA703CF62E2379E94CFCBD057912770E6F5DC4F9E6D7FA1F9C6AF7BBDE8D3C55DCF6D3DBA
0EFA00E6A55F637EF64FAD9794FF2CE3A21BF24CFD26DFF09D7A7EC427B74AAEF01B3CE1397925CFE83CBA6F58ADFBDCCED6
757B54AD7BDDF42DF67AED8FB88FD7F33AFE9EDAB8302336B6B8653FA5FBABCD05E9B1CEF98FEA19FFC0CFF09BF23515B7C8
E728E003BE092F189F7E9F79E01BFB3A78075FEC9C4CFD947DDBDD233363A5FFC0BFE3DE2C3925BE83CED1A1D5EBDA801BDC
751B1C6D5FAC78820F7851BFCDF7EE8E6F927D1AFA4FDEBEF4F5A4A4F296FD8BC6471FDD26B2A2B59B8FE8A377CFC2391237
7E469EC92F3E850F111771C20B7C2BF9B5F2F1BB779F377EC32F788E1FF46B971E2BFB60F4C3FAE9BBEDFD331DCFFA565D07
3E811FE42EAC48932BFD4EE72E3B3A3ABF4BAA3AB0FAD56BE349F91E3878C0173B57503F42CFE499BC479ECD1AE4F8B7FA20
7C0487BEA7A36A659F891F131FF1322EB8D87E47F189D9523847C68537F8AFD5A9A8B05AF95E25F1249E5879243DBDAA659F
A2F903CF3E05B9BF4A5DC28FB9BE776BD620E983C023BA5DEEAF322FF3111F7C054F780BFFED1C437540FEADFF541CF0FBEE
0165938577E8121EA24FF8089EB6DF523F225EFA6BE28E2C4B5A213C85AFF8BAF5F70F2DCC70F711EC3FFE9EB3ADABAC13FD
812FB8A207F80E2EE0043FC01D1EA057EA4FFDD3C59FC4C4149A7FE05B9553CA7FCECCAC365E926FF899F0ECDA80ECEC5AEB
2F899B7DCF8323F31F97FCE15BE481F8A92BE00CEEF0115CA813E0732629B154EA05FD26FD13EBA14F645DD61FBFB0F47577
5F401FFDEF5D3BEF105EC077F2059FE137BC8507F097D79127F80F2FF0A1F4D41D1D9D3AA57C86C7E00A6EF09CD7A17B7CA5
E8BD8AD672452FAC83FCC11FDBCF2A8FE8A7DA6FD9112C8FD101F8DB399FD695F8CD7EE3E577C8E8E7C86BE2B0B2C9E273D6
FF29FF888F78F1477CE5BD179356C8FE0A5E46ED8A9BE47CBF5579407D8677768EA23E082F1887719987F5D8FA58BFF68FC4
01EEF4F1E06FFDF45D2B8FB8FB0E5BA7AE8FFE009EDB7E5FF90E1FE96B8EFB64C6CA159D77AA9FDBCDAD7770A10F031F70B6
7DA7FA0EFED473687147E123FE14D729A8AFF463F47FC41D5593D6C6F9FC80FA053E64AFD73ECECE970F470E74C749FFC3F3
F658FB2774457CE82BF6AAD05EF2B91BFA30E6A1EE75E9B56DADF3B94CD5B7BD5EC7E53EABD73BD706B8EB3679B5F363ED4B
EE783F759ABCEF797B4ECACDA20FEA80F9E1D37193E4FB09E89CFC90177C8A3CE33BE47B4B5071C7DADA1F8CF7E816DE900F
F8832EFB3E99E67CCF1F1EA3337C77EF1BDBBAA6A656B4D4B57E7EE3DDF58D7DC663AD4AAE947CA16FFAEB677E2C1F253C82
D7E4D9FA0DCD373CC29FB27BEE08962B3C7AE7AFF9CDC26B7C0EFDE143DB9EDE79875CE9D7DAFDBEC3F3B93FE52FFA24CFE4
3DE6E6CC83922FFC0E3E920FF2C3FDE8983C9377F3E7072207CAF7822EF5DF773282FACAF790C963DAF9F9BB62638BAC5FCF
78646EB7CD9B732C6EDB4FA87FF17AFA71EE837FE00C0F99FFDD9F437BC9F70659077CA4DE1EF2C97F3C2222DB2BEEBFF9CD
729E086FF0797CBCF7F9CC83C25F78030EEF2ECB992C3CE675B19372263BE7927A3F79850FF4EBF0027ED9B98BFA2C38D2E7
8167CA85C5231313CB5BF6239A0FF0B7FDB1FA546AFCB6B5E2C3F831F891579EE775D15B136F73CEB9556FD401786EFB6BE5
3BBABEAD5B62A9E4057DB34EFA74D6CBBEE7DE37B38608AFF1697808FFD03FFC375F551DA043F0478FF01B1EC24B700677EA
6BF2C32BEE91FED0F611BA5EF240FCD41BC35FF9803EA8DB1D1AB73BBA46DFE8035EC133FCF5EDF62551121F7EC67AE85B59
17FA216E74449CF0FC52DE739FED6394B7E411FEC25BEB7B140F78CDEB883B71EEEAA3D9D9352D7DB5AE97F8E03FE3C2477C
AEE0F21D339C735AF5953E1D126F93CF4F831F7DDFCB2F9797BA75409EA9C7E43B2624E6BCE0034ED449F0617CE6A32EBE75
F3FC5DEEFA6F7169DF4A7CF0D6F6EFCA5FF48A1F575C533E4AAEEC0FFEF9E4B6C765DDE81EDCADCE2AFEB67E5D3775DDFAE3
9B444F2DF59CFC914F7C33E116DF82FCFC7AEB0B8903BE53DFF637975C2957FA76FAA3D36FA4DC2CD7E4EC922877DF068F6D
7CEDE39827A626254D5E876EF10F70B47E5EF164FFF1D0EA82A7241FF84AAC4FC42E77FCF41DE4953CE377F1658129D2BFD0
57B12E1B47EFA78FC1E7888B38AD4F6858718FBB5F201F761EACFD8DCDABFD1AF383BBED13157FF4DDFDDEC265A23FFC0CDC
F0E3AD0F6D77AE71C19BB6CA7A59377D1ACFDB639D9FFBACDFD738C1D9CEA5B4DEA0B3DE8362CECBF71CF01170043F7C10DD
F77E32254D3E578AAEE82F5EEFB67DBAF3FF925497569767AE3EEAAECFB7DEB5A583F415ADA6C45F257E42DD409FAC87F5C1
13DBB72A5F129B0A97C97EDBFA3EF59DB807327AC9F3F8177C64BDAC9FF9D0153E80BEC0997EA1292B6BC8A64DE9E6EB9D7E
2D8D90F5A063AB47E77C0BDC7509DCCDE7157FF3ED9111BBE4FB48F835BCE85193D7417C179F86EFF868494979A95CE92B9F
1FBF7D9E5B67E0627DA0E203BF891B9EE3EFEF1CCEE825F880B7D59D170353DCFDAFD597C19BB6CABE9AD791273B8FD4BA9B
F0A4C4F3A3E1445D88F7CEAA95E7790C6E56E7C767D50A3FAC4EF96EF9DADDA7E0F7E8183EC00FF88E4FEFAED9F6B85CC91F
F9C44FD997DC3926ED2EC109BFA7AF7F626A4977C927FA860FD42178010FE11FBA420FE0802EE03D7C801FE06BFB08C5993C
103FF5C3700347E53BBC827FF08B7C901F7C161EB35EF89CDE34B75EEA08FD387D3BF9C3CF7293B64F976BC695A511EEFE12
9C527F5BB4203171BBF5AF8C076EF004DE9017F806FFE017F5AD716CC15372657CFA50E6615DAC13DFE4CA3A5937BC271FF0
1FDEC257F08627B61F54BEB06EC6432FF01A9E831BFD4E9B90ED8FBB758C6FBD3B62CBD7CEEF0D296FD00BB8325E4AC7659D
D3D22A5BFA70CD07BC601CC68D1B1A3DC8CD27EA30385CAA5BF462FB7AD58DCDABFD2CF3E3975D2E2FAA103EE28BEC03EEEF
915D2A7FC7EFC1D3CEC3B49EB14ED647BD4FCA5EF54056568DF5A3D6DFF3BCF677D6AFE6253F2C9FAB409FF80A3AC077CBBD
B7CF932BF5ADFD8AB2D6CEFFAF53DD8397F5AB8A1B7D5FEB88B258791E1DC7BF18FE827BDDD4F9C48F7CAECDCFAFB3BECDD6
CBF3DA7FD9B92AE3E8FDF407E78A63A64BFDA2BF6B5E127F955CA9739DEF9D5BEF7CFEE3C857CEFB2CE062754DF1216EDB5F
AA6F263C1B707571F121EBB758D796BF05DF29BFBBC33AE84F989FF5D057A5FE5E5421F9B77E53F56FE36B9FC53C36BE8ECB
DFAD3FBE62D102773D377D04450F92EF3FE023E00EDEF8293CB4FDA9F211DCADCF51FCE12375A676714977E7BC4875FCCE6B
C90FCBE76BD11B7CEF7A3E3F5FF48A1F5B5DCE5CD6D95D9FE18FED239447F893F9A9E6113D132FF1936FFA9A936BD2EE0A0C
4CB2BEFB91A68283CEF782D4B7582FEB67DC2D77A4C53BEF6FD067A96EE00575A068D0F6C7E54A3EEC9C43EB0D75A55BFB55
0FB8EB21F8F42A0B7F41BE87836F329EED7BD4C789C3CE1DB5FF48B9393F5FD6697C52FFCA9E58D6DAEDF7E491BEE9DF67B7
BDEDD631F3F07AEEB7FA18E873ADBB4EE2BB3D4FA5C5CBE70EE1013877CBCE59287E894FC70F8D3B27F3C33F7C97B87B1E08
BE53BE2F67FD17F56867C0D5EEBA94E0BBF549E97BF065F0816FE819DE317F8F115B9F947E817527B597F5FD60F9623C7485
6FA02FFAF75BCF260C9471A903F0047EA037F866FB3EE51D7C05777880BF76C82F0E155FC29FC099FA7F2421BB342222C7F4
C37AD1117926BFD40FF888DF57BE5F162B57F4047FE01375F8D5ABCADA38FF1756FD00FE581FAE3C027FE383EA09DEB35EF8
4F5FF9F70B25AFCBFDE89EF8EC1C4BEB22EBB5FDB1FA3579843FF0091EC203F8088F59177CA6FFBDBB677A9DE41D5F846FF0
0CBCF1895E4171E764FF06AEC4CDF3E627574766EEDAD58207F5907CC3437809EF592FFC274EFE8E4EE08FD56BE5117C87FF
E003BED46F70462FA65FBD8F75A113F2822ED8E7A20FE2254EEAF7D6578B43DD7D27FC4B6B5E1020F5C5CECBE9BFF579FA5C
FE0E7F6C9FA23C8217F8F4B65BB7BD2D5778441DABBFB7E0A05CB7AC4E9C29F3A35BFC29B5FB9231A9A9BB5AFA6AFA789ED7
BE8FBFE397C4459CCC6FFB30F5F72D79A14B77EC386979A18F40F79D1A0A1A85BFF85972C3CA9E59597BADCF637EFAAC679F
296D70EB92F8EC5C44EB13E3D027321EEB611DF41FF4A70F2EC8392479C6CF8883BA91FF73591BB9A223DB57AA9EF8BBED4B
345F49FE6B3BE5E5D559FFC7BAD053CF9F2233657F88CEC10B9CF03B1B47FBC79673DBA011DBB635593CF4513C6F8FB50F4B
9CEBD753AEF4792DE3E9F3FA3AFE8E9E7A8E4D9C29FB0AF47F61EC66E7FB24F4A7F485F80271112775F8E5D1A5F96E5D822F
FD1038673C5CD02879403FF808F3B31EE6C3BF5F1F5C7295D3E7A94FB02EFAA6B3B7260C942B7CA16EECBDA9E475B9C24FAB
53CAD3841F532ECAFCD637A99EAC4F6DB520C05D57E1F55B0DB947A48EE2D3CC6FE791DAF798CFBF16BA54BECF830FC24FDB
4F294FAD3E262E19E3AE93C461E7015A57E8071F7BAFB09DE4175FB1BEEDE1953DDD750C7D937FF048B839F6CDDDBBCF192F
F17B7CB4C73529179D730A78839F8F0B1A21DF2BC3D7E183F569CA0BAB5F1FAFEDE4AE63690FE71E91BCC303FC9AFCD1471C
DF925E27D7A44519A39D7E517D1CDCACEECCF4EBE9AE3FF0C7FA01E511B875CDCECA703ECFAA75033EC003740BEFC1115CB9
AFDBB48CD1F2B930F295D23E2B43D669FD9AAE973C935FEA0DF1D9F99CF613E04C7DD85951DA20D75B5E89992EF3DD5C15FD
82D43DEA18F8185EAA27780BEFE12FE3D93E527D3F37B8E42A773F0D6FE09BED479477E8A453FBAF0EFF2FBD288FF93B7CA6
6FFAE77F8AEBE5774AF009FC967959073C8527F0063CE021BCA4CFBD3D38C9399FC6D7C92B75F5D0E7398762628A2C6EEA43
C9DAD27CE77323EA433DCEC7BE299F478707FC9DD7B36EF446DCE80E1C188771E1333C867F09BE9BAFDDB5EB8CF1837A8BFF
B50F2B5C2171E167E88AF9D117794117B62F519E1A6F952FE8C2FC4DF5419EECBC47EB34E3323FBA63FDAC9B3A0D6EE8029D
704527E4057D994E75BDF40BCF179464BAFD26B1C79612795F04DDE267F4C5F77E9CD957F88BEFC243EBE7958FE8EBA6D59B
03A54F4467D4B736878B4FB97D25F1D94D35121F7152B7C91375697F64613BE9C3E0A7EDCB94A7368EDE4FBD0777AB838A7F
7AE1BC48F91C33FD26FD27CFF33823FBABC3F239B5B4365F2F9475585F4AFFCAF3DA27F2F7EC01852BDC7D213A44073D466C
BE56BE9F800EC93F79C7EF88DBCE1BB49EA51C5D3E6DEBD66AEB0B993FF1446074595993C54FFFC2EBE907B90F3E583FA9BC
B0712EE9939233577F9E9777A0A59FA44FE579EDCBF83BFAE89EB0A5443E378C9E93B27D9714151DB471E9BBEC791EEB78F4
A10F35E546C8FAF05174DCA13AAF4E7C0A9FC62F2CCF9A7778687D97F211BF042FF083FFB64F541D100771311FF8507F76CF
2EAE976B524E725B392FB0BE89FAA8F877DFB9A946BE6F83EF701FE3506FE833E9E3CEFF2BFA05F795E7E9AB8A1E283EE5EE
D72DCE7FE4D5B9FB18FC2531327A7155D539E31F3ECAFCF87D79BB12E77B69F44D4F1C2EEC24FAC70FF079E6651DAC8B3EE5
F4E0A4C91B3624993F757B23B9AD7CCE12DE581FD9765EA4BB4E32BFEDB7743E7862FB17E50BBCE8743CFB5BF1237CDAEA6C
D1D70BDDF5D6FCFFAEC068F97E177E0DDFAC8E28EF88C3CEEDB42FB13EF2BEE5D3DCF50D3EC003F49F3238CDCFF179FA35E5
0B3C8627C61BEA487BDF25FFAB4FA4AE755CFDB9BBBE812FB8526FF0E357BB176E97F8F0CB8CBBB3BF95755ABFA9F586BAD2
3538CD4FFA26F023AF6FD566EE17FFE775E809BEC1BFB42732F7BBC7256EF247DFD1E49FD957DE8F8317F004FDD1073D7D6F
F1A36E3FE33E3B67D27A4CDD7BB1777167B71FC067DB1728AFC93B7CB0FE5D75DFFDBDE8C5F2797770A2AF6CD533D657F841
BD62BDF67ABD3FE9A6F0EE8283F52F5A0FE1ADF99BF2179EC00FF8447CD4CFC67B7323E40ABFD1113C0777ABCB8A3F71C22B
78C6FD1DE7FFA7D1AD47F8CD3CF0DC7056FEC027E2254EEA373C325E29CEF4A1772E4BE92CBCC48FC19DFE03FCE13DEB6AE9
EF3CF1F07774429EEC1C45EB6E7E62E176773F0D8F93BE89BF5CF48A6EF14BF86F3EA33A806FD4577807FFE13DB8A21BD391
E683FCA117BB7F49F06DE5E5272C8FD47F7C97F5130FF151AF6A1F28EC24FD01FDDDFD115BE7C83AF151FCE0F5F0FC11B21F
C2FF9897F9E81B18CFF6DF5ACFD047B70BE1DDE57B0EE80ABDD8BE4C75439EC92FFE95FCB0FFD1B2B2A3361F7D44594C7167
F73E84F8999F3AB6ABA9F851B9DA387A3FFD077CB37D90F28EBAC3F8CC074FBBDD187FB97CFE189D338FED17F5FEF41B17FA
CBE762ADFFA4EFE479ED2BF93B7E41DEC863C6AC395BE4F718B89F7E35B572E93AC93B7D21E32537ACCB2A2A6AB438E9AF78
3D7D22F7D9EBF575F479F45F8F7E9477AFE80EBF4CC9F3DE2479B63E933E95E7F57EFE4E9E2C6F1A676EFFFC1182B7F5D1EA
2329E31247CB15BF437FF0D3EA97F2143EDBFE45790DEFBA4D09BECD79FF40FD8E3AF07243E102B90F5FA4AFF9C77F0A7E96
78F11B74D0FE444EA9F83B3E9DFC4464615555734B3FA5BE4E5D7A36BBE82EB71F7C531FEA7C0E8FFE8E7E0EBD9207F2F2ED
81E023EED7733FF7D167356F89F5952BF1BDDDE07FD4ED83F84AD75589A3E5F3A3F098FBEC7C4EFB1B70A7AF027F7869F9D4
FC821738A167F2435EA807E6CF0FAFCB72FB74F65D39A5F23D4CEB5BD5075917FDC2C949299D434232BCD2EE49C9747C927E
4AF96275F09B85FEEE7A66FD59EBA5EBDC75CCEAC86BDE9BDCAFB37E74CE9C2DEE7AC8FC762EA5FD07BCE878DDD60BE273D4
07F8667D88F20EBED17FC03BF031BC5437D489B78EA564CAE701C135F3FBAD17246EEBE3B45EA47F977ECAFD983CB15EFA88
23EF6F9D23EFBFC05BEBAF94BFF0DFEA82EA80F83A5F9F7E4ACE4B585FC9C3850BDCFB01F80D6FE0A7F5DDEA076FD74616CA
B9393C635D76BEA2F53B797CE8DC8A8AD3C627EB2B344F8C73A9CFC02B7886AFB34ED6CDB9C64DC323FA499CD42BE6A7AED6
05E6DDEBF8AEF210FEA107FAC1DBAE8DCF15BFA2CE103FFC470FC4453CD45DE6B17302ADBBF856DB81F933657C7C115E184F
9447CC07CFE1FDCEF64577B9F77BE485D7A1974BF583BED051CA8CD8A7E5F328D637A92FA20B78833E8883BA5D3DBBE067B9
322EBA4387F4837747A7EE111DE1F7D45FD64F3C8CC7F8BC0EFCD10BFA49DE1ED459F26D7D96F61DE808FD80377C37FE6BBE
D10571F277FAACFFB9396BB0E81C7F655EE6A35F615D97EA197DBCFDFFD83AEFF0AAAAECEF2B30A8888A224D447A13105140
8AF4624361444606C4916207A49700092524404209044209218626BD444A1A1D11A94AEF08A107E9A0A38CEFB3F65D9F7537
797F7FE54972EF397BAFF56D7BDF73CE8D993342EE978057D4DFF293EA2FFE53E3934D69BE0EAD6A363D9FFB1E253D1F7901
5CA153E00BFF2FF7E786BA327E746BF39C8D83659D69395AF9068F9E2C3DC75D67079FE005F90F7E301EC6412E81CF365F9D
3FB87E7BE4929765BF18FEE3DB2F55DAD4D6D7B355676333376F3E69C725A7802B5B7728BEECF5FA3AF20DE7E1BC8C6FF5E2
F19BDCF3E1F47596B3F83BBF6B0E5BF3ED98A5F2BC19FB9D3CA87FE777F222FCAEFFC3FA97A46EE8F1DA1F46CE97E72DF13A
CB99FC9DDCA839129DA35FF46F4D8EE5510FE42FE5B1E9768D596FC87D38E81778266F816BF0867F80BB55E3E75F93FB91C0
25FA0DAE2C672ABEE823BE483F39FFDBA9D3F3F93ABAB1F6FA970477962B95E7E086F1331FF8FAE6FDE551EE7A7AC5D3F5E9
F33F7C60BF4D7312E367DCE8F54FB336B9CFD76D5DA2F5CD9ED7382EC723D75CB9B36C7D42C2F7567FEA0E9FF9BFED3F69FE
303D7F3E36D3D77570D1A8707A92AC53D15DF0637EA138BA5D7066A68C87BCC9B8390F79E1FC07ABF7C97D8B6BD39377385D
227F295ED05DEA405DCC475A8FDFE4FB07C7B3FD1BF57BCB671F8F59FA408EE367ED91F37D7FA06FD647E525C7C3B74F5D4B
FF58F6EDF9C9DFF15774F1B53F37BAFBDDD1BF9FF76F6A2B3A60EB3C9DE7BA0BE949D22770865EC32770082EF181379A24EF
70FB9E8A3F706978261FA8AEBCDD63FE35D9D706C7F81AE3617CCC035F3D926B435DF949FF9BB4483925FB0DF812FFE7F51C
775572D252C1B5E52CF5C3D4D9AB2FBBFB2FC995DA7FEA6DE3D471739E661FADBE2CD79131FFF4E494538217EA65FEAAE7E5
7CF8227501E7E01E1C823FF8435ECB7771C136C95BF80CF9ABC0FF56CC117F43EFE10B7C80374F94FDF68CEC53E5BD9CB85E
8E878FE1272F9FDCD8C5D739C607CEC1FDEACB8B5AC9BA0E3EA37FE497A2BFAF7D56F40D1D05178613C535C78517F0047C81
4370099ECDF715D7C62BE511BC228F949A917E49BE4F183D5E757D6677F73D6BDA0FF202FE5F71C186DCBE4E73DC66994397
FB3CFD31D7C6A7FD9C0D4FECF87A5C7200780367F41B1EBCDD2269A9EC8FC31F704EFFC03BF5A7EEE83878347C6ABDA91FF8
04AFABDF9D5A43C665B949730438377EE9F1781FF8663EE82275A02EE84D9D31EB8B094ED1697865F3D5F9331EC6415E019F
967F14A786A3028B5AC975D5F0165ED83A45F9018F6C1DA47C02B796C315BFABF74ECA29F711320E7210E7B1F3EA7CECF5FA
3A7212B82207802FF049CE00A76B068F3D26F771F27EF214FC7DFCD2CC4CE9333CB6D767CB61A6C34FCDEC2EF7B7A03B6B3F
1CB54EEEA3E57596B3F4EFF63B796BF0D299821FCB53CAFB3D253776F1D727E020A5D388A5729F34EF27A7F1777E4F2D3A74
79E0395EF3F2496EB33CA53ABDB5EFFA62C227CBE3AA0BF803E7651CE00D5D0077CCFBAD95536B3C902F147FD4877A81D37A
B9D3CFB9E71AAB8E1EFCD786DCFEBA8E71307EC68DEEA3A3AFC66C282EF842B7F0298EC371E1EB1B2396CE94FD45F04D5E22
BF5C7E78C51C794E353F6DFF46F3C7B5420BB6F9B98CF7D347FA878EF03EF245E6236B9F95FBB3B2E73B8ECBFF793DF9E344
97F44BFE3E02F3E4FFFC9DD799BF549D94D3F7995B0512D7FBB98D71A02BE0017C988F8C187BCCF7137220C7E1B8A98D56EC
713A46EE535C9B8F7C376A9DEF27F481BEC0FB4D7FA49F937529B844672DCFCD1CB1D4F70D7066382407A87EBC55665E3EB9
CE1E5CAE3E99B8CB5D7F4D7E51BFC2972A77593FD85F4780D3065FA47677DFE3A5BEC178EDF87A3E74BDD9BA157B647F0EFC
815B7416FCA6D7FCFEB8CB39E426AD1F7C816FF086F1336E7C111F2EDB38A384E81DBA485DC003F8D81097DA5DEA4A7DF121
F2C50B5D5372CB71D0517250E132C94DE4FFE83AF568B468CD3A59C7E36FE4A967BE5B745FDE876FF0FF26DBBF3F2EFB4FD4
6BDD076BD6C9F92DF7693DC01578826F3BBB6F282EBA68EB2AC5F1DAFA0BBA086E2C4FA9FE91CB9EAC3ABBBB8C0BBFA22EE0
1CDCC30370082E57FF35235ACE47DDF179FA6C7D476FF4B8E01F3EA0AF8C9FF98043FC0C3CF27EF00F1F6C3C3A0EFC1E5DAC
596F5DB4AFDFF0E0AD9289BB645F1EFE18AF3A872EF6F905FEC11F3CA0CED417BD0687F82E7804E7E0DBFA5A774A2BB90EDB
7293E602EA0D3EC12BF8244F80537064B8D23EC04F9BAFCE1F3C9363C035E3611CE48F7D33D60F96DF6D9DA7FD03E7B68E50
BCD3FF37362EE822FBCBE8057C313D51DEACF936E639B9FF8FF39153E0057CA08EF0C878A7F3A64EF088F7737E1B8FD6059F
645ECCD3C6A3E3200FC1235BEF289FD6BE127555AEBBE775E420D3E79C33A2651F1F3DB2D7939734276DDFB42EDACFBBF010
3CA39BE07AEDCD88ED721F83E52CCD474733324AF8EB5BE6997274F16CA90F3A824ED871C857E4ADBC73CACAF92D4FA96FE0
0F55CBAE4BF4751D1FE6BC8C23E5CCB0E5E2D71C979CF6EBA494DCFE7E03EF87AFB58EA52F135CE103E4045ECFFBA9EB9B5B
A6B47AC0FF15F7D49F7E5C3896DC449EF36AFB2B7A5CF20F7FE775E8137DA14FD487BAE04F573F5C747FD6AC55C17D293D2E
7987BFF33AFA62FD507D317DFE38E6395FA7E177F3B28B67CB3E2B7CC3F7CAFF95914374011DBDB17776773FF7300EF20B7F
E775A6DBBBA2AEFAFABDAD4CFA32D9D70397E8147840AFC145DAA2655BDDBE0FF94871C7786DFCCA4BCB494F446CF7F59EFE
597FB5DF778A4C3B21E32597311FFBFC2047CA7CA9037A0CEFDFB835A7ACEC8F838F35D712CECAE749968FD4DFC05DDDCF52
C7895FE20377F64D39EC9F9771D83C38BE9E8F79342F366CB9EF33E01C3E83F7BD47D725FAEB10F047EE782EE4FBA9321F7C
80F1336E7C115D6CFAC1B2ADB2FF076EC823CF262DB92D7DC207A82F78037FE00A3F035FE825E364DCF0039C809B2DD352C7
095EC00D7A9DB168A5FBBE25CB838A1772CA536BE65E91F1E233F4A1DEA1B595458FF113C60D6EC13138047FF02D65C0FC6E
320ECB4DAA73E859B50B19BDA52EE8EE9AABD392E4FA1EEA8B6FE3438D3F58B951F6A5A8BBF147F10CBE3755585B59D64B96
43D5A7D61F5C35C1FF9D7A30BF8615574D90753AE703E7E400F0CE7AFEF14BF11BA47EF815FD070FF0CFE6A5F32107C08F37
F3259C95FB22E0A1F153F10CBE0FE7CCC8E1AF3FC1073C02DFF0093C933FC0357DA11FE82C7EC2F139DFDA2763BBC8754B96
A73447803BF006FEC0BFE573E501BCB5F9EAFC393EC7251F805BF00ACEC00F3A088EA877F3C8F9DDE4B975E8C5C9FB29F3E5
3A74FC8F79AEFD707C25B9FE9BF3912FC01778A39FBB0B67F4F6F3343C844796FF954FF8C38B5BD3631EC8DB3A2E1BA7D6C5
C6A3E320C7C047D349E525FECCBC98A7E9F933D3927CDF5EBB6BF41F723D14C725D7D037F8027FECF5E41FCD31F008FED0A7
7323BF9FBA74E956F37BC6079FE017754DE919BE67CD9A3D765CF20B7C649EBC2FA574D2CB8237CB4DEA27E847F5A4B4DFA5
EEE87D5AB18509EE7A607299EA10B98371326E789DE7FCB413927FE137F5C8AE238CDFF6C93487A1D3E004DC5CE9B0E4B67C
3F33798571D09F376EC47679208F284FE8237DC56F4B7F9F72D75F079177383EE7A36ED40B3FBBFEE6DC2B7E5EE1FDE0161F
02BFA6FFDF8DAFE4FB00B986E3705CEA4FDDD11DF8D4F4CCC204D99785B7F8EDF3F75625CA79D1DD1D1DD27EDFB5EBA2F10D
DD31FD7F65F41FBE0FDC2E18BFC1CF2B8CEF409DF4187F1D429FA8AFD55FFB01AED017F0055F9B1F4F7A59F6C1E9D3DACCF8
EBEE3E2AF28EFA18F985F1303ED3F971E17B1ED07BAD87D54775C2D6DB2D535F95CF8BD0FBF47B4B96091ED147F04EEE28B8
77D93CA92B3EC0F1183FF341EFA90FF5026FF81FB863BECC139FA38EE004DCDC4D98B8C9E539CD77D4057ED05FFA8D5E5689
4DDBEDEB3D39E5E9C1DF55137FC3672C1F3CBA6499ECFF81B31F97A7BE2A78C99E930CC7E047F1040EF14BF0983A796E57D9
8FB05CA33AB47673DC32B96E86F9E393C7DF49B99B9676C0FC86FAD2B7DA17D75690F53ABA4E8E78A2DAAC10F93F7E427FE9
2BBC5A9F63F9343FDFD17FE603AEC0998D53C787EFE21B8C937183EB378AC65F977D7678020EE10378046FF83DB843BF1BDE
5F3E4DF6B1B2FBAAF50BBCEBB89BA5843D90EB7F28B4B6827C7E6639577599BED00FF48EBABEFEF9AAE3A2E7F8C9DA1F26F6
90EB7B2CD7A8FF6FAC93DC46EA89BE53D73C95A71C96FD9F3CCF4FDE25EFC71F385EFDADF2BE7D364F708B6E82DFCD53571D
97BE316ECE73F6D15589F2BC07FC927E183FA93F3AA3E367DCE4067803CEC1BDEDDB0E4AB9E1AF2F7EF93A6DB7BFAE8087F0
8B1C03CFF049C6C9B8D7DE1C5B43AE4F621CE40BFAD9ECABB95D65FF1ADDB1D7934334678067700CEE799F1D07DD45B75F8F
5BE6EBF7A5AACBE689AFE2AF8C3BA5C5A8875CEE23D768EE808FE0185EDAFE4BE7D414DF07C02F7806BFF097DC008FD151EA
4CDD6D3CE41FCD13F836E3673EA9C5125F917D10F08DCFA4EC1FFEB35CA76579447308FD872FF0075F2D1EB7BAEDB973E74C
D7D31B7F37CE5DA74A6E529DB3E3EB71539B872D95EB564CE76B4FECE1EBFD6F23BEAB9698B82AB84FA3F3815FF08A7E836F
FA485FC12D7E037EE1357D82DFF0CFF8A8FDA16ED40BDF22A7304EC64D9DE029BC45071EEB387193CB89AA07E800E3E0FCBC
0F9E337FF39D27C6D6F0FDC73EFFA8B2229FF405FF39129E72C35F8F81A79BBB6785F8F986F9D047FA873EEE894F4DF1F328
7A01CED129F06E3E923CEA21DF4FC82F9C9771C0832619DF8D93FD60F86FFC04278A1B7C867931CF948CE977E4F327CB53EA
ABF0BED999C457649F1DDC80677C1A5CDBFE42E5940F45C7F02BCB5B9586FFECFB0C3922FFA085B545B7F131E6C1BCD0A13B
6727EFF2730C75E1FF8C9371DBBE61CB9469BE3FAC2BB7689CE0CC721CFCD33A307FFC1BFDA08FF4D5EA0BAE146796BB749C
8CFBF417ABDB8ACFE037F483DCF1647AD23FA47EF800F8C4D7C129E70327E066E7BE940FA57FE00C3D4A193279F5B66D676D
3EF870DA94D91D653FC27296EA1C3EC3381937FC6E7464D138D9EF047F767C3D2E7E6938D6BE83037051E3FC9A85326F7C80
BED00F7805EE9A379E7E47F6BBC13778833FE06E43DDA59FB95CA17A479F6D9D59726A8A9C17DF0087F82E78043FE820384A
E934A1BFE881E51AF561DB9F6DBD7ABDF4091DA5CED4177D3CFFCB8A7CF2B9017E49DDA9ABD55979707079CA347F9D60BEC6
78C803EAF7C607EAA675FCA9C89A857E5E417FC135F8014FE86B832D4B3F93FD3670804F327EE663EBA5E855ABC52FF08D94
9E5175DDFE16794473C1A6D0154F49DF2D0F6ADFE085E573E5871D47DF8FAF838BE6C327AFF6F5EFD16BE31365BC8FEE1937
4FF645F013F48F7A525FD3D32F677794FD6ED317FD697FD7FA6EEBBD6AB5D48F3AE243968387AD784A7C8C3ADAFE54A5B5BF
F93E60EB93C2C943659CBC3F2B7C61EDB973D3CCD7A97BCAD2889CF2FD46CC9F3CB2E542F250E10DE3A0AEF67A728BE60CF8
024FE0213C8557F0D5F6EB9357BE2438C707F073C6C9B853FF4CA8ECBE2745718FCFA4161DBAD77D3EA7E320E7802FF0060F
E0177E03CFE873F39913FAFB7A0CDFC90DF01E3E184E1437196FCFEB23FE8DEEA08BF000BEC09F1B4D92FEE1E700E6CF7C99
27BE7262F1EAF5A2F7E83EF8333D1F1755D7D7757C9EE3733EF8024FE81FBCA34FF0CF3E1759B6E86FA91F7EB26FFFDADFFC
7C0CFFC129B8A50FE0161D04BFD403FCC307EA491DD141F3857F46E47CC01FB41FD647ED2BFE43DDA8E3ED935353FC1C40BD
E011F3874FF0023E5047F08CBF826BF8C7FBA9476AC8D45BF2B994E526F555F4A059AE84CAB2CF0ECEA8077C641CE0BAF1AA
797D647F1ADDB1FD97DC6B6EF9FE66F92973A8FB1ED366BD862C949FB62F5976F53EDFAFC81DD4877AD9E74F23E7B493F3E0
7BA66B3A7EE6C37C99277ECF78193FEF4747C11538A3AFE01CDCAFAF29FB95D782B94CF997D962E54BF21C6E7C111CDC7B79
DC3C9987E52C9DA7F1075C29CEC03FBE0E0F528B4E5AEDFC93FCA27E8E8F715EC6712F296AA87F5EC6B1F78F35B7FC5C885E
D8E7555FCF08937AE363E903BF7D4FF862F94B7591F1300E7C9DF9D077E38BE2ABD9C8A9B7649F1AFE984F6F5FD0459E5B0A
CE6C5D7D6CF5025FD7C11B7C0677E00A9D055FF48B3EC173DB476EBFF27FE7CF9F379D4E39336EA07CCE643945FDFF72AB45
7F2F5CB8D1FC8CBA1F3DBA7ADFFAF5474D1FC11375B37AA9CE8267F201B8E6BC9C8FDCB0B1CDE28A327E7413DC19EEA9A7D6
97BCF05874CC5AD113FC84BA52677880EF312FE6090FC009B849D93FBA9E5C97C3F8F0F3DD65562F101C812774D65E8FEFAB
AFD3FF6699931ED051F4927A525F74A242AED5B1BEAE5B5E59B8B8A23C87123C05D775AB9E157FC16752EE84E792CFAB2C5F
A89FA3074D22BE7D4FF6BBA9AFE99DFE9DD7D971C8039A0FE011FE0D9FAE85CC699798B8DA7C97BADBFEFB334B42057FE8FA
E66F9766B8FD40D57DFA9F961E5F5CF22E7842D7E10B7C8037E01C9F00EFA66FC5C60DF4756E478755CFFAB9013FC0CF193F
F3019FE012BE814BEBAFF6DBF2E2C74B33DC7381B57FB6BE6A9D5C4BF6C9F0C31F0AAC28E23E3F561FA21EEBAACF79DFF76D
F489FA5017FCC0F2FDE51545E4FC8CE3C7C5C9B5C4EF9927E70147E00ADEC003700EEE7F9DB3F27FFEFA131C9BCE571A5DCF
D77B786A395FF97A2B7646989F0FA8FBA13F57C7FAB91F7EDAE71CA3E7F5913EE31BF0051D8437E63779C273F93A8DDF725E
C661FD03278A1B7C86F9327F7807FEE01FFDA24FE8233807DFE029ADE894CBF21C3DCB2FEA7BF0055F8437F0B56993F8E272
FF24F8B37DCC12C943FC1C0EEEC0217DA69FE01F3ED8BECC95D5397DFF619C8C0F5F343FDB31E77DB9DF0C1DB9DB3B66AD9F
03A8BBD5179C286E4C779897CE135D040FE0031E5177F8C4F9190FC7855FF08A7EE0FF4FFC915058F0855FC11B7844FDC817
CC8B795EB8BA2474D9B21FCC2FC10D758577F0109E326ECE637551DC8263DE076F994F6A7CCC4AF91CCB7293FAF986F2F35F
96FFA353E8073CC28FE1137EC8F8990FC7E7B8F8B7F12A257497CFAF8C06B35EDDBDFB4A3017A87E1F7876754E3F77A22FE0
B769E694CBBE4FDEFB6FE4472ECF69BEA3EEE01C7D07EFA9CDA307C873A32DD7A86FDBE75531B14B641D865F8173F402BCDB
BEF3ACA52BE5F8F8891D9F3CA07E6FFB053FAF2AE9FBD2C953C943641F047D04C7F4857E98BE6BDDA817FA6E7E7F78BEBB7F
0BBE5DED37AFCFDCB9E9E6AFF40FFE583DB5BEA945235F97CF992CA7682E00FF96879507F67AF283E600F8471FE92B7D6BD6
3966A5AFA3F827E364DC9BDE5B70D6ADFB5407C129F363BEF018DDA59ED437B5C4F01CEE73231D1FF9C4F6FBFE4C5EE8E7E4
473B440D151E3F3A604C84AC23F1258EC3FB539B87EE92EB597E796955493F17A29BE84AE30DB35E95E72C5377D341FD3BAF
4BEF383DAFBCCF728DEA3D756A96123DC0D7455B67AE48BE21FDC1676C5F7ED377EFC9BCF0879BB9130AFB3E4ADD2D7FAC58
70563E67064F860FED17FD03FFF80D3C804796B7954FE01FFEC003E6CB3CF10FFC967132EE2DB3177772FB63AACBE0029C80
0FF8439D9A6546BEEEEBE5FA87922A88EFC22B74E8DC274B574A3FD165F0B4E7BDE41B7EBE41AF6D3D367C655FC9B9F80C78
070FE0C3F2D6478B3BC9E7A2D9FD96FA836BD3CF53C3733CA0A3F4033C283EC0337A07AE8FE74A5E28EB28F4057E6E2BB974
81FBBC537D82BADEE91EBBC4F763FA81DE531FEAB573C8CABEB2AEA32EF8907D2EB2357188E0113F499B157B42F6232CD7A8
2F51BF9AF597AF101CE17B96834F2E5D2039913A1A4FFF11F991E417F84A7FE92BBAB97DE3F215F277C6C7BCE1659384E979
659F1DDCE3DBD481BAA0F7C55A2E9B26B8C767E01D7E09FF982FF3C42FE99FE1331BAFE109FCB67DA2CCE491BE2FC18BEC7C
307D615E3A4F74113C800FF38FBF932AC87E393A057EC13F7C00BFE019FCF23E8EC33898277806DF6939C62FDCBAF5543047
A83F5FA9F7DD7B0B176E32DF027FF722C744F87E4F9FECF3A11371A1F2B9107EC5F1392E3E0D7FA92F3CDEF0DF39774417D0
1174C2F2869E9771C023F8433FC051D34F624FF8BCC5DF9817F3CCB83833BFF008BCE203D41B7EC13778002FE803BCC3A7E1
5FEAD831BDA5CF9653D4E78F144D1EE9E766F40B3E521778C971783F39005EE007F003DEC137EA05AF392EF3607EF094F7E3
63855E5C70468E8BFFC00BF40E7E504FEA88DF181FB42FF4E9CCF265D3C4DFD17DCB0FCD47D6913C6179477384EDBFEC4F7E
C2F7375BC7D6189B2075C30FC1157842A738BEE513CD19D7EB240E91EF37C5A7C10D386A7A7FFC42DF07CC3FFE31E78E3CC7
163EA7BE34F421B9EED9F28BE616784FDFC1013CB23CAC7C42EFA90FF5C29F1927E3E6FC8C073D402F4B3CBFBC9A9FCF3757
981FEBD65FAA5FF00FFD685C68667ED9D7A67FA69FBDC6F4F675D474475FCFFBD32FC4DD927D10CB59EA3F8603EAA9F53D58
29F9093FBFA2838F9C0A2F2FFB218F5C1C51D3ED5BA8FFD8BEFCCEA405321F7C00DCE207E0D7F43065641D5F171927E3C33F
6E978A0BF5FD98BADB3A7CCBCA66BE6F58CE38343F563E97061FE0169E805F708B8F82DFF54B66DD13DD01AFE8137EC87818
DFC5030BCEC87310D06B70437FE92B7CB0BAD347EDEBD6F0057FBBFD2ED55370919E7BE24ED95F305F55DF00BFF48BFEEDAB
BBB299D491FEA1A7E0109D028FA7CE2EAF26CFCB87FFF084FA514F700D1E9B148EBB25CFA1064FE83A75A02EC1796ADF1507
D499FAA25FB66E5CB2A2AAF8147E451D983F3E878ED619B4E06FF93C161CD8FAADDB74F7DC11FCE1EE4F63137C1FA5AFB60F
9BB1A8CA850B174CD78D873A5FE6FF639EC555DDE788AAEFF4CDF8B269E843321E78035FF027786375F979E8433E2F6CDFEA
7FCBFBFBF976EFFB2BAAFAB904FFC09F9917F34C7B7B6C823C87C57287FA217A46DFC101F57E2572D94EE1037E08EFE00FFC
B3FC7A677155F93C967E980EBE3FEB9E3C2F189DB29FFA77EA0BFEC13DFCF929796909E121BE48BD9917F3C1A7196F8D164B
4BC8B818DFAE01CB764A9DA817C7FBEDD5A405F3E66598DF8063F8048FE015786FF2C7C49DBE8E823B70086F52D78EEA2AFB
239647D4577F2F34A2A6CB55E42CEDDFBA5ED3B3DCF7E9281FD06F5B67D69D102FB8C6AF362E4E4A96DC80BE98CF725E3D1F
BE8B5F315FE6FF7BBDA10FF9E3617CF0917EC34BF8855FC2B3D40F46D494E7D7588E501F36FE5037AD233C823FE0E0C4C3CB
FB4BBDD023F40E9EE22BF095BA512F7CC5C6431ED07C40DF9AAE1A9BE0EB317D865FF0CDD6398FCD9922E7C5DFC03BF8A7FF
F0179D85C7E71B2FAAF2FDF73BCD27E0217DA47FE8237CA4EEF0D2F671AE2DEFEAAF93E00B3CA1BE37274E2F23DF0F818F82
035B4FB68FE822F9143F04BFF55B2725CBE7E7F0193EC12FE64D9DE0177C33FD7F635457DF07E029F3E178BC0F7E321F74A2
D1D8E959F23DE5E009DD0227E006FFA1CED43DA3D5E483C267CB1DEAAB864BC503F860DED401DDC21F9EEFBDB8B2E00F5F25
175067EA0E7FC919F098E3312FE6B9F9F4DC06EE7A0C7284F2DBFC62D1889ABE6F302FE6835F1ECDB7BCAB9F77D1355B2FD5
8D9F26F8C197C039FE04DE8D27F451FB7A67CB8478DF4799FF86BA33378AAF9BCFAB9ED97AFEC2F2A77C3F01CFF0045CA7B7
9AB0E6871FCE047390FA98E96CF1B90DE4F36DF07AE5EE9C298B166D321D0707E0DC72ACE29D3A535F78887F322FE669F3A7
FEDA0FF0D878D9E483B2BF0C9E183FE3C6B7C00FFA058E7E68327FB28C031DB4FE6B3DA9237A77367A7165D10F74043EC003
FA0E0ED07BEA43BD6C3EA786DCF771455DA9333C385478F953FF57EE341C5307AD4B5ADFA818D98F30DF563FB475D7A84903
058FF809FA5AB0EFBC065217F4FEDE9C882EBE5FD2277054E9D4B28EF27EFC8AF39AFFAA2FA25BB5D3E74F96BC069EE0057E
063FE077F1E68B4FF8F9D9F471CBCC8D727D80E51DFDC9DFA90B3A44BFE89FE501F67B749EE0A8C9B2096B7C7D4AFD21A28B
DB8FC0E7D56F7FDCBB60A9E017DD0747F00E5EC13F3B0E7EAD7EB9BFC4B28ED24FFA8ABE834F70097FD69D9CBA49FE6F3EAC
FA7A7D4BFC34791E3D7E03FEE8CF4B7B9776F1FDCD7CA5EA82A5C257FA94DA7CD82B6EBF003F571FDFD420F16BF7DC24FC53
F506BC834FF04A9FE93BB8DE317451A4BBFE41FD86FAFD5E6290BB2FF1F77A617FFBB9C8F0AD7DA14F3F575DDAC5CF0DF899
AD330B8F1E2178C707CC6FC644C5F8FA4A3DAAAD5CF29CEF93FCAC1EB628528E47BDA83F75C70FF01FEA4F3F76BFBBE439DF
A799F7BDFF0E98FDC07C75FEF014DF82AFF0149F80AFA753169F101DCAAE47F014FCC1575B17354B3824C7C17FE01DFA08FF
4C876B4774F1F5F8D298790DE4B9D2E83E7C0327E0033D852FF0047C8317F00F1F8C87DA2FFA074EC12DF882F78D4A4EDDE4
AE2B52DC64C44F5A23F76F581E517FB3BCF8D2927ABE1FDEEE3F69A0EFAFF415FCD6DB90F8B57CCE8F1E05FD69D82BBEAEB3
2E7DF41F61F964FF033F349DD271326EFC817A525FCEC379793F7A46DFC1017CA7DFF09E3A307FFCD2F27DDEF9DF5EBC78D1
7C129EC24FFA873F531FEA659F6BAADFC33F78038FE81B7D86BFF0D9F809DE147F5B062465BAE795AB0EA243F09779C26378
07DF980FBCE6758C8BF3C353DE9F3E78DCB7721F8BE50BF5D50D0B6744CAFFC1377A6CEBB42253AA08CFF0379B0F3851DCC0
47FC0F5E9EF879493D3F67A3B38CC7F61BD477EF5E183DC2F769FA040F1A779EB4C6F7B7DF52130EC9F7A5594E51FCC10BF8
9F9D1F65CB2D69E1AE7B527FA32FF403FEC3CBBA914999F2B93AFCB1FD9A53830E073E670BF0837A524774167FB6FD249D67
5AF8E811B21F61B943FDD9F8AFF365FEBCDE7C5CFD18FCA3A7F0E0FCADF9DFCAE7F1E81A7CFBE1F4DC32EE733EFC59F188DE
534FEACBBCA9033C86AFF4977EA7AE0DFFC8ED239017D487CDFF468CFBD6D7457BBDED0704FC1B5DCF1F9E7858E6856FD8BA
A878F418C1293E73F4C892167E7E4597EFFD372C9FEFC7F4031D6AD06646A45C6740DD4DA7F4EFBC0E5EE067F003BD29BAED
BB143F9FA37FD49F7E58AE2ABBE431DFDFD6579A122B7E0C6ED07BF389E273CBC8E7E1E0D2F286ED0F05E6697E3368F4085F
878C3F5A7FFA012FE009FC0027E0039EDC3C3FA58A3C1FCDF282E2667BEFF9E55D1E529D055FCC8BF9A0E39B56241477CFFF
51BEA143878E2E79CCCF61E8ECBDA43E6181796AEED0F9077536FC235F9FA8137503BFB6AE1B30BCB1ACB3D17BF00B7EC093
E58FD28BFB4AFFF1257CC3F64BB42EE871CDE8F9E5451FE89FE15EFB429FC039BA09DE7756FCEEA2BBBE101FD3BADAE767EA
4FF0807A5247F428B3D677296BD7FE62FC87B719DB6326CBF3CACD57D51FF61F5FDCD7CF31F80A3EF0E4A5A96E3D8ADEC39B
86FBA7C4CAF305C099F977D14593649CF85BD6A0C4C3F2BDC5E8353CC1CF5E3DF8DD45D153EA08BFD03578B6FBAF859D4567
F139EA643E9A736167E92FC78577E01BFEFD9CB968D2CE9D176C9E1C8FFAD8BE82FA163807DFF004BDA00ED4C574F3BD84E2
F2F93FFA62FC04B78AE33BA7A3C7F8FE07CED09517662D18EEE7DB47FF3160B67CFEF2E86BFD93DDBE85FA12E7E1BC8C0F3C
824FF801DEE11D3CC407E817FD43CFC01538B3F9648655F5F9887E17B89D14277DC4AFD2C7448F91FB432C1764CB03B66FA4
75097EBE18F05D78BCE5B7C470D111740A5DE1F8B6AFA03E0AAFC101FCDEF0C7B4B70557E6E7AAC7F0111E8203F3E99A3193
7D7D05BFE0997A5357F805DF6CBDD4757C2F198FF90C9F4BA9DFC0BFB4BE917D44AF2C17A8DFFEFAC982E17E9E46D7787D70
3F20E0BBF7228737F67D91BAC35FEA028FA90F7541DFE1233CA45EA60B5A1FEA052FE12975A28FF40F5DA34EF014DEDE2838
75A53C0FDDF28BE2151D60DC9C2775ECF0366EFF82DCA179C378D4A14F98F0013E711E7481799AAF3C9D182ECF5B82E7F014
9D82AFF873705F2A50771B0FFB1645C3AA0A9FE077C9E8057FF939D6FCBB6FF418DF4FE017FE01CF2E3F9E1427FA6A794479
CB78193F3A818F514FEABBAD7DD232A77BAA77F0CC7252EE696FCB730AE8377A000EC1A5E991BE9EF7AFFF77ECA772BDA3E5
05F5A5A05F46F6F1F5D2D66935467E2E7DC60FEFFDD83FD9CF07D41DBFCA777B460EA90F7E6BBCD57A525FFBFC4CFD155E9E
ECBDE02F3FEFA29BE859E17173F2B975A8FA15B9C5C6A3E363BEB6BE57BF415FE917FD335DED35BC8DAFAFE84AAD3949CB04
BFE07273CE197B05AF960B54CFC05199D10BCBFA7E787BD2F85EBEDF8003FA4EBFE1A7E14FFB429FC02DB807BFDB6FCE9921
E330FF541C651C1EDF4BF4D0FC567DC0D6571D42B2DCBA5CF59EF9315FF08B3FD8FE87CE073C369C17FBA9FB9E15ED37F5B1
75B0EADAB17E0BCBCABAD4FC50F5DAF047FDB51FE01A3C800FE665EB6FF5A1E0E76701DF02BFE0137D07A7E61F4FCC9921F5
A0DFF4ADFCD805F77D1F331EEA7C993F7E92F7FEC4EAD27FFCE442CF39F956ADDA65BA006FAF3D3E23873C67193D0597E01C
7D01EF3B17CF6B2AB913DDA7CFF4E7F5BF66EC95E730C1FFF4C1637AB9FD00FC597D8CFFF37AE6F3E89E9EFF96E33C7AADC7
7F24D7A3D3E00A3C81D7BB3F8DFCDCF733F071B8D782FB7EBEC157184F705D1EF04BF8028EE18DF956EB794DA53FF4893EBC
18B320DCF74FCBAF2DE636F473ACF962F9F1BD7C1D4277A93FFDD8133E7FAECBBDEA43D41B5C8013F891F6F6C8CF65BFC0FC
5F7D6F63C5B83CE2A3E01EDDC59782FB31813A1A1FD517E1E581EE0BC2FD1C835FDA7AA6C6CC42325F7C093D036FE06F6BE3
8412EE7953F8BFEA16FEFDF2A0F9735D9ED7BA53D72A63BF4B77EB0B7D1DFC83BFF0F9975EDFA5CBFCCCCFB57ED4C7F280F9
F1D07FBAFD087E575F361DD17ED13FEA491DD1775BA7D51DD55FF635F113F8024FC0B17DBEA57E035F83FE34A697AFBBE800
B8440FEE258564F9FE475FA933F54507C12F78064799C9731B4A2E4117D03BEA0A8FE0D5AD1C13ABCBF3B0F00F708CEFD978
747CE0B44E464209B92E12BD305F5935F2735F77E111BA039FD089FA07E3F2C8E7E2F409BEC013EA0E4FA9177C5DFFE5C4EA
8263F367F525E3A7D69F7E305EC68F6E59FE2B3FAF949F03C12F78A6DEBC8FF1331F78841FC027EA0DEFE0A1F953CAD07F06
AE9F1FF292FCCCFA696621F91E62741C1EE293F48BFE6D5B9958D3DD7F87FF2B2FEFE5EBF19F40FF743F46FB0A1FA923BC84
77F08D7EC0575E473DA8BFED4FA88FDA7A6F4F68295937E09FCC1B5E733E7405DC82635B47ED8ADD2EF5C46FEF2674ABF7C0
BC749EC1CFE1023E0D8F8D878A4FF0BAB9FAB450F7BD45E40BD557CBE5F312B22E5DBA64FE89AED317FA94B17D6C4BB97ED1
7284FAE4AF87E795CAC838643C44EF8C27DA7770607EF06E624DF97C1EBED15FFA8A2EDCD932AABFEF9796A3C81353265677
D72F281F1867701D1FF0515B476DEF1D239F77E11B960B8BCC5BEDAF6FA81BF5424FE1113C874F3FB5487ADCE50CD5CDECBE
1FDC4709CC87795307F86A7AA1F3629EE9AD467D24CFC7B45CA0FE693824972A1E79BDADCBD557C01B3A0EEEC0237DA4AFB6
DE1815F5B5BCCFD61DE7E7ADDEB0E198F51BFD45375EDF312D54AE97643EA647FA775E979663447B793EB8F99EFA5BD0CFC6
B6F475EEC62BB1DBE5F947E8387845B76A24273D2ED73382B3CBFF4EC892E71EA12FF090F366F7577881AEC10FCB5B45E6FD
D7F7B77B2F8796F2FD8C7E6FAC1DBB5D786F3EACFABD2BCFEC78773D277EAB78013FE0063E186EC9998A5FF353F64B741C5B
FBCE4873DF5FA7FC81F7E82875A38E969B1B270E9679E263C62BADFFFFE737CB467DE4EB0D3C820FF0893AD87A5AF5F4D8F9
79FF15BC801BFCC3D61B47E22AC938F019FA497FE10D3AFDCA9DD9F1B20F070ED02DFA0E0E2C97949897E5FBED9E7373327D
3FA03FF0009C837BD3D9FB23DAFB3A67EBA8E861EFCA7A1D9DBE57B3778CEF4FF42FF8B95140C7E1D9E153F3B2FC1C83AF98
3F3D372753D66DCC9FBED36FF4CBF2448579E37C3FFCE5D0DCD3D217FC89791BDFB4BFF4DB7CB4C2DCD3BE4F1E38346F9CEF
F71CCF7C99F9EAFCC17F76DCAFFFF7849CE2AFE6B7EA33E847FD1F62B7CBE7FAE0F276FFA8AF7D3F00DFA6A76366A4C9F58F
E80EFC0597F0F84246E260F9BE49F367D5417004AEC09DE99A8E87F131FEE0FE41C087383FE3E1FDC61FC50F78BAB3EFCBF2
81E720A9FF6B1DC11D7C813FF0179D82C7F0087D874FF84F707F2550475B475DEBB7D3ED07A89F6CBB30D3E10A3D42572CEF
E6F876B39C0F3FBC562EAE923C97C77284F2107EC12B70C0FC9B170B6DE9F328233EAA8A5C77687EAEBE077FE91F3C8607F0
82BED9E76AEA8BF0121F7BA2C8B8BA5227FC7073D3B84AA2D7F007FD663CB64FA0FE8AAE801FF0447FE02F7C0EFAFA849CF2
FC23F8061FE121F582BFCC131EA3DFD499BAC34B784A3D2CAF3F37ED03E917FE969E3BE20D79AEA5E505F5E1C73A76AB27FB
B18FA5756DE6AEBB51FF611EF01DFE8313F0811EDDED3DEC5DDF17C11F75A6BEF8043AC33C9907E7415FA817E3B77D05F575
F389C233D35CCE545D301DD1FAD38FCCFBDF6E96EF57B31CA4BA9F1A3FACB5E421CB0B9A13F05BDBD7D179F27A5E977226B4
A5DBCF875FE44FE599F96BE7A82ABEEE5A3ED6F1A31F3F8D9835D05D87435E505EC2BF62A793D6FBB9149EA20BF0D572585A
5C25B9DE93BA984EE9DF795DF073A0803FC11B5B87C4847F2EFFC7C736569C384AFCDEFC5F75DD7CE58F88377CFDBB797E5C
5DB92FD1728AE2FB5EBE7E3B7DBFA4EEC61FAD1B75847FB64FA07D05CF792A7F595EAE1F01D759E7A67D20DF9B627941F90C
5FD04D78637E3078D640793F383B533C69BD5B07911754AFA943701D1FF009CB0D9D87B5F6F5CF7202F3D5F96F8D98F641E0
7BB2D4E7559FC083AD8F9587F4BFC4BED9DD7CDFDBD524F186BB7E52F5171CE103C1FD92403F2C0747CFA87FF9F265F319C3
ABF691BEDABAE2D284D3F23EFC873A5137707DB769D766BECF317FF4CCF609B44FB65FA93806EFEB2B8DDD27BA643EACBE81
8F3DD63124B7E8083E403D6D7DAC7A14F4F589A3DC75048A33F395F4C41BD27F70805F337EE6639F03A9DFC08F5395677713
BD323F579FA08EE0017C302F5B7FAB0FD1CFD20766EFF0FD10BC834370B9B744520977FF9DFA0AFD37FEEB7C99BFE9CEC069
1FC8F592E8C29DEEE10FE838B8B1BCA5AFE7B8E00D9C81D78C56A39F101F355F55FF41EFAB9E4A2AE1F2ADD6FDF88BB377C8
7ADE7C557D085F7A6CCA37FD9CDFABBE528F72A7661FF1FD93F3DA7A5AFDE372EF19F5E5FB00CCCF5547F6E59F9DC3F71BEA
47FFC103B8C6076CDF42EB051FC12BBCC4E72A67CDCE21F773337FF3A7FD63F705AE370CF06DDBEF33EA3BBEE3B7AA07474A
CC3EE2E70CFC927957BC347B8AEFABC1CF5D023A0E5F6D7F5CF50EDE6FAE1E7B58E66FBEADBA7EB0E0EC29F279A1F9B3D629
6D56784D796E85F9B3FAA7ED53A8FFC1CB1B05279C96E70A59AE51FE731CF35FF557700D8FE095E5E08AF1E764DCF89EADAF
4E447EE6BEF74E7D2CB568584B79CE23C7C5AFF39CFFAC94E0F8F16A9F16917D1A7C031EC11FF06D7AA1FDA27FD42DB81F10
F003F4CCF627147FE01DFC832FF3B965A39FF075F76E42486EDF17C11F38B1FD00D53B7403FCA11FE87D707F25D00F5B3FDC
9EE4F0801F9A8E3F32A3BED405BD303EAABFC24BFA037FE1B3F93BE3D7F9C03BF846DD2D377C125ED3D7CB605E893D2CD747
D20FF804BFA8ABE5421D3FBA75E160FC39B93FC37284EA297A42BDD0950D7F4C784BF06F3EAFFE633AA8E3617CC667ED237D
FD6972FC39A793F8B9F29B3AC16B780EEFE01BF532DC5EEA5454F418FCC27B74807AA00FCC87FFDBE743EA5BF006FFCCF3FC
905AF23AFC8D71A1038CEFEE97DFF4F37D9D7EDBBEB3EA4B76BE3E7739A1999F57A9B3ED13A82FDE2E15F999EF43E0D87285
9E97715077F0039EB6F69D7255F065B94075EEDAE3933E58B06063308F286F83FBF501DD841FF0D8F615144FE6135FC59F73
B951F9B3EE64D43CF175F379F543F060EB6FE5B3E13E33ACA58FFFDBBB3F2D12F83E32DD3FD0F9E33FB65FA2F5821F0D724F
784B9E37042E6F17EC54D43F0EC73D5720A1993C3F0E3EA0BF8CDFD6EBEAAFF8F06369BD56BAF5B7EA3D750EAEB303FAB8EB
5F09CD5C2E212F280F6CBDD16D4C53A9377E65F95BEB098F2DAFD7981C2378C2C72C2FFE38ABBDBF4EB27D01F521700D6FA9
0FF5421783FB1C017CA45F88B82ED7F359BE605DAD7F371F543FB3FD56C53D78373D1A33E5AA5C2F898E981F2C4868267505
4FC17C16783DE34E2B3ABCA25C57C7F9F0377C2CCFC62F5E95FD5474FACED921B57CFF0037E0A251C9A879BEFF707C5BEFAA
8F9E796D567B3F77A2FBF000FC80A7BDD567B577F79DE1ABAC4F140FB6DE553E6CCE39314274D67C5575D772557A6269DFDF
F0139B97CED3F0A9FE044EB75D98BAC97D2F23FEA97A90F5D3E418799EA6F9B6F21C7C501FEA659F97A86E825F5B3F681FC1
AFF97AE188EBBE3E996FED98D55EEEFB0107D481F9A38FC17DF349BACF13E081ED7B951CDA4EC685AEDF9C38A6A9DC5F67FE
AFF83ED524B1B4F48DFEE143FB2A2496DEB1E3BCF90D7D835FE0189E59CE589F18EFFB27F36B9A39BCA2AF67779BF65AE9FB
04FD03B796C714BF9653E7C5B979E227E67387124B4BFDA823FD0DAEA7037A74BC6162BCCC8F79E267E0143C83EF839512E3
A5DEF81CF5B07CF053E216DF57C179767C9BEFEE4F8CF7FD193C8107F0617E4A7DB45EE83A7DA4AFE853709F2380E3233512
B7FC5F79C2F2CA5F132302D75106F0BD61E1B8C7C577CD87D5AFD0E3A7EA46FF11D8770EF8A1E972E1A99B64DF055D03BFE0
191E987EE9791947305F068EC3EB180FE3C087D001C39FEA81F14AF106FEEED4FFE255DF2FA92F7E95A772DFEB6E1DAFFE80
1F9ABFEAFB7F1A31FD7DA757F8B6EA87ADFFD557E01F38825FF08DF1351F3EA4A5CF0B7807DFC0D3E57FC7C5C8F72858EE50
9D5DD76BF420B9CECFFC597D091EC00BFA70F3E8C7C502DF6711F07FE6B9B5716CA45CC760BEADBA7BFBE4D0760FE40BC51F
7CA7EEF09EF1983FAA2F06F7A3033A959DC7851B4F6BE9E743F0D8A0CDB8C703F735047842FDE02F7C86A7F0933AC2337847
3DF0C33CF15F7776CF87537F48EF187ED8F7737C15FCA3EFF0001D60FEE8013809AEE3037A747D4BF41FF3E6A507F385F216
BFB23A6BDDE98BADCBD50F18A7AD9BD5A7CD6F064F7F5FD667F013FD617CD483FAA10FA617FA13BDA0BEA617DA17FA84DE04
F75702784D0D195A52F06F398275B8FE9DDF53860C6929FA6BFB05EAD3F0121FCBFBF5884F653FC3D6F7BAFEC9DBBD6301F9
FC1CDF08E6EF401DD01570D468ECE841BEEEEE6A32639ACBA5E411E5F7858C692DDD3A51F9864FA07F4F371B1F2678C587D1
63FA4BBF2DE765C44606AE330DD4D1F44EFFCEEBEEECEB7BDDF745F0013F8AE698BECDCF939B564C784DEE73031FE837F36B
92107ED8D745C3BFFA2B3C08EEFF06F8092F8D875A37EA882FD93875DCA6FBE933A6B91CABB8B41C307268495FFF183FE346
F7B7AD9CD25CF4C8FC5CF5897E05D7D301BE818727CA7E5C4C3E8F05177B4BC44FF6F5943E83D3E07E40A0DFB61E561D0777
99F7A76F93CF7B2C2FA8BEDBBAE2FCC00A721D253E6038D3FAD38F3B9DBFEEECFB10750CE6FE40DDC0EF6FA9E3C3E47928E6
F3CAB75BB1233E4D48480EE602C59FE5B3D333EEFA3E096ED12FF0BBA16ED4CFC277EA8E9F986FEA3819377DB1F5B4EA14B8
AEF7DE84D7E4398FE0C67CE854FC64976F151FD4053C800FC6C338F095E0FE410097F0E0564CC702BE9F316ECBC7E1317BAF
5CB962FEB32FFFCC26C213F35BC5053C64FCCC073DA6CED4FD4CF1197745D7CCE7D59FCC67753C8CCFF4E8DD29CDE5FA4DF8
6FEB873F225A497FF09960BE8E77D72BE1B7E80A78001FE03CB8CE0EF0C4F2999E9779E27B7936764F74EB63D5DD8C8B916D
64FF8EBAE387E6A359339BF8F902DE04F70F023CE138B63E56DF32BEA0EFCA9B83056726F97E437F82EBA2007EE0C7A9CAF1
87FD5C859F592E3930B3ABEFCF3FB588DBE0BE47037F565D313FDB12F5B35C87076FD3D2876F91EB81CC6FD5C7CC472FCD4C
F2F3C596DF269691F19A6FAB1E731CF341F53FCB0DA7662EF3FDF9488999CB1EF061ADC7F1176776F5730FFF4FC9087DC15D
DFA6C7C54FF1C3BC0F75A9E4BE0F5CF53EB8EF16E00F3A71A7C8C00ABEDF80FBAC41317BE5F91AE6F3AA9BF028B8BE0FE0CF
F8AFF5A4BED487BAA0FB4F966EDF42AE137C725ABB4A824FF41B1C05F71B027C83D7E817FCA66F8D0B45B6F175175FB29CA5
F30CAE1B0378065FE0C4F60954EF8C2FEA97F006FEE213F0D8F2EEED49E7647CF8EACDDC11AD9292524CAFE1B3E97E72DC06
F12DF422B89F15E8137A60BEDE64F8165F476D1D7269706B9793D4C7E057707D1FE81F7A50F7E9896502F7CB04FA47BD83FB
0401BEA137B65FA2F8B0DCA9E3476F82EBDE005FA8B7E1537D1A9C6E6A30F64AE07E7DCD05EA4BE089FED2EF3BF5BB27FA3E
4ABF393FF3629EBBF24CBDF6408E503D08EE2B05EA8B7EA0F7B68FA2FDB3F5C0D6D105A4AFB63E563C3F35E9C34E723E700D
AF83FB0A817E50D7E0FE4480E7C1FD85803E503FCB0FCC57E74FDD6C1DAF3E1ADC9709CC87E3719EE03E41607C971F9F744E
BE57C0F291EA3B75076FE06F5BFB49BB657D62B940F5EFD6DF5D2AF97ECCB82D5F8E9B9C26BCC0CF396EF3C6A12FF878B67D
5EE5397C323FB833F59AE006FEE0638F57EB79C6ADCBD567D62F193D43F284E50BF537FC9971326E7062EB7ED505F8856FC1
B31B9FB573EFC3FF392E7CAAB761EC1577FFA9E298F1D8BA59FDF27AD70F3BF9C7E1B8B70B0E6EEDEB32B8A4FED41D9DDD73
6EDA2F81FD3BCD058A77F86AFB0ADA6FCBEB5A4F787CBDCEE802723F183A05FF331A8CACEEF6C3F179F549DBC7571D841FA6
2F5A07EA72A1E7E434798E3938B3FCA0C737FF551FC5C71E3FF6793BF7FDB4AAC7E8BAE535AD57EA9F43D7B87D727DBFAD3F
F5EFE69BEAA3962F5BC4BDE3E74CD3B5399376BBFC468EE3A7FE9D79A2732F3F37ED173FBFD0FF86EF8F9EE1AEBF565EA1D3
4F540BFB52EA879F585EDF19B544708FEFD93A191F52BCE39F4F6476785EF4187D45B7E8237DB5FD18E515BCFCE5D074773D
B7F9BCE268CB80F15B4497C1193A0DDE6C3DADFCC94C8E7B47D6F9E6E7EA2BDB6F4E7ED9F10BBF55FDA01E8D378CACEEEB0D
B8A16ED4F1F6EE9E677CDD0707E0DAD6F18A33FA405FC03B7D7EA1FCD449FEFA89F9324F74163DAE5261FAA77E8E009FF819
38359FCE35748DAF73E837E3673E070E4D5FE0FB047DF8F5F0D4496EDD813FAB7FD077CB51DAFFDF5E8D5AB268D1E6A06F2B
FFA9BFADB355479877C922D3C607AE2F0BF821BEF16285E9EEFB3A98BFE158FB951DCFE823B8B6FD3EC5193C036FE00C3C1D
3E353DDEF727EA012FD017F881DE3F7EA96F4977BDA2EAF4C6C551B3C517CD27D507E06BDDC8F15B7C7D3A797EDA78B92FD7
FC59FDF2E6EEB02FE5FA39CB05CA2BEA54A6C87477FFAEAD3FF467F912D3E37D3FBF5DE6F376BE2F8203F058F389C92FBB1C
A8BA73ECFCF47A7E3EA01EE00E1C826BF00EEFE021FAF1EC63E39EC9CACA325F35FDD23A5017FA0A6EC1B1F9BB8E9FF95067
EA8BCF310FE6C5F938EE1B450797F679813E815B707CB36887E77D7F651CE8BDF545FBB473F19437C417CDCF556FF0558EC3
717FDBDA3AD23D875A7D98D7A1C74FA60FEF2CB8C3DFD00D708C7EAC7F28F22BF9BC9CF9E393F01A3E837B5ECFEBF0BF1F4E
C7B4159D859FE82EBE97F7A1AFA6067427E003F0129E820BFA8F6EA023E0BD7EEBA8D9EEBE34E5B3F1487D053E5DB93BEE19
B9DFC67284EA6CFAC011CF094F18377E8B9E50177485D79B2FAB1FC377784EBD528B85AEF07DDBD6BDFA777E5F9B39B8B4AC
7FE93BFD46EFE01F7CA44ED403BEC37FEA4FDDF1037CF589FE9FB8E76BA2F7F0EFD5D653DE70B954F90CAFA83375BF5DB06F
49DFB7C0317A42BDD015F88BDE65E771C1BE13EAC9EFF839BA815E505F7489E3520FD62DF99AFEEB0DD937C24FA807FAC371
C14BC3BF23BFF27D8079530778BF273C2E5FE0BA2FCD0BCA4B5B3FAB4FC033F377AB4FA05E379A0CEF2CD7B159AE51BD8017
F82EFC403FC02D3846876A178F69EBEE07D73E315EFECEEB985F938811CFF9BA88CFE4BD3CA08DFB9C427D6FF3E9B15F0A7E
2D47A8FFF0BE66674257F8FA7769CC847AF2FC3ECB2FAAE3C64FAD3375371CAADF8347FAF8749DD6917E3FD1AD7C2347D612
7CE09FE825F5A4BEE62B83E2F2493DC1F7ADBFBF9AEAFB167D62BECC133F008F45F24EACE8E7CCEDBD27CE111C581E519D03
27E0033E9B6FEA791907F8870FF0807AD12FFAF74BAFB831BEFE8247E6C77CC1AFF55BFB45FF6E8EFAA4EA033EA4E3DBF0DF
D11D8567CC1F1F387F6B624579AE83F9B3EA327560FEE88EF99D1E9FF3815F70026EC0EFEBC5C77EE9EE6353DCE00F79BB77
3BE9F687581FDB3A36A0EFE0171FC4FFAEA5FDEB0DFFA7AD77757CCC9779E267E64363E3C6C8F9C193E5BF6D937EF173E0B5
9091B5E47E27F37FE52D3C645ECCF35681016D7CDD04C707BAC725CBF70DE01FF419FCA32FF0C0742A7AE21CD99782FFEBAA
8F6C287E6CFEACFEC6FF793DE30357E0091C731C5BD7AAFF65D69AF48BF81638C087E00B388637F8CC13AD3B77109CA273E6
673171C9526FEA9E3665D8236E1F5BCF877F5A5E9935F955A903FE8C1E3C333AF28C9C179F41772D6769DD39BEADA7F1E3D2
8393DDF7B6E9EFF832FD070FF0EC993EADFE90CF5BF2176E99269F83A25BE0B7C13F467774F775299F0FF78A8BF2FD8C3EEF
AC3869B6D41BBD4057D0E327AB0EEAE73E67C167D4DF9E9AF8EF579C7EAA8EFED064DC45D16BF37FD56F74117C8257F807AF
E1F9AF9F4C7E55EE5BA3AFF8A8F151FD095E9A1F8F8D8BF2738AE5A7E8296FFBB9E058BFB8677C3FA61E961B46C73DE39EBF
A5C73DD97BCADB725F04E3E07DF08A3A67F79B463B463694EBDDD141EA435DD07D7807DFE0037C07C7F0FE564CB793BE7F80
33FA48FFD047CEDFF4CB618FF87A0CDEC13FF832BEA82FC21BDED7FCF8E064F9F966BE909703DF0B197946EE474277D059FA
69FB19CA5BF389839366CBBA1F1D0157F82EF8C2AF982FF38777E811FC43576AA78FBBE8EEDBCBD627FA929D17F0817E807F
F41D1E305EC68F9E6DAE10955F7C8BBEE263963B378D6927E3C36F6DBF4CFB8AAE983EEAF899CF8DBD83FAB9EBD854AFE115
7807FFF433B8BFD5B374605F21E057BBFF8A6DEDF63B5567E13DFD8177F0F0E6E2CE1D7C9FA61FCCD7D6AFEA67D417BC813F
7C98E3705CF400FC8027780A3FE913BCA65EF0FB4AF2BBB5FCF568D6859669F2131EC36BEAC4F8C033F866DEF01DFE733C8E
8FAFA31B8C83F35CEFF6EF577CBF67FE1C0F1D619E1C8FD7F37E7CC3FAAE38408F9F5A133A4470843FFFB877C2F3A21FF01F
3DDE503EB29CEC9B5BEE507FBE526F4CBBE5CBB71B6FD165F4A05ACED8D62EC72A6F390EEFC7A7C11578428FF0DB27FA7F16
2F9F63A3DFF0F2F54351F9DDF7AB283FE147A117A3A75FBD7AD5FC9C79311FF434E3ED112F483E601CF836F8C727E0017C44
2FE0E5CF99B1C9EE7BAFC923CA0F8E6FEB6FF56FEA417DD083D4C9A17FB9FBAEF179F573D3239D2FF3E7F5D97300BEFA54B9
8F1ACBE7A8E8313EFF4CD9F77BBBEBFB552F6F1EED59DAF70370039FE82FFDA6DE0D0E4796F37511FD79ADEA84E77D1D627E
FC9DD75D7F337488BB4F4DF5083E5F3C103D5D9E1F079EF003CB11AAA7E017DD7AA9686CB2CBBD8A3FF33B9D17F3345CA93F
812FCB7F19E312FC1CB8353C7A88DC87036ED069E6DD78D58817E43A7574013C8003F0BDFF78EC269793F16DC50BFAFA64E9
3E9FB8E727A86FD02F721A7D333FFE2AF42F5F878C0FDA17FA04BF9F1E3CAC93CC07DF339CA85F82971D436372393E920B54
0FD00FCB65DA3F78071FE01F75A35EE83EFDA15FF0E0E6A8CFE27D9FA07FE71B8F4B90FE5B2E50BF418F2B958EDDE4E711FA
F67CEF091FFBF907DFE0F89C0F5C8337F077E8686C98EF1FF48DFAD277707069FE5BA3DCF1D44FD07DEA491DD1A34DEF8DDA
2AC7A52EE8EBF5231F35F6759F719F8D9EF0B1F0DBFC56FD86F3E003578FBEDFDBF997FA4985B2B161BECF336EF0039EE06B
9D41D1437CBDE1756FDC0A59EFE3CD7C4AC7C9B8393FE3E175E01FDCC313FA53BC79CC11DF9F8F1E9994E5FB1C7D002FD5C3
6272B91CA8FA823F3C91F9E539B9DFCED6DFCA2F741C9EFD3662582777DF91EA0E7A475FE807BE857E511FEA856E807BF403
5F2F5B6E52969F234EA7C41C71EB14F282FA3AF32EF5D2A4C7FC5C00AEC139BC39F1F3A4C7DCFA825CA075A21ED48771A303
E8267A70E3789F4FE47B0ED07D780E0FE01D3C5C5F33FCA69C87BAE09FE84AFE41E1C5A40EF8E4EE77636648BFCDE755579E
6DF16E2D59871678B1C5C372DF153AC7F16D1DAF3E6AFB6D39FED34DAE5F46BFE129FC0467E8BDCD4BE799DE78582EB77FAE
C7C56F7FCC33365EF0687941F53E65C0E0FFFAFEBA366FC87AF99C07FDB31CAA78E2F8E6E3EAD3F8EDD3FDDAB495EB82D15D
8ECFEB380FFCABB762D4565FD7D1377086CEE11B4F7DDA3FDAEDBBA84FC27B74001C195FD4FFE00D7D4667D01DEA435DD077
FA4B5FD13BFA5BB0ED5BA3E4FE13FA0C7F9817F384EFF09CBEA20FCC139D30FF581933C3CF53F00CDE31EFACF0F06272DF8E
E51AF28CD6B5E1F6F09BEEFB79541F99373C85B7378B7E79CEF7277066F8541F05A7F0147E324FF8CE7CB2F39EDF9907F821
078023C6850E703EE6D12463582E5FEF2D37440EFEAFAFEBD48F7AA21FF01ABDC8CEEF02CB2216C8EF9687ABC69C093CB733
A0EFF0DEFC9DBA691DD1A7D7EE8C8D77CFC9507C987EE9DF79DD96D9639E137E5BBE50FFB4F5B6EA3B3C0067F4977E5F9FDE
3FDA3D075FF50E3EE3E3F8EAE5032D1E76EB37E537F8014FF8CC93D3BA55739FBFABCF5CBFFF9F6EBE7F337F8EC7F1F14BE6
C57CF021CB037A1C8E6B7EB037E68C9FA798277D0707975B452C58B1E227C33DFACEF9F1E9DF46B769EBFFE4EF8CE3A7E471
5F0BDFE1397AC9F9C01538030FE000BD409FF25D0CF9D67D3EA2BE844E5B5FB44F8536BFD142E655786AB32BD24F7405DC15
7E66D4DF82437C787F899883BE1F80C78D6D22C264DFDC7C5E7D8CFA514F7800CEC13D7807A7753F1AF39CFB1E3CC519BEF4
54B94E4BE5735DF495F3DA3A587D8E7A524774163F7CF6F5F7960BFEE1BDAD4B55AFC19DE9FBA998837E5E58576EF88F729F
18E7C3F72E5C1DF5B7DC576F3EAFBA0FDF9817F3C4AF9ECEDBB6B3D337D5397800AEC019E7B575B0FAE2DAFA8372BAE783A8
EFE1A769C5426FBBFD737D9DAD87F5EFE68FFA3EFA53A3C5B8AFE57948E805FD7F6E5654113F8FDDF8AC5B355FAFC1D9A1C2
31037DBF012F1C8FE353976B8542BE95EF51432FD005F006CEC0BDAD57558FC03B78A9BF3022CC7D6FACF2D67CEB42CC40DF
B7D14B1BBFCEE7DC275145E4392BF40FDF023FF82E38327CAA3F8153F0027EE0C1B69251A564FD062ED15D78071FE0DFAE01
E32E3BDD579EA313F4A158CBE81BBE9F1FCD376187EC07E18BF4013D78E6BB2193DDE723EA3FE8EB53937AF477FB01E83DF8
545F04A7E7733549973AA1EB173E6F76C5DD8FAFF56E7464F88F727F5D76FD7A63E3A09CF2F3AD3203F6BBE7812B1FE817FD
439F2C0F2A3ECCBFAF4DD8E1E79133CBA36F08FEE913FECDF99B9E09BDEDFB337C82A7F096FA95787EDC48A91339827ED127
7CE5C4C3E3770B0FC907D49B7A501FFCC972C9FFC6EF16DC329F5367C78D94FB211937C7BB7EA4D352DFE7C02BFC8257F004
9C800F74D07C538FC37119177E7565F37BCBE527FC802FE0171D0097E801FDA3EFE000BCBF1239EEB25CD780BE701ECE6BFE
78AB6D67DF1719B7F9A6FE9DD7A11FB54E4695F2F5F5EA8743262F5CB8D174043DDDFC6DE47AF93E54FA875FD17FF80EFF19
2FE347B79827B8CA8E33700EDE781FE364DCF8D853133FFDD37D7F95FA8FADABD587E025BC86CFE0005EA383F09BF9324FFC
129F2FD4E3CD76E243E8D9DEF7C7170A7C6F94E60BD51B7806EFC0C5F5AE3DFACF9AF5BDE93EBC4597C007FA84DE3C9B34B4
93FC1DFFA7DEF0173EC32F78C5FC37D40D6F247ECF7CF061FC33DFA5F6E3E4FA62F4189EC13BE601DF191FBC2717E49FF0CF
F1D237F40F1DB2DCAA38DBBE31EAA8D4CB7287EA3AF380D7F09CF1DBBA5C7D1DFFB17A6A7DD10DC6C73C381E7A41BDD0EF7C
3FF75EEAF627D427E93B798FFEC39BD73F8E5CEFEB71FABDA1F3DD7EBE8E8F1C917274F0353F3758FED0BF932FF83FC7B175
BCFE1F5DA9BA647C213F278243F087AE195FD45FE1CD950E433BB9FB7994E7E83B7DA11FE8F8D9E90D1BBA7D50F531F41A7E
5207EAB2AFEEB89BB27EB03CA2FC60BC8C1F9D00D705F78E78F8B7DF7E333F0757E406F0056EF157F07BBDDBA77FFABA0F3E
E0113C874FF4ADC196F046EE7B5A55A7D09F9AF5A38EFA3AC4BCF1818BE3DF749F07320F5ECFFBCD8F743C8CCFFCFED1A1F3
7D5D34DF2D3BF89AAF8B962F54D7C0297CB2FCAB7D456F2A6F1977D3CF53D75EEABDD47D9F9DEA0EFC641ECC0BBFF9A1C0A8
15D23F70800F5CAA3AE261773F89E2063F0007D493FA16F9B349BAAC778BE669DC5978049FD0EF7C4DBF7CD75DAFA1FA7DAD
60FB71BE9F514770059EC0FDC14AE306F9FA0EEE98077E98F5CD3FC73FE09B7A7CCE07EE8A5419D95570886FEF1C12DD48CE
03AFD02D3B5E363FA66EB69E565F41579E69D9CF5DBF86FFA05BD617ED137DA38FF086FA821FF0046EDE2A30B0B08F1FF805
EEE119BA5E71FFB8417E7EA15FCF7FDAB0A1DF37CEC779C0EBF95F467675F75D280EF0277CB260E43B0F0B7FE0F1A6D091DD
DDE70E5A177C893A535FF4085F7A3AEFC77BA50EE8D991A26393C517CCE7B5FFF4AD687264BC9FBBE071EDCBA35688AE8057
78039EC137E3B475B3FA10FE96BFE20771C27B74083CA3B3E01AFF29973936D9CF05E016BC826F748471326EFAFFEAF0E846
EE79F6AA3B992D22E3E57A2BEA8E8F32EF17DA8F1EEBAF7B8EE78A8E77399CBCA0F55B9F63F8936E5F5DE7898F711ECECBF8
D60E0EB9E3FBD9EACB030BCBF5325797F7DBEBAE97527D41A7AEA57DF9AEAFCBE01EDC8143F08E7F97FE333ADEE57FAD23E3
349F55BF643CF82AE3FA75CEE8B1EEBE62AD0B39007CE2EBE0943A952C11555BCECF38D007F8834E9C3C15555BFCD9F285D6
D5D6DBAAD7F03E6D51688CEFEB9603F4EFCC87FFE3BFF84FE69DC64E67F6BC17DD4674191D449FE00BBE086FF01FEA4F3FE0
41BD6123BBBBEF27577D810796F795E7E8D6B37906EE773EADBEFAE3E2D1AD653CF016DF60BC8C1FBF311EA9AFC327780D9F
E109FEF0749DAE33DC7E86FA217A697853FCC1337807BED007FA874E30EF86F7873FE9EBABF9F188903BBEBEC24FEA461DA9
03F3C70F7EBBF5F15EDFCFE807F5C0AF2E0D78E761BFBFF4157DB4DCF041688CDC6FCFF8E80F3A802E989FEA7919073C7B79
45741B3F97A10FE802FD605C8C13FF8397F0947AA33FD4171DBA7267E07E77BD94EA82E5D1831FC4F93ECAB8B3FB36AF2347
14B9DC749BF0171E336FF80EFF991FF345CFE02FBC65FEF01E1D607EE803F3E1FFE8E26BAD47B7F6F51B7C802B70F6CB4BD1
0545EFD17D78CC78D113C607FEDE1ED9BF968F435BC7AABEC0CB2D1722DE91E3823FFC9679334EC68DAF3EDDAF93FBBC1A1F
8377F819FC83F705FE3778A5F0887CF1DBD6AE33DCBEB9EA3BFC04AFCC03DC5217FE0FBE19BFADCBD58FE93BB984FE932F0A
9C7AF7A87B9EA0EA0A7A607853FCE1DBF9737FB841F0805E9EFEFBF5A3EE7A07D55DF46B4787D1DFBAEB8E9587E81CBA52E5
E7E8827E5E435FA9037541B79EE9F3CD49B70FA1FE001EC001FCDF58674479F77981CE1F7FBEFCF0E095EEFE16E50F7ABD26
C7C0477D3F5E15D2BF963CA793E3D87A5A7D94D7E39FBCCFF0A0FE012E0E3C1B35D2D765F00BDEEB168E78C7C73DB8285C26
EC53C109FE99DA68F028DFFF383F7DC48FE967C6A2B00BEE73041D373E49BDA93FFCA65FB6AE555DE6BCF82DE7E7F8E6BFFA
7F5B372AAFC03BB840BFCF1768BACDF91AB8523F035FF08DFA50AFDF46772AE9FB01B841475FBC1235D2CF47E4F2624BEBDF
94EBF0E001E7673CE8F58563619F2627EF347CE013F0C0F2AFE289FA554F1AFDAD9FFBCCBF749C8C9B3ED7DF3AA2BCAF73BC
8FE3307FF8FFE6FD818FFA7A76E8CF31B30597F80DF8BA1AF5CDC93973D24C47E0333C323D523E81AFE7428665FA390DDC82
57F8469DF099CB25DE75FC671E8D3F08BB20CF2F409FF8D96CDDE051BE9E65FDF1E106DFB7A8D7B6DE1107447FC11F7A4C1F
8B3FF4FA51B9FE9E7EA2BBCF94FD7460E073F780EE322EC6890FF2D3FC52C7815F55C83566B69F83769719ED9E4B82BEA067
E00BBCC1A773238765BAFB1C1437F82B38A6EFE0003DCB9FD53BAFF0C67C497500FEA007478F8EEEBD7973D077E93F7D7BFE
5E78253FA7A1A38607C50775B6F5B7FA07B9A048EE373E93CFEFD1157CBB6CD9D1BDFDFC020F2CFF2A7FCE3E1A5EC9DD5FAD
F3C7A71957F1B8885764BEE482138B23EFB9EF9953FF667EFCBF54EBC87B7EBEA08EE00DFC9DFE22E215F77D347A3E8E07FE
E01138C4C7F2576C9FD3E9B6EA377DA7DFE823BC836FF004FF2C10D2EAA4ECEBA3678C8FF3C3077C135D3CF3CFFA37DDFDAD
9175F3C8FD0A25B7D78E169C83F7CD53475C93F3D027FC0D9C825BF0082FF021F881AED48A8E3820EB0970097EAB1D1B1DEF
BE4F4C75CDF0A0B9C170A1BA012ED18FACFCBDF3BAEBF35477D04DC66FFB04EAAFE81DE3617C9C9FF1F0BAAB473F1D28FBD7
E83EF85ED5AD7F05F161FC191F078F96F7959FD49BFAE307B6FE547D8497F0CED65B8A83BD7F8C6EEAEB233AB1FEE0D02CF7
3982CE131F5E9B1E12E5FB37E385F705BEED57DCE543F54F7806EFE8B3E156FD12FC1A0E35078047FC87BA51477483FEA11F
8C9F71E3F78C9F5CC03CD26B0ED9E9E700FECFDF793FFFA70FE8037AF15391C8BCA247F00D1F40B7C0153843EFF317FE628F
7BFEA8FA1B7A62EB33ED1F787FFDF311D7E4FA2C7484FEE36FE7FF78C33DD710BD417FE8033A46BDD033F3E989FD2BC8CF16
4FF6FD4C7E323F74099D42DF783FE7E1FFB6DED2F9302EC6693E77B07D4EDF47E937F80497E81A39E8F9C98D5AC87D47F01B
FDA89A7B74533F27D2175BAFAB7F1C4F79ED27390F7A76B266ED68DF57190FE3BBFC71BFE2EE7A3ED505FC03FED317FA149C
5FC0B72F8F6CE5AE47B3FC547168965C778FCEF27FF37B1D87E58F262151BE0FF0BE26DB87ECF4757ADFFE51D37D7D87DFD4
83FAA04FF0B5D0D8814F5DBB76CD7C1B5EA323F09BFA9363E803F5A39EE814BE947FC2C76DA53FE83AF847DFE1017A59E37C
645EC9D7D497F1F2775E9775E18B3DF21C37F414BEE1FF853BBD5DC7E989EACA0F85C23BBBE7D12A0EF025FA475FE8137CB5
BE2B0EE8077CA12FE85CE54AA3A6FB390B5F2DD0F68347E47A58F4EF62AF814FB9EB99144FE83EE3B4F5B1FA10E3A3EFD971
C0B8B2E3C1D69DCA2B7494BA51C79DFB22AE495DE01B7A863E3DDBA26B31B70FA13E7370796482EF2BE00E1C15593468BFE0
CAD685DA17FC98FE8017F4ECEC578D5A3C801FF389008EC8172586D59D2A38048FD49BFA83775BD7298EC1FBA60AC36689CF
515F7C6FCD7303E2E53E43F365F56974B462CBC8043FBF305EC68FDE59DFD5CFE83F75A6BEE814E331FF553F8307B52F8677
F6F5E6FC0783F6BBEB751437F8D391DDD593E4FCE81CFC645EF821F3CBFAE6E3B6BEFE828F23E1116DC4FFF021FA4CFDC97B
F4017E8067F04DFF8BFE1EFA85D41DDF5EF741E83AF169E689DFA6CE1EB4CEF763C6CBEB6D1DAF3ECDEBED777D1F7A41DDA8
A3F98ECE97F98397572B475C13BD4177C0B9AD2F9557F419BDBF30F3ED3A6EBF54FDB0DCA088367E0EBA92DCB5983CAF0D9D
425FE0173A08CF321F09FDC25D7FA37DC58F1917E3649EC7DF09DF2075C077E91BF57EA1EBD0AD7E5E62BC8C1F1FBA3CFF83
477C1FA25EE01DFCC333F32F7D3DEF271F94FE3E7C83CC87BAA00FF0109D00D7F50E0D9BE5EB31BEF1ECEB1DDDF7A1E11BBF
4E1ABAD5DD3FA075217F30BF92398667FAEBB693F787670AFE791D75327F3D37205E9E4387AE932F8AAE697256F40ADDDA13
1F11EAEB297A892E1658F74D5B773D8FFAE78FCBC387CBEBE039FE53BAF96B3FC9BAA9CCC69A6D7D3EA1C7E004DCC0275B1F
28FFE1373A7A2AB4AEBBFF8CFA355A14BACEF70DCB0D1F0D5AE7FB067CB7F5A2E21EFED017FAC4F8CB55AB9EE4CF83F9DA7A
5DFD0CDF2B10F2E13B727D31BA0E3FE00BF8A2FEE40CFA809F17FAFBDDBA4E57555FE95F8B467D3FF4F32CFA03CED021EA44
DDD071FA6BEB69D56BEA6DEB00E527F8C49FC029FC78B97344A89FBF2E37FAA6ADBBAE43F98F6EC25F5BB7693FA83B7DA4AF
F48F79D3477806EFA82B3AC0FCD1832B9B3BAE96FD6B741FFCA12BAFB50C1FEEEB283847A7C03BF363BEE813E366BCE06FCB
B461EEFBB6C1073E4C5D6DFDA13C87BFB64ED27A8107F204B8E0FC8C9FF9C0CB425FF70A71D723692EF8E5EB91A9BEBEC36F
F4033DA1AEE80C75E4FF8CD7D6433AEE551FF74BCFC808E689E4757D3F9475A4F142F301FCA03EB60FA1F9C0F286FA4FE69B
4DCECA4F7C1F3D3856BF665BDFDFE82FFD663CE402C685DE183E15AFE4ACE23DEA5D123EC24BF4B8C08B9F1670FB82EA638C
8B71E2371BE2C23A88AE301FFC7FEDAD81EB7CBF675C3B3A84BF2B7D8087E828C7E1FDE4038E432EE07897477EF88EEF67D4
05BDA9123B32D5CF8FD48F7AA21717637B852C5DBAD5788B2EA7270F4E10BE320EF2087FCF9E57C09BAD775577CC97D5572F
3EF46E5DF949FDC95BF4C1FC9AF5AACE13BED6FD6CD87D5FBF395E76BF37FFFBB65FBAAFA3F0A348CBBE834437F0ED0375C2
E7F93E016FE897ADBFD50F6C3DA97A01FFA837F5473FE02DE3673EF870D142CD57C8F34EE03D7E52E054BBD972BD307A0C3E
D165704A3D1A7C11D6C1D757F3C5BC03D7C973FDA8CBE5039F1690FD68F4119EC00FC383E2833ABEF758EFC96E5F45EB891E
BFB8357C9E9FA7CE2FEF3B48BE77073CE137E6DF2D0447C1F151A7EA49E1EFFA79909FFC9D3AE293854EB6EAE47446F586DC
517A44ADB63EBEC03B7A737A7CBD4BF2F3E0EA6A937DBD8127E0A2689701ADAF5FBF6E7E7B38E78826827FFC0CBC581FD5D7
E9E7B632C3730B5EC00DFE00DE6C5DABFC41CFE80B7D42870AB6FDA29BFB7C4FFD0A3CD17770401DC831D483795307740D7C
D06770B2BBF0882DF279367C46B7F08DF27F8D68E2E71DC3ADEA3EF8B5F5AAF204FE8117F0036F32670C68EDAEE351DCE0A3
D4D3F615D4579247F5C9129FC28756DCED3D79EDDA5F8238511F052FF4B3D8FC903FFD1C753463D814F75C5DF563FA4C7DC1
6D761C533FF00C6F0CFFCA07C6895F315EFCB66CE36153A4AED4F74CDB903FDDF5313A7FFCF8728976B37DDD0737CCA344BE
21DD7C3F3FD1252CC9DDF7AAFECDFCE83F3A9D79B1F90AF949DF83EBFB800E928FCB9F7BB5BECC177C719E5233C292641ECC
67D31F6173C517E9133E79EADA906EEE7EE16C79087ED93A5579627EA7F365FEE841AD63C373FBFABAA642FFDDBEEF52EF4B
F3BFE8E63E17C9A677CC9B3AE02B8C9F71E3D3E804BA727C702DF73B3C4597E12BF9E285690D07CAF3F9D009F05EEDC2882D
7EDEB4F5A4F2041DA2FE15DFAA36D9EF03FC838FF0893A9023A8C7C592AD3AF97E497DCD4FF5EFBC0EFDA5CED49DF1327E5E
B7EEC2900AEE7E0EAD1739021F2B18F9F123F279383E663EDDBF4F96EF33A9B743A2FD9CC13CE093E983EA05E7B5DCA27985
E3D8EF7A3C74B4D0E6AF57B9FCAC7EBEF7E8F08BBE2EA37FE831F8014FE8067A018EE93BF986FE8317F0832EDBBA4EFB8D3E
A107E8037DB63C943B6CAEDCCF874E51C796977A9E953A524FF2C5736BDE5A26F735A28BDB370D4B95F9C173FC8ABAA23FE8
117E5BE8EFD6F3E53A62F4DE72C6A1FEBB7D3F417F6C5DA8F5C27FE1D3E1E75EAD2F3FD131E68F9EA12BE80CF5605CE809FA
6279A8F0900ABE7EA3571C97E3592E793C24DAD7758E67EB369D0F78C6E7C0F5C5D7BF5EB568D126D317FC00FED02FFA07DE
6CBDAEFA4B5F6C1F42FD0F7E572D3BFCA29F37E92FF9833E832FF4FED7CF1A0E949FE01F9D820797067CFC88ECA7E32BE09E
BE822FFA4BCE2A35A84EA28F77F854E48BEEE75DDE56FF67BC8C1FFD305CA94F83AF7D3386A5F97E002FC10BF881572BA6F6
AE22EB6DFC963CC0BCA903BAFFF3AD2AEDFD3CB8AFEF4B737C9FA10ED4059DAB596F58AAAF8BF415BC813F74C2F0A0F8A06E
E433EA473D383EE70397F4017C6EED1B3655EA013EF05BF4B5508F8EDDDD7E9EFA12E3E378E0FE7C5CF7F36EBF5BF18ACE72
7CF371F545EA876F9C7BF3AD65F27355E3BED77CDFA50F171F6A3DDFF71BEA8B9E55EE322CCDCF831C8FE3E347E65FE65B81
E3725E7C8FF39BAF7CDEBB8AAF8B3BBB0F2B2CF3816FE81F382DDAA4E73D97FFD587C1A1AD7795B7E005FC80CB83FF1AEA9E
CB831F805FFCBFD8FF9A3472DFF3A17CDD587B48079737B4BEF82739ABDCCA1AD5E5F9C2E0C5D69FCA5B7896526440943C57
D17C5FEBC1F1CD97D52FE10D3A7122BC8E7B8E357DA02FF095E3E3AB9C071ED4191336D5E703FE50B8D38739E5F352F432A3
F0A05CEEFE081D07FECCDFB3FB37784667C135F5B13CA475CA4CEF794F9E430D8EF133FA4EBFD1597869F8573EE03315170C
6DE8E7AC8BE33B76977D6DF4119EE3ABCFC5B53820F735A267F49D1C40FFC157B1F7FAFC2A78C3FFCD4733FA5EF375175DA0
CED4FD48AED041EE7E75F55770473FE92F3A483F5F8D1956D8CF71E0D6D6A98A7BF48FF9327FDEC77118DF99157D7E75D7E1
69DD2D3FA80E159EDA65AF7B2EBAFA0FB9B652F5971BC97CE927FE5CEECFD0417E0E025FF57F18D2C1D75DEA573CA2FF4C3F
47A177CDCF0F88127C505F5B972A7FE0FDB199837BB8FBF6D5EFA92BFD404FE80B38C7B7C03BF3AA3CE6A539F23A9B9FE208
5CC133E6D1F8C2A05CEEFB4874DC862BF575F0458E29D369700F3F1F9D1ED87FA65C8745DDC931CB627A86C9EFF8E4F2823D
CFCAE700D4A9646648193F779D2C1A52C6DD8FA4EFA71E3FE60AABE1E714FC07FD43578EBC5BA3BAFB3CA2E4D032BE6EA2BB
F0CC78AABC35BF569F39F3709346F213DCD37770F052DE2AED45D7ABDEA9E4BE87897A93474A2EAC57417807FF2ECCFC30A7
EC23E343E09BF3705E7C88F9324F7C11DED9FA55F903CEC00B78B3759DE20CDE9B3FEA78185F72E33E2564FFDDF28CF6EDC2
E75DF6BAE77CABEEA05FF003BE802FF0DAAA7BCF30F9F9CF5F7B34967AD11FFA855EC0537416BE5277F0095ECDAFD59FCF7D
D1E280DF77DB5750FDE5FFE6EB3A7F8E4BFDC03FB8303D517DA1BFE427FA4CDDC805D40FBF06A7FB77BCECFA0DFF5E3E39B4
8C9F4FF1AB22B9DB7515DFC17FD0AD227F7E51C6ED6FA90FA31BB62E547C6C9E33A48ECC1BDCE0F3E8E56B7F86D5F075DDD6
818A0FF40C7EC377FA8C1E8307F0C13C98173A4DFDD007F482FF331EC6C7F819377982FC522CB1795DA74BAA4F6B4AF59BBE
61C3B1601ED33EA027E808F5427F98273AF4F3FED053BE0FA11FB62ED4D7510FE661EB363D1FE321BF30AEDD535F8CF573FF
DE3C957E76FBC0EADBCFAD695554F41C5D5FFFEF4165DCF5A53A7FF28BF974469F12BE3F90D7CA26BFD65B700D4ED32A0DB8
EEE71BC6051EC127FCE4BC9C8F9CC471F89DE3C16374F4649B7A1502CF47F8A28CDBDF57BD407FC139BE05DE7F9A153A58EA
0D6FD17BEBBBE618FA0F7E5F6F3FA48EAFD3D4817C463DA817BA44DDD0AD972A859EF27328BC438FE01FF363BEE838BC2C1A
F7751797E7D5E7C195ADBF55A7CEFFD1AEABDCDF8E1FC02BCB0D27FA4DF775DD70A27E0C5EC011B8A2AFD4D3D6F7EA2BFB87
0E6921FDC33FC03B3CA49ED9FDAFE1BC41657C3DC637183FF3410FACEF8A03EA4A9DD1098ED774FF80EBBE1E93CB5ECCFB4A
15F99D7EA28B452E77382E7E816FEC28543ECBED2B2BCFA80BFDC137CEFCA7B95BE7657EF17517F7FD068A4B747F45C55E5D
84C7F81F39867E8023FAC2F1383E7E73EECD56457DBFA15ECC8B7CC6FCD0C31A9F840E96FBCCE80775E2EFBC0E3DAE1436A4
85F8273886C7F0F2688BD77AFBE3623C8C0F9C16FBA2FB48B72ED0BC407E29F96C832C798E26FCA33EF81E750237E004DCFC
506AF006C121F8C3970EBD30A887BBBF457D0B3CDA3A53F9069FC0077D0027E0191CC33772D94B8BAB74F2EBCDF8381EF83B
13D77DA43C6F1A1CE05B8C9F71E34FBB1A0FE9E1F319DDA2DED41FBDB07A6A6EA0AE8C933CC978F181A2855A87884EA297AB
F2F54997FB602D6F681FF0A50ABF0EEAF17FE5327002FFE0A5E15FF960FC525D8667E0A2F8173D87DDB871C3F2023E5C6CD8
5B9F385D527D3A5FA0C371771D96EA1ABA50EDF31763250FBE1A51B1983F6FF3B38382AFA04E1F5930F03D59F7E3DFE0853A
E0B7D483E357BF583ECB3F0FBA083F0EDC7AA58A5BAFC6F51CE6EE17D17EE3EB1B578664B9EB5DB5DFF83039A0DCBF06BEE7
AF87A847C9497D5AF9390A9C800F74F978D17E59EE7A17CD0DCC8FF3723EFC98E395CEEC97E5AF9FA8D7FB63BF39EBE7EE94
52FDAEFA7E4F9D4E76EDD3CA5DFFA7F3E4FCE84DED138337F83E80AED347FA8AEE16CDD3F14D773F2B3EC9FA56798BAEC077
780E6F1927F982F166B41BF0977B0E87CE9FDC6178D31C00EE787DF6BC028FF04BF8042FF18793571A38FF801FAF640CE9E1
E759F497FA502F7256994BB56FFBFC867FF011FEC133D319D51DFC0A9CFEDCBA8A5BAF580EB8D627DDF71FC6C538191FB806
E7E811F363BEF8D6D2DF7BA4C9754CF8DC925EDF9C95CF619817F989F9655E6C1D22DF5F829F810FFC9779EC1C58D17D8FB7
AD3F1507E887E5B67743B27CDFD8933938BFE01E1D475FD101F80F9EF031C6C3F82C4F9CE877D5F70D5BCF287FE90F7D301D
505DA0EEF08CFA67DEE9F8A6ECB3A36BE8E98F8B075DF7798EEF594E501F3E13FAD627EEBA5372D5DC017FF9BEC4FFB3E709
FA850FD337F8687AADBCB4759DD61D7D427F4C57B5BEE00E1C8227F406FDA17ED4C9F444F5059ED05FFA8D8F155DF3AF4BE2
4BF80CFAF1FCA75D16B8CFA5D5FF6D9DA6E3E6FCD499FAE2D39C1F5D639EE49712298DCFCAF360D00F7265C542D53B4BDDC1
3B38047FE8383C7EB9E8E0FC7E9E4517C103F8007FE0081C1A2F55878E15AC7D5B7EAECCEC75DDF775FA0C9FE0177D0297FF
7CA4479AEF03CC8B7CC9FCF0E762892DBE93FB47D145FA040EC1E5B65E65BBB8E73D28BF1987E578AD033AC4FB993778FE65
68C897A2CBE83EFC66BEF839F3A65FE424FA66F9B8F5A0EBBE1E33FF0F36756F2FF3A70E5B5E08F9C85D37ABF8C0CFCF4EEF
B2C0EDBF2B6FD165E6C1F1391FFDA15FE8A9E1597D055C9337ABC4546DEAF781F158AE50BF07EFC5EE7D7ED7E570F54F7801
9FE1C7EAFB7DEE8B0E5AFED1FA715E7224E7CF7CF35F97643F1DBD06AF8607F55170812E56090BF952F00EEE392FBECDF9D7
7F3B6099CB2D3A1F7C7BF1F26FCACAF547F8CAA27ADDDB8B5EECA817B2CBE7157A094FDE2DDAEBBAAFBBE810BC3978B1BABB
AF175C995F286F198FE507CD07E4CD575ABC38DBAF535AFD7EE7FDDCC0FCF013EA461D793DC7E57DE6B3AAF7A79A370E7C0E
AC7A00DEC0DF81A706F473F70BA99FC10FC6497E65BCE4DA5A63CB76919FAFBD5E6E94FC9FBFF33BEF3BF3E8E777E5396EE0
1B9F23A79469FBFA4AE1173C438F9F9FDCEE13F7BD16EA27F0ACEEAF211FF9BACEFC982FBE45DFF12BFA0F8EC015FA41BFE1
317D3FF39F16DFF9FE473FE047F18CAFCBBBF588E60B7065FEAB3E617EAAC7E1B8F8D88BD707F4139E817BF005DEE03FBAF6
568E3EF77D5D472FA80FF542A7C1EF2FDDAB36757E4AFD3557D10774009EA30787C6F4AB2F3CC3A7C18BE5A78F072C93F7A3
BFA71B7F5DDEED5F6BBFF169CB191BFB9DF7759A3E54DF14B2CBCF95F4A375CB6FCAFAFA4ABD4BAEE85E5DEA4FBEE07D1C87
BA8073F00D6FF1F30A7DFBD59779D30FFCB9C48CE6B9E5FDB6BE54BF2AF6BF567964BCE8EEF2593D6BAD5B17F44B7007CF0C
FFCA077C8CFEEC4A7ED17D1EC23C3FCCDBF582CC93F91EB9DA7B9BBBAE51F3047D802FE83EBC011FE01A9C9C7CAF7B75E101
FD20AFE0B38CEBC7CDE54605AEA36DF789ECB3A3D7E808F52EBDA767213FD7F1B3DC33BDB7F979095C8373F8F7C3E481FD7D
9EE093C75FEE59C85D37A9E363DED419DFA6DEE01E1C82CB5DBB07B673394A7519BD8707F002BE91BF2AB4A97942F0093FCC
7F555F8FCD7F7DA55B5FB0CE549EA007E09B3E8073E6CB3CF163C39BE60C7087FE155BDAE1A07B7E9BFAB0F9A9F6853E317F
CE0FDF16B6EE3647723AFEF7DD2DC1D76EE3357C86278673F563F04E7DC873D4095C1BDF95FFDF67F42E299F0F910FE81BF5
A39EE8B8AD8B948FE0057EC017FA477FE92BFA0BEE5B7ED2B396AFEB8C871CC1B836FE31E071771DAFF6831C61EB46ED2B7A
B37961A9FBEE7942CA1BF8429DF11FEA0DDE5FA936B09D5CDF81DEF13E78C73C2D5F7E35B0BFAFEB671E6E9547AE9BC787C0
01FA60FAAC7DA50FC677E53FF3B5BCA479857AA143D4CDD65DCA6BEB8FE6CACA73ABF5F47161BEACBE77AA4BF3DCF2332557
9F4C3FDFD007FC9079314FCBF75A37749D9FD4877AD93A53FB843E9DF9678783F25C00F40B5D665C8C9371304E7216E34597
4C87B5EE193DFA45B8FB9BB48EE432EA6D3AA6EFB7F5A48E8FBA721CCB757A7E74007D470FA83F399B3EEC6DDB7F9ABC1F9F
43FFD039F4ECD0C29A275C1ED0FC55E6E106F39CFEA80E8183368BBBCDF17D095D7BA771EF92BECFA04FC51FEAD844F0428E
C0275F98D63ACC3DDF44FDAADEAA520DDC7307EB9758E3BEB7B24DA9FBF2935C5CAD5EE541FEFC9636F92656AE77C35FE199
E5C2DC031EF77D607BFE01F57C5DC067D04BF0003E382F7997F39377195FDD42650ABA1CA2B9A044CA5B2D9D4EAA5E5A3D35
575157EA4DFDD1418ECF71B3E7E89ACF54C8E58F8BBFF33AFE8F5E569DDF7F9A9F732D6FFDD927D3F725FA07EFE923F831DF
549DB59C37BE5F84EF33F4057FA63FD4DFFC5DFDEFF4DF1D9BC8BE3F3CC43FE0337DA4AFE83DF8DDD7AE5A4FE70307FB8E77
D71DAA6F8177F0866F813BEA417DD015705A22B24B82CBF3EAE7EBFA96F887FB9E6AD5FBF51B4BAC71F7F7D177F569FA6FFB
5ECA4F78F0EB67ADC3E4FB60D06BF037FF74D7DAA20FF806BE4C9FE80FFD022FE8EBB1FF3598273FCDE7D2BF89F5F57B65E7
9E6DC50FF05DFA837ED6CC1A50CFD751FC877AEFDE547990FBFE50CD41E5B7D79EE2E39A3E8107F0415D393EE7C337AA54EC
3B5E700A5EA92F7EB8E5629982EEE747FDF3F9B8C19FE93F788057E68F5A67EA0EBFC11BF863DCD41B3C5347748F7A526FEA
8F7F9C1AD025C1EDA72B8EF15DE6435DB75FAD90CBEDB3695FC823F4E740A3DE8DE477FC1B5CA397C57BFCEBBEBB9F597D00
FC837B784BBFC867F40D5C976AF4F9C99B376F5AAEA1CE967F341F30FFB6C5BBD6969FFFFEEAEB48D19B53CDDF6A9990F0BD
F92575373FD5BFF3BA1D31FD927D1D418FF1BD17D7F56EE4E7BED5CB7ADD9579930FA893AD1F143FF0F8C4BACF4FBA7D73AD
3BBE0BAFE121FC3E34A9C76AB7EFACFE4DBDD727F7DDE57291CE1F3F673CE402C6B5686FF77FC8F537F8103C8637F008BC52
EF3215BEBA2FF527A7988FC6F76CEBEB31E3B15CA1F9007F2E7DB6F176A73FAA43E4890A5D7BACF6D761E8163C835F6925FB
F494E359CE21F7E8DF391FFF3F76E8ABFB323FEA6B7946E751AE56B79B7EEE3AB2ADDB4DB70FAEAFA3DEE4C44A575E3D27F8
04A7967767F7CFE7FB007C855FF02DE8B72D5689CEA3F7E00D9CE143F00F3ED227F8854FC033EA048F793F7DFF7F6DDD77BC
55D5B52FF0E4AA79263189051105B1A1882828204540AA3401E9E550CFE1D08B806241BD18150BA244C5822D46B0A0888A05
3B6010258A5D63A2295E4D62347AF5AA31565EC6DAF33B5997F7FEDA9F5DD65C738EF16B73EDC6F7F43FFB7BD2E9375AB7BB
AAD867265E679F4ABCFDF39261DFC6FFF1F0033A050774081ED44DCE533F3AAB6EEA881FAD669D725FE8315DCEF93DF10E9F
E04A5E81AF2E8BF6DFA9F89F899E0DAF2CE7357A825774E59651D38ACFDBF2B75BAF9C7E417C8F2DE78F0173FF59F6B79CEF
D33CAD27E7B9BEF3B6C47CE8B5BAAA335DA68BFBFFBCEAACE23A4EF2613EA0CEEA0E4774838EC89B478E6E7653B9DEF0C6EF
E12EE7A13F9E3CA7EC27FA37A4F9AC9DCA7E427FB2CE277C3CFFC393DF2A3EAF95749C5ECA9DEADE69F70336C5E7ADCD534E
355FE7A3DBCE7BF7C727DE12E7CD392DE137EFEB926EC0C5A64DF33694F9C387E9987ED333F3343F79D63CAA564C7B276E47
BD39AD43B18F4838807778D0073A4377E4DFA3776DF24A79DD7426EB79AAAF7C6C1EEDBE39F882381E4FE00ACEE8341EBF5A
A765F17F41CE2B673B7FCE03C9F7DE6AD07573F17BBFF675490FD5D57CCCC378F4531D1DF7C70FFA3E10D7FDF92B1CEBA31C
A19FEA94F536D5E59173F6FDA89C5F1F7DA861F13D1EFDCDF92AE500EBC8F926F978F6FD341FF393071BCFEB7061593FFEF4
9F5567DD75D7C6AC6B741CCE731E487A8D3F7081472B7BCEB824F613FC95AED0ADA3FE79F25BE55CAFAFC6C34FFA40DFE984
3EEB3B7DA11F076C1EF379711D28E5083E73C01D039A15FF539F7CE3FE4E73F68DEB0439D7A5FEF019F87A6179B3E2F359B9
6E299FA99F7558175DD43FBEB4FEC3033615DFE3FE8F931F2FF74F7EC9F9BEDDBC0D657FC0B7013F3BF196B2CFBC78F2DCC7
8B3C99FC890EA93FDFD3073AAABFFAAD0E72937ADCFCEBA92B82D7FC6145A3691DE2F7B1E49746DFEFB147595FADDBFCADC7
7A734E4B39893FABF76FFEBBC92BE57A797ED30E075F10B7D64F2FD5E1A107E7CE2EE721EB78FC84935694CF271F79BDDCE4
B83FB61EF3795CF7877B7E639C9CAFD27174B7F9C5731F2FE7DC8CB7E4377067DEC31E9A714959EFE10BDEE07BF391271F50
E61B5DA707F00557C6A5DBC697D79A4E6A7D6D19BF3987A6DCA09FC79DB5EF47C5FF645C52BFF85FC95E331A5C1DB7F87450
8FEA31A1777244F6F7A4AF6F2CEA50E42FF8844B3A92F3C4FA39FB96FDE4E52DB3AB22A7F357B8C66BFCC76FFCE8F8EDC98F
97FDE18F4307348BFFB3A1BF707FD7E4593B87DFF051FCB05EEB949B72DD526E533FEBB35EFA9BF739090FF400AEE15CBFF5
0DDEE15CDFB38F27DD7FEBE1EA3171BD1E2EF92E1EC32D1CD37B755017FE76448BD955E55C79EB9DD39B462EA6D7784F977B
F69E3BBBACF7FA8797FA98F3DF3D27AD28EB3A5DD02FFDCBEB4E7C572FFDE29FFA06778D3EADBD2BFC57AEA1EB079E31F8FE
E27794936F647F4C7EF8E6773DF688CFA1BFF6B3998F14D7DF53FE803BCF7BBD3ACA73F0D1ED88FD9615DFA34B79B3F99D47
FDB48C17FD6FFDC2C90794F3F29BBBD4DE1575D53F39408E38ECE3998F94F35DCED9691C75B4EE439E9BFC5639D7E151F685
C443F3376FB9EFB72DA6BD535C1F4FF9433DE0819FC0051DC0437AA00E72AD7A3C38B7FE57E51CB3F68A06C5EF31CB1D8D47
765AF9BF7428ADE3D02DD3DE29E7ABDFB59CFC565C97572FF3A59774E8D565ADAF2D7E4728E1154EE1165E075D336BE7B29F
6C1C76D2C832DFF833BCC2337CE7FD58E211DEDFFBC08955F1B95479027EF1127FF1997F361AD1FB8CB26FC8C1FA74ECCE07
350B3D7F76E9DCA34347E9297F90AF5B1DD3F493721FCC5B3DF0D1FAF9863A38AFF3C9BFC697A79D674DAF7D9A06AEA2CF71
1BFD8D5B7CE47F78A93E72873AA9BFBACB0D23874C6F1AE71BD361EA8AF2BEC4BC3DEFF18CAB94CFE04B7EB79E763B37EE17
C7FDE1BCC1F7C7FB0A799F99F434BF3EBDCEBAE10BDEF8013DC9FE957848AFF0846EE5FD5BE23B1CE1097F78ECE5FD96557E
77B88297AC1B49471E6C3EFBFCF89C74CE6FA9AEEB2E987B76F1FDFD544739087E72AE487EC397F0E9C52147FDB4B80E9674
F4A07387772C7EA724F970C66DCA0DF0CB67D4537DEF787AC639A1E77C97AE2CFF66EADBC1677E23775A97BC667DD695735D
CA673937AF3C6964D967E85BD6F9D4A7BC1F4EFDA0A378DC6AFADCA3CBB99E8EF4EF7362555947F2BE31E9957ECA4FFDD656
F8D177D53E4F14BF23B465CE65E51C697D1EB71ECFEB73D6BDA48379FF99E6EDFCF0C8F7DEB8B5D3CAE277A7D37159E7535D
F88ABAA923FCF1E70D5F1CD42C6EF33E27D5C7F17C1B8E9E7DAAE927C5FBECA9EFFC5CFFE5D0A67BB6BBB0CC57FE8007F0FF
E66DBDCF58B1E291ECC770667DD6CBBFF589CEE9D7960F667FAFB8FE9E7C918E5AA7756CFAA271F1FF18FAA87F72227EE00B
5D707EF3346F38E7ABF0FED699C33BC6FF48D11D3E91F3D68BB3CF2FFB58CEBBA7CD3DBBEC4BB99E29A7AA2BFED25F3CDE74
E79C5BCAFCE77BA37798FA76E8DCB8BF4C39337E3F8EEE3ACF712DE65C56F613F51DDA76C639653FA1138D3654AD2F3E6F94
F2027F6BF4F7BE5F057EF8031D777EE7D57F788443E7739CF3D2D11675667FAF9CC7E58BC6F3BA7E16F8A65B77D799B93AF4
9A4FC329FCC2B37EC34FF6DFA4A772217CF49CB4EFE5E1DBF44F7FF51BDFF4515FF9291FBD6F700577FA95F349F24BF88077
38F955FB2977C7EF2F859EC6EF3ADDB4CF9433E37D26F9B479D3964F94F1A2AE7800FF2FB6987951F8093FC31BCFE38B7AA9
279F5757F551173974F54EF57E1E8F5BB7F5EA9B3ED29B378FAD5A1FBFCB80277C31FB72D2BF57DF6F57DCDEBEFBF43AA173
FC06FF721D524E558F8C93941BE0459FE01F1F329E93FFC1B5D7E90FBEE0C121478D3D3AF45A2EC8F97EC89C5BCA7A6CDD1E
F73A3ED37CCBCC8B228FC3FD86A7E6EC5CC607FFC747FCC74B7552377A2A77EA53D756FBFF2B7ED7E6BE4B66AD89FECB19FA
FC66DDBE5FC5F719E832DCC361E85F7C4E1F1EE5D39697351B59EE43F6D70F66AE2EEBF7CBD5D35EABBCCF5AF16F78344FF3
935BF59DFFE8BFFAE45C94F2C4EF9E1F7B74BCAF004F7C977F0FFCAA82D3813FD8FBC7054E531DE414F5C0173CA117D62B37
5B77F6B5543775E483F4EFA165FB5E5EF95D8C8A7EC03F3E8CDD38E5EEF85D8A9A1B27CF08DDA6DFF0D5F8A5EAA59F7EFA69
CE3B72B07A753CB0D13F423FF59F8FBDB1A8EB67BFFAD543B91F74575F9EF9EAC4E59B366DD3417E228F1C71E3B4D7CAFB5A
E319DFFAE5EBB61D9BD496EB649EE6B77D1E77DF718F7D34BBB84EADAFF2149F39F8BE011714FFAF94FCF0A166272E2DE730
FDE443EAAA9E79BF9778479F721E483AF4E2AB2D9F28EE9F34F983E27D8E94B7E0D73CCD4FFE339F7C3FCDEB8D66D54B8BF7
57122E730E4CFD6C725DED2BD15F79D06DD3C5933F28E7C38CE78463B9465F877F38BD4ED957E5B2C3BEEAB0B2AC8BAFD7D6
BE52BC6F91E6617D8E573FE3D075BE16788EDB1BA74DBE32FEA78C3FDD503D7946FC2E159DA453F4922E1E7BCC9C9DCBBA0C
8739CF24BFA4BF7DE7CE5A53F649BCE27B8F3FBBFFBF8AEF238C9AB177F4470EA01FF427FB60E2B91CD178738F71A1CB7C00
2FE9097DC9B84838A1BFEA2CB7A937DDB75EEBE7B770F7DCAC6623E3F6963E539F887CCFE7E4CEBCCF4C7CDC1EBF6D763A71
7939D7E30B3FC31BF5EEB6DBECAAB21FF2717D5EFD65055F392FBE74E2D2B21FFEBEEF800B56ADFA75F609BAA64FF446BF72
3E49B9E3C93F342A3E8F6CBEE6AF5EFA288FE827BCF125B8A39370462F9D4FBD9F7EB249F13962390F7F8E7FB4FE3DE1432B
57D52DDE57550F75A8BE72F295A123136B279D1C3A42A7E943F6B5843378C177FCE703F00BB7198729AFC0E3B31D6756857E
F0213A4D0FE8035CD0C543D60CD9BBB8DE9FF2025FD22FFDDBB8C3899DCA7C933BE0900FC323BD525FBAA59EEA2827D25FF5
5337FA4377D4519DE892D7D143E7B3EE9C7F92EFBDB6538795F1B9667CC22FBC95979BFDB1F51DF13E34BEE2275EE0C9758B
265D197987AE5D7BDDA493E3FFBAF001DEF102EE07AF98B177D96FEF79654671FD23E7EA846BB964D8E00AEE869EB2D782B8
D517FD907BACCF7AF9B1799B17DEC217FEE19DF9D17BF3AC7A60EA1371AB6F7C40DDE46CF57BA3758F71F1BD083900CE8C6B
3CE3E4FC905EEF7838876F7E80AFF8ABAFEA29C7A92B5D6CF5E4CCAAF2FE430E86CB1EFFD9F0BDC89BEA41279CEF77FD86EC
1DEFC7D0357E401FF9DFFDDDEBDF53CE97FAC79FF191FEE2E5BA3AB39E2FF755DEA99D37A9E0C994B613AF8AFFC1A5377048
FFE0D1BAAC478ECEFB906F4EEC54F69FBC1F4BF5A1BB37CF9E7276FC2F9B9C4F5FB79C38ED85E27D88E49774E8819366DE1D
7D92B7F4C13EA065BDA3762CE385AF1E3AB4F73B655FA24F8DFF32FCB6E2FA6ACA0B2B2F997643E829FFA643EAA65E72985C
D374CF4E9F17FA9974D4FAD4415DE8BA7EEAA375C971D697D79572BFF53DD266D6F1E579C87D1EDF3E17D28B130E9FF1CFB2
5ED0637883BF6B6E98B83AFA4E77AE7E7AE255F17F9A398724FD7BE980D677C42D9D8523F8B14F808F2E630EBC3BF898FD2E
D5E38E8B2AF8F57AAF93DF6F3976CF9BE279C779BD7EF16D7DE34B2D7E31ED85C07DE667DA1FB4B9A1E9A6725DF5517ED2CF
37F6197E5BBCEF8287FC156EC3AFE2F72CE1D7F8F615CEA3EEF44FFDED43ACB3C38443BE29BE47D661E661D1477AC427BDDE
EB9C07CEF937BCEB8F7ED171FCCAFE9E7430E7E30F663D5FF625F981AE3FFCF386EF15D7B92E9EF2C3E27D8B943F72BE4D7C
6AD27B5493627F91F291DC43A70634D9FBF3FBEFBF3FFB479FC533EF2EFBDB982553CE8EDF65842FFE30B966E2EAF85FD0E9
7BD5BE187A45B7F09ADEE177CE9DCFCC3ABEEC5739CF269CD027BE58B5A182BB915FD7AD1FB7BFBDA3F73B2B576EC87E833F
FCC67CCD53BF87CD9D7643D9F7F08BBEE259CE0F8967CFFDEDA8CAE759121E725E4A3E2D8F1C79F2941FC63AE1DEF9CCCB79
5F7F705493789F048EE51E3AA5EF70409FF84EF4BD9C5FF54FFEBA7245EDDB715E7AB1F4BDDA17A31E393F241F7EF5FD4E9F
C7E7F7B38FA77AAAEF4B174FFC495CFF93CFE06BF5F7A60F0A7CC931780F5F87FD7EEC959F7DF659CE89F40C7FF1F9970B27
8F8EEF8FF21BFE06DFFCF9899B0FBCBBC29F34DF347FF38657FC865B3A08A7F02947363B79E24FCA79F2B583C75E1939573F
E44F39B1D9C6769BCBFA647D872FAD7E37D66B5CFED664AF7EA796FDFF95E9D5EF46AE36AEBACA1D70F74C4DD3E27A391EB4
FBF5CCC3CAFB24F59783F421EFA3D2EBE1893E641F4FBAA3BECEFFEBEB0FF9A6F83C52E21DBEC935B7DE39F598C0033FB64F
C00F7CA177708CEFF84F97E8027DF23AFDA40F1BDACFF89FE2BA68E29BDC346D54EDDBD1D759DFAF6DB46EDDBAAC43F8878F
7425F731E557FDE47F7876F7EB7B7F1EB77244D349DD26441DE93ADCF21FF8952FF1B0CFBB0DB686DFC23BBEE3FF33D3A66D
29DE8748BACF67E83A1EC0FF4D4FD49916F3E2C3B7EE58E98FFAE43C9972121DA073F4E0F29B6A9BC7BE992E5CF65D6DA3D0
257EAE1EEAB0E697D3DF89F9C97970377EFEE4D1713C9EF107F3339E79C28FBC0047783B70EBF441657FA63F74078EF18BBF
E219DEE0115CEA8B7EC8CBD625F7591F1D834B7AF6FA7BFD4EBDEBAEA7B2DFD035FA414FE08EEEC0297CF227F553B7C78E9B
B1AADC3FB952FFF1154FD58BCFA81B7CD131BA268FB778BEE58705FE130F32CE53DE8277BE3BAE4B057763AFA8F3A7B835CF
9C7B53DE5D3B61C690B89577D5D57CE98C3E5DBF68D2E0980FDDE7C378917357F2DD99E36A9B471D677F3AE1F2E03BDED355
7A0A377488EFBDD4BEDDE6F83CB5799AB7F9D273FDDE5ED7DDD76775A5E7E691F75D1B67FC4FD9CF470E997A4CF447FFE502
7E0B5770C6370EDB65E06EC5F7C2539ED26FE319479EA3D74356D7AB895CA0EFF62DFA4F17DB5C396D4B795FB8E48A09D7C4
FE922E2CD965C2E5F13B71F60FF873DC83FB17FF3BA2AEEACCA7F9903AA80BFDA14F7428EF4BD37C1CCF0FFA8D9FFE4ED9E7
E9B53CF440FD065BE3D63CCDCF3EE6CEBDA63E13D77FE5153AF4EAB26E13E2FB2AF2847EFCE68229CD63FFC72FE9F4F6B9C9
F1F6316D0E6B36AB5C6F3CC21FBEBAFC81C99FC47AEC5FF8893EEA9F7C4747F000FEF33EE19119ABCA7E8B5FFC0ACFA67CBE
FB9C187FE65EBB3D16B75366ED7E5ADC4E98376970F4816ED1F5BCAEB48FB23EFB25F5EDBCB2D143C5FF0CA4F3F7BA7EC690
B29FBFF6E9C0DDEEBD7773D64DBEC867F23E3CF5F1FAE67B8C895BBEF7AB1995FAC101FD8307F3310FFB337D979BF4DFBEAD
FD734D4694D7F7EBD1D33617393FE98EDCE1F5C6751C9EF2217C95139BFDB163DDB2BE5E7EE66E8FC43AAEBEAC52F7CBDFAB
F4815FD1AD2D47B5FC303EC72A7734FDA0C73FCABE470F9AB61A7247F13E61CA537C453FF5511DF9B17A2E3B67626DD48FFE
CA0B7CECE8D3A6342FEFCFE8031ED309780E9F0F7D816BFC874F78C5AF9C7F92EECB0D7071E7C07A35651F9AF062051713DA
EFF11CFE9773D223BDF72FFE0FCD3A8D671CFDCD7932E527BA0B9FF09A7D67C6846B9E7AEAA9EC3FCF5D3A7158E04A2EA113
B9EF295FEA3F3EC02F5EC0F190F7A63E53F673FD87073A3BF3DC0A7EE66CDEAD4999C77227BD38A1B67EC7471E79E47B77AF
997664D4576EC303FDD45FBA3BBACFE44FE2FC70C4AFF4CF3AF0FCD56787DC11DF3BC4E79C6FE5B184EB675E6B362B6EAD93
7F6DEF67CEEFBC19FF2997E0C125EFEE7658993F97B6AED4030EE48E75C31A3D14B793CE9A581BE3D1513A8737872F1AD1B6
785F23E53E759313D56FF182090DE23A011EF047B9A0E59C89C3CAFB0F7AC6A7029F71CBDF3A2C9FB6B9EC93FC9B4F8CEC55
F78BFBEEBB2FAF475D37B66C52FCCF06FED27B3C7E62F0B4CAF77712BEE54A7D765EF3801FB89197E90C5DA337F870D2D913
1AC4FFD29FF4E1847DD6AF5FFFBDE7CFA8E95DFCCE68CA61F891F99B78CB0F6FBA7BD293C103FEC7AFE8231FB87E63E5B857
EBF4F847FCEE3C9FC65BEBCDF939E5D4AB66D62E89EB82744DEE985BBF829F53CED9757D994FAFCC1BD1B6D8FF271CCB79F8
C7FF5F3AA063DDF87EC6FD274F3B357444DE8597DBF7993222DE979033EC7FE806BDA08370D76C97516B3EFFFCF39C43E5F4
160FB4D9A5AC2BF41DAEE1594E3DEABC9ADEE5EB1C728ABAA9E3453D77DD58E6CF456755EA617DD6EB783CC66B3A452FE900
5D501FB95B9D5EEC3D6ED7D8FFC9CDF0F2F0F0692BCAFD93D32FAAAD7934F6AFF8B778F709FB04EE5EFA74D49AF87C903E19
4FFD9A3F386ED7A8A7BA187FFB7DCB53274D1959EC3B920FF17F7E23274DFE64F7FF2AE3928FC267E66FE2AD1C47F707F49B
766459F7E904DDA0C77C8C9ED1317C3BE5BA9A478BDB7FD43C1639572E38FC963EE797F3C8290F55FA7B46E75D7F5EC67BC6
95FD40C2977C8FFF7D9A36BC3DC6E7C7F4EC9EEBEA778CDBF356EC7A7E1937E7ADAB9C07FEE52F3CC0777901EFA75E5EBB24
FACB6FE9349DA37BF0472FBBAC9AF678395F9C7F5CCDBCE003BC5CB447CD63A153678CAACC73C1D65D6794EB41FFE91A1FF8
F7EBE696D7E3387DA0C7E6A7AFFC4D7FE50EF55337FD3BFD919A79E11FA77F5C737AD4594E6CDEB4F30565DF834FB8946372
BF529FEC7BC60D98F464E81C9DE0E7F09AAFCBACADFB45DCD26D7A4DA7E8279FA1A3F84D27E9A671ADDBF857F49BF0DBB84E
467FE54779E2F88BA79D5ACE47E71E54D337E6AB9FE7FFACE6F4C003BFE70FD7FCB4C2473A47F7F0497DCF7CABA6F81FDD33
BFA81958BC7F93FC0F7EE049FE503F75B39F6BFDEE5143627E7803773D6E9FB6A29C8F166CAD393672B1759CBB73CDC098BF
7978FEDFB7DD229FC3CDF0BF4C1951CE47F4E098C5534696F7C774891EF1553A2027D1033E8167F8E57CD6EDBCFC081FF852
CEA5C90FB7F469B34B7C5E5E7DEC0FD5E9866B2616BFE7CFD7E514FB46783DEE67075E177579A5AACFF991BFF93A5FB16E75
D03FBA2C07CF58B3FB81C5E7FC93BFF015F886D3ECEBC93FE42DB964CC3FF77C20F41B1FE578BCE427FC05EE56AD9E7275F8
BF9C4DFFE4727E30F89BBD9747CECA3E9FEA6BBC5C9F5417FBDA196B26FC36FA277FD14BEBE07FD68F57F409BF72BE4ABE85
474FCF9FB832E62F97F043FBE3F6E39A7E53EE333F906B1F78B5E1ED710B177092F344CA1DFAABDFEB474F39A7F8DDECA46B
722B1F3BE2B2FE4DCA3934CF27EDC3CDCBFE5CBD3A3768BC6BE8CC2F1A4D3826AE7FD01DF93DBF3EBDCE787401AFF15C7FED
0FF4F9969D27EF10FFDF64FFCBDFE81E1DA57FF2163DBBA25F05B7CEA76FF4047EE0464E5737F59283F18EDFE25FEE43CA73
2FBEDAF982F8BE975B8FEB135EE0835CC3C7F45FDFF1509E9E73D66E8BE3B666F2C4BBC337E8105FE2E7746573FDA386C4E7
8BE914BCC3F9896F4E3826F29F1C8407700BAFFC923ED0057AC91FDA2E9CB8B2BC0FBE77C694FBE3BAB31CAFCFFA1BB9267E
974E9FF187AFDEFCC33D1FC087726E7BE4E303AF2B7E573FE5F716CF1F5355D6959767F5FF77BD9ECD7AC75FF176F0C02957
C7F9F8C4E24F6B1E0E3CC3B15CCCB7D45B9D1FED3B6579191FF2B5F5CAE9D68D77FC03FFE43CF8BDF49C4A7FAF3DA9766D5C
EFA1EFF2027CD31738CF3924E5037593779ADFD9BD4ED9C736CFA87922F6DF72097DE2039D964F39A7ECE7F2079CAEDA61EF
E571AB0ED66F3F907366C2E3C69B9A7E139FB3830B38E13FFAC13FF5853E355B35F0D9E2739A2957C281F3AC7BA7F1AE718B
2772F92927EC7640F1BDDCFE531E8C79D99FE88F7A599F75D1377A43E7F009BFE8B2F1ADDF794EDAA5E6E1C0917C8A67D6AD
CEEACE2FE58F9AFBEADC50E431384FFB0A78B7AFA05FFD1F6DF0E7B86E50F5C5E41D429FE0571E515FCF7B9C6FE00BFE5CF4
58CDD0C03D3ED807D011FA4137E47BB8BEE89E4AFDE9195F0BFE96FB64FF70C621BBBE1DB7FCBFF515354F44DFE94BD6C3A4
8374E8A5C1039F8DDC8CE77289FAA88B5CCF0FFA5F31E5FEB29F9FD2AD6668AC47CED537FB15EB3BEF779579667EA51C8367
D6E3B87FDF16FD9CB8B8766DD4931ED3EF954B27EF1DFC936FECF7E800FFA107CF1E37EE57C5EF6FA51C86C7F2105F19B1B2
5ECB871F7EF87BE79F5DF393E8B3FED9073C5935E9FBE5FDBABC80F7DDEF9BB2BCECDBBFFACBC4CDA1DBFC981FE265F3CF87
D616EFF7A59C4B67E80B5FA177F48BEED109BA816FF4192FE9F4E90B6A7E12B7F601FA66BEBDEF9DF260D9FF738E4ABEBDE5
A863AAE2FB46F61FADB7B49A5BE6199FA13774064FF1939F9FDBB9FABFA2CEEA6BBF42C7E83E3D934F5B3D32EE57E5EB4C2F
0EE95E27BE272937E0E7D55D26BC1978E0077250F6F1A43774864FF0D51BFA56E6FBE28F86D6C6F721F044DEC533BCA3D367
AEABFEAFD03FFB2CF5560F3AAC7F7845D7F16BC1D6EA57631C75311EFCCB8378C077F1065FB67C6FF488787F5A2EC703B83B
6AF2880FFEF9CF7F96F659953EEBBB7C06FF702FFF651F49FE416F72AE4B7DD1273A8E2FF8633F46777A5FB6FFADC5F5ED34
AF165B478F8879C241C661DAEFC1E3F3D78CF820BE2FAE4FD64DA7E91E1DA403F82FE7E69C99F4F8DEEE0DFE5CE487C46BB9
0CBF733E49FCA4E7D917122FF9039CC91FF0C6CFE4C5C90BF6B8BBF85E4CAA8FBAD89FC929478EEEF940FCEE91DC94FD34E9
347DA6877490EEF073EBB65EBA00A7F421FB51F221FDC7838EB74CFA7E79DF6CDFD6BE61B33F17FD4A7DF3BCD7C3C7942726
BC1975E39F7C006ECC8B7E8DDD67E2E6F0693CE0E778356CFAE4BDCBB929FB6FEA53F6A9A4CF5977924EAFEE3BF927A1FB72
39FEC9D12DEB75DC52F637FA44AFE8447EBF3AF5AFFBE2461F05CF737DD27E529DE043FEBA6D58BD96C5E753532EB16EEBB5
4ECF5B6FF6B58447FEF6788749C383C7F82CA7D34FBE4A47CD035ECCA7E4E3CF957191F3C0EFF2FDE27975B07EFB60BCA4B7
F437E78DE463FC6BE98E131AC4755CFA2E87DA4F775ADBE49E725D732E4DFEB6B945ABB9F1BD10FDB52FD167E3989FF136F6
99D0295ECFA7E502F55447B99F7FF3996BCEAEF09A3ED365FCCCFE98FAC6276FBCB7F6CEDFFF7EDBF56C79870EC855F4801F
E0131ED17FF8E6033927241FE35F5977930EC375D6ABA4537C3EE7D8945F1F9CB5FFAD95EBFF15BEC8E5782397D3DFC113EA
778E1C9A7D30D5831FC2059DA7FBDBFCB7721CDF7C6179CF07E2FBCFF2013F90FFE49231D3EAFE237C80EFA88B753B0F5FD1
BFFB9A4E7A27FA6CBF0447B74E9CF850FCAFB2FD3BDF803B38E4BB745F5E9F317CF7E2FFB71E3A74D28A329EEC9F9CD7BECA
F9BDDEEB3C3FFDEB090D020772107D854B3E049F700EDF72003F574FF5E51303EF9BFC9372AE9177F06FE3DBCDFE1CD7E7F8
D5511D8F1F577CDF34E53E3ED0F5D7938697F342C649C287FD0ABF823B78FBC5E335A7873ED005B9188FF8393EE51C9BF2DA
A327352A7E5F957EB77F6042A7F2BE9CFED04B3AA4EF91D7CAFDD767F9E8B9BF75DC72C71D4F6ECB87A94EEB7B35B9A7F85E
617A3E5F2F4B759733E8CA15B7577052DDBFF6CED07DBCE78719E709DFF20EDF556775D73F7EA28FD6AD0E8EB70F68FD6E9B
0D657E6F3A647CB778BF59CEA097F4081EE190AFF77D75D23BE5DC810FF2E59C76BB0F2E7E5F2FE1A2E76F27AD28E702BEA8
1EEA7062D79AD3A3FFF2111EC0237CE267E663E2219DC313FE7CF395758BE79F7FF2F871F1FD6DFAC2CFF8259CAEBABE7EE7
B8BD6ED6843B0287749CDFAAFBC86B273E54F673F7ADCFE397CCACFE2C6EE15BDEC42F7E8067724ACBCBBA8C2BFB86FC8D57
FD9B35FC5DE4467CA2E7C18BB8CD7D4FF95BFFD5453E501FF7D5431DD61D386163F1FB63495FE4117AD0A2EF09BB14EF07A5
5C96FD34E914FCCB73F871E9A60A4EEEDC73E2CF8BFF8D4CB9CD3E8A8FB5FBDDF86E312E9ED30DFCA61F7889A77003AF7803
B77C852FD6B4DFF385871E7A28F30D6EE138D73FD55DDE877B39FB940377FB6BF13B3F895F7845AFE75E5EFD59F45F6E8277
7A466FAC877EC32D1CEB8B7EC8374FBF3FA66EEC1FE415FCCEFE9278865FCB1BD516F8E07F74575FE42FFD5954AF7A71F97A
B3DCAA5EEA8DEF5BEE3B6197F8DE3BBEC90719E72967C0BB3C9AF7697FA8D4D3BEA4FDCE477D5BC6359EF20D7C9557E8E988
FDF7FE4D7C4E0F8EE08A5ED65E36E18EE029FDA3D359DF92AED11BFD96A7CFF8D96EB7C46DCE0FC94F36D76FB321BE47851F
2DCF1F744FF1BE5ECA77F24BDB3DC7D42D5F4799F7B7EAC5816BB90F5EEC4FE0AF77BD038BFF1DC9F549FB137592EBD5F5BC
8F2BF32CED039EABDCEEB6246ECF7FBCBA4BF4575FE5DCCD978D7C22FEB74DAE82AFEC7789E778CD7F3AFF61C2C6B28F653D
4C3CA68BEAA46EFCF5E1EF6A7B96F12D67D11FBA49879E9BD5655C7C6F961FE3CF73A70FBA277E4700FEE4C5BB3F9A7844F9
7D0F78CF754EF5B51FCABA91F4826FF33338821FFC923BF0EC9A86359BC22FF88A9C028FEAA68EF0D2EABD61DDBFF8E28B9C
7FF9161FB8F7A586BF2BAE3FA7E75BCF1AF944BC1ECE4EEF5ADD2570258FC3D7B37B0DEB5E5C874A75D1E79C73523DD5979E
0E797FE2CFCB7984EEF3E11B36EEF942F13F5B699FD7E9F4C32F2FE3547FF5553ED617395A7FCE1D575D27721C5CDA7F64BE
279ECBA1F48EEED25BFB467DED7E44E38B63DC3CCFB49F345FBA4A4FE93B3DA66FF4EECC9BAAEB449DED77D43B9F379DCF79
FEFDBA9DE2BEF538EEFEBFD73E59AE4BDE4FA4C7DD57A7EC47C987E8963EC929FAC5E7E1087EE48E56C71CB753E4133925E7
AB941BF477F49BB573633C3C9083E8197DA3BBEBBF1CFF49F08C2FCA17FC589E9B5CBFCE1FD6AE5D9BE7ED3CC6CFB922E90E
BFE053F4865FC9BF6D0E6BDFA6EC3F396726FDDEF8C551DFC6F739E8367DA1DF93DFAED914E3F255BE0167F49B6EE77C95EA
7DDB9FF6FECDFF2F7F799C5EF6F87E6DCF721EA11B7C8F7EE0BFFCF7E0DF0EFC2C742AFB54F227B8E0D303769B7844391FD1
79FA0EF7704D3FE969F6EB5477FEF2CB9727F489F1E409FE7CDBE6DA3DE27F4EEDBFE8227DA257F00227F021F7F0297DE657
4B9FACDE1AF5A1A772A2BECBC1FABFCD072BFEA74E72361D19F4CEBE4747EED387E3EBD63E59CE3BF4193EE8341CC1159FCB
BE9CEAE5FC39C7263CAE3FFFF0CBE3F3B6F82827E165CE0949D7E9391F8253F894339DE7D1971B5F1CFB2AFA4FF7E1409FF9
82D7F127F88457BE63FDFA291FF0876BDEADF387B87DF6A9E3768AEF81F3697A4A273AFD60FC27E57D6DD687A40B7C58CEE3
C7A33FA9B7207C20FB45C20BDF807BBA08FFF019F9B78C53BA8867F895FD34E93D9D375FF3572FB88227FB89ACFF69BE7CA0
E423C5F37C813ECAA333FE4F9DBDA21E782777E01FDF5667759FDEB17A6BF0834FD2337A7FF4AEBDD646BE92B39E7C624CF1
BF64FC8F7EC3B5BCF0CC6BEDDBC4F717731E483EA1EF70834778E538E398F7A3BFA92EFC0FFFF3F5F1C40F7CC137E3EA339E
DA37B46F7874DFF89D2D7A30FE88097D4237F0868FD113BA4B572E1B553D3C72017D9167E110FEE411381BD1BA768FB20FC3
11DD0E3EC6EFB6AFBA76C2D2589F3C6A3FC6FFCCD73C8DEB79E3AB9B7AE5BC95784F07D4918FE2EB15FFAAB357B1FF52B7B4
AF503FFB1038EEB7EE8083427F72AE48BCB9ABC1BE4717BCE59B0977F066BD729F75FFE6BF7BAD8DEFC9D3053E74FD9C9A4E
310FBEC20FF946C72E63FE5ADEE7C12DBCF215FEDFE686639F2BEB2BBCF2E7E53FAD577C8E7FD68AEAE1316F3E8F3778286F
CDFE668F33A38F7844AFF189CFAB93FAA8A73ADA7F64BE279ED3BF9C1FD278E6F9E08A9AA14F3FFD4EEEB7FC88C7AD77EF3B
33F28C5CB37146D5CAE23A60F27BFCBEF7B29A6F631CF9527F2EF970FC73817FB897B39CD7F9E45479BFD3C4E6ABCB385A31
AD6643E8AC1C40B79CD7F1CEAF9FFABBBDBE763FBA7A45D997E81CFCC053F6A3A447709EE79972BFF966BE273DC07BFB077D
EAB5E1E04B43A7F9107FAA7E65AFB1F17D7F3909CF96ECB0C799C5FFB7A57598BFF5E4FC90F475E3DB47F78DDFEDCBB84D39
067EE9DDE089139696FDCD3CCDCFFE63F3877D67469EC32B7E8F5FFC04CFF879FB2BAA5606CFF0CD3EA6DB67874D2FD76FEE
EEE39F8B3AC929788327F2E8BC2FF7F849F17F26691CF3339E1C464F87AFAEFF687C5E8B8EE0253E4EB8B4A653F489DED075
786FD3614083C829F21A1CF27F787C6AEBF0AEC575B49487F063D1A9E3AB635E702F87E69C90747D4DE7030E8A75E99B3EF2
397DA2DFFA2517E6FDCC0FF6F849D9377B8FAA195AF64FFCC55B7A27871DF3BDE15DFFF5AF7FE5FE3DF3EB010DE2770AF45B
AE7BA6E6D8E78AEB25C977F1E7F18BC66E28AE3726DDE4C3748FDED127F3EA3FABE6DBB20FD31F3A4887F4A1ED25831F8D79
9A379FD34F7D7C7AEEE0478BEB5069DEFA034F72FCFCAF777F363E479273455A9775CEBB707C75E0468E83D365DFAF1E5CBE
BEC8FFE110FEE4B87C1D28F19CFFCBF5FAB970C7DD9F2D3E8F91788347F48CAED21BFA6A3DC65BB075F7FB035FE77F3B7E8F
B8DE097F726EF6F1A46B7426E7A5A42BEBAF6DBE3A3EDFCD7FE0176EE9257DA017F03AEACA9A0DE5DCE1BEBE799CCEF2D51B
0FDF6B6C992772CFDA630FBE34F412DEE04C8EA18FF08F0F7CBB6DC7AEA795FDF6F4FF18BF47E435395A9FB38F277DA24BF2
60FB716D9F8ADF6FA3AFE72E18F797A8B3FACAE92B8FAAA917AFB75FA01F395FA57A3FF6E3C3A6C7E7B6CE3C7BDC5F623EF2
BC79A90BFF561FFCA11378B460EBB8D7621CF3301EBDA2537481DED2597CA587784817F9969C316963DDA9F1FB817206BCDC
3EB0FEA3D13FFD305FFC872F78A38374882EE1C14D5F565F1AF59303F8CFA4EFAA07078EE83D5D572FF3321F7CEA7ACAD80D
E57D58D6EDC433FAAD8F9143CBFD84773A4937E9039FA413F41ADFE876C655CACDF0655EE609DFF2347D19B460BF33E23A52
F6D3A45BF48A1ED24138A0ABFA4A5FF11BFFF0EEA1A1E36E2EF607896FF2237EC919996749A78DAB9FEA449FCD8F3ED36578
BC72E1F885F13E3C9D95B3E8BEFEEBFBD34F763D2DFEC7923FD1A975AD47DD565C2F4AFEC7AFE8AF7ED061EBB54EF935D73F
D55D1E5FBD7FF5D0679E7927E760FB96ACC3095FF418EFE9003EC0B9DC01EFF20DFF187DC93ED38BEFAB66FDAFE8BEFAE91B
FD8CF3C5E7C5F984753A2EFB5DD2457A485F863D5FE1017DB75E79DCBAF90A5D5FD6BEEED4F2FBA8F4821FE23D1FDE7853DB
A756AF7E2AEF033A9DDEE2BB8217891FFC597FF5FBBEE6E37F1EF3D1A7BC3F488FBB6FBE7CA0DDCE3D5A947387F39B8FF365
9D4FFDA4F7EAC2DFD427FB4BD2757A3EEE07D59716DFCB4FF8E587D3E68F5F18B7FC879EE53AA4FD867AD05FB96AFA863D6F
8ADF638573F8E697F621F0DAEFC4831EA9F843F28BB43EBE613DE66B9EDB7CAA721C7FA14F3DEF187773D987E1806FC1035D
EFBC79D46DE12F7439E33CE15B9EA203F4951EF01BBA75D7D9FB9D51FC9FB5F5A675DADF5CB675DCEB310E1D91A7EC4BBABD
7DC405F1BC7AF3AD817FAA1E5AF65BBA6CBDD66F1DD6E575F242FBE73A8C2DEBE50DA78D9FFFA73FFD29FB19DFD8F4458F16
91ABE817BF59BF65C4A8E2FA5AF26F7A96E79FF645D661BF64FDBD161F3A23FACC0FD54FDD6EB96FFCA3B13E79850EE967DF
17C7FFBCECE7F4177EE1369F379DCFBCF08E4FE01F5FA41F4B8FDDF3A6D031B8E1E7CBE7EE333DFA8CAFEDBEE9FD71F179DD
9453BCCEFABC3EE78DA437EBCF6FF15DFC4FF2ACEF8D7B3DE6CD3FF14FBFF15BDFF10EDFE8191ECA4FB39FD97387D05BBED7
A9C58851E1B370AFCEF2857AD3337D87838C9F940FE0E8E1AF46DF5B5CC74B3CE7CFF8C007D6FCE2A047E27DD04B678FEB19
B8807B79070EF90F3C6EDAA1F7C7F1FD67F8E3DF4FF619F25EE44879021EF1014EE046BDD59FCEC2331CD32F3EC727AACFDA
BB577CAF53EEA10BC387EEFB5DF8305ED309FC8627BA0357F0D27EAF7E2BA33F729079EB3F9CABB3FACA0735178C9F1FD703
F080AE677F4FB87BACE111C5E701E50C785FD266CF1DA24F7252C70786BCF7E5975F66DCE81FFF5E7BD2A1C5FF23C0115CD1
9B8DEFF55B19FFA3AB4FF2D5AFEF183824BEDFEE7EEE5BCB0E63E37B717C0EFED4A3C3D08143623EE6977D2AF1065FE62C19
D73372AE1C814FF853D5B78277BE927D369DD73CF88DBAAA67E67BD20BBCC7377971DEEB75AE089CDF71CDB82DEFBCF34ECE
ED7846E7F0032F32AE129EE43BF3362FF35FF4FBB1BF8DF74FF047BE5BF6F7710B8347749C5F994FE4A6F2BCE49A4E6B5B9F
5EF03AF11BEFF10E0F332F121FF8E4FD4F8D591F39072EE520FED363A7D1F796FD2AEB55D2297A2617C2E9A22675AE887AC1
9771B2AF251FEB30A1F3AE655DA70F7401CFF1002FE882FEC9AFF3DFAF3320AE9B5897F5C873F30E1EFBDBC09B9C057759DF
921ED0B9BB7F35F6D9725E955FE9029FB9F19CBD7BC57ABDDEF91C87B774047FE568F55BB8679D01515F3CC66B7CBDA0D3D8
0B023F7023575ABFF1166CADD32EF89F7192F21FBCDC7CD0B8B302477C9A7EC335FDA01B740F7FE85FCE03891FB7DFB1EF77
711C9EE3375E1B8F2E455FE3FA977C893F835A1FF8B3E277BBF978F263757E74C9C80B8BEB3949F7F918BCD331BA76DAFAB1
17C43CE450FD37AE75983FDDA3777089A77C0F5FE9FED0C9E3B694FDFFDC9BC60E8F7EE993BC3CA9EEB885A15374910FD04B
FAA9FFFC4C6E987466BD01F1BB2EF4563F1C875FFC1ECFCE1C377678D4410E570F75DA5EF7176C1DDB33E66FDE8ECBFA9974
931EE43EA6FEC9D1F4E3F863C6AC2FFB3C1CCA35F04827F1925ECAD7DD3E3BF2B7655CFFFAFACEBB3EF0C096EC437430FB66
E21B5DD2273AAD5FE63560ECD867CB7E9E7388FCD1ABF5E9F17D0838843F7E9FFD3DF5597FAFEA37F69EC00D7D973BB20F26
7DA24B79BD29975BB71C43C7471FDB60557CAE5E5E57F77E9D0EB92AFAC627D547BDE8B2BEE893F559AFE3E950F7D9232F0C
3DA0EF598793FEC2055DEF7860B717E3FB767C38CF33CD6FFBBCDF6B7ED3CA7AD3BAF569CC5B15BCF055F7D5DBE3799C34AE
F1B24F251CD0DBACAB0977F435F32EF18D7F5AB73AA8237DE70FCBCEAD578C9FF343D2D7BB361FF8B3E0CB36DDAEE875AE5F
C201DDA257D92FD2BC1DF7F84F87352FAE4F25BFA48BB77E33B63672831C8367F08D07709E7195F024D74C5D33F69EA82F5D
A75B70633CFCA17BEA4AFFE494CE0D8EA95BD6A7070754AD2FDE5F4EFC9147E89D9C31FD82BDCE8BDF2DCA3C4A7A8C4F59E7
93CED1B7ECD7C93FF425E78AC4EFC77E7CE46FE373A9379E3B765CD4892FD27BFEA09FFAF8E41FBABD18DF63A43B7C48FDE5
387DC8BA9AEA445FF1855EE3CD36DDAE1C47872FBF75CCD7312F3A22BFA8A73ACA43D9DF930FAF597FC855910BF94CD74F86
358FEBAF788EAFC7EEDCF3AAE2F396C9FFE91A5CC043F6F7E4EB6B1736BD20FA95AF4FA7FCA1FF700BAFF4549FF8F3F20D0D
56C5BCB35FA73E785CDFF4917EDDF9CD98D6EFBEFB6ECECFFAF0C4BC41AB8BEB20C9D7F1862EEAA3BEC203BF828B7BEB8FFA
F69967DECDF595FBF80D3D587ADA5EE705AFCD27726239D71BC7F172247E8EDC616C6DD48D0E6EF8A2E755916FE08C7FEB9B
3ED227B90A9F87D7D9FF17F1FED4CC9163BE8EFEF11FFCA667BDEFAE5A5FF687CCC7C443BA287774593468F5575F7D95F104
2F9D763F7E64DC9753F0596E987D43DDBFC6F761F0094EE173DD86137E1AFFD32097E81BFFA3CBD5EDEA9F5FFCEFC199957E
C1297F5047BEA19ED9EF92CFAD7BE798BAF13BA3EAA70EEAE279AF879BF51F1E3F32F2937E98AF75773EF6849F7EFDF5D7B9
5E97FE73CCA4E02FDECA3B784DE7F0DBBCADCFBA727F535FE52FFCA523788C97781AFC89FAC953DD8F6835A4CCD7EC6B4937
E1597FF0DEB8F8855778CE2786EC30A675D97FE8083CC2219FC6A7253575FF1AB880C7FEEF8EFAB6EC3F5987132FF0E4DA57
C65447BDE838DFE03F9D57767CAEACBB6B7F33FCA2E2BA4DD223BEA43E72983ACDF9E19849316FBE8D4F78801778B9BCDF98
FAC17B7E4C57329E138EF9309EC855F31EAE7B72B16F4A394FBF07FDB8D1FBE10B17771AFD3FB13E78926BF0110EF032FB5A
E2DFED1FECFF8BB2EFE227DEE3231EEAB7BEE299FEE679A6F9C97DEAACBE7C4BBED1FF453DEA9E1CF3C2437D87037DC63B3C
941F7B2D6EF646E8B5BEC1319FB97153FDF3637C7DE1DBFA83A77889AF788CD7FACC97D4553DF54FCE9BFF66DD23833F1987
C917E1F1E4F5A3FF277025D7C0173DD13FE7CFEB4D79D7BAE19F8FE1813AD1297DCDFE9BD6A13ED9D712FFF151CED6D77ECD
9BDC5B397F7ADCFD34AFECD7F6512FB71A123A7AC139A37F19F38257F951CE858B858DEA1E19F58403B8A303D9D7122FF151
DD8DB7606BDDBD8BFFC33CBC521F3AC437F0B5D7D1C32F2AE7FD75C33A3E17DFCFA1FB74C7F35E0FAFFC895F4EEAB6CFE2E0
F769678DFE659C47BED367F5501F7E43F7CDD73CCF7D72F4CCB8552FB9D271E669DEF8397ACD98FA653F5CB368C4C7C575B3
C443B9E6F6AEA336446E919FF023EB52D2233A7166C7CA7CE455EBBBBAEEE885711EFACBDFCCC3FACC475EA183A377DAAF43
F8DB82ADA3C714EB4EEB741EFADD654C972BC2EFF85ECE15C9AFEEFAACD1FB31DF875B0C39A8B87E94FC89DE67FD4CBCA7A3
EA20C7A987BAA9973C97FD2BF1848E65FF4DBEBBF6A4666FC4756478C57B7CBFABAA6A5439373ABFC7E37CE5E7B32FA7F3AC
79B1C9BD813B7EDE7D719BCFCA7A93F19970C97F72FD53DDE58BAC6349BFF49FCEA9B73AFFF292519FC6BCF81C1DA223F478
D963FB2C8EE3F5453FF8047DD01F3A41777B6C197250F0156FB38FA4BAD38D276EEE7245E403BCE51FF0437FE168CADF2BB8
A52B7402BFFBCD1BF171D997E09CEEC2BB3EEB7BF03CEA4CBF863F3E6A43D93FB3FE279DA00FDB74ACA25FFA90FB9B72893E
67FD4CEB771C7E746D755C75EC1BE48247CE1970757C6F9EAFE347F6CDE407FA45D7E483E9557BAF8BDF63D1077DA18B6ED5
4D1DB38EA53ED333F8E1CF70A4DE836EA91A55F6812B4E18B5399EA72F7C9BFEE93B1C64BF4EF85DFED57E1DA29E599F939E
D131FA0B17F000EFF88777395FD8C7A6F1F59BFEE9BBE3F0176EE43C3C1CF6E14153E33A61D6C3543FBA28171C77D680ABCB
F9E4F1678FAB8E9C04F772C036DDAE8C4757E1A2DB11BD0F8FE3E58E471FEAFB8FC87F7203BCF02FFEF3E8496D3E8BDF6373
EB71FD90CB7A6D6831BC8C5FE7E9DEB3EF3FBEF9E69BBC9EC75EEE7D78E415F3767EF551177973FCDC519FC678FAC667EEFB
60E85DC5EF93267DE17BFAE138AFA7E3746BE92D7BAF0B1EE10B9ED0AF19F7547048D7E907DED117FCCB7548B94D3DF004CE
E19E2FE67DEFCE0DF7287E17101F130FE9BAE3E10E5F569C5C352BFC83AFE319BCC00F1EE3179DC333F55447BE4817E483D9
F3F7DEAFF85E4ACA8BFA34F09B43772C709B7CE9B89FB5BFA7ACC7F80AB770BCE4AC51CD8BFFBF4A3CE7DBF266BF4E47FC47
51CF54D7EB1EAFFA30D6CB87E8B4F998C7F6B9D57DE3A99B5CA07ED9D7923EC209BD800B78C8BE99F463E5EE074D8D7AE117
5EC159EE7BEA377FC6D38CA7C4D7D5770DEF1AF38B7EB8DE14F7F109BFF4993EF7AD33F4AEB24EF3513C58B270EFFD2AFBFF
8A3EF1B3B5C7B6185EFCEE5CEA0F9EE22D1C19DFF9669F53E923DFC31B78E06770B1F29C116F063FE4217DB05EEB9493E889
BA58B779D1117537BEE39D877EEB9B7EC115BF822FB89793E62DAFB72AEA75F18AAAB5E5FDA55C839FA32EAE9A55F60778E2
27377ED1708FE23A6BAAB7FA87CEC6E7D81EF9B8FD3DE17374167F321F130FE944F6A9B43EE3675F4BBEB47A8743778CEB14
0FB41BD4BEB83E92F48E8EE2239DC0CBEC6B092F6BD61FF11F71DEDAAE551FC6FAF0834E67BF48F3B71E75A1CFEA23DFC0E9
A2D1F55645DD4F1E55A9333F872FC73BBFF3EA935C35FFC57A9323AF65FD4F7C867F7C86DF8CDBA45BC73DD8E9CCB2DEE3D3
C041C3BB96759D8ED3EB49FBEC5BAFF89DEAE4F7783EEA9D03F7087FBB605355315FB8919B8C37ECAC116F96F55EEE53A785
CDEB4D0E1C58A7E7176CADF8C7D57FAEFA61E0993ED27B7DD76FBACFB77BD73BFAABF7DE7B2FF3153FE0087E4E6B5799BFDC
A13F74A8CFA641EDE33A07BDB8E9CA911D637CFE8467700457F870EEBB55CD627C75918FE88EFAEB873AABAF1C91F99B788B
37F80BEF787C66FDCA79E528EB531FBEAE4E0BB65635283E2FF0E4FFCE5FD6AD0EE6FD48EF4E67C6EF21E2339D55573C5397
8CC3A49BF0F8E0DCFECB9F7AEA4FD957F0186EF00C7EB2FE27DCE357C66DC22B5FCA7A9EF88777FAC83F73DE48FCE8F19FDD
0616FBB8E46370485FE151FFCC17EEF1935E2EFBCBBEF5027F53F6ABFA619C175FE9CDDDFF33B8C051E0D87590B84F6FD459
DDE515F818F6E743CE2DAECFD18DA417F0C2677A5FD27F79D967E92C7EAC6870E01E719C9CD2EFC423CF8DF5C24BE67BC21D
DE67DD4BFDC0EBEC0BC90F1EFCDBD15F455FB28F247D579F877FDE6D60F8343CF1C5BCDEB44EF968ED157D3E2AF27EF24B38
524F75E4C7CE67DDD69BF53FF5475FB26E279DA00FB93E2987A993F595FD2FF60370D47352CF838BEBA4C9FFE92CFF9BDEAE
41EFA24E8917F840EF6EB96AF859313F39037EF164DCB4911DCBBE74C517237BC6F374931F6ED3DB8ACEC28B79F59AD1E7A3
6FBFFD36E326EB76EAABE31E5AD6F3E0C82BFAA40F72AA3E0DFA71D36181138FE7FBA97EE6AD0FE66FBD729075675D4D78A4
AF700DE774EB8E87863E18C7C3813C85EFF408EFB32F241DA59FC671BC7965FD4CF3A2A3FCA7F765ED1E095DA07FDBF4B972
1CBDA5AF037E32F8D3B23FF05DBA36FE83FDDF8FEF0DD1037CC11FEB368E71F19C4EAFDCEF9073C347AEBF6544B358073FA3
7333761ED9331EA7A3748BFED2B9A59B1AF42EEB3DDEAFF9C591E7C6E7B1F0143FE13BE33FE19E7E67DC26DD855F3A5B3575
F859657F8053B8D57FBCD357FCCB7D49FEA73FD9BFD23AEEFAACE9B058C72F761C715BCC037FF8ADBEE807DFC6673E3DBBA6
FEED5107FDC117FCD1E7A13D873E58F613CFE3953AE2A9F558F7BD8D079C221FC53E9E9FD039FDD7F7EC47A9CFFA9BFB96FA
48D7E0847FC24BF6BBA4E3F0A04FF8AD5F0FCE6A57BC1F471F03C7C5EF2DA4F378DCEBACCB7AE8FB895F8FB8ADB84E99741A
4EF9109C2EB9A1FEEDD1577C506775E70BEAF4CB3AFBBF1FAFB71ECF7B7C42D58866811BF5A35BD6E9798FE37D9FA61D1797
F5461FE00EDEDCC2A17EE4EB31EB5A3529709BF06B3DEA0D3F8B7F34A269E0165EF9365CCA07F32EAF7F745C679727E078D4
6F0F3E22FC36E33FF1120FF0ACFF1B034E09BDA13B19FF09F7F8FDABAB87758DF763F9161C39CE38FA94D79BFCD5BAF33CD3
FCE4826B6E1BDE355E4FBFE847BE2E95F465D2F7F7FB28BE0707F770043F7C0E8E16CDAC5F7CCEDEFCF99275ACE9D5F7D5E0
5BE85AF810FEDDF95F83C614EF8B24BCF29F933E1FD1348EE733F09B7535F1076EF5894FCF7F7A9F8F031F7045BFE0EB8157
3B2E8EBC8457A13BF17942F3310F7E76E18F86BF1CF7E183AFF2E741AD9BBD11F3D707FEAF4E0BDBEEF371E0E9D6BB87F42F
FBB53AC901FA34ECC326AF459DACC7780BB6EEF3BBF003E3989FF1F0A6CFBB5DA614D7F593AEE779A65C61BEF9BCE97C9EA7
97FDD6F67DB5EC17A77E3EFCE5380F7FD61FB8A27FF0957195F0C407F0959EADE9DCAA49F17BC9898FF88497590F53BFF42F
F334F1139FCEFBC1F01BE2787D923BE010FEE866D6FFA4372B0E3DF888E027BD183B6558D7B24EDF37B8F79EC5FB45FFF69B
F85C235C9FF1E5F01BCA79479DACCF7AADE781FA5DA6C4E768E02FFCACF27EE9F08B631CF336DEE411C3BB06BFF192BEE227
DEE00B9D18B2EFA031655DCFFA9FF883BFF072FCA3C7758B7EF34F7CA3F7CBBEDBEFA3A893E7FBAEEABDE777DF7D9771A35E
745EDDF8837AAB337CC11B9DBBBFFB71DD220FAB8B7A1B573DD4C1FCB23FA5795AF7C80143FA977DC079D447BD328F925EE0
53F69D84DFBB36377B23F472E98DC3E687EED1533E907529E9119CC299F156EEDEE4B5789CCEE001BDE163FD9BB5F9B0F81E
6ED2337ACD3FA61FBC5FC3C055D6E7D457FA031F788CBFF804877009E7F08D3F375C38B469AC973FE13D5FC4F3F16F37DA21
C6D9A66315FD528F55DFF47FADB8DE9D78C807B27EA6793B0EBEE00D5EF415EEF457DDF8A2FAD107F5A613FAA9BF780B1FCE
034FF8021F78A3BEF4800EC0BBFEEBBB7AF24F75D507F8C707B7F8A04F70488FE1114FCD4FFDA6570F9B1F75A11F7425F731
F54F8ED02F7DE297F44A5DE956C667C2251FD8A69395E3E81B3DA2374B7FBF5FC3A25EA94E39FFA53EDDFB529B0F8B5C4C4F
138EB3DECAC7E975BF5834F4B5A80F7EF2193AA21EEA70F3E5837F54CE33FAC3E7072D6871CEDFFEF6B78C2F3CE443B30737
3CABF89E4DC2E9E01DFABF56CE71F05A73EAD0A61F7CF041D6B5DB1E1BF041F1BE49AA2F9F779C71E032EB73D23B75337F79
C13A8C6F5CCFD395FE8FB6DF58D63FF35327F39457E063D43B879E1DFBF4ACF309EFBF6CD86887E85BAE5BCA1FEA67DE8EF3
FA3C7E1A574E59BDD3F18795FD9C7E9C386FE86B312F3A07D772CDF03A478C2CCEEBFCE971E37A9E6EC2E392550DCF8AF92C
9B3EA44DCC97DEE1195CE11F7CE57A26DD51D7C5270E2DBECF017FF4FEDEEEED3716FA9B708CFF7467CCCCC13F2AEB9CFE8F
E836E083B23E651C26FCE11F5DC6C3895F1EF483381FBCF2A7791736FC61F027EB61E2D35D67B73827C6517F75A70B7465E0
57C71F56CE8F5967527FF405BE4EA8ED746BB13F4A3A9A7532F17E458343CF8EF5E0B1FAA897BED051FD59FD65CF9AC8BDF4
115EB2CE249EE8AFF559AFF3E19B3AABFB49BF183A34D643EFE18E6FC1D1A2531BFE30D6A17FF8A48FEAA46ECE97752FF5E1
F60F8E1859BEBEA57F74FA9EEB3ADD1A3E01377C425FF9EEFC27F7DD18E35CD86FC8975137F8E3877C62E00F7AD66CDDBA35
F7557F0634E97661CCD3EBEE7EBDDB85C14FE7536F3CD64F7DE4E7EAB4B0E3BE1B037F93960E69137581177AA3EEEA07E7D6
63BC055BF7BD2BEA94F995788967EAEE3CC6A7DF83DE69DDBBCCD753D70CF9321EE757FA7CFB63FDAA623CB8A7D3F08CE770
9DF995788597F047FFAEFDC1413F883E9DD778C893D117FDE0AB374E1CF8413CCE87E01DCEE11B3FB3AEA63E183FEB64C22F
DCC217BCE143D6A5842B383BE38DCA3CF9B83A65FD4CBC827F75532F3EB460EB907B621CEB34DE955D074F8DC7E92C7D8523
FC30EFAC3F09B75987F437F98D3EAB279F54577E8657E33F683CBCF85DB2A49BF4747AC303DBC4E751969F7BC28D719CF5F0
537C8103BCF17AAF737E3CC21F7DCD3C4DB8C6577DC03BC767BE279EC31DBE9A873AE5F5A675F25B3A34BC5BBFAAB2FE657D
4B7DC5DFAC878937F8A28F745A3F33BE92CEDED5A075EFD075FE3C7CE8513B16FD4A7DF3BCD7C3A53E9BA779D38BEA6B077E
F0D1471F65FE6FD38D8A5EA8533E6FCA01CE2F1FA8CBE89D0EDFA9B22F4F3A96D66FBC698F0F9E1AE7C1337A661CC73B4FC6
55C2135DCF3A90FA4F0FAC878E5917DD1A3CE19883CB3AA46FA3CF3CE1C6B28E667D4B3CC6DF95AB7ADF5A5C8F4B3A42E79C
CFFAAC0BCFE9D9D2B70F6C13FCCCFA9FF8FACB3A8D8757FEEF22F13DAD0FEF2FFBE9A075711FCFE9F7369DA91C471F9C8FAE
39AFFBCEEFBC784B7F67F73FE0CAC053D6BD84C7DBEF386AC7E8873A583FFD52673AA2DED73619F0493987C03B5CAABF7EA8
8BF32FFFEAF09D625EABAE3FE6E0D01DF8C6737C1D36B8F7ADE5BC605CB8D05FEBB02EFC99F5C9A075D13F7C85537819FCCD
B1BBC47CE996F338AFF5E8177EEB9B3AD057F5B8E3A2E35E2FFC29E91A7CD159B89FF4FD43AA42CF32DF531DD485EEC0DD92
7B0FB832FABB624ADF9DCAFE4C8F17BF34A87DAC0F9EE868C643C201FED1C7A1A71CF77AE123A9EEAB7638769798977ED05B
F8D24F7DC413F55667751DB2BACBE7711EE7BB736097CF431F8CAB4EF04AC7E72D3CE080D89F5BAF75F215F55447BAFECB76
FD4F2DFB9DFED0D7E1AB5BBD10F3846BEB804F78328EF3190FFE26BE3EE093B23E657D48FCD4AF8C8784037C73BCFA1927E3
33F10D4E4F6A36A8FDFBEFBF9FF50CAEF3BA928E5B1FBD86A345F30F3820EAAACFFA0EEF70A7EF7070E192817F88F3C0139D
D627BE327FC3FE6F15FF039C7833EAEABE3B95F520EB55C2D1B2EF0EA9AA7CCEA9E237FA57BDF3616F45BDE019FFF2FE733B
FD33CE957F1ED0205E4F8FF0D5F9C76FEA7F6AFCBE9DFE39AFF36DEF57A38F6DBEA6F83D825457BEA7AE0B8FDDFFADE273DE
A97ED66BFD700EDF708017FA8A1FEAE83C0BB6EEBF31789BE793E6675E992F8927F87DEAEC817F089CF01378C13FE7D37F7D
C647F3CEFA9970076F596F936EDD3EB0D50BA127E7F51AB822AE7FC307DFF3BCD7FF3FFA9678029FF4914E4D3FF8E08B038F
FA4B77F5F98CB50357C4FCF8E6361F1D7859F1FF7B8DFFB7AFC201FEC1C32DC7F6EA5D7C9E23F193BED28911FBB7BDBBCC4B
38A61F377E71D85B511738701EF8CB3C4AB8C6A7CCA3C41FFDCCBA94EABD7C43F335F1F8B4FD063488F9EA0FFD307FF3A69B
B99EA98EF432EB49C20B9C645D4AE7317ED693D427F8BFF6CFFDE6C5717436F727F1529DD4473DE994BAC263D5865EBDCB3E
BC8D17153EA813DCD083A5BF3FF8E2CAFFCA249EA6F939EEB24B4E283E774017E84AC66DEA0FFC5A171FB03EF3324FFDD61F
FD82FB9BFEA3CFEB719C7EE47DD49FDADE5DE86CE2031E5A8FFA5957E675C23B9C1BDFB8E67BEB8EDDE716BF37917407AEE1
77C4CA0ED7C5BE8D0E99AFF9985FE651AA133EE9279C5A3FDED2A3D9831BFD3DF65DDB785A192FF331EF3B2BB89CB44F93B7
432FE9D4C8AFBBCF2DFB55EE63D203FDA4D3A32F697976E00ABE32DFF7EB37AFCCFBDB8675B82EF44EDDE9DEAD6B3B2F2AFE
2F29DD573F751AD9ABF3A2B28ECE9A7BC2AC78FF1D9FF032AF2BAD87CFC0CBB86FFBBC5ED609FDD2277CCEEB4ABA6F7DD683
A7D695FB9B78A5CF5937D27CCD33F72BF5894EE0399E2D59D5E8EFD16F3E643DD5ED9A2D0E5C5F7254FFBF465FF08ACEE99B
3EE24D1E271D6F7D9E571FF5CA7C49FDD23F38A39BF32E6C744AE0177FE074D95F9ABCFDBF7423F1CABA327F130F96CF6D79
768C030F70A04F579DDEF7BA582FBD8347F3D61FFD851FFDD81E47EEAB93F95987F95BF7F6789DFB7CFFBF465FE9043C5A87
75A9E34D4FF4DCA378DF27E18C3EE53E269EE9273D858B45A7363A25E693EB90F432EF8BE843AAEB8D9B9A2D8EF35FF842FF
F3625EF04177CDC73CE89ABED2FDF94F363A32FC0CBF477FD2A6367E5F04DEE15B3DD481BEC0F7F4764D768FFA65FC24DCE0
C17597F47EA1ACFFD6E771F3F33C9F51A7851D1B1D59E4B123FB9F177E47FFF4C77A1CB7606BA3BD8BEB22093753CFEF7B5D
AC076FCF9BD2BF4FAC57BDE871E64B979E7B947993F99BF0AD5F70A43FFA0277C631AEE7F5119ECC4FBDCD33F725E994FEE4
FAA7BAD3C533AEEEDF27EAC317B6F944FFD671BC757A1E7EE0148E96FFB44D6DA10B09DF78F2AB19DD8AEF4BE12DBCAB53ED
DCDE2F041FD53DE32DE10C2EACDBF8CE77F9CEC7BF17F3C17FBCC97C4CFDD777381D33ADFD90382F3DA07FD567B5181DAF53
3FFDA3074B3735D9BDF81F48FA90FAE971BAACBE93BA1D7179AC37F337E1003EE9CED82BBAB52EEB749E4FD267F3CAE3A771
3DAFCEEA8BCFDB705EC1B7BADE7C65FB2181037DA26337FFB053AFF005F7F52DE33FF51F0FD46FCC3F3BF58AEBBED693F998
D6E9FCFA45B7B26E27DCC3C3F6F8807738A123F46276CDA173621D99EF09D7F0ACAEEAEC7CF080277091F19FD68D0733BF38
FEBD783D7CE1DF361E558EC397CCDF84DF1BCF69313AE66B3DC6334EE651EA8B3E65DD48E32D7BEC88E277DBAF6FDE63A722
87A67ED3C5AB7FD4EBBFCBBEA0DEB96EA95EF420F73DF51B9EF2F8695CFA7BC925C78F8E75E21B3DB03EEBB51EF371BC79C1
0BFCE84FC655EA87E7F1072F97DC70E89CA8873A6D8F63BC543F75CBBC4ECF1B079EE8CEBCCB1B7F57FC1F51D283EA575A77
2BBE0F90F00247135EECB153CC0F5E9D47BFD445DDF040FD1C37E5F35EFF1D9F97878FB9738F1F1DE7816BB85327E7350FBA
A3AFD3AB0E7FAAB87E61FE7424ADE3A25DFAFC2BC6D53FFCC65BF55E34B3F177B11EF383DBED716CBEE699E793E641B7AEDF
D8F5A0E2FDB08433F8CDF54F758753FDA05FF39F6EFC70E0EAC6C35B772BFEE72FE10C0EE8E2A4338FAAADFCBE495AB7C7D3
3C3C9F719BC631AEFA5AB73E9DF2699F7F45FDF05E7FE89CFA2D6CDBF8E118175E6ADA1FD336C681EB8CE7843F75B64EE32D
D8DAF8E6E0ADBAA9177E5E7E66AF36655D85AFF32EE8F358CC4B5FE9193E4D68DFF5A0783F139E325F127E97DE72F85351F7
8C93845F78B9A1EFB14716BFC791F8A69F679CD6E7B1C823F4509D6ED8784CDBF02F7D729CF9D309EB5097ECBB693E0BB6F6
B93DF0653DCEA3AE35F71D7B647C1E6D7B5DB15E7DCD3C4DFD5F76EE5145FFE1589FF5F7EACB8EDB3FE6A5EEF42CE333E132
EB46EA3BBEE9BF711C6FDD99770907F093F372E2C3ECF9871F5CFC8E4EC2A97AAA6FE6E9B9BDDAC4EF01E0A7F559AFE3F531
EB7EC219FCE83F1C19579D8CBF0DB715BCAA47E6455A8FE7F547BFF4479FF5C7BC33EF12BEE02AF73DD5378F87BFB38EDB3F
BEAFA32E1987097FD69FEB9CEA4BB732DE523DE08EAE4CDA7874FFE27A66EAF325EFF6BC3D5E878778B38D1795F1E0D6BACD
D3BCE1000F962C3CBCF89C663E6FD22FE7CF7C49F5F37A7A673DD32F387259F1B9E6743EFD549FFCFAF43AE7C9EB4A7CB5BE
6B7EDAA55DF1B9B5841BFCC8FA9DF8306F79D33551F7B9F57BDE1EBF43A89F7094F19970B9ACFDD1FD8BFC936E335E53FFF2
79D3F9F02AE3339DC7F897BFD7BDB6AC33D6A16F93EBB7EB11AF87A3FCFA34AEE32EEAD9B34BF44F5FF123E336AD63E96947
2E2BFE9722E171F2275DDA9575171FF479D1E8A66BE2B85CB7542F78749C71F4EF9AB33B169F7780637DB8E6DD763D0A7CA7
7A795E3FF073FE8B4D67864F5BF7E4051DBBC6B8E60D47D6675DEA067FEA97719B8EB3AE531EEAD925AEC7E381FE671CEED5
BD367E9F43FF729D13EED4DBF9D5193EE88ABA2E6CDEB4F81EE8361DAA3CBF606BE57399E7ADE8B963CC43FFF0DB781967A9
8EEA9CEB9BF097F19CD6675DF8A67FB36F68FE551C679D8EB3DE3346F5DC3174943EA853C64FEA8BFE2FD8DAE3BF03D7E6ED
B8DCDFD45738729CF3AB7BEE7BEA9FF3E0FDF40DAD1616792EF539E32AE1495DF27AD33AB7D70FF78D97F19FF0A15F97B6EE
F6725967F43FF72BE142DFAEE8D7F9F4E27FC3523FF0755BDF2BFD568F8CC3546F75CEE74DBC77FE8CAB343FE3655E243C2E
A969FE553C6E3EE661BC8CDB84AFA5C7B65A58FCDF415ABFF18C93FB9BD6ADCF8E338E7E6CC34FE53878D0EF399BBBBD1C3E
0977F030634DE7D3CBFCD37F8F7B1D1CC3D3BC879B5F58E4E9848B19FFA7ED82A82B1C399F7A3BEF15B77758507C4E2BE146
DFAC5B9DACF38A7FB55D107EB43D0F9D67C6F00E0BE27A41C66D9A0F9C99579E7F9A379C5E7456B7A9F17A38832FFDD69F45
3D9A5F18F35117EB53B7BCAEB41E38B00EF3B71EE7853BE7C713F39BFD4CCBF78AF70B3CEE7E5A5FC64FAA9F79EA13DCCF7F
B379F7E2FB08FA9BEA663DA79CD36D6A5987F42DD739D5D171B92FA91FD67FE9399DDE8DF9593FFC673CA4F318DFEBBD4E3D
CE5BD7EDB0C0B17EE001FEE9CFC246CDBB97D7E9F9055B9B370B3DCAB84A7D58D2A6E57B651CAB0B5C9CD1B9DB61312F3CDA
C6AB6E7BC6EDB679549E37DEF678CD3839ABD3BB512FEBCFB84A7D32FFBCDE840BEBF63AEB703E7D30BEF3C1C5EC6FDAFC25
E6AB6FB9BF697DD6B5AD6E957AA987711D67FEB9FE69DE8EBB7453FBBF063EE140FFE1587DE7BDDEF2B4E8739E67C2B3F9E6
3AA77AA8F7921DDAFCA5C059AAAFF358DF9C76EDFF1A3EA62EF9BCE97CCEB3ADBF95F3E8639E7F9A375C5E744FA7E27F19E0
527FF4C17CCC2FBF3E1DEFB8DCAF344FF5CDFD4DE32D6AD2F2B4A8A7FA7ADEE359CF52BD8C97EB8387693DB98F27743A2D5E
97F997EEAB87C7F5090EE6BFDFB2558CE3FCE6E97CE7FDAE53DBE27A52AA8BBE9997F398471E3F8DAB3F7064DD0BF76CD9AA
3C0FCF2FD8DAB261B1BF48FD9CF7659B53E2FCEA911F4FE37A3ECF33F5C77CCF38A453DBD03778508F8BFED0FEF4428F521D
F533F7459FD23C176CED745019B7C6CBE36C87A7453F68734AE17BA98ECE93FB92E6653ED677CA81ED4F2FFC35D537F7258D
63DC6D75ABD44BBD735FD238CE97EB9FD753392EAF2BD5CFFACC477D8DB7AD9E95E3D4C1BCAC535D9CCF3ABC6EFBF1CC535F
AD67FED76D5A16B9D7E3EEEBFBC7EDDB16FD4E75CFEBF0789A87E7F37CD23C16EED8A665E57DB04A3D3CEF71B75E6FBC3CFF
9FB56F5BE02AADC3BC3CEE75D6B3ADEE6DF62DF42ABD4E3DB7D5B7FD81655C7BDEF8EAE0F86DE357C6F5BCF99BB7BAE4F1D3
FD6DE74975C875D97E9E95E3B6F5BFF2BCD71BCFF9BD3ECF333F5E396EDB7895FBDBCDEFFF02702A31AF
}

\begin{document}
\begin{abstract}
We realize generalized Markov numbers (GM numbers) for arbitrary $k_1,k_2,k_3\in\ZZ_{\ge0}$ as $\lambda$-lengths of essential simple ideal arcs based at a distinguished cusp of a hyperbolic orbifold constructed from local monodromy matrices. The McShane identity of Tan--Wong--Zhang on this orbifold yields a summation identity with multiplicities for the gaps between the associated Lagrange constants and their limiting value, together with a summation inequality that counts each distinct value once. These results extend the summation identity and inequality of Gaster and Loustau to generalized Markov numbers.
\end{abstract}
\maketitle
\tableofcontents
\section{Introduction and main results}\label{sec:intro}
\subsection{Generalized Markov numbers and the main results}\label{sec:intro-main}
Fix $\kk=(k_1,k_2,k_3)\in\ZZ_{\ge0}^3$ and put $K=3+k_1+k_2+k_3$. The equation
\begin{equation}\label{eq:GM}
 x^2+y^2+z^2+k_1yz+k_2zx+k_3xy=Kxyz
\end{equation}
is called the $(k_1,k_2,k_3)$-generalized Markov equation~\cite{GM}.

\begin{definition}\label{def:GM-numbers}
A positive integer solution $(x,y,z)$ of \eqref{eq:GM} is called a \emph{GM triple}. We write $\mathcal G_{\kk}$ for the set of all ordered GM triples. A positive integer occurring as an entry of a GM triple is called a \emph{generalized Markov number}, or a \emph{GM number} for short. Let
\[
 \mathcal N_{\kk,i}
 =\{n\in\ZZ_{>0}\mid \text{the }i\text{th entry of some GM triple is }n\}
\]
be the set of numbers occurring in the $i$th position. The same positive integer may belong to more than one of the sets $\mathcal N_{\kk,i}$.
\end{definition}

When $\kk=\mathbf0$, equation \eqref{eq:GM} becomes the classical \emph{Markov equation}
\[
 x^2+y^2+z^2=3xyz
\]
and the entries of its positive integer solutions are the \emph{Markov numbers}. Listed in increasing order, the distinct Markov numbers begin $1,2,5,13,29,34,\ldots$.

Markov numbers have applications in various fields, including combinatorics~\cite{FZ02,Propp}, algebraic geometry~\cite{Rud,HP10}, and symplectic geometry~\cite{Via16}. GM numbers were introduced to extend the theory of Markov numbers by generalizing the positive integer solutions of the classical equation and the rules that generate them~\cite{GM}.

The following conjecture asks whether a Markov number uniquely determines the triple in which it is the largest entry.

\begin{conjecture}[{\cite{Frob}}]\label{conj:classical-uniqueness}
For every Markov number $n$, there is exactly one Markov triple $(x,y,z)$ satisfying \mbox{$0<x\le y\le z=n$}.
\end{conjecture}

The author and Maruyama extended this conjecture to the generalized Markov equation with $k_1=k_2=k_3$.

\begin{conjecture}[{\cite[Conjecture~1.8]{GyM}}]\label{conj:GM-uniqueness}
For every integer $k\ge0$ and every $(k,k,k)$-GM number $n$, there is exactly one $(k,k,k)$-GM triple $(x,y,z)$ satisfying \mbox{$0<x\le y\le z=n$}.
\end{conjecture}

GM numbers are also related to Diophantine approximation, which studies the quality of rational approximations to irrational numbers~\cite{Gyo}. To describe this relation, we first introduce Lagrange constants and the Lagrange spectrum.

\begin{definition}\label{def:lagrange-spectrum}
For a real number $u$ and an irrational number $\xi$, define
\[
 \|u\|=\min_{p\in\ZZ}|u-p|,\qquad
 L(\xi)=\limsup_{q\to\infty}\frac1{q\,\|q\xi\|}
\]
where $q$ ranges over the positive integers. The number $L(\xi)$ is called the \emph{Lagrange constant} of $\xi$, and the set of its finite values
\begin{equation}\label{eq:lagrange-spectrum}
 \mathcal L
 =\{L(\xi)\mid \xi\in\RR\setminus\mathbb Q,\ L(\xi)<\infty\}
\end{equation}
is called the \emph{Lagrange spectrum}.
\end{definition}

The Lagrange constant $L(\xi)$ is the supremum of the numbers $A>0$ for which the inequality
\[
 \left|\xi-\frac pq\right|<\frac1{Aq^2}
\]
has infinitely many rational solutions $p/q$. Thus it measures how accurately $\xi$ can be approximated by rationals infinitely often on the $q^{-2}$ scale.

\begin{theorem}[{\cite{Mar79,Mar80}}]\label{thm:Markov}
The part of the Lagrange spectrum below $3$ is completely described by Markov numbers. More precisely,
\[
 \mathcal L\cap(-\infty,3)
 =\left\{\frac{\sqrt{9n^2-4}}n\mid n\text{ is a Markov number}\right\}.
\]
\end{theorem}
We denote this discrete part by $\Spec_{\mathbf0}$. Its first values are
\[
 \Spec_{\mathbf0}
 =\left\{\sqrt5,\sqrt8,\frac{\sqrt{221}}5,\ldots\right\}
\]
and its accumulation point is $3$.

The author extended this correspondence between Markov numbers and Lagrange constants to GM numbers. For $a\in\{k_1,k_2,k_3\}$, put
\begin{equation}\label{eq:Ldef}
 \LL_a(n)=\frac{\sqrt{(Kn-a)^2-4}}n.
\end{equation}
Then the following realization theorem holds.

\begin{theorem}[{\cite[Theorem~1.2]{Gyo}}]\label{thm:spectrum}
For every $i\in\{1,2,3\}$ and $n\in\mathcal N_{\kk,i}$, there exists a quadratic irrational $\xi_{n,i}$ such that
\[
 L(\xi_{n,i})=\LL_{k_i}(n)
 =\frac{\sqrt{(Kn-k_i)^2-4}}n.
\]

\end{theorem}

Define the set of Lagrange constants obtained from this theorem by
\begin{equation}\label{eq:spectrum-set}
 \Spec_{\kk}
 =\bigcup_{i=1}^3\{\LL_{k_i}(n)\mid n\in\mathcal N_{\kk,i}\}
 \subset\mathcal L.
\end{equation}
We call $\Spec_{\kk}$ the \emph{generalized discrete Markov spectrum}. Each real number is counted only once in this set, even if it arises from several pairs $(i,n)$. For $\kk=\mathbf0$, this is the set $\Spec_{\mathbf0}$ in Markov's theorem.

Each $\LL_a(n)$ is strictly increasing in $n$ and converges to $K$ from below. If the distinct elements of $\Spec_{\kk}$ are listed as
\[
 s_1<s_2<s_3<\cdots
\]
then $s_r\to K$. We study the series $\sum_{r\ge1}(K-s_r)$, which sums the differences between these values and their limit $K$.

We also consider a sum with multiplicities defined by GM triples. Write $t_i$ for the $i$th entry of $\mathbf t=(x,y,z)\in\mathcal G_{\kk}$ and define
\begin{equation}\label{eq:arithmetic-multiplicity}
 r_i(n)=\begin{cases}
 1,&n=1,\\[3pt]
 \#\{\mathbf t\in\mathcal G_{\kk}\mid t_i=n>t_j\ (j\ne i)\},&n\ge2
 \end{cases}.
\end{equation}
Thus, for $n\ge2$, we count the GM triples whose $i$th entry is the unique largest entry and is equal to $n$. For $n=1$, the initial solution $(1,1,1)$ contributes one term for each $i$. These counts are finite. We denote the multiplicities for $\kk=\mathbf0$ by $r_i^{(\mathbf0)}(n)$.

Gaster and Loustau~\cite{GL} studied sums of the differences between the classical Lagrange constants and their limit $3$, and proved the following identity and inequality. The first identity rewrites equation~(6) of their paper using the multiplicities defined above in terms of Markov triples.

\Needspace{4\baselineskip}
\begin{theorem}[{\cite[(6)]{GL}}]\label{thm:classical-sum}
The following identity holds, with absolute convergence:
\[
 \sum_{i=1}^3\sum_{n\ge1}
 r_i^{(\mathbf0)}(n)
 \left(3-\frac{\sqrt{9n^2-4}}n\right)=3.
\]

\end{theorem}

For the sum in which distinct Lagrange constants are counted only once, the following inequality holds.

\begin{theorem}[{\cite[Theorem~2.2]{GL}}]\label{thm:classical-GL}
Let $L_1<L_2<\cdots$ be the distinct elements of the discrete part $\Spec_{\mathbf0}=\mathcal L\cap(-\infty,3)$ of the Lagrange spectrum. Then $L_r\to3$, and the series of positive terms
\[
 \sum_{r\ge1}(3-L_r)
 =\sum_{s\in\Spec_{\mathbf0}}(3-s)
\]
converges absolutely and satisfies
\[
 \sum_{s\in\Spec_{\mathbf0}}(3-s)
 \le4-\varphi-\sqrt2,\qquad
 \varphi=\frac{1+\sqrt5}{2}.
\]
Equality holds if and only if the Markov uniqueness conjecture~\ref{conj:classical-uniqueness} is true.
\end{theorem}

The main results of this paper are Theorems~\ref{thm:arithmetic-sum} and~\ref{thm:GL} below. Theorem~\ref{thm:arithmetic-sum} extends the identity in Theorem~\ref{thm:classical-sum}, and Theorem~\ref{thm:GL} extends the inequality in Theorem~\ref{thm:classical-GL}, to the generalized Markov equation.

\Needspace{4\baselineskip}
\begin{theorem}\label{thm:arithmetic-sum}
For the multiplicities in \eqref{eq:arithmetic-multiplicity}, the following identity holds, with absolute convergence:
\begin{equation}\label{eq:regrouped}
 \sum_{i=1}^3\sum_{n\ge1}
 r_i(n)\bigl(K-\LL_{k_i}(n)\bigr)=K.
\end{equation}

\end{theorem}

Theorem~\ref{thm:GL} gives an upper bound for the sum in which distinct Lagrange constants are counted only once. When at least two of $k_1,k_2,k_3$ coincide, it improves the upper bound $K$ obtained directly from Theorem~\ref{thm:arithmetic-sum}.

\Needspace{8\baselineskip}
\begin{theorem}\label{thm:GL}
For arbitrary $k_1,k_2,k_3\in\ZZ_{\ge0}$, define the constant
\begin{equation}\label{eq:Cdef}
 C_{\kk}=\begin{cases}
 K,&k_1,k_2,k_3\text{ are pairwise distinct},\\[7pt]
 \displaystyle\frac{3K}{2}-
 \frac{\LL_k(1)+\LL_k(k+2)}2,
 &\text{exactly two }k_i\text{ coincide; the third is }k,\\[7pt]
 \displaystyle\frac{7K}{6}-
 \frac{\LL_k(1)+\LL_k(k+2)}2,
 &k_1=k_2=k_3=k.
 \end{cases}
\end{equation}
The series of positive terms over the set $\Spec_{\kk}$ of distinct Lagrange constants converges absolutely and satisfies the bound
\begin{equation}\label{eq:GL-main}
 \sum_{s\in\Spec_{\kk}}(K-s)\le C_{\kk}.
\end{equation}

We have $C_{\kk}\le K$, with $C_{\kk}<K$ whenever at least two of $k_1,k_2,k_3$ coincide. For $\kk=(k,k,k)$, equality in \eqref{eq:GL-main} holds if and only if the generalized Markov uniqueness conjecture~\ref{conj:GM-uniqueness} holds for this value of $k$. In particular, for $\kk=\mathbf0$, we have $C_{\mathbf0}=4-\varphi-\sqrt2$ and recover Theorem~\ref{thm:classical-GL}.
\end{theorem}

These upper bounds can sometimes be improved further by examining repeated values that are not explained by permutations of equal entries of $\kk$. For example, when $\kk=(0,2,4)$, the first component of the GM triple $(1,1,1)$ and the second component of $(2,9,1)$ give the same spectral value $\LL_0(1)=\LL_2(9)=\sqrt{77}$. Accounting for this repetition improves the upper bound from $C_{(0,2,4)}=9$ to $\sqrt{77}$ (Example~\ref{ex:spectral-collision}).

\Needspace{5\baselineskip}
Define the difference between the upper bound and the actual sum by
\[
 D_{\kk}:=C_{\kk}-\sum_{s\in\Spec_{\kk}}(K-s)\ge0.
\]
This quantity can be viewed as a weighted measure of repetitions of spectral values not explained by permutations of equal entries of $\kk$.

\Needspace{5\baselineskip}
\begin{question}\label{prob:exact-spectrum-sum}
For arbitrary $k_1,k_2,k_3\in\ZZ_{\ge0}$, express the exact value of the sum
\[
 \sum_{s\in\Spec_{\kk}}(K-s),
\]
in which each distinct spectral value is counted once, in terms of $k_1,k_2,k_3$. Equivalently, determine the exact value of $D_{\kk}$.
\end{question}

For $\kk=(k,k,k)$, Conjecture~\ref{conj:GM-uniqueness} implies $D_{\kk}=0$ by Theorem~\ref{thm:GL}. When some of the $k_i$ differ, repetitions of spectral values can arise even from GM triples with different largest entries: distinct GM numbers associated with different $k_i$ may yield the same spectral value. Such repetitions are not explained by permutations of equal parameters and can make $D_{\kk}$ positive. In the example $\kk=(0,2,4)$ above, the largest entries $1$ and $9$ give the same value $\sqrt{77}$, so $D_{(0,2,4)}\ge9-\sqrt{77}>0$.

\subsection{Outline of the proof}
Through work of Cohn~\cite{Cohn71} and others, classical Markov numbers are related to traces of closed curves on a once-punctured torus.
Gaster and Loustau~\cite{GL} prove Theorem~\ref{thm:classical-sum} using McShane's identity, with the sum taken over all simple closed geodesics on the once-punctured torus, together with the correspondence described above between Markov numbers and traces of closed curves. We replace the once-punctured torus considered in the classical case with its quotient by the elliptic involution, a hyperbolic orbifold on a sphere with four marked points. By generalizing the types of singularities and ends on this sphere, we obtain a framework that also applies to GM numbers. We then replace McShane's identity with the Tan--Wong--Zhang identity~\cite{TWZ}, thereby generalizing Theorem~\ref{thm:classical-sum} to Theorem~\ref{thm:arithmetic-sum}. We describe the construction more concretely below.

We first define three initial matrices $L_0,M_0,N_0$ in $\SL_2(\ZZ)$ from $k_1,k_2,k_3$. Taking the quotient of the upper half-plane by the group generated by their projections to $\PSL_2(\RR)$ gives a hyperbolic orbifold (Theorem~\ref{thm:orbifold}).
Let $a_1,a_2,a_3,a_4$ be four marked points corresponding to the singularities and ends of this orbifold. The points corresponding to the three initial matrices are labeled $a_1,a_2,a_3$, respectively, and the remaining cusp is labeled $a_4$.

We consider essential simple ideal arcs based at $a_4$, that is, arcs from $a_4$ to itself whose interiors are embedded, avoid the other marked points, and separate the remaining three marked points into one point and two points. We refer to these simply as \emph{ideal arcs} below. Let $\EE$ be the set of their unoriented isotopy classes. If an ideal arc $e$ encloses the single marked point $a_{i(e)}$, pushing $e$ slightly toward the side containing the other two marked points gives a simple closed curve $\beta_e$ separating $\{a_4,a_{i(e)}\}$ from the other two marked points. This construction gives a bijection between unoriented isotopy classes of ideal arcs and of closed curves enclosing two marked points (Proposition~\ref{prop:slope-bijection}).

We next specify a horocycle of length $K$ at the cusp $a_4$. Using Markov--monodromy triples, we prove that the $\lambda$-lengths of ideal arcs are GM numbers. We also show that the correspondence between Markov--monodromy triples and GM triples intertwines the Hurwitz moves corresponding to exchanges of ideal arcs with the GM exchanges. If $\ell_{\beta_e}$ is the length of the geodesic realization of $\beta_e$ and $k_e=k_{i(e)}$, we have
\[
 \lambda(e)=m_e,\qquad
 2\cosh\frac{\ell_{\beta_e}}{2}=K\lambda(e)-k_e=Km_e-k_e
\]
(Theorem~\ref{thm:lambda-GM}). We apply to this orbifold the extension of McShane's identity (see~\cite{McS}) to cone points and geodesic boundary components established by Tan--Wong--Zhang~\cite{TWZ} (Lemma~\ref{thm:TWZ}), and use the formula above to express $\ell_{\beta_e}$ in terms of GM numbers. Indexing the sum by ideal arcs $e$ through the correspondence above gives
\[
 \sum_{e\in\EE}\bigl(K-\LL_{k_e}(m_e)\bigr)=K
\]
(Theorem~\ref{thm:weighted-sum}). Grouping this sum by $(i,n)$ using the correspondence in Proposition~\ref{cor:multiplicity} yields the summation formula in Theorem~\ref{thm:arithmetic-sum}. Finally, we account for the multiplicities arising from permutations $\sigma$ of the indices satisfying $k_{\sigma(i)}=k_i$ ($i=1,2,3$) (Proposition~\ref{prop:orbit-sum}) and derive the inequality in Theorem~\ref{thm:GL}, namely \eqref{eq:GL-main}, for arbitrary $\kk$.

\subsection{Relation to previous work in hyperbolic geometry}\label{sec:prior-work}
Chekhov--Shapiro~\cite[Theorem~3.6]{CS} gave a method for constructing hyperbolic orbifolds with cone points, cusps, and geodesic boundary components from a graph with a cyclic ordering of the incident edges at each vertex (a fat graph), prescribed cone-point orders, and a real coordinate on each edge. The mutations in the three directions of the generalized cluster algebra associated with the GM equation correspond to the Chekhov--Shapiro flips of arcs on a sphere with four peripheral components, so this construction can be used to obtain a geometric realization~\cite[Lemma~3.1, Remark~4.6]{CS}. GM numbers are obtained by specializing the initial cluster variables to $(1,1,1)$~\cite{GM}. To realize these numerical values as $\lambda$-lengths, one must specify a hyperbolic metric and horocycles for which the three initial arcs all have $\lambda$-length $1$.

Our construction, on the other hand, differs from the Chekhov--Shapiro method in that it starts from the GM equation and its initial solution $(1,1,1)$ and constructs a hyperbolic orbifold using the corresponding triple of integral monodromy matrices. The three matrices, given in a common coordinate system, act as isometries of the upper half-plane preserving its standard hyperbolic metric, and hence determine the hyperbolic metric on the quotient at the same time. We also specify $\{\operatorname{Im}z=1\}$ in the upper half-plane, thereby assigning a horocycle of length $K$ to the cusp. The corresponding three initial ideal arcs then all have $\lambda$-length $1$. Under this normalization, we identify the $\lambda$-lengths of all ideal arcs with GM numbers (Theorem~\ref{thm:lambda-GM}). Our realization also falls within the scope of the general Chekhov--Shapiro construction, but we adopt this method because explicit matrix calculations and side identifications of a fundamental domain provide a more elementary description.

\subsection*{Acknowledgements}
The author thanks Tsukasa Ishibashi and Shunsuke Kano for helpful discussions and comments.

This work was supported by JSPS KAKENHI Grant Number~JP25K17224.

\subsection*{AI disclosure}
The author used ChatGPT 6.1 sol (OpenAI) to search the literature, assist with proofs, and organize the manuscript.

\section{Markov--monodromy triples and Hurwitz moves}\label{sec:monodromy}
In the GM equation \eqref{eq:GM}, fix $k_1,k_2,k_3\in\ZZ_{\ge0}$ and put
\[
 K=3+k_1+k_2+k_3.
\]

\subsection{Generating GM triples by Vieta involutions}\label{sec:Vieta}
We first describe how the positive integer solutions of the GM equation are generated. For $(x,y,z)\in\mathcal G_{\kk}$, fixing $y,z$ and viewing \eqref{eq:GM} as a quadratic equation in $x$ gives
\[
 x^2-(Kyz-k_2z-k_3y)x+(y^2+k_1yz+z^2)=0.
\]
The operation replacing one root $x$ by the other root $x'$ is a \emph{Vieta jump}. For a GM solution, the values obtained by exchanging each of the three entries are
\begin{equation}\label{eq:Vieta-jumps}
 \begin{aligned}
 x'&=Kyz-k_2z-k_3y-x=\frac{y^2+k_1yz+z^2}{x},\\[3pt]
 y'&=Kxz-k_1z-k_3x-y=\frac{x^2+k_2xz+z^2}{y},\\[3pt]
 z'&=Kxy-k_1y-k_2x-z=\frac{x^2+k_3xy+y^2}{z}.
 \end{aligned}
\end{equation}
We therefore define
\[
 \mu_1(x,y,z)=(x',y,z),\qquad
 \mu_2(x,y,z)=(x,y',z),\qquad
 \mu_3(x,y,z)=(x,y,z').
\]
Integrality follows from the polynomial expressions in \eqref{eq:Vieta-jumps}, and positivity follows from the fractional expressions. Since the other root satisfies the same quadratic equation, each $\mu_r$ preserves $\mathcal G_{\kk}$ and satisfies $\mu_r^2=\mathrm{id}$. We call these maps the Vieta involutions.

\begin{theorem}[{\cite[Theorem~1]{GM}}]\label{thm:GM-enumeration}
Every positive integer solution of the GM equation is obtained from $(1,1,1)$ by repeatedly applying the three Vieta jumps $\mu_1,\mu_2,\mu_3$.

Every solution other than $(1,1,1)$ has a unique largest entry, and exchanging this entry strictly decreases the largest value. Taking the resulting solution as its parent, the positive integer solutions form a trivalent tree rooted at $(1,1,1)$, in which every ordered GM triple occurs exactly once.
\end{theorem}

Consequently, repeatedly descending by exchanging the largest entry always reaches $(1,1,1)$. Reversing this descent enumerates all solutions. The three children of the root are
\[
 (k_1+2,1,1),\qquad(1,k_2+2,1),\qquad(1,1,k_3+2).
\]
Thereafter, keeping the largest entry and exchanging each of the other two entries gives the two children. We denote this tree of GM triples by $\mathcal T_{\mathrm{GM}}$ and label an edge by the index $r$ of the exchanged entry. The sequence of edge labels from a solution to the root is uniquely recovered by repeatedly applying the Vieta involution at the largest entry.

\begin{figure}[!htbp]
\centering
\begin{tikzpicture}[x=1cm,y=.75cm,every node/.style={tree example vertex}]
 \node (treeRoot) at (0,0) {$(1,1,1)$};
 \node (treeOne) at (4,2.8) {$(\mathbf2,1,1)$};
 \node (treeTwo) at (4,0) {$(1,\mathbf3,1)$};
 \node (treeThree) at (4,-2.8) {$(1,1,\mathbf4)$};
 \node (tree12) at (8,3.5) {$(2,\mathbf7,1)$};
 \node (tree13) at (8,2.1) {$(2,1,\mathbf9)$};
 \node (tree21) at (8,0.7) {$(\mathbf{10},3,1)$};
 \node (tree23) at (8,-0.7) {$(1,3,\mathbf{16})$};
 \node (tree31) at (8,-2.1) {$(\mathbf{17},1,4)$};
 \node (tree32) at (8,-3.5) {$(1,\mathbf{21},4)$};
 \treeexampleedges
\end{tikzpicture}
\caption{Tree of GM triples ($\kk=(0,1,2)$)}
\label{fig:GM-tree}
\end{figure}
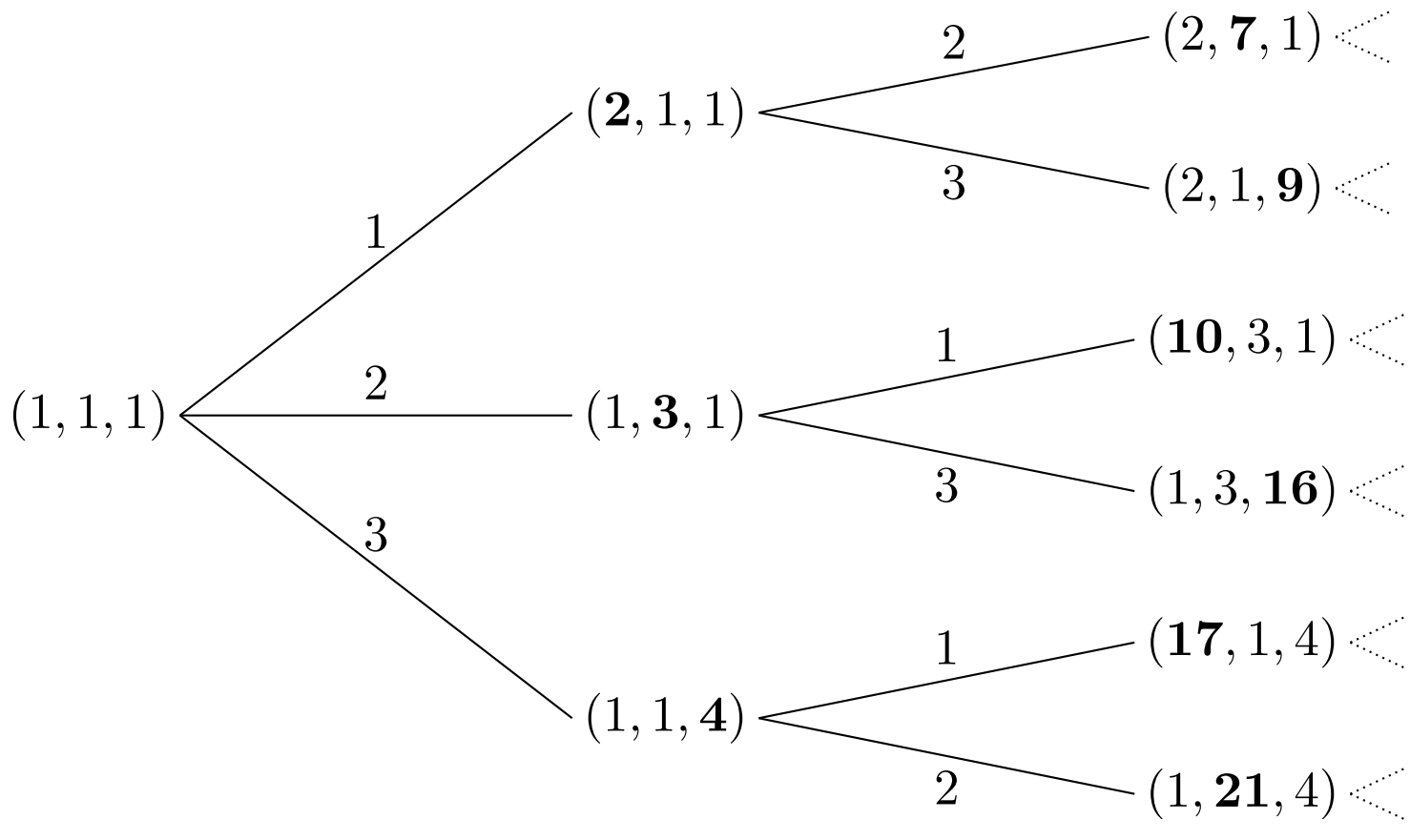

For example, if $\kk=(0,1,2)$, then $K=6$ and the equation is
\[
 x^2+y^2+z^2+zx+2xy=6xyz.
\]
Figure~\ref{fig:GM-tree} shows the triples obtained by at most two exchanges from the root. The exchanged entry is shown in bold.

We now realize these three exchanges and the enumeration tree by matrices.

\subsection{The initial triple and its monodromy representation}
Take the following three matrices as the initial triple:
\begin{equation}\label{eq:seed}
 \begin{aligned}
 L_0&=\begin{pmatrix}0&-1\\1&-k_1\end{pmatrix},\\[4pt]
 M_0&=\begin{pmatrix}
 k_1+1&-(k_1+1)(k_1+k_2+1)-1\\
 1&-(k_1+k_2+1)
 \end{pmatrix},\\[4pt]
 N_0&=\begin{pmatrix}
 k_1+k_2+2&-(k_1+k_2+2)(K-1)-1\\
 1&-(K-1)
 \end{pmatrix}.
 \end{aligned}
\end{equation}
The $(2,1)$ entries of all three matrices are $1$. Direct calculation also gives
\[
 L_0M_0N_0=\begin{pmatrix}-1&K\\0&-1\end{pmatrix}.
\]
Put
\[
 B:=-(L_0M_0N_0)^{-1}=\begin{pmatrix}1&K\\0&1\end{pmatrix}.
\]
We then have
\begin{equation}\label{eq:seed-data}
 \tr L_0=-k_1,\quad\tr M_0=-k_2,\quad\tr N_0=-k_3,
 \qquad L_0M_0N_0B=-I.
\end{equation}

\begin{figure}[!htbp]
\centering
\usetikzlibrary{decorations.markings}
\begin{tikzpicture}[x=2.8cm,y=2.8cm,
 every node/.style={font=\small,inner sep=1pt},
 >=stealth,
 based return arc/.style={black,line width=.8pt,
  postaction={decorate},decoration={markings,
   mark=at position .50 with {\arrow[scale=.85]{>}}}}]
\shade[ball color=gray!10,opacity=.24] (0,0) circle[radius=1];
\draw[black!60,line width=.55pt] (0,0) circle[radius=1];

\coordinate (basedAfour) at (0,.70);
\coordinate (basedB) at (0,.48);
\coordinate (basedAone) at (-.55,-.23);
\coordinate (basedAtwo) at (0,-.52);
\coordinate (basedAthree) at (.55,-.23);

\draw[black!40,densely dashed,line width=.5pt] (0,.67)
 ellipse[x radius=.35,y radius=.19];
\node[text=black!65] at (-.43,.77) {$\Delta_4$};
\node[text=black!65] at (.45,.76) {$\partial\Delta$};

\draw[based return arc] (basedB)
 .. controls (-.26,.37) and (-.76,.08) .. (-.76,-.23)
 .. controls (-.76,-.48) and (-.34,-.48) .. (-.34,-.23)
 .. controls (-.34,.10) and (-.16,.33) .. (basedB);
\draw[based return arc] (basedB)
 .. controls (-.14,.28) and (-.13,-.19) .. (-.17,-.52)
 .. controls (-.17,-.77) and (.17,-.77) .. (.17,-.52)
 .. controls (.13,-.19) and (.14,.28) .. (basedB);
\draw[based return arc] (basedB)
 .. controls (.16,.33) and (.34,.10) .. (.34,-.23)
 .. controls (.34,-.48) and (.76,-.48) .. (.76,-.23)
 .. controls (.76,.08) and (.26,.37) .. (basedB);
\foreach \basedPoint in {basedAone,basedAtwo,basedAthree,basedAfour} {
 \fill (\basedPoint) circle[radius=1.5pt];
}
\fill (basedB) circle[radius=1.6pt];
\node[above=4pt] at (basedB) {$b$};
\node[above=4pt] at (basedAfour) {$a_4$};
\node at (-.81,.11) {$\gamma_1$};
\node at (0,-.05) {$\gamma_2$};
\node at (.81,.11) {$\gamma_3$};
\node at (-.55,-.50) {$a_1$};
\node at (0,-.80) {$a_2$};
\node at (.55,-.50) {$a_3$};
\end{tikzpicture}
\caption{Peripheral loops}
\label{fig:based-peripheral-loops}
\end{figure}
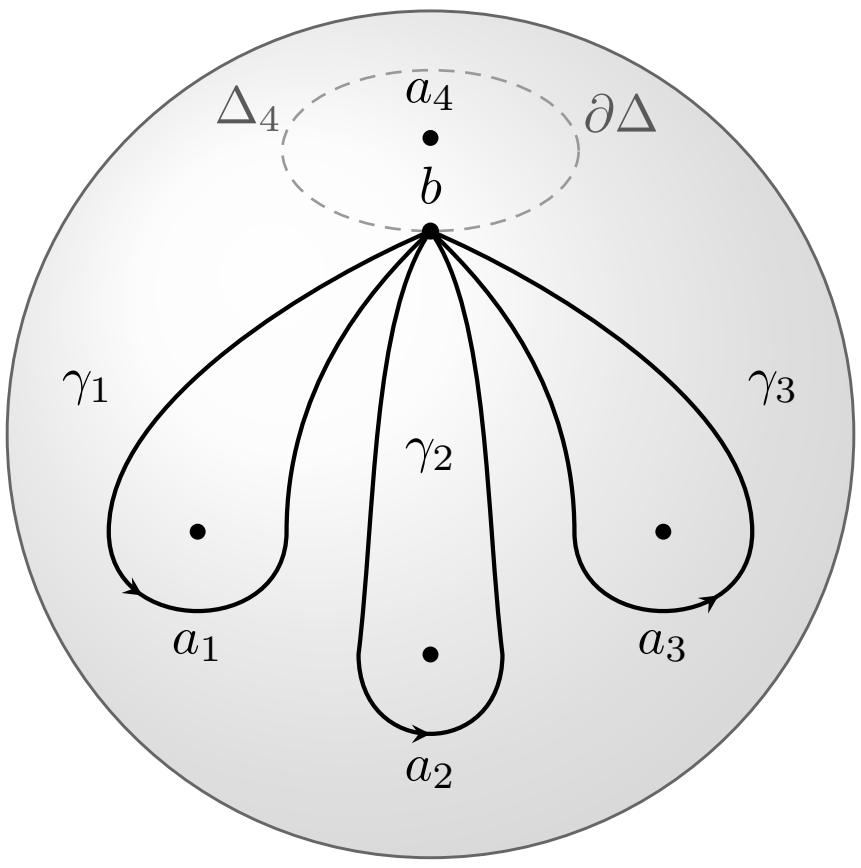

We interpret these matrices as a representation of the fundamental group of a sphere with four points removed. Let $\Sigma_{0,4}$ be the oriented surface obtained by removing four marked points $a_1,a_2,a_3,a_4$ from a sphere. Choose a small disk neighborhood $\Delta_4$ of $a_4$ and set $\Delta=S^2\setminus\operatorname{int}\Delta_4$ and $b\in\partial\Delta$.
Choose three connecting paths from $b$ to small circles around $a_r$ ($r=1,2,3$), with pairwise disjoint interiors. Let $\gamma_r$ be the based peripheral loop that follows the corresponding path, goes once around $a_r$ in the positive direction, and returns along the same path (Figure~\ref{fig:based-peripheral-loops}). Order the paths so that $\gamma_1\gamma_2\gamma_3$ travels around $\partial\Delta$ in the positive direction. Its inverse $\gamma_4$ goes once around $a_4$ in the positive direction. These loops give the presentation
\[
 \pi_1(\Sigma_{0,4},b)=
 \langle\gamma_1,\gamma_2,\gamma_3,\gamma_4\mid
 \gamma_1\gamma_2\gamma_3\gamma_4=1\rangle.
\]

By the product relation above, the assignments $\rho_0(\gamma_1)=L_0$, $\rho_0(\gamma_2)=M_0$, $\rho_0(\gamma_3)=N_0$, and $\rho_0(\gamma_4)=-B$ define an integral representation. Each matrix is labeled by the index $r$ of the marked point $a_r$ encircled by its corresponding loop $\gamma_r$.

For a triple $(L,M,N)\in\SL_2(\ZZ)^3$ with labels $(i,j,\ell)$ given by a permutation of $(1,2,3)$, define the Hurwitz moves
\begin{equation}\label{eq:Hurwitz}
 \begin{aligned}
 U(L,M,N)&=(L,MNM^{-1},M),&\text{labels }&(i,\ell,j),\\
 V(L,M,N)&=(M,M^{-1}LM,N),&\text{labels }&(j,i,\ell).
 \end{aligned}
\end{equation}
The inverse operations are
\[
 U^{-1}(L,M,N)=(L,N,N^{-1}MN),\qquad
 V^{-1}(L,M,N)=(LML^{-1},L,N).
\]
All these operations preserve the product $LMN$. Each new matrix is a product of the original matrices and their inverses, and the same holds for the inverse operations. Thus the generated subgroup is unchanged by these operations.

\begin{definition}\label{def:MM}
A triple of matrices in $\SL_2(\ZZ)$, equipped with labels $(i,j,\ell)$ in which each element of $\{1,2,3\}$ occurs once, is called a labeled integral monodromy triple if it satisfies
\begin{equation}\label{eq:MM}
 LMNB=-I,\qquad(\tr L,\tr M,\tr N)=(-k_i,-k_j,-k_\ell).
\end{equation}

In this paper, a \emph{Markov--monodromy triple} is a triple obtained from \eqref{eq:seed} by the operations $U,V$ in \eqref{eq:Hurwitz} and their inverses.
\end{definition}

\begin{remark}\label{rem:BG}
Markov--monodromy triples for $k_1=k_2=k_3=k$ were introduced by the author, Maruyama, and Sato~\cite[Definitions~5.1--5.2]{GMS}. Banaian and the author~\cite[Definitions~5.1--5.2, Proposition~5.3]{BG} extended this construction to independently specified parameters $k_1,k_2,k_3$ and to Laurent polynomial entries.
Let $P=\left(\begin{smallmatrix}0&1\\1&0\end{smallmatrix}\right)$ and denote the initial variables in the latter paper by $\xi_1,\xi_2,\xi_3$. Specializing the initial triple to $\xi_1=\xi_2=\xi_3=1$ and the free parameter to $f=-k_1$ gives
\[
 (X_{\mathrm{id}}(f),Y_{\mathrm{id}}(f),Z_{\mathrm{id}}(f))
 =(PL_0P,PM_0P,PN_0P)
\]
The $(1,2)$ entries in that paper correspond to the $(2,1)$ entries used here, and the fixed product agrees with $XYZ=P(-B^{-1})P$.
\end{remark}

\subsection{Half-twists and Hurwitz moves in three directions}
\begin{figure}[!htbp]
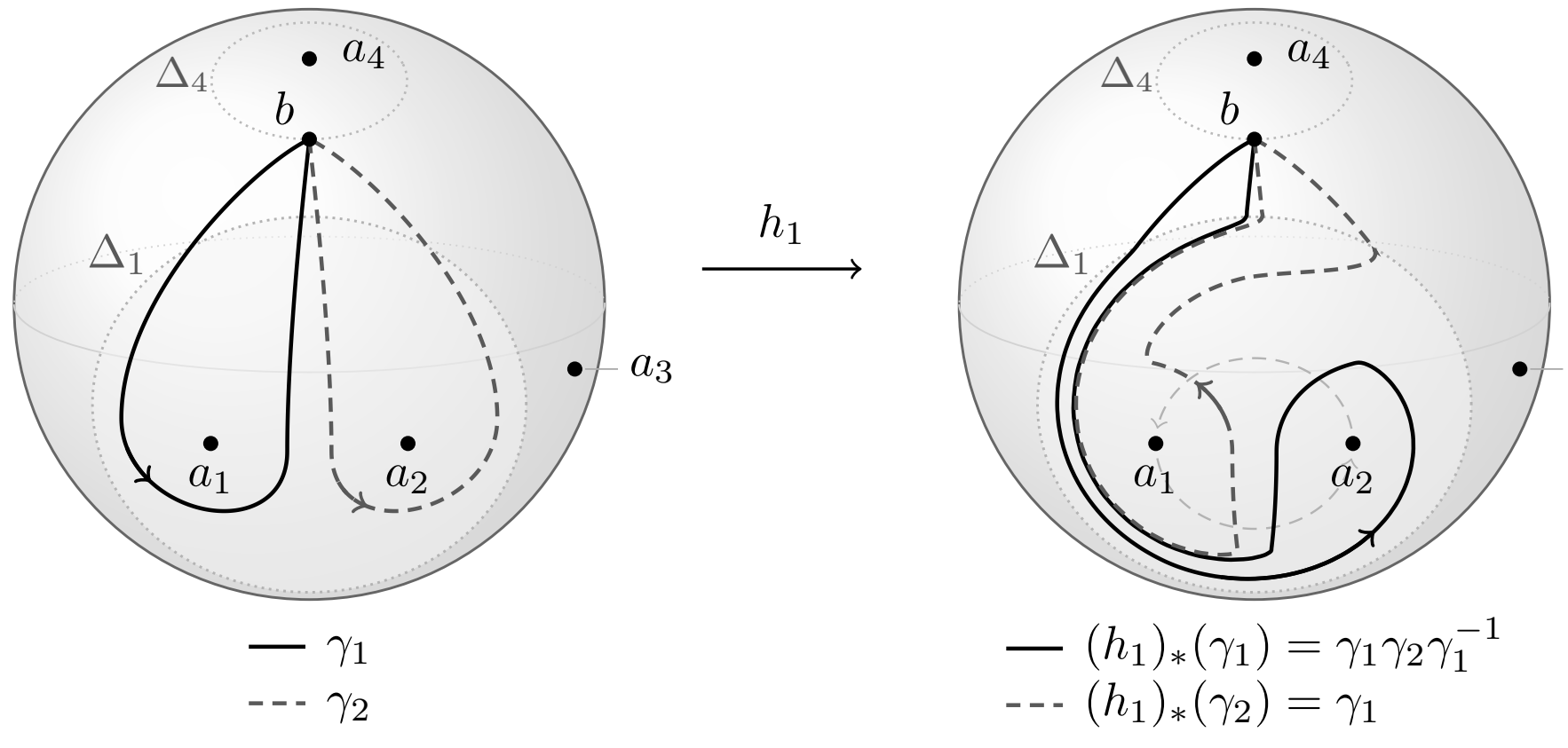

\centering

\caption{The half-twist $h_1$ on the sphere}
\label{fig:half-twist-loops}
\end{figure}

We realize the matrix operations $U,V$ as transformations of peripheral loops on the sphere (Figure~\ref{fig:half-twist-loops}). For $r=1,2$, choose an arc in $\Delta$ joining the adjacent marked points $a_r,a_{r+1}$, and let $\Delta_r$ be a disk neighborhood of this arc. Choose $\Delta_r$ so that it contains neither the other marked points nor the basepoint $b$ and is disjoint from the connecting paths used to base the other peripheral loops. Take an orientation-preserving homeomorphism that exchanges the two marked points by a half rotation in $\Delta_r$ and is the identity near its boundary and outside it. Its class under isotopies preserving the set of marked points and fixing $\partial\Delta$ pointwise is called a \emph{half-twist}.
Of the two possible directions, let $h_r$ be the half-twist that sends $\gamma_{r+1}$ to $\gamma_r$.

\begin{lemma}[{see~\cite[Theorem~15, (14)--(15)]{Artin47}}]\label{lem:Artin-action}
The actions induced by $h_r$ and its inverse on the based fundamental group are
\begin{equation}\label{eq:half-twist-loops}
 \begin{aligned}
 (h_r)_*:\ (\gamma_r,\gamma_{r+1})
 &\longmapsto(\gamma_r\gamma_{r+1}\gamma_r^{-1},\gamma_r),\\
 (h_r^{-1})_*:\ (\gamma_r,\gamma_{r+1})
 &\longmapsto(\gamma_{r+1},\gamma_{r+1}^{-1}\gamma_r\gamma_{r+1})
 \end{aligned}
\end{equation}
and they fix the other peripheral loops. In particular, $\gamma_4$ is fixed as an element of the based fundamental group.
\end{lemma}

Let $\rho(\gamma_1)=L$, $\rho(\gamma_2)=M$, $\rho(\gamma_3)=N$, and $\rho(\gamma_4)=-B$. Moving the peripheral loops by the half-twists above gives the representations defined by precomposition:
\[
 \rho_U=\rho\circ(h_2)_*,\qquad
 \rho_V=\rho\circ(h_1^{-1})_*.
\]
By Lemma~\ref{lem:Artin-action}, their values on $\gamma_1,\gamma_2,\gamma_3$ are $U(L,M,N)$ and $V(L,M,N)$ in \eqref{eq:Hurwitz}, respectively, and both take the value $-B$ on $\gamma_4$. A conjugated matrix retains the label of its original marked point, so the new labels are $(i,\ell,j)$ for $U$ and $(j,i,\ell)$ for $V$. The same construction applies to any system of peripheral loops by moving the connecting paths and the arc joining the two marked points simultaneously.

The Vieta involutions are indexed by the labels of the entries they exchange. The formulas for $U$ and $V$ use the positions of the matrices in the displayed triple. To exchange the entry with a given label, we first move its matrix to the first position by a cyclic change that preserves the product $LMN$, then apply $V$. For $T=(L,M,N)$, put
\begin{equation}\label{eq:cyclic-Hurwitz}
 R(T)=(M,N,BLB^{-1}),\qquad
 (i,j,\ell)\longmapsto(j,\ell,i).
\end{equation}
The product relation $LMN=-B^{-1}$ gives
\[
 R=U^{-1}\circ V,\qquad
 R^3(T)=(BLB^{-1},BMB^{-1},BNB^{-1}).
\]
Thus $R$ also preserves the generated group. Moreover, since $(BAB^{-1})_{21}=A_{21}$, the three $(2,1)$ entries read according to their labels remain unchanged.

\begin{definition}\label{def:three-hurwitz}
Identify $T$ with $R^n(T)$ for $n\in\ZZ$, and write $[T]$ for the resulting equivalence class, called a \emph{cyclic class}. The matrices and their marked-point labels are moved together under cyclic changes.

For $r\in\{1,2,3\}$, choose a representative $(A,D,E)$ whose first matrix has label $r$, and define
\[
 \mathcal H_r[A,D,E]=[D,D^{-1}AD,E].
\]
This is the Hurwitz move that exchanges the entry with label $r$. If the displayed labels are $(i,j,\ell)$, the three moves are explicitly
\begin{equation}\label{eq:three-Hurwitz}
 \begin{aligned}
 \mathcal H_i[L,M,N]&=[M,M^{-1}LM,N],\\
 \mathcal H_j[L,M,N]&=[L,N,N^{-1}MN],\\
 \mathcal H_\ell[L,M,N]&=[L,MNM^{-1},M].
 \end{aligned}
\end{equation}
The displayed labels on the right-hand sides are $(j,i,\ell)$, $(i,\ell,j)$, and $(i,\ell,j)$, respectively.
\end{definition}

For each cyclic class $[T]$, let $\Phi[T]$ be the integer triple obtained by listing the $(2,1)$ entries in the order of the marked-point labels $1,2,3$, independently of the displayed order. As verified after \eqref{eq:cyclic-Hurwitz}, this value is independent of the representative. We use this map to relate the Vieta involutions $\mu_r$ of Section~\ref{sec:Vieta} to the Hurwitz moves $\mathcal H_r$.

\begin{proposition}\label{prop:Hurwitz}
Each $\mathcal H_r$ is independent of the representative and satisfies $\mathcal H_r^2=\mathrm{id}$.

The map $\Phi$ from cyclic classes of Markov--monodromy triples is surjective onto $\mathcal G_{\kk}$ and satisfies
\begin{equation}\label{eq:c-actions}
 \Phi\bigl(\mathcal H_r[T]\bigr)=\mu_r\bigl(\Phi[T]\bigr),
 \qquad r=1,2,3.
\end{equation}
Thus the Hurwitz move with label $r$ realizes the Vieta involution in the $r$th position on the $(2,1)$ entries.
\end{proposition}
\begin{proof}
Two representatives whose first label is $r$ differ by $R^{3n}$, that is, by simultaneous conjugation of the three matrices by $B^n$. Since $V$ commutes with this simultaneous conjugation, $\mathcal H_r$ is independent of the representative. If the first label of $T$ is $r$, then $r$ is the second label of $V(T)$. We therefore choose $R V(T)$ as the next representative whose first label is $r$. The product relation gives
\[
 VR V(T)=R^2(T),\qquad
 VR=R U^{-1},\qquad
 VR^2(T)=B\,U(T)\,B^{-1}.
\]
The first identity gives $\mathcal H_r^2[T]=[T]$. The right-hand side of the third identity means that all three entries of $U(T)$ are conjugated by $B$. These identities also give \eqref{eq:three-Hurwitz}.
By \eqref{eq:three-Hurwitz}, every cyclic class of Markov--monodromy triples is obtained from the cyclic class of the initial triple by repeatedly applying the three moves $\mathcal H_r$.

Let the displayed labels be $(i,j,\ell)$ and put $(x,y,z)=(L_{21},M_{21},N_{21})$. The product relation gives $(LM)_{21}=z$ and $(MN)_{21}=x$. Using the identities $(A^{-1})_{21}=-A_{21}$ and $(AD+DA)_{21}=(\tr A)D_{21}+(\tr D)A_{21}$ for $2\times2$ matrices of determinant $1$, we obtain
\[
 \begin{aligned}
 (MNM^{-1})_{21}
 &=(MBL)_{21}=Kxy+(ML)_{21}
 =Kxy-k_i y-k_j x-z,\\
 (M^{-1}LM)_{21}
 &=(NBM)_{21}=Kyz+(NM)_{21}
 =Kyz-k_j z-k_\ell y-x.
 \end{aligned}
\]
Applying the second identity with label $r$ placed first shows that $\mathcal H_r$ realizes the Vieta involution in the $r$th position. Since the initial triple maps to $(1,1,1)$, induction gives $\Phi[T]\in\mathcal G_{\kk}$ and \eqref{eq:c-actions}. Surjectivity follows by following the exchange sequences of Theorem~\ref{thm:GM-enumeration} with the same labels.
\end{proof}

\afterpage{\begin{landscape}
\thispagestyle{empty}
\begin{figure}[p]
\centering
\begingroup
\setlength{\arraycolsep}{3pt}
\resizebox{0.99\linewidth}{!}{\begin{tikzpicture}[x=1cm,y=1cm,every node/.style={tree example vertex,font=\normalsize}]
 \node (treeRoot) at (0,0) {$\left[\begin{pmatrix}0&-1\\1&0\end{pmatrix},\;\begin{pmatrix}1&-3\\1&-2\end{pmatrix},\;\begin{pmatrix}3&-16\\1&-5\end{pmatrix}\right]$\\[2pt]
  labels $(1,2,3)$, $\Phi=(1,1,1)$};
 \node (treeOne) at (8.1,4.6) {$\left[\begin{pmatrix}1&-3\\1&-2\end{pmatrix},\;\begin{pmatrix}5&-13\\2&-5\end{pmatrix},\;\begin{pmatrix}3&-16\\1&-5\end{pmatrix}\right]$\\[2pt]
  labels $(2,1,3)$, $\Phi=(2,1,1)$};
 \node (treeTwo) at (8.1,0) {$\left[\begin{pmatrix}0&-1\\1&0\end{pmatrix},\;\begin{pmatrix}3&-16\\1&-5\end{pmatrix},\;\begin{pmatrix}16&-91\\3&-17\end{pmatrix}\right]$\\[2pt]
  labels $(1,3,2)$, $\Phi=(1,3,1)$};
 \node (treeThree) at (8.1,-4.6) {$\left[\begin{pmatrix}0&-1\\1&0\end{pmatrix},\;\begin{pmatrix}1&-1\\4&-3\end{pmatrix},\;\begin{pmatrix}1&-3\\1&-2\end{pmatrix}\right]$\\[2pt]
  labels $(1,3,2)$, $\Phi=(1,1,4)$};
 \node (tree12) at (16.2,5.75) {$\left[\begin{pmatrix}5&-13\\2&-5\end{pmatrix},\;\begin{pmatrix}18&-49\\7&-19\end{pmatrix},\;\begin{pmatrix}3&-16\\1&-5\end{pmatrix}\right]$\\[2pt]
  labels $(1,2,3)$, $\Phi=(2,7,1)$};
 \node (tree13) at (16.2,3.45) {$\left[\begin{pmatrix}1&-3\\1&-2\end{pmatrix},\;\begin{pmatrix}20&-49\\9&-22\end{pmatrix},\;\begin{pmatrix}5&-13\\2&-5\end{pmatrix}\right]$\\[2pt]
  labels $(2,3,1)$, $\Phi=(2,1,9)$};
 \node (tree21) at (16.2,1.15) {$\left[\begin{pmatrix}3&-16\\1&-5\end{pmatrix},\;\begin{pmatrix}53&-281\\10&-53\end{pmatrix},\;\begin{pmatrix}16&-91\\3&-17\end{pmatrix}\right]$\\[2pt]
  labels $(3,1,2)$, $\Phi=(10,3,1)$};
 \node (tree23) at (16.2,-1.15) {$\left[\begin{pmatrix}0&-1\\1&0\end{pmatrix},\;\begin{pmatrix}16&-91\\3&-17\end{pmatrix},\;\begin{pmatrix}91&-529\\16&-93\end{pmatrix}\right]$\\[2pt]
  labels $(1,2,3)$, $\Phi=(1,3,16)$};
 \node (tree31) at (16.2,-3.45) {$\left[\begin{pmatrix}1&-1\\4&-3\end{pmatrix},\;\begin{pmatrix}13&-10\\17&-13\end{pmatrix},\;\begin{pmatrix}1&-3\\1&-2\end{pmatrix}\right]$\\[2pt]
  labels $(3,1,2)$, $\Phi=(17,1,4)$};
 \node (tree32) at (16.2,-5.75) {$\left[\begin{pmatrix}0&-1\\1&0\end{pmatrix},\;\begin{pmatrix}4&-1\\21&-5\end{pmatrix},\;\begin{pmatrix}1&-1\\4&-3\end{pmatrix}\right]$\\[2pt]
  labels $(1,2,3)$, $\Phi=(1,21,4)$};
 \treeexampleedges
\end{tikzpicture}
}
\endgroup
\caption{Tree of Hurwitz moves ($\kk=(0,1,2)$)}
\label{fig:Hurwitz-tree}
\end{figure}
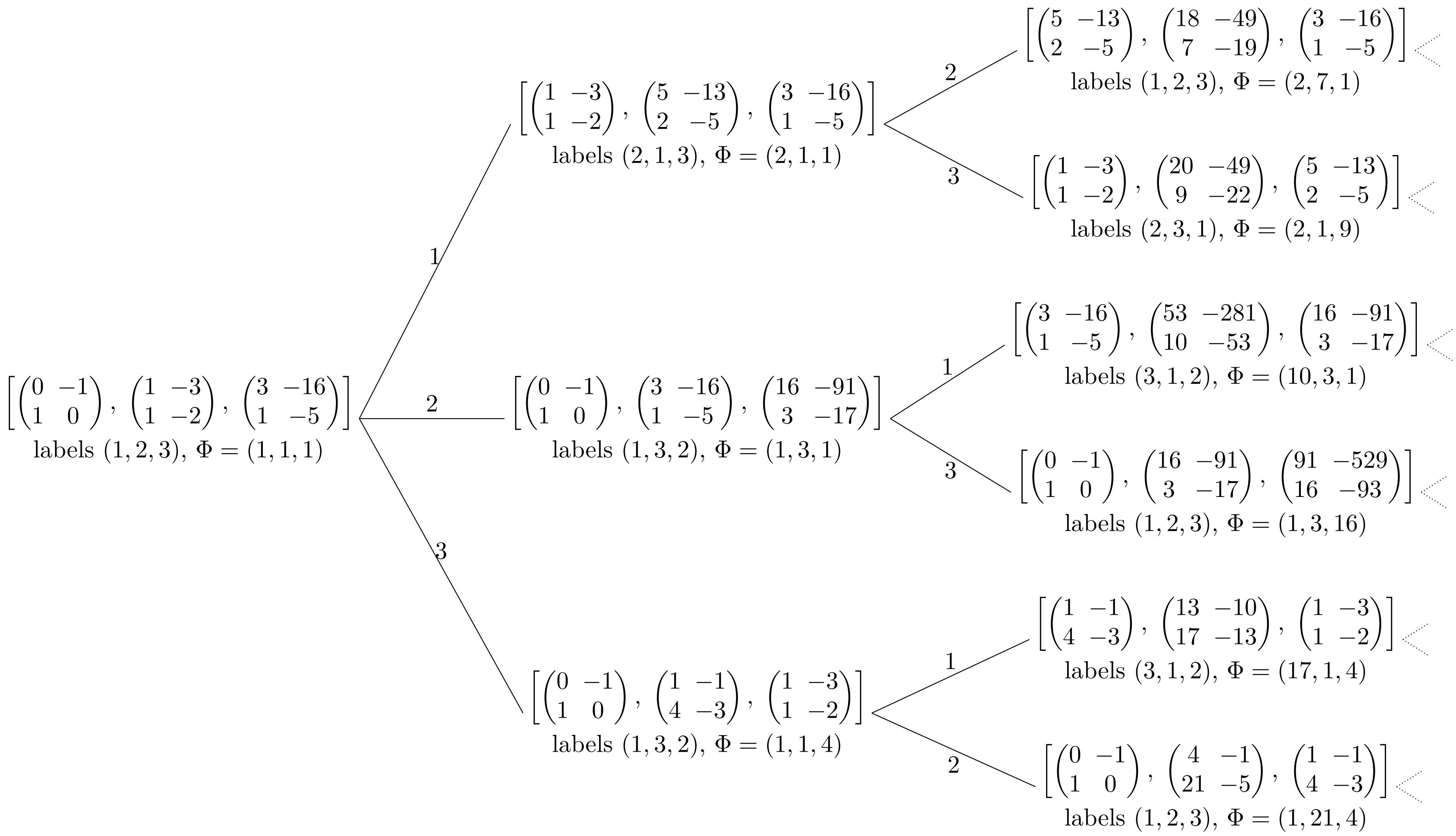
\end{landscape}
}
\subsection{The isomorphism between the tree of Hurwitz moves and the GM tree}\label{sec:Hurwitz-tree}
Let $T_0=(L_0,M_0,N_0)$. Define $\mathcal T_{\mathrm H}$ to be the graph whose vertices are the cyclic classes obtained from $[T_0]$ by repeatedly applying $\mathcal H_1,\mathcal H_2,\mathcal H_3$, with an edge labeled $r$ joining $[T]$ to $\mathcal H_r[T]$. Its root is $[T_0]$.

\begin{proposition}\label{prop:Hurwitz-tree}
The map $\Phi:\mathcal T_{\mathrm H}\longrightarrow\mathcal T_{\mathrm{GM}}$ is an isomorphism preserving the root and the edge labels. Thus cyclic classes of Markov--monodromy triples are in bijection with ordered GM triples, and $\mathcal H_r$ corresponds to the Vieta involution $\mu_r$ in the $r$th position. In particular, $\mathcal T_{\mathrm H}$ is a tree with three edges incident to every vertex.
\end{proposition}
\Needspace{9\baselineskip}
\begin{proof}
By Proposition~\ref{prop:Hurwitz}, $\Phi$ intertwines each Hurwitz move with the corresponding Vieta involution. Consecutive moves with the same label can be canceled because they are involutions. The remaining sequence maps to a nonbacktracking path in the tree $\mathcal T_{\mathrm{GM}}$ of Theorem~\ref{thm:GM-enumeration}. Since the sequence of edge labels from the root to each GM triple is unique, $\Phi$ is injective. Surjectivity follows from Proposition~\ref{prop:Hurwitz}. Thus $\Phi$ gives an isomorphism $\mathcal T_{\mathrm H}\cong\mathcal T_{\mathrm{GM}}$ preserving the edge labels, and $\mathcal T_{\mathrm H}$ is also a trivalent tree.

\end{proof}

A sequence of edge labels $r_1\cdots r_s$ from the root is represented by a sequence with distinct consecutive labels. Two consecutive moves with the same label cancel by $\mathcal H_r^2=\mathrm{id}$.

For $\kk=(0,1,2)$, as in Figure~\ref{fig:GM-tree}, the initial matrices are
\[
 L_0=\begin{pmatrix}0&-1\\1&0\end{pmatrix},\quad
 M_0=\begin{pmatrix}1&-3\\1&-2\end{pmatrix},\quad
 N_0=\begin{pmatrix}3&-16\\1&-5\end{pmatrix},\quad
 B=\begin{pmatrix}1&6\\0&1\end{pmatrix}.
\]
Figure~\ref{fig:Hurwitz-tree} shows the first two levels below the root. At each vertex we display a representative of the cyclic class, with its ordered labels and the value of $\Phi$ below. For example, following the edge labels $1,2$ from the root gives the value $(2,7,1)$ for $\Phi$. In Section~\ref{sec:farey-tree}, we relate this trivalent tree to the Farey tree.

\section{Construction of the hyperbolic orbifold}\label{sec:quotient}\label{sec:euler}
We use an overline for the projection of a matrix to $\PSL_2(\RR)=\SL_2(\RR)/\{\pm I\}$, and set
\begin{equation}\label{eq:quotient}
 \Lambda_{\kk}=\langle\overline L_0,\overline M_0,\overline N_0\rangle,
 \qquad\OO_{\kk}=\Lambda_{\kk}\backslash\HH.
\end{equation}
We equip $\HH=\{x+iy\mid y>0\}$ with the metric $ds^2=(dx^2+dy^2)/y^2$.
The group $\PSL_2(\RR)$ acts by orientation-preserving isometries. Geodesics are vertical lines and semicircles orthogonal to the real axis, and hyperbolic angles agree with Euclidean angles (see~\cite[Theorems~1.1.2, 1.2.1 and~1.3.2]{Kat}).

The Markov--monodromy triples in Definition~\ref{def:MM} are obtained from the initial triple for fixed $\kk$ by $U,V$ and their inverses. By the invariance of the generated subgroup established in Section~\ref{sec:monodromy}, their projections all generate the same group $\Lambda_{\kk}$.

Write $d_{\HH}$ for the distance induced by this metric, and define the \emph{translation length} of $g\in\PSL_2(\RR)$ by $\ell(g)=\inf_{z\in\HH}d_{\HH}(z,g(z))$. The projection $\overline A$ of $A\in\SL_2(\RR)$ with $|\tr A|>2$ is called hyperbolic. The geodesic joining its two fixed points on $\mathbb R\cup\{\infty\}$ is called its \emph{axis}.
The hyperbolic distance formula and the normal form of a hyperbolic element give the following relation.

\begin{theorem}[{see~\cite{Kat}}]\label{thm:trace-length}
If $A\in\SL_2(\RR)$ satisfies $|\tr A|>2$, then
\[
 2\cosh\frac{\ell(\overline A)}2=|\tr A|,
 \quad\text{or equivalently,}\quad
 \ell(\overline A)=2\arcosh\left(\frac{|\tr A|}{2}\right).
\]
The infimum defining the translation length is attained precisely on the axis of $\overline A$.
\end{theorem}

Let $P\subset\HH$ be an open convex hyperbolic polygon with finitely many sides. Its sides are the geodesic segments, rays, or lines forming its boundary in $\HH$; intervals on the real axis are not sides. We also allow a straight side to be divided into two sides at an interior point, creating a vertex of interior angle $\pi$. Vertices in $\HH$ are called finite vertices, and endpoints of sides on $\mathbb R\cup\{\infty\}$ are called ideal vertices. We assume that two sides meet at each finite vertex and that one or two sides are incident to each ideal vertex.

Assign to each side $s$ a distinct side $s^*$ and an element $g_s\in\PSL_2(\RR)$ satisfying
\[
 (s^*)^*=s,\qquad g_s(s)=s^*,\qquad
 g_{s^*}=g_s^{-1},\qquad g_s(P)\cap P=\varnothing.
\]
We glue $s$ to $s^*$ by the map $g_s$.

Choose a finite vertex $v_1$ and an incident side $s_1$. The point $v_2=g_{s_1}(v_1)$ is an endpoint of $s_1^*$; let $s_2$ be the other side incident to $v_2$. Continue by setting $v_{j+1}=g_{s_j}(v_j)$ and choosing $s_{j+1}$ to be the side incident to $v_{j+1}$ other than $s_j^*$. The sequence $v_1,\ldots,v_r$ up to the first return to the original pair $(v_1,s_1)$ is called a vertex cycle, and
\[
 C=g_{s_r}\circ\cdots\circ g_{s_1}
\]
is its cycle transformation. We use the same procedure for ideal vertices, stopping if $s_j^*$ is the only side incident to the image vertex $v_{j+1}$. If the procedure returns to the original vertex-side pair, we call it a closed ideal-vertex cycle and define its cycle transformation by the same composition.

\Needspace{6\baselineskip}
\begin{theorem}[{see~\cite[Theorems~13.4.7 and~13.5.3]{Rat}}]\label{thm:poincare-polygon}
For the polygon $P$ and the side-pairing maps $g_s$ above, assume the following two conditions.
\begin{enumerate}
\item For each finite vertex cycle, there is an integer $m\ge1$ such that the sum of its interior angles is $2\pi/m$ and the cycle transformation $C$ has order $m$. For $m=1$, this means $C=\mathrm{id}$.
\item The cycle transformation of each closed ideal-vertex cycle is parabolic. No such condition is imposed on a sequence terminating at an ideal vertex incident to only one side.
\end{enumerate}
Then the group generated by all the side-pairing maps,
\[
 G=\langle\,g_s\mid s\text{ is a side of }P\,\rangle,
\]
is discrete, and $P$ is a fundamental polygon for it. The space obtained by gluing the sides of the closure of $P$ in $\HH$ as specified is isometric to the complete hyperbolic orbifold $G\backslash\HH$. Moreover, every relation among the generators $g_s$ follows from the inverse relations $g_{s^*}=g_s^{-1}$ and the relations $C^m=1$ for finite vertex cycles.
\end{theorem}

We refer to a cone point or a surface end of the quotient $\OO_{\kk}$ as a peripheral component. If the holonomy of a positively oriented loop around a peripheral component is conjugate in $\Lambda_{\kk}$ to $\overline L_0,\overline M_0,\overline N_0$, or $\overline B$, we call it the first, second, third, or fourth peripheral component, respectively, and use $r$ as its peripheral label.

\begin{definition}\label{def:marking}
For an oriented surface $X$ equipped with a hyperbolic metric and with four ends labeled $1,2,3,4$, a \emph{marking} of $X$ is the isotopy class of an orientation-preserving homeomorphism
\[
 f:\Sigma_{0,4}\longrightarrow X
\]
that sends the end corresponding to $a_r$ to the end of $X$ labeled $r$.

When $X$ is obtained from the quotient $\OO_{\kk}$ by removing its cone points, choose the common basepoint $f(b)$ and one lift of it to the upper half-plane. The marking is said to correspond to $\bar\rho_0$ if the holonomies of the four loops $f\circ\gamma_1,\ldots,f\circ\gamma_4$ can be written, in order, as
\[
 g\overline L_0g^{-1},\quad g\overline M_0g^{-1},\quad
 g\overline N_0g^{-1},\quad g\overline Bg^{-1}
\]
for a single $g\in\PSL_2(\RR)$.
\end{definition}

\Needspace{7\baselineskip}
\begin{theorem}\label{thm:orbifold}\label{prop:geometry}
For the three initial matrices in~\eqref{eq:seed}, define the representation $\rho_0:\pi_1(\Sigma_{0,4},b)\to\SL_2(\ZZ)$ by
\[
 \rho_0(\gamma_1)=L_0,\quad \rho_0(\gamma_2)=M_0,\quad
 \rho_0(\gamma_3)=N_0,\quad \rho_0(\gamma_4)=-B,
\]
and let $\bar\rho_0$ be its projection to $\PSL_2(\RR)$. Set $\Lambda_{\kk}=\operatorname{Im}\bar\rho_0$ and $\OO_{\kk}=\Lambda_{\kk}\backslash\HH$. Then $\Lambda_{\kk}$ is discrete, and the quotient $\OO_{\kk}$ is a complete hyperbolic orbifold. It has four distinct peripheral components that admit the labels above. The fourth peripheral component is the cusp corresponding to the fixed point $\infty$ of $\overline B$. The types of the first three peripheral components are as follows.
\begingroup
\setlength{\topsep}{1pt}
\setlength{\partopsep}{0pt}
\begin{center}
\begin{tabular}{cl}
\toprule
$k_i$&$i$th peripheral component\\
\midrule
$0$&Cone point of order $2$ and cone angle $\pi$\\
$1$&Cone point of order $3$ and cone angle $2\pi/3$\\
$2$&Cusp\\
$\ge3$&Funnel end with entrance length $2\arcosh(k_i/2)$\\
\bottomrule
\end{tabular}
\end{center}
\endgroup
The surface obtained from $\OO_{\kk}$ by removing its cone points admits a marking corresponding to $\bar\rho_0$.
\end{theorem}
\begin{proof}
The initial matrices belong to $\SL_2(\ZZ)$, so $\Lambda_{\kk}$ is discrete as a subgroup of the discrete group $\PSL_2(\ZZ)$. We determine the topological type and marking of the quotient directly from a fundamental polygon.

Let $A_1=L_0,A_2=M_0,A_3=N_0$, and set
\[
 s_1=0,\qquad s_2=k_1+1,\qquad
 s_3=k_1+k_2+2,\qquad b_i=s_i+k_i.
\]
Equation~\eqref{eq:seed} gives
\begin{equation}\label{eq:polygon-pairing}
 A_i(z)=s_i-\frac1{z-b_i},\qquad B(z)=z+K.
\end{equation}
Consider the open region
\begin{equation}\label{eq:fundamental-polygon}
 D=\left\{z\in\HH\mid 
 -\frac12<\Re z<K-\frac12,\quad
 |z-s_i|>1,\ |z-b_i|>1\quad(i=1,2,3)\right\}.
\end{equation}
This is a convex hyperbolic polygon obtained by intersecting the exteriors of unit circles with a vertical strip. Its sides are the circular arcs and two vertical rays appearing in the boundary of its closure in $\HH$. If $k_i=0$, the two circles coincide; divide their arc into two sides at the fixed point $s_i+\mathrm i$. There are thus eight sides for every $k_1,k_2,k_3\in\ZZ_{\ge0}$. Figure~\ref{fig:fundamental-polygon} shows the case $\kk=(1,2,3)$. The light shading indicates the exterior of the circles, and the thick lines indicate the boundary of the region.

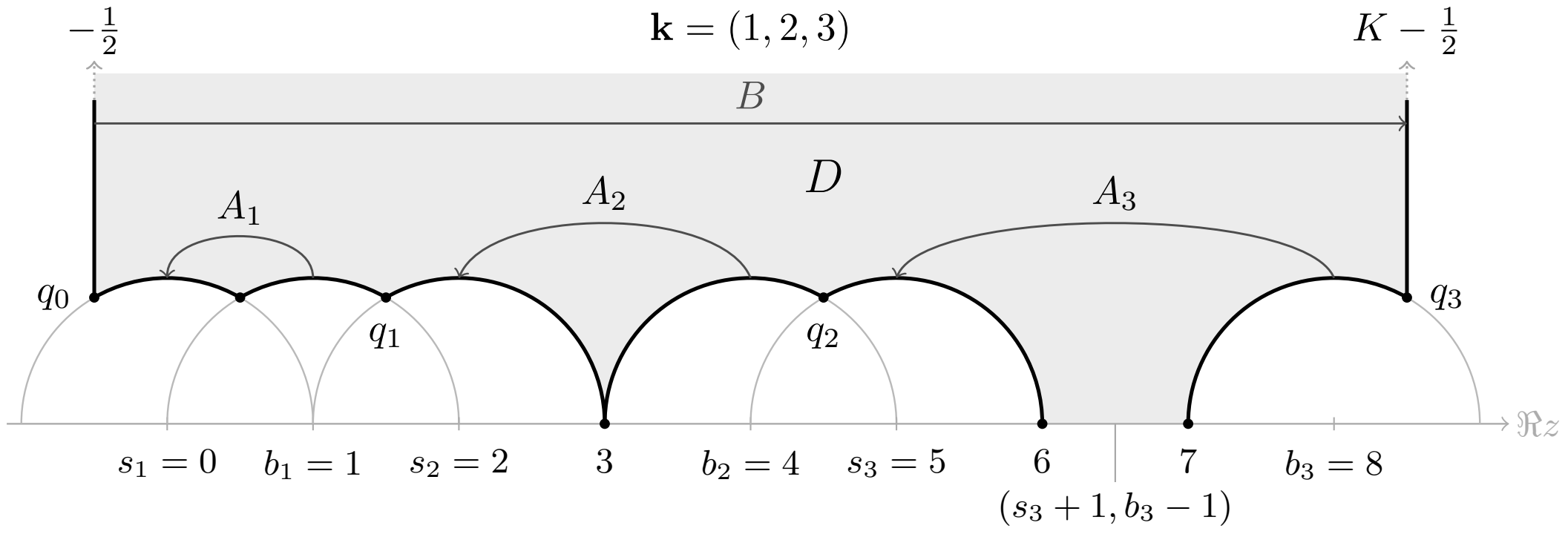
\begin{figure}[!htbp]
\centering
\begin{tikzpicture}[x=1.3cm,y=1.3cm,
 every node/.style={font=\small,inner sep=1.5pt},
 polygon side/.style={black,line width=.95pt},
 circle guide/.style={gray!55,line width=.45pt},
 pairing arrow/.style={black!70,->,line width=.55pt}]
\pgfmathsetmacro{\polygonHeight}{sqrt(3)/2}
\begin{scope}
 \clip (-.5,0) rectangle (8.5,2.4);
 \fill[gray!15] (-.5,0) rectangle (8.5,2.4);
 \foreach \polygonCenter in {0,1,2,4,5,8}
  {\fill[white] (\polygonCenter,0) circle[radius=1];}
\end{scope}
\draw[gray!65,->,line width=.45pt] (-1.10,0)--(9.20,0)
 node[right,font=\footnotesize] {$\Re z$};
\foreach \polygonCenter in {0,1,2,4,5,8} {
 \draw[circle guide] ({\polygonCenter-1},0)
  arc[start angle=180,end angle=0,radius=1];
 \draw[gray!65,line width=.4pt] (\polygonCenter,.045)--(\polygonCenter,-.045);
}
\draw[polygon side] (-.5,\polygonHeight)
 arc[start angle=120,end angle=60,radius=1];
\draw[polygon side] (.5,\polygonHeight)
 arc[start angle=120,end angle=60,radius=1];
\draw[polygon side] (1.5,\polygonHeight)
 arc[start angle=120,end angle=0,radius=1];
\draw[polygon side] (3,0)
 arc[start angle=180,end angle=60,radius=1];
\draw[polygon side] (4.5,\polygonHeight)
 arc[start angle=120,end angle=0,radius=1];
\draw[polygon side] (7,0)
 arc[start angle=180,end angle=60,radius=1];
\draw[polygon side] (-.5,\polygonHeight)--(-.5,2.22);
\draw[polygon side] (8.5,\polygonHeight)--(8.5,2.22);
\draw[gray!70,densely dotted,->,line width=.6pt] (-.5,2.22)--(-.5,2.49);
\draw[gray!70,densely dotted,->,line width=.6pt] (8.5,2.22)--(8.5,2.49);
\node[above] at (-.5,2.49) {$-\tfrac12$};
\node[above] at (8.5,2.49) {$K-\tfrac12$};
\node at (4,2.70) {$\mathbf k=(1,2,3)$};
\node[font=\large] at (4.5,1.70) {$D$};
\draw[pairing arrow] (-.5,2.06)--(8.5,2.06)
 node[midway,above=2pt] {$B$};
\draw[pairing arrow] (1,1)
 .. controls (1,1.38) and (0,1.38) .. (0,1);
\node at (.5,1.48) {$A_1$};
\draw[pairing arrow] (4,1)
 .. controls (3.8,1.50) and (2.2,1.50) .. (2,1);
\node at (3,1.58) {$A_2$};
\draw[pairing arrow] (8,1)
 .. controls (7.8,1.50) and (5.2,1.50) .. (5,1);
\node at (6.5,1.58) {$A_3$};
\foreach \polygonPoint in {(-.5,\polygonHeight),(.5,\polygonHeight),
 (1.5,\polygonHeight),(3,0),(4.5,\polygonHeight),(6,0),(7,0),(8.5,\polygonHeight)}
 {\fill \polygonPoint circle[radius=1.3pt];}
\node[left=4pt] at (-.5,\polygonHeight) {$q_0$};
\node[below=5pt] at (1.5,\polygonHeight) {$q_1$};
\node[below=5pt] at (4.5,\polygonHeight) {$q_2$};
\node[right=4pt] at (8.5,\polygonHeight) {$q_3$};
\node[below=5pt,font=\footnotesize] at (0,0) {$s_1=0$};
\node[below=5pt,font=\footnotesize] at (1,0) {$b_1=1$};
\node[below=5pt,font=\footnotesize] at (2,0) {$s_2=2$};
\node[below=5pt,font=\footnotesize] at (3,0) {$3$};
\node[below=5pt,font=\footnotesize] at (4,0) {$b_2=4$};
\node[below=5pt,font=\footnotesize] at (5,0) {$s_3=5$};
\node[below=5pt,font=\footnotesize] at (6,0) {$6$};
\node[below=5pt,font=\footnotesize] at (7,0) {$7$};
\node[below=5pt,font=\footnotesize] at (8,0) {$b_3=8$};
\draw[gray!70,line width=.4pt] (6.5,0)--(6.5,-.40);
\node[below,font=\footnotesize] at (6.5,-.41) {$(s_3+1,b_3-1)$};
\end{tikzpicture}
\caption{The fundamental polygon $D$}
\label{fig:fundamental-polygon}
\end{figure}

$A_i$ maps the side on $|z-b_i|=1$ to the side on $|z-s_i|=1$; for $k_i=0$, it maps the right half of the divided arc to the left half. Indeed,
\[
 |A_i(z)-s_i|=\frac1{|z-b_i|},
\]
and the endpoint correspondences below show that the entire sides are paired. Moreover, $A_i(D)$ lies inside the circle with center $s_i$ and radius $1$, and hence is disjoint from $D$. The left and right vertical sides are glued together by the translation $B(z)=z+K$. Moreover, $B(D)$ is disjoint from $D$.

The side-pairing maps are $\overline A_1,\overline A_2,\overline A_3,\overline B$ and their inverses. Since $\overline B=(\overline A_1\overline A_2\overline A_3)^{-1}$ by~\eqref{eq:seed-data}, the group generated in Theorem~\ref{thm:poincare-polygon} is
\[
 G=\langle\overline A_1,\overline A_2,\overline A_3,\overline B\rangle
 =\Lambda_{\kk}.
\]

Set
\[
 q_0=-\frac12+\frac{\sqrt3}{2}\,\mathrm i,\qquad
 q_i=b_i+\frac12+\frac{\sqrt3}{2}\,\mathrm i\quad(i=1,2,3).
\]
Since $b_3=K-1$, we have $q_3=K-\frac12+\frac{\sqrt3}{2}\,\mathrm i$. The distance between the centers of neighboring circles from distinct pairs is $1$, and
\[
 B(q_0)=q_3,\qquad A_3(q_3)=q_2,\qquad
 A_2(q_2)=q_1,\qquad A_1(q_1)=q_0.
\]
The interior angles at $q_0,q_3$ are $\pi/3$, and those at $q_1,q_2$ are $2\pi/3$, with total $2\pi$. By~\eqref{eq:seed-data}, the cycle transformation is the identity in $\PSL_2(\RR)$.

The remaining finite vertices are $s_i+\mathrm i$ for $k_i=0$ and $s_i+\frac12+\frac{\sqrt3}{2}\,\mathrm i$ for $k_i=1$. Each forms a one-vertex cycle, with interior angle $\pi$ or $2\pi/3$ and cycle transformation $\overline A_i$ of order $2$ or $3$, respectively. The cycle transformation at $\infty$ is the parabolic map $B$. For $k_i=2$, the two circles are tangent at $s_i+1=b_i-1$; this vertex also has a parabolic cycle transformation, namely $A_i$.

For $k_i\ge3$, the two circles are disjoint, leaving the real interval $[s_i+1,b_i-1]$ between their endpoints. Only one hyperbolic side is incident to each endpoint, and $A_i(b_i-1)=s_i+1$. The endpoint correspondence therefore stops there and does not form a closed cycle. We have now examined all vertices and verified the angle conditions at finite vertices and the parabolic conditions for closed ideal-vertex cycles. Theorem~\ref{thm:poincare-polygon} thus shows that $D$ is a fundamental polygon for $\Lambda_{\kk}$ and that the quotient is obtained by the side identifications above.

We next determine the topological type of this gluing. Take the closure of $D$ in the closed upper half-plane including the extended real axis $\RR\cup\{\infty\}$, and collapse each ideal boundary interval $[s_i+1,b_i-1]$ with $k_i\ge3$ to a point. This gives a topological octagon. After the identifications, $q_0,\ldots,q_3$ become one ordinary point, and the remaining four vertices become four peripheral points. There are five vertices, four edges, and one face, so the Euler characteristic is $5-4+1=2$ and the resulting closed oriented surface is a sphere. The four edges form a tree joining the central ordinary point to the four peripheral points, and the octagon is obtained by cutting the sphere along this tree. Matching the labels and orientations of the tree edges with the connecting paths in Section~2 makes the holonomy of $\gamma_i$ equal to $\overline A_i$ and that of $\gamma_4$ equal to $\overline B$. Thus the resulting homeomorphism gives a marking corresponding to the specified $\bar\rho_0$. Figure~\ref{fig:polygon-gluing} illustrates this construction for $\kk=(1,2,3)$.

\begin{figure}[!htbp]
\centering
\begin{tikzpicture}[x=1cm,y=1cm,every node/.style={font=\small,inner sep=1.5pt},
 glue one/.style={draw=blue!65!black,text=blue!65!black,line width=1.05pt},
 glue two/.style={draw=orange!80!black,text=orange!80!black,line width=1.05pt},
 glue three/.style={draw=teal!70!black,text=teal!70!black,line width=1.05pt},
 glue four/.style={draw=violet!75!black,text=violet!75!black,line width=1.05pt}]
\begin{scope}[x=.72cm,y=.72cm,yshift=-1.05cm]
\pgfmathsetmacro{\glueHeight}{sqrt(3)/2}
\begin{scope}
 \clip (-.5,0) rectangle (8.5,3.90);
 \fill[gray!15] (-.5,0) rectangle (8.5,3.90);
 \foreach \glueCenter in {0,1,2,4,5,8}
  {\fill[white] (\glueCenter,0) circle[radius=1];}
\end{scope}
\draw[gray!60,->,line width=.4pt] (-1.10,0)--(9.10,0)
 node[right,font=\footnotesize] {$\Re z$};
\foreach \glueCenter in {0,1,2,4,5,8} {
 \draw[gray!50,line width=.4pt] ({\glueCenter-1},0)
  arc[start angle=180,end angle=0,radius=1];
}
\draw[glue one] (-.5,\glueHeight)
 arc[start angle=120,end angle=60,radius=1];
\draw[glue one] (.5,\glueHeight)
 arc[start angle=120,end angle=60,radius=1];
\draw[glue two] (1.5,\glueHeight)
 arc[start angle=120,end angle=0,radius=1];
\draw[glue two] (3,0)
 arc[start angle=180,end angle=60,radius=1];
\draw[glue three] (4.5,\glueHeight)
 arc[start angle=120,end angle=0,radius=1];
\draw[glue three] (7,0)
 arc[start angle=180,end angle=60,radius=1];
\draw[glue four] (-.5,\glueHeight)--(-.5,3.50);
\draw[glue four] (8.5,\glueHeight)--(8.5,3.50);
\draw[glue four,densely dotted,->] (-.5,3.50)--(-.5,4.02);
\draw[glue four,densely dotted,->] (8.5,3.50)--(8.5,4.02);
\node[above,font=\footnotesize] at (-.5,4.02) {$-\tfrac12$};
\node[above,font=\footnotesize] at (8.5,4.02) {$K-\tfrac12$};
\node at (4,4.45) {$a_4:\ \infty$};
\node[font=\large] at (4,2.40) {$D$};
\draw[glue four,->,line width=.6pt] (-.5,3.30)--(8.5,3.30)
 node[midway,above=2pt] {$B$};
\draw[glue one,->,line width=.6pt] (1,1)
 ..controls(1,1.50) and(0,1.50)..(0,1);
\node[text=blue!65!black] at (.5,1.75) {$A_1$};
\draw[glue two,->,line width=.6pt] (4,1)
 ..controls(3.8,1.55) and(2.2,1.55)..(2,1);
\node[text=orange!80!black] at (3,1.80) {$A_2$};
\draw[glue three,->,line width=.6pt] (8,1)
 ..controls(7.8,1.55) and(5.2,1.55)..(5,1);
\node[text=teal!70!black] at (6.5,1.80) {$A_3$};
\foreach \gluePoint in {(-.5,\glueHeight),(1.5,\glueHeight),(4.5,\glueHeight),(8.5,\glueHeight)}
 {\fill \gluePoint circle[radius=1.3pt];}
\node[left=2pt,font=\footnotesize] at (-.5,\glueHeight) {$q_0$};
\node[below=4pt,font=\footnotesize] at (1.5,\glueHeight) {$q_1$};
\node[below=4pt,font=\footnotesize] at (4.5,\glueHeight) {$q_2$};
\node[right=2pt,font=\footnotesize] at (8.5,\glueHeight) {$q_3$};
\fill (.5,\glueHeight) circle[radius=1.45pt];
\node[below=4pt,font=\footnotesize] at (.5,\glueHeight) {$a_1$};
\draw[fill=white,line width=.5pt] (3,0) circle[radius=1.45pt];
\draw[fill=white,line width=.5pt] (6,0) circle[radius=1.3pt];
\draw[fill=white,line width=.5pt] (7,0) circle[radius=1.3pt];
\foreach \glueCoordinate in {0,1,2,3,4,5,6,7,8}
 {\node[below=3pt,font=\footnotesize] at (\glueCoordinate,0) {$\glueCoordinate$};}
\node[below,font=\footnotesize] at (3,-.62) {$a_2$};
\node[below,font=\footnotesize] at (6.5,-.62) {$a_3$: funnel end};
\end{scope}
\draw[->,line width=.65pt] (7.25,0)--(8.37,0)
 node[midway,above=4pt] {gluing};
\begin{scope}[xshift=10.4cm,scale=.82]
\fill[gray!15]
 (-.12,2.30)
 .. controls (-.16,1.65) and (-.38,1.18) .. (-.82,.97)
 .. controls (-1.04,.86) and (-1.19,.69) .. (-1.28,.62)
 --(-2.18,.18)--(-1.28,-.28)
 .. controls (-1.24,-.60) and (-1.11,-.83) .. (-.82,-.98)
 .. controls (-.37,-1.20) and (-.16,-1.68) .. (-.12,-2.30)
 --(.12,-2.30)
 .. controls (.16,-1.68) and (.38,-1.20) .. (.82,-.98)
 .. controls (1.14,-.81) and (1.28,-.54) .. (1.47,-.44)
 .. controls (1.98,-.48) and (2.42,-.80) .. (3.00,-1.02)
 --(3.00,1.38)
 .. controls (2.42,1.15) and (1.98,.83) .. (1.47,.80)
 .. controls (1.28,.93) and (1.08,1.00) .. (.82,1.02)
 .. controls (.36,1.20) and (.16,1.68) .. (.12,2.30)--cycle;
\draw[black!65,line width=.6pt]
 (-.12,2.30)
 .. controls (-.16,1.65) and (-.38,1.18) .. (-.82,.97)
 .. controls (-1.04,.86) and (-1.19,.69) .. (-1.28,.62)
 --(-2.18,.18)--(-1.28,-.28)
 .. controls (-1.24,-.60) and (-1.11,-.83) .. (-.82,-.98)
 .. controls (-.37,-1.20) and (-.16,-1.68) .. (-.12,-2.30);
\draw[black!65,line width=.6pt]
 (.12,-2.30)
 .. controls (.16,-1.68) and (.38,-1.20) .. (.82,-.98)
 .. controls (1.14,-.81) and (1.28,-.54) .. (1.47,-.44)
 .. controls (1.98,-.48) and (2.42,-.80) .. (3.00,-1.02);
\draw[black!65,line width=.6pt]
 (3.00,1.38)
 .. controls (2.42,1.15) and (1.98,.83) .. (1.47,.80)
 .. controls (1.28,.93) and (1.08,1.00) .. (.82,1.02)
 .. controls (.36,1.20) and (.16,1.68) .. (.12,2.30);
\draw[gray!55,densely dashed,line width=.45pt] (-1.28,.62)
 arc[start angle=90,end angle=270,x radius=.14,y radius=.45];
\draw[gray!55,line width=.45pt] (-1.28,-.28)
 arc[start angle=270,end angle=450,x radius=.14,y radius=.45];
\draw[gray!55,densely dotted,line width=.45pt] (0,2.30)
 ellipse[x radius=.12,y radius=.045];
\draw[gray!55,densely dotted,line width=.45pt] (0,-2.30)
 ellipse[x radius=.12,y radius=.045];
\draw[gray!55,densely dotted,line width=.45pt] (3,.18)
 ellipse[x radius=.20,y radius=1.20];
\foreach \glueSign in {-1,1} {
 \draw[gray!65,densely dashed,line width=.5pt]
  ({\glueSign*.12},2.30)--({\glueSign*.06},2.88);
 \draw[gray!65,densely dashed,line width=.5pt]
  ({\glueSign*.12},-2.30)--({\glueSign*.06},-2.88);
}
\draw[gray!65,densely dashed,line width=.5pt] (3,1.38)--(3.65,1.68);
\draw[gray!65,densely dashed,line width=.5pt] (3,-1.02)--(3.65,-1.32);
\draw[black!75,densely dashed,line width=.8pt] (1.47,.80)
 arc[start angle=90,end angle=270,x radius=.14,y radius=.62];
\draw[black,line width=1.0pt] (1.47,-.44)
 arc[start angle=270,end angle=450,x radius=.14,y radius=.62];
\coordinate (glueQ) at (0,.40);
\draw[glue one] (glueQ)..controls(-.62,{.40-.62/sqrt(3)}) and(-1.20,.17)..(-2.18,.18);
\draw[glue two] (glueQ)..controls(0,1.02) and(.03,1.60)..(0,2.30);
\draw[glue three] (glueQ)..controls(.65,{.40-.65/sqrt(3)}) and(1.26,.18)..(1.61,.18)--(3.20,.18);
\draw[glue four] (glueQ)..controls(0,-.32) and(.02,-1.55)..(0,-2.30);
\draw[glue two,densely dashed,->] (0,2.30)--(0,2.98);
\draw[glue three,densely dashed,->] (3.20,.18)--(3.82,.18);
\draw[glue four,densely dashed,->] (0,-2.30)--(0,-2.98);
\fill (glueQ) circle[radius=1.65pt];
\node[above left=3pt] at (glueQ) {$q$};
\fill (-2.18,.18) circle[radius=1.6pt];
\node[text=blue!65!black] at (-.78,.37) {$A_1$};
\node[text=orange!80!black] at (.34,1.40) {$A_2$};
\node[text=teal!70!black] at (2.13,.43) {$A_3$};
\node[text=violet!75!black] at (.30,-1.48) {$B$};
\draw[gray!65,line width=.4pt] (-2.18,.18)--(-2.02,.66);
\node[above] at (-2.02,.69) {$a_1$};
\node[right=5pt] at (.12,2.74) {$a_2$};
\node[right=5pt] at (.12,-2.74) {$a_4$};
\node[right=3pt] at (3.65,1.65) {$a_3$};
\node at (1.1,3.35) {$\mathcal O_{\mathbf k}$};
\end{scope}
\end{tikzpicture}
\caption{Gluing the fundamental polygon $D$ for $\kk=(1,2,3)$}
\label{fig:polygon-gluing}
\end{figure}
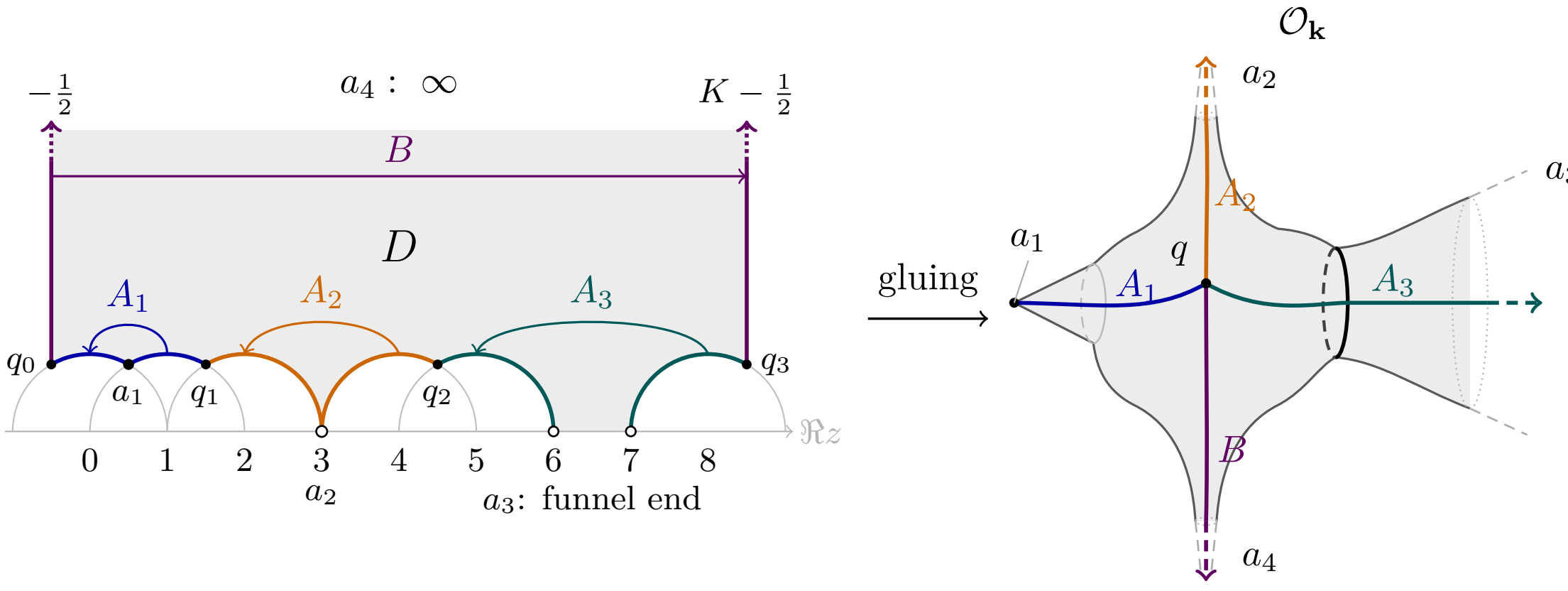

\begin{figure}[!htbp]
\centering
\begin{tikzpicture}[x=1.1cm,y=1.1cm,
 every node/.style={font=\small,inner sep=1.2pt},
 local region/.style={gray!15},
 local guide/.style={black!45,line width=.45pt},
 local angle/.style={black!65,line width=.45pt},
 local side/.style={black!80,line width=.8pt}]
\pgfmathsetmacro{\localHeight}{sqrt(3)/2}
\pgfmathsetmacro{\axisRadius}{sqrt(5)/2}
\pgfmathsetmacro{\axisHeight}{sqrt(5)/3}
\pgfmathsetmacro{\axisAngle}{acos(sqrt(5)/3)}
\pgfmathsetmacro{\sideAngle}{acos(2/3)}
\begin{scope}
 \node at (0,1.86) {$k_i=0,1$};
 \begin{scope}
  \clip (-.94,0) rectangle (.94,1.42);
  \fill[local region] (-.94,0) rectangle (.94,1.42);
  \fill[white] (-.5,0) circle[radius=1];
  \fill[white] (.5,0) circle[radius=1];
 \end{scope}
 \draw[local guide] (-1.66,0)--(1.66,0);
 \draw[local guide] (-1.5,0) arc[start angle=180,end angle=0,radius=1];
 \draw[local guide] (-.5,0) arc[start angle=180,end angle=0,radius=1];
 \draw[local side] ({-.5+cos(110)},{sin(110)})
  arc[start angle=110,end angle=60,radius=1];
 \draw[local side] (0,\localHeight)
  arc[start angle=120,end angle=70,radius=1];
 \fill (0,\localHeight) circle[radius=1.5pt];
 \draw[local guide] (0,1.24)--(0,\localHeight);
 \node[above,font=\footnotesize] at (0,1.24) {cone point};
 \node[below=3pt,font=\footnotesize] at (-.5,0) {$s_i$};
 \node[below=3pt,font=\footnotesize] at (.5,0) {$b_i$};
 \node at (0,-.62) {cone angle};
 \node[anchor=north] at (0,-.86) {$\begin{aligned}
   \pi &\quad (k_i=0)\\
   2\pi/3 &\quad (k_i=1)
 \end{aligned}$};
\end{scope}
\begin{scope}[xshift=4.65cm]
 \node at (0,1.86) {$k_i=2$};
 \begin{scope}
  \clip (-1.2,0) rectangle (1.2,1.42);
  \fill[local region] (-1.2,0) rectangle (1.2,1.42);
  \fill[white] (-1,0) circle[radius=1];
  \fill[white] (1,0) circle[radius=1];
 \end{scope}
 \draw[local guide] (-2.14,0)--(2.14,0);
 \draw[local guide] (-2,0) arc[start angle=180,end angle=0,radius=1];
 \draw[local guide] (0,0) arc[start angle=180,end angle=0,radius=1];
 \draw[local side] (-1,1) arc[start angle=90,end angle=0,radius=1];
 \draw[local side] (0,0) arc[start angle=180,end angle=90,radius=1];
 \filldraw[fill=white,draw=black!80,line width=.45pt] (0,0) circle[radius=1.6pt];
 \node[below=3pt,font=\footnotesize] at (-1,0) {$s_i$};
 \node[below=3pt,font=\footnotesize] at (1,0) {$b_i$};
 \node[below=3pt,font=\footnotesize] at (0,0) {cusp};
\end{scope}
\begin{scope}[xshift=10.15cm]
 \node at (0,1.86) {$k_i\ge3$};
 \begin{scope}
  \clip (-1.7,0) rectangle (1.7,1.42);
  \fill[local region] (-1.7,0) rectangle (1.7,1.42);
  \fill[white] (-1.5,0) circle[radius=1];
  \fill[white] (1.5,0) circle[radius=1];
 \end{scope}
 \begin{scope}
  \clip (-.5,0) arc[start angle=0,end angle=\sideAngle,radius=1]
   arc[start angle={180-\axisAngle},end angle=\axisAngle,radius=\axisRadius]
   arc[start angle={180-\sideAngle},end angle=180,radius=1] -- cycle;
  \foreach \hatchX in {-2,-1.85,...,1} {
   \draw[black!45,line width=.36pt] (\hatchX,0)--++(1.2,1.2);
  }
 \end{scope}
 \draw[local guide] (-2.65,0)--(2.65,0);
 \draw[local guide] (-2.5,0) arc[start angle=180,end angle=0,radius=1];
 \draw[local guide] (.5,0) arc[start angle=180,end angle=0,radius=1];
 \draw[local side] (-1.5,1) arc[start angle=90,end angle=0,radius=1];
 \draw[local side] (.5,0) arc[start angle=180,end angle=90,radius=1];
 \draw[local guide] (-\axisRadius,0)
  arc[start angle=180,end angle=0,radius=\axisRadius];
 \draw[black,line width=1.4pt] ({-5/6},\axisHeight)
  arc[start angle={180-\axisAngle},end angle=\axisAngle,radius=\axisRadius];
 \foreach \localSign in {-1,1} {
  \draw[local angle]
   ({\localSign*(5/6-.105*2/3)},{\axisHeight+.105*sqrt(5)/3})
   --({\localSign*(5/6-.105*(2+sqrt(5))/3)},
      {\axisHeight+.105*(sqrt(5)-2)/3})
   --({\localSign*(5/6-.105*sqrt(5)/3)},{\axisHeight-.105*2/3});
 }
 \fill ({-5/6},\axisHeight) circle[radius=1.5pt];
 \fill ({5/6},\axisHeight) circle[radius=1.5pt];
 \node[font=\footnotesize,fill=gray!15,inner sep=1.5pt] at (0,.46) {funnel};
 \node[below=3pt,font=\footnotesize] at (-1.5,0) {$s_i$};
 \node[below=3pt,font=\footnotesize] at (1.5,0) {$b_i$};
 \node at (0,-.62) {entrance length};
 \node[anchor=north] at (0,-.86) {$2\arcosh(k_i/2)$};
\end{scope}
\end{tikzpicture}
\caption{Local peripheral geometry}
\label{fig:peripheral-geometry}
\end{figure}
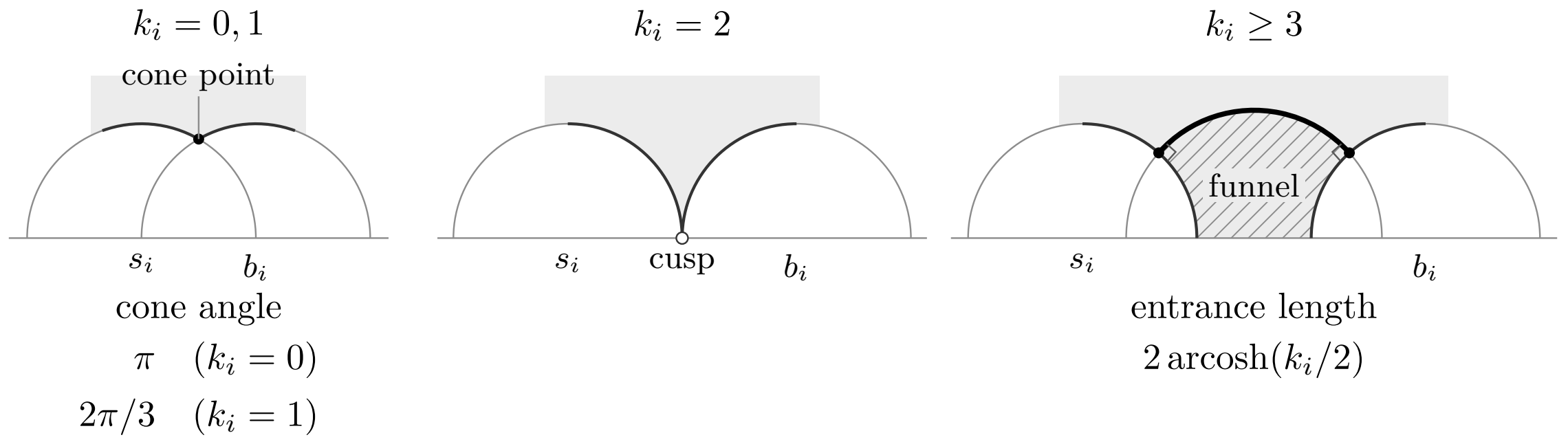

Finally, we verify the geometric type of each peripheral component using the trace classification (see~\cite{Kat}). The finite vertices for $k_i=0,1$ become cone points of angles $\pi,2\pi/3$, respectively, and the tangency point for $k_i=2$ becomes a cusp. For $k_i\ge3$, the axis of $A_i$ is the semicircle with center $s_i+k_i/2$ and radius $\sqrt{k_i^2-4}/2$. It meets both unit circles orthogonally, and $A_i$ identifies the endpoints of the segment between them. The side of this closed geodesic facing the real axis is a funnel. By Theorem~\ref{thm:trace-length} and $|\tr A_i|=k_i$, its entrance has length $2\arcosh(k_i/2)$. Figure~\ref{fig:peripheral-geometry} illustrates the cone point, cusp, and funnel using $k_i=1,2,3$, respectively. The light shading marks the exterior of the two unit circles. In the right panel, the hatched region projects to the funnel and the thick arc to its entrance. For $Y>1$, the region $D\cap\{\operatorname{Im}z>Y\}$ is a vertical strip of width $K$, and its two sides are identified once by $B$. There are no further identifications in the interior of the fundamental region, so the transformation corresponding to one circuit around the cusp at $\infty$ is $\overline B$, and the elements of $\Lambda_{\kk}$ fixing $\infty$ are precisely $\overline B^{\,n}$ ($n\in\mathbb Z$). After removing the cone points, all four peripheral components become punctures, giving $\Sigma_{0,4}$.
\end{proof}

\Needspace{6\baselineskip}
\begin{proposition}\label{prop:orbifold-group}
Let $E=\{i\in\{1,2,3\}\mid k_i=0\text{ or }1\}$, and set $p_i=2$ for $k_i=0$ and $p_i=3$ for $k_i=1$, for each $i\in E$. The representation $\bar\rho_0$ induces an isomorphism
\begin{equation}\label{eq:orbgroup}
 \Gamma_{\mathrm{orb}}
 :=\pi_1(\Sigma_{0,4},b)/
 \langle\!\langle\gamma_i^{p_i}\ (i\in E)\rangle\!\rangle
 \xrightarrow[\ \bar\rho_0\ ]{\ \cong\ }\Lambda_{\kk}.
\end{equation}
Moreover, $\Gamma_{\mathrm{orb}}$ is a free product of three cyclic groups: the $i$th factor is cyclic of order $p_i$ if $i\in E$, and is $\ZZ$ otherwise.
\end{proposition}
\begin{proof}
Let $A_1=L_0,A_2=M_0,A_3=N_0$. Apply the presentation statement of Theorem~\ref{thm:poincare-polygon} to the fundamental polygon $D$ established in the proof of Theorem~\ref{thm:orbifold}. Every relation among the generators $\overline A_1,\overline A_2,\overline A_3,\overline B$ follows from
\[
 \overline A_1\overline A_2\overline A_3\overline B=1,
 \qquad \overline A_i^{\,p_i}=1\quad(i\in E).
\]
These are precisely the relations obtained by imposing $\gamma_i^{p_i}=1$ ($i\in E$) on $\pi_1(\Sigma_{0,4},b)$, giving the isomorphism~\eqref{eq:orbgroup}. Eliminating $\gamma_4=(\gamma_1\gamma_2\gamma_3)^{-1}$ gives the asserted free product presentation.
\end{proof}

\begin{remark}\label{rem:holonomy-realization}
Goldman~\cite[Theorem~A]{Gol80} proved that a representation of the fundamental group of a closed oriented surface of genus at least $2$ into $\PSL_2(\RR)$ is discrete and faithful if and only if its Euler number has maximal absolute value, and realized such representations as holonomies of hyperbolic structures. For compact connected oriented surfaces $S$, possibly with boundary, the formulation in terms of the Toledo invariant was later given by Burger--Iozzi--Wienhard (see~\cite[Theorem~3]{BIW}): a representation $\bar\rho:\pi_1(S)\to\PSL_2(\RR)$ with $\chi(S)<0$ has Toledo invariant $-\chi(S)$ if and only if it is the holonomy of a complete hyperbolic metric on $S^\circ$. When there are no cone points, our construction gives a concrete instance of this general theory. When cone points are present, we restrict the orbifold holonomy to a torsion-free subgroup of finite index in $\Gamma_{\mathrm{orb}}$ and obtain an instance of the theory on the corresponding covering surface.
\end{remark}

\section{Realizing GM numbers as \texorpdfstring{$\lambda$}{lambda}-lengths of ideal arcs}\label{sec:lambda}
Fix a marking $f:\Sigma_{0,4}\to X$ corresponding to $\bar\rho_0$, where $X$ is obtained from $\OO_{\kk}$ by removing its cone points and carries the induced hyperbolic metric. The symbols $a_1,a_2,a_3,a_4$ label the punctures of the topological surface $\Sigma_{0,4}$; the end labeled $a_4$ corresponds under $f$ to the distinguished cusp of the quotient. Let $\EE$ be the set of unoriented isotopy classes of essential simple ideal arcs based at $a_4$ on $\Sigma_{0,4}$. The interior of each arc is embedded and avoids the other marked points. We identify arcs under isotopies fixing the marked points and reversal of orientation; tangent directions at the cusp are not fixed. Each $e\in\EE$ separates the remaining three marked points into one on one side and two on the other. Write $a_{i(e)}$ for the point on the side containing a single marked point, and set $k_e=k_{i(e)}$.

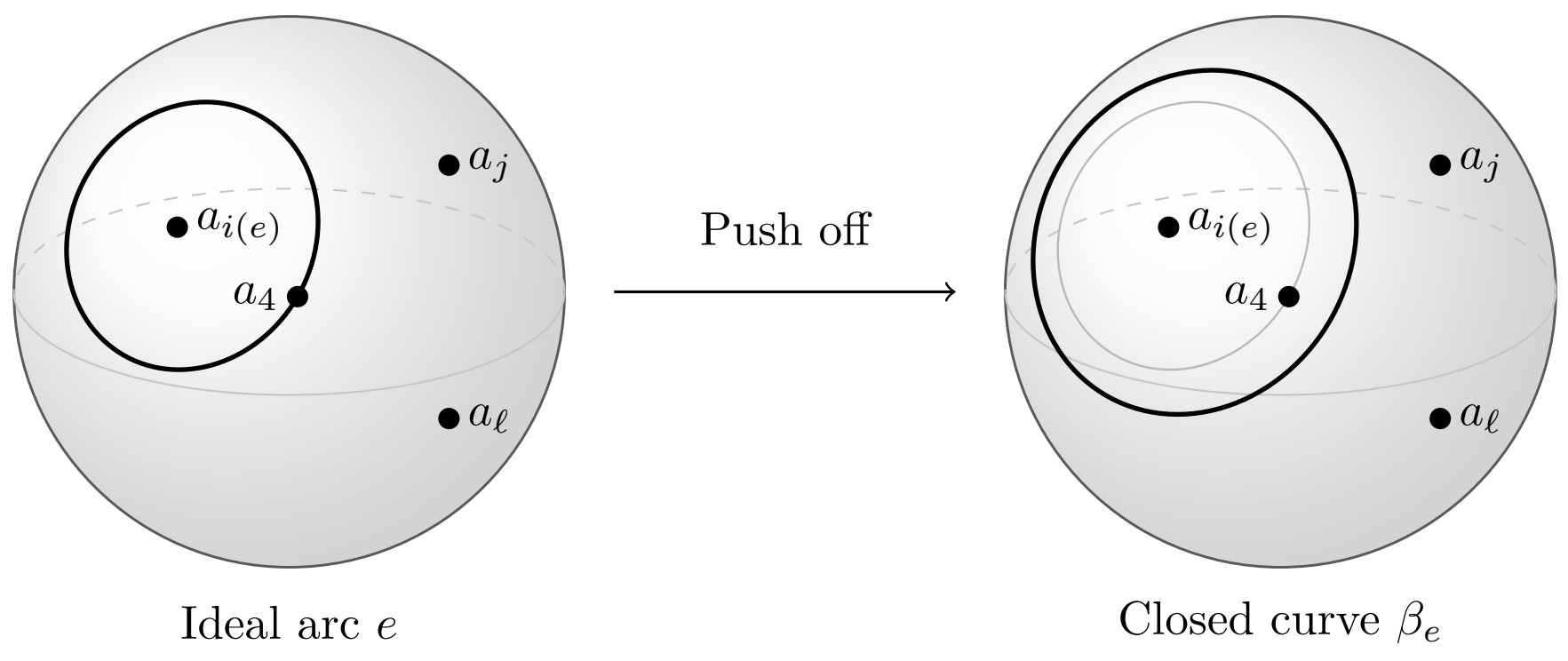
\begin{figure}[!htbp]
\centering
\begin{tikzpicture}[x=2.05cm,y=2.05cm,
  every node/.style={font=\small},
  marked point/.style={circle,fill=black,inner sep=1.6pt},
  point label/.style={anchor=west,inner sep=0pt,outer sep=0pt,xshift=4pt}]
 \pgfmathsetmacro{\returnNx}{sin(28)*cos(150)}
 \pgfmathsetmacro{\returnNy}{sin(28)*sin(150)}
 \pgfmathsetmacro{\returnUx}{cos(28)*cos(150)}
 \pgfmathsetmacro{\returnUy}{cos(28)*sin(150)}
 \pgfmathsetmacro{\returnVx}{-sin(150)}
 \pgfmathsetmacro{\returnVy}{cos(150)}
 \foreach \returnShift/\returnRadius/\returnPanel in
  {0/30/{Ideal arc $e$},3.6/40/{Closed curve $\beta_e$}} {
  \begin{scope}[shift={(\returnShift,0)}]
   \shade[ball color=gray!10,opacity=0.28] (0,0) circle[radius=1];
   \draw[black!65,line width=0.55pt] (0,0) circle[radius=1];
   \draw[gray!45,dashed,line width=0.4pt]
    plot[domain=0:180,samples=61,variable=\spheret]
    ({cos(\spheret)},{cos(68)*sin(\spheret)});
   \draw[gray!45,line width=0.4pt]
    plot[domain=180:360,samples=61,variable=\spheret]
    ({cos(\spheret)},{cos(68)*sin(\spheret)});
   \ifnum\returnRadius=40\relax
    \draw[gray!55,line width=0.45pt]
     plot[domain=0:360,samples=121,variable=\spheret]
     ({cos(30)*\returnNx+sin(30)*(cos(\spheret)*\returnUx+sin(\spheret)*\returnVx)},
      {cos(30)*\returnNy+sin(30)*(cos(\spheret)*\returnUy+sin(\spheret)*\returnVy)}) -- cycle;
   \fi
   \draw[line width=1pt]
    plot[domain=0:360,samples=121,variable=\spheret]
    ({cos(\returnRadius)*\returnNx+sin(\returnRadius)*(cos(\spheret)*\returnUx+sin(\spheret)*\returnVx)},
     {cos(\returnRadius)*\returnNy+sin(\returnRadius)*(cos(\spheret)*\returnUy+sin(\spheret)*\returnVy)}) -- cycle;
   \coordinate (returnCusp)
    at ({cos(30)*\returnNx-sin(30)*\returnUx},
        {cos(30)*\returnNy-sin(30)*\returnUy});
   \coordinate (returnAi) at (\returnNx,\returnNy);
   \coordinate (returnAj) at (0.58,0.46);
   \coordinate (returnAl) at (0.58,-0.46);
   \node[marked point] at (returnCusp) {};
   \node[anchor=east,inner sep=0pt,outer sep=0pt,xshift=-4pt]
    at (returnCusp) {$a_4$};
   \foreach \returnPoint/\returnText in
    {returnAi/a_{i(e)},returnAj/a_j,returnAl/a_\ell} {
    \node[marked point] at (\returnPoint) {};
    \node[point label] at (\returnPoint) {$\returnText$};
   }
   \node at (0,-1.2) {\returnPanel};
  \end{scope}
 }
 \draw[->,line width=0.6pt] (1.18,0) -- (2.42,0)
  node[midway,above=6pt] {Push off};
\end{tikzpicture}
\caption{An ideal arc $e$ (left) and its closed curve $\beta_e$ (right)}
\label{fig:return-arc-closed}
\end{figure}

Let $a_j,a_\ell$ be the other two marked points. The closed curve used in the McShane identity is $\beta_e$, obtained by pushing $e$ slightly toward the side containing these two points (Figure~\ref{fig:return-arc-closed}). The curve $\beta_e$ separates $\{a_4,a_{i(e)}\}$ from the other two marked points, and this correspondence is bijective. Below, sums and the injectivity conjecture are indexed directly by ideal arcs $e$, and $\ell_{\beta_e}$ denotes the length of the geodesic realization of the marked curve $f(\beta_e)$ in $\OO_{\kk}$.

We define the Farey relation by the Farey tree in Section~\ref{sec:farey-tree}, and the correspondence between fractions and ideal arcs by slopes in the plane in Section~\ref{sec:slopes}. We write $m_e$ for the generalized Markov number attached to an individual ideal arc and $(x,y,z)$ for the triple attached to a vertex of the tree.

\subsection{Farey triples and the Farey tree}\label{sec:farey-tree}
Let $\widehat{\mathbb Q}=\mathbb Q\cup\{\tfrac10\}$, where $\tfrac10$ denotes infinity. Write each element as a reduced fraction $p/q$, normalized by $q>0$ or $(p,q)=(1,0)$. The corresponding primitive vector is $\pm(q,p)$, up to sign.

\begin{definition}\label{def:farey-tree}
The vertices of the \emph{Farey tree} $\mathcal T_{\mathrm F}$ are unordered triples of fractions such that any two components $p/q,r/s$ satisfy $|ps-qr|=1$. Two vertices are joined by an edge if their triples share two components. The root is $\tau_0=\{\tfrac01,\tfrac11,\tfrac10\}$.
Two fractions are \emph{Farey-related} if they occur as components of one vertex of the tree. The three components of a vertex are called a \emph{Farey triple}.
\end{definition}

If the shared components are $a/b,c/d$, the triples joined by the edge are
\begin{equation}\label{eq:farey-exchange}
 \left\{\frac ab,\frac cd,\frac{a+c}{b+d}\right\}
 \longleftrightarrow
 \left\{\frac ab,\frac cd,\frac{a-c}{b-d}\right\}.
\end{equation}
Thus two primitive vectors $\mathbf v,\mathbf w$ forming a lattice basis are kept, while the third component is exchanged between $[\mathbf v+\mathbf w]$ and $[\mathbf v-\mathbf w]$. Here $[\mathbf v]$ denotes the class of a vector up to sign.
This graph is a trivalent tree, the dual tree of the usual Farey tessellation (see~\cite[Fig.~4]{HT} for adjacency and mediants in the Farey tessellation). There are three directions from the root, and thereafter two directions other than the exchange returning to the parent.
Fractions are the \emph{components} of the triples forming vertices of the tree. The same fraction occurs at more than one vertex.

Label the root fractions $\tfrac01,\tfrac11,\tfrac10$ by $1,2,3$, respectively. When a component is exchanged, retain the labels of the other two components and give the new component the label of the one it replaces. Give the exchange edge the same label. Since the edges incident to each vertex have labels $1,2,3$, there is a unique isomorphism
\[
 \mathcal T_{\mathrm H}\cong\mathcal T_{\mathrm F}
\]
preserving the root and edge labels. The slope description in the next subsection will show that these labels agree with the peripheral labels on the sphere.

\begin{figure}[!htbp]
\centering
\begin{tikzpicture}[x=1cm,y=1cm,every node/.style={tree example vertex}]
 \node (treeRoot) at (0,0) {$\left\{\dfrac{0}{1},\dfrac{1}{1},\dfrac{1}{0}\right\}$};
 \node (treeOne) at (4,2.8) {$\left\{\boldsymbol{\dfrac{2}{1}},\dfrac{1}{1},\dfrac{1}{0}\right\}$};
 \node (treeTwo) at (4,0) {$\left\{\dfrac{0}{1},\boldsymbol{\dfrac{-1}{1}},\dfrac{1}{0}\right\}$};
 \node (treeThree) at (4,-2.8) {$\left\{\dfrac{0}{1},\dfrac{1}{1},\boldsymbol{\dfrac{1}{2}}\right\}$};
 \node (tree12) at (8,3.5) {$\left\{\dfrac{2}{1},\boldsymbol{\dfrac{3}{1}},\dfrac{1}{0}\right\}$};
 \node (tree13) at (8,2.1) {$\left\{\dfrac{2}{1},\dfrac{1}{1},\boldsymbol{\dfrac{3}{2}}\right\}$};
 \node (tree21) at (8,0.7) {$\left\{\boldsymbol{\dfrac{-2}{1}},\dfrac{-1}{1},\dfrac{1}{0}\right\}$};
 \node (tree23) at (8,-0.7) {$\left\{\dfrac{0}{1},\dfrac{-1}{1},\boldsymbol{\dfrac{-1}{2}}\right\}$};
 \node (tree31) at (8,-2.1) {$\left\{\boldsymbol{\dfrac{2}{3}},\dfrac{1}{1},\dfrac{1}{2}\right\}$};
 \node (tree32) at (8,-3.5) {$\left\{\dfrac{0}{1},\boldsymbol{\dfrac{1}{3}},\dfrac{1}{2}\right\}$};
 \treeexampleedges
\end{tikzpicture}
\caption{Farey tree}
\label{fig:Farey-tree}
\end{figure}
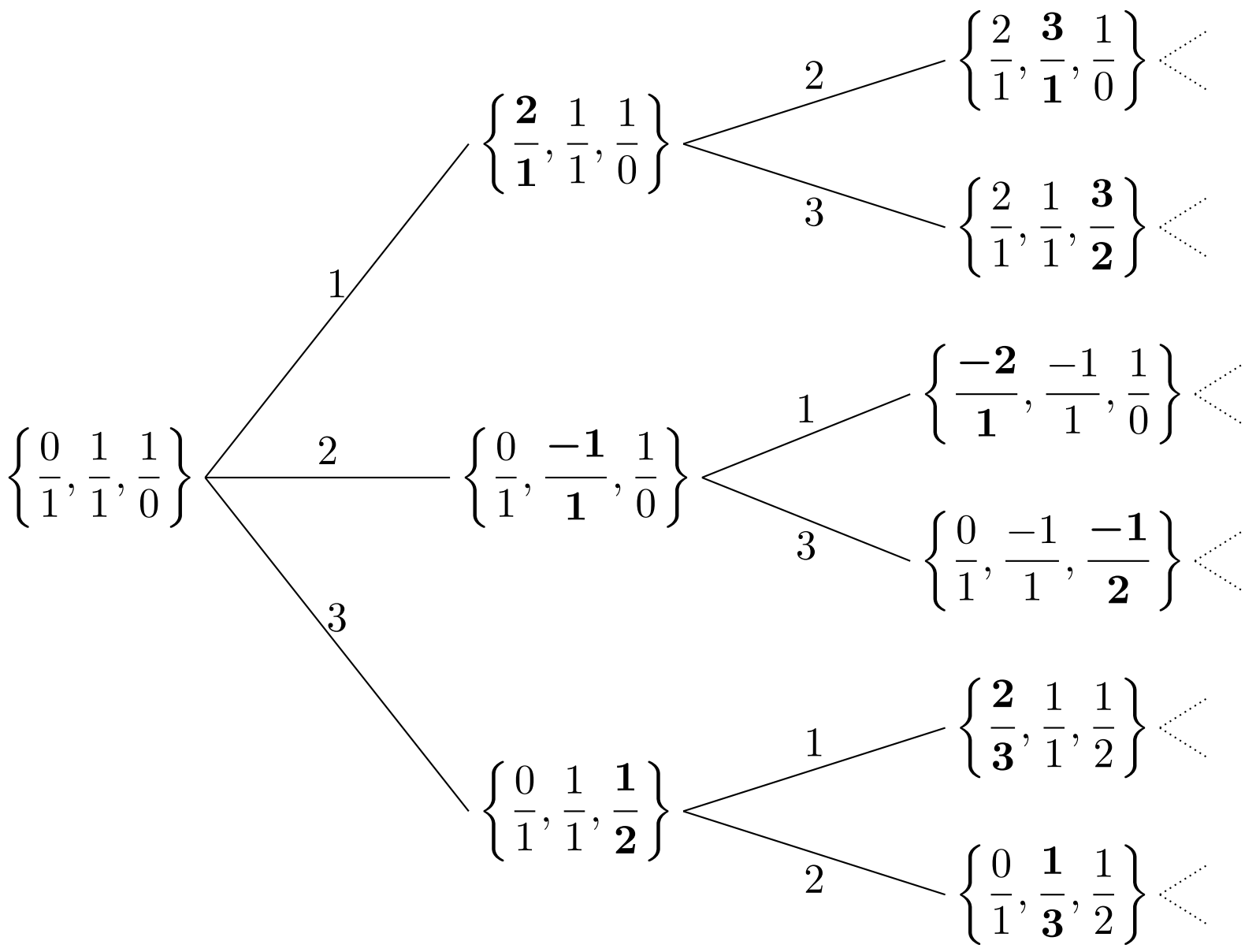

Figure~\ref{fig:Farey-tree} shows the first two levels from the root. In the figure, the components of each triple are listed in label order $1,2,3$, and exchanged fractions are bold. For example, the path with successive edge labels $1,2$ is
\[
 \left\{\frac01,\frac11,\frac10\right\}
 \longrightarrow\left\{\frac21,\frac11,\frac10\right\}
 \longrightarrow\left\{\frac21,\frac31,\frac10\right\}
\]
and corresponds to the path with successive edge labels $1,2$ in Figure~\ref{fig:Hurwitz-tree}.

\Needspace{6\baselineskip}
\begin{lemma}\label{lem:farey-fan}
The Farey triples containing a fixed fraction $[\mathbf v]$ form a bi-infinite path in $\mathcal T_{\mathrm F}$.
\end{lemma}
\begin{proof}
By B\'ezout's identity, choose an integer vector $\mathbf w$ such that $\det(\mathbf v,\mathbf w)=1$. Then $(\mathbf v,\mathbf w)$ is a basis of $\ZZ^2$. Normalize a representative $\mathbf u$ of a fraction Farey-related to $[\mathbf v]$ by $\det(\mathbf v,\mathbf u)=1$. Writing $\mathbf u=a\mathbf v+b\mathbf w$ gives $b=1$. Thus all such representatives are $\mathbf u_n=\mathbf w+n\mathbf v$ ($n\in\ZZ$). Conversely, $\det(\mathbf v,\mathbf u_n)=1$, so each $\mathbf u_n$ is primitive and satisfies the Farey relation. The sign normalization implies that the classes $[\mathbf u_n]$ are distinct.

Since $\det(\mathbf u_n,\mathbf u_m)=n-m$, the Farey triples containing $[\mathbf v]$ are precisely
\[
 \tau_n=\{[\mathbf v],[\mathbf u_n],[\mathbf u_{n+1}]\},
 \qquad n\in\ZZ.
\]
These triples are also distinct. Two distinct triples $\tau_n,\tau_m$ share two fractions exactly when $|n-m|=1$, so together they form a bi-infinite path.
\end{proof}

\subsection{Classifying ideal arcs and simple closed curves by slope}\label{sec:slopes}
We construct a bijection between reduced fractions and unoriented isotopy classes of the closed curves $\beta_e$, so that each curve in the McShane identity is counted once. Using arcs between marked points and ideal arcs based at $a_4$, we describe this correspondence by slopes in the plane. Take an auxiliary closed torus $\mathbb T=\mathbb R^2/\mathbb Z^2$ with the involution $\iota[\mathbf u]=[-\mathbf u]$, and let
\begin{equation}\label{eq:pillowcase-cover}
 \mathbb R^2\xrightarrow{\ \pi\ }\mathbb T
 \xrightarrow{\ \varpi\ }\mathbb T/\langle\iota\rangle\cong S^2,
 \qquad \Pi=\varpi\circ\pi.
\end{equation}
Here $\pi$ is the universal covering of the closed torus, and the composite $\Pi$ is given directly from the plane by
\[
 \Pi(\mathbf u)=\Pi(\mathbf u')
 \quad\Longleftrightarrow\quad
 \mathbf u'=\pm\mathbf u+\mathbf n\quad(\mathbf n\in\mathbb Z^2).
\]
The fixed-point set of the involution is
\[
 \mathbb T[2]=\bigl\{\pi(0,0),\pi(\tfrac12,0),\pi(\tfrac12,\tfrac12),\pi(0,\tfrac12)\bigr\}.
\]
The map $\varpi$ is a double cover branched at these four points. Set
\[
 a_4=\Pi(0,0),\quad a_1=\Pi(\tfrac12,0),\quad
 a_2=\Pi(\tfrac12,\tfrac12),\quad a_3=\Pi(0,\tfrac12).
\]
Removing the four marked points gives an unbranched double cover $\mathbb T\setminus\mathbb T[2]\to\Sigma_{0,4}$.

Set $o=\pi(0)$. On the torus $\mathbb T$ with only $o$ marked, consider arcs with both endpoints at $o$, whose interiors avoid $o$ and have no self-intersections, and which cannot be contracted to $o$. Let $\mathcal A$ be their unoriented isotopy classes under isotopies fixing $o$. Tangent directions at the endpoints are not fixed.

\begin{lemma}[{see~\cite{FH82}}]\label{lem:torus-arc-slopes}
For each reduced fraction $p/q\in\widehat{\mathbb Q}$, let $a_{p/q}$ be the image of the line segment from the origin to $(q,p)$ under the covering map $\pi:\mathbb R^2\to\mathbb T$.
The following map is a bijection.
\[
 \widehat{\mathbb Q}\longrightarrow\mathcal A,\qquad
 \frac pq\longmapsto[a_{p/q}].
\]
\end{lemma}

\begin{figure}[!htbp]
\centering
\begin{tikzpicture}[x=1cm,y=1cm,
 every node/.style={font=\small,inner sep=1.5pt},
 slope main/.style={black!90,line width=.9pt},
 slope reference/.style={black!45,line width=.45pt},
 slope arrow/.style={black!70,line width=.55pt},
 slope point label/.style={font=\footnotesize,inner sep=0pt}]
\begin{scope}[shift={(0,2.2)},x=1.15cm,y=1.15cm]
 \draw[black!12,line width=.35pt,step=.5] (0,0) grid (2.2,1.25);
 \foreach \slopeX in {0,.5,1,1.5,2}
  \foreach \slopeY in {0,.5,1}
   \fill[black!35] (\slopeX,\slopeY) circle[radius=.65pt];
 \draw[slope arrow,->] (-.1,0)--(2.3,0) node[right] {$x$};
 \draw[slope arrow,->] (0,-.1)--(0,1.5) node[above] {$y$};
 \draw[slope reference] (1,.5)--(2,1);
 \draw[slope main] (0,0)--(1,.5);
 \fill (0,0) circle[radius=1.5pt];
 \fill (1,.5) circle[radius=1.5pt];
 \fill (2,1) circle[radius=1.1pt];
 \node[slope point label,below left=3pt] at (0,0) {$0$};
 \node[slope point label,below=3pt] at (1,0) {$1$};
 \node[slope point label,below=3pt] at (2,0) {$2$};
 \node[slope point label,left=3pt] at (0,.5) {$\tfrac12$};
 \node[slope point label,left=3pt] at (0,1) {$1$};
 \node[below right=4pt] at (1,.5) {$\mathbf v/2$};
 \node[above left=3pt] at (2,1) {$\mathbf v=(2,1)$};
\end{scope}
\node at (1.3,1.05) {$[0,\mathbf v]$};

\begin{scope}[shift={(6.6,2.8)},x=1cm,y=1cm]
 \begin{scope}
 \pgfmathdeclarefunction{slparccontourx}{2}{\pgfmathparse{(1.15+#2*.40*cos(atan(4*sin(#1)/3)))*cos(#1)}}
 \pgfmathdeclarefunction{slparccontoury}{2}{\pgfmathparse{.60*(1.15+#2*.40*cos(atan(4*sin(#1)/3)))*sin(#1)
   +#2*.32*sin(atan(4*sin(#1)/3))}}
 \pgfmathdeclarefunction{slparcx}{1}{\pgfmathparse{(1.15+.40*cos(2*#1))*cos(#1)}}
 \pgfmathdeclarefunction{slparcy}{1}{\pgfmathparse{.60*(1.15+.40*cos(2*#1))*sin(#1)+.32*sin(2*#1)}}

 \node[inner sep=0pt,outer sep=0pt] at (0,0)
  {\sourcestorusimage{plain}{3.30cm}{2.50cm}};
 \foreach \slparcsign in {1,-1}
  \draw[black!55,line width=.45pt]
   plot[domain=0:360,samples=121,variable=\slparcphi]
    ({slparccontourx(\slparcphi,\slparcsign)},
     {slparccontoury(\slparcphi,\slparcsign)}) -- cycle;

 \foreach \slparcstart/\slparcend/\slparccolor in
  {115.175385/180/black!90,244.824615/360/black!45}
  \draw[\slparccolor,line width=.65pt,dash pattern=on 1.5pt off 1.7pt]
   plot[domain=\slparcstart:\slparcend,samples=91,variable=\slparcphi]
    ({slparcx(\slparcphi)},{slparcy(\slparcphi)});
 \foreach \slparcstart/\slparcend/\slparccolor in
  {180/244.824615/black!45,0/115.175385/black!90}
  \draw[\slparccolor,line width=.95pt,
   preaction={draw=white,line width=1.5pt}]
   plot[domain=\slparcstart:\slparcend,samples=101,variable=\slparcphi]
    ({slparcx(\slparcphi)},{slparcy(\slparcphi)});

 \fill (1.55,0) circle[radius=1.6pt];
 \fill (-1.55,0) circle[radius=1.6pt];
 \draw[black!65,line width=.4pt] (1.55,0)--(1.82,.28)--(1.82,.72);
 \node[font=\footnotesize,anchor=south,inner sep=1pt] at (1.82,.72) {$o$};
 \draw[black!65,line width=.4pt] (-1.55,0)--(-1.80,-.28)--(-1.80,-.77);
 \node[font=\footnotesize,anchor=north,inner sep=1pt]
  at (-1.80,-.78) {$\pi(\mathbf v/2)$};
\end{scope}
\end{scope}
\node at (6.6,4.15) {$\mathbb T$};
\node at (6.6,1.05) {$a_{1/2}=\pi([0,\mathbf v])$};

\pgfmathsetmacro{\slopeNx}{sin(40)*cos(150)}
\pgfmathsetmacro{\slopeNy}{sin(40)*sin(150)}
\pgfmathsetmacro{\slopeUx}{cos(40)*cos(150)}
\pgfmathsetmacro{\slopeUy}{cos(40)*sin(150)}
\pgfmathsetmacro{\slopeVx}{-sin(150)}
\pgfmathsetmacro{\slopeVy}{cos(150)}
\pgfmathdeclarefunction{slopex}{2}{\pgfmathparse{cos(#2)*(cos(#1)*\slopeNx-sin(#1)*\slopeUx)+sin(#2)*\slopeVx}}
\pgfmathdeclarefunction{slopey}{2}{\pgfmathparse{cos(#2)*(cos(#1)*\slopeNy-sin(#1)*\slopeUy)+sin(#2)*\slopeVy}}
\foreach \slopePanel/\slopeShift/\slopeHeight in
 {1/12/2.8,2/9.45/-.7,3/3.75/-.7} {
 \begin{scope}[shift={(\slopeShift,\slopeHeight)},x=1.20cm,y=1.20cm]
  \shade[ball color=gray!10,opacity=.28] (0,0) circle[radius=1];
  \draw[black!65,line width=.55pt] (0,0) circle[radius=1];
  \draw[black!22,dashed,line width=.4pt]
   plot[domain=0:180,samples=61,variable=\spheret]
   ({cos(\spheret)},{cos(68)*sin(\spheret)});
  \draw[black!22,line width=.4pt]
   plot[domain=180:360,samples=61,variable=\spheret]
   ({cos(\spheret)},{cos(68)*sin(\spheret)});
  \ifnum\slopePanel=1\relax
   \draw[slope main]
    plot[domain=0:60,samples=61,variable=\spheret]
    ({slopex(\spheret,0)},{slopey(\spheret,0)});
  \else
   \draw[slope reference]
    plot[domain=0:60,samples=61,variable=\spheret]
    ({slopex(\spheret,0)},{slopey(\spheret,0)});
  \fi
  \ifnum\slopePanel=2\relax
   \draw[slope main]
    plot[domain=0:360,samples=161,variable=\spheret]
    ({slopex(26+34*cos(\spheret),10*sin(\spheret))},
     {slopey(26+34*cos(\spheret),10*sin(\spheret))}) -- cycle;
  \fi
  \ifnum\slopePanel=3\relax
   \draw[slope reference,dashed]
    plot[domain=0:360,samples=161,variable=\spheret]
    ({slopex(26+34*cos(\spheret),10*sin(\spheret))},
     {slopey(26+34*cos(\spheret),10*sin(\spheret))}) -- cycle;
   \draw[slope main]
    plot[domain=0:360,samples=161,variable=\spheret]
    ({slopex(26+40*cos(\spheret),15*sin(\spheret))},
     {slopey(26+40*cos(\spheret),15*sin(\spheret))}) -- cycle;
  \fi
  \coordinate (slopeA4) at ({slopex(60,0)},{slopey(60,0)});
  \coordinate (slopeA3) at ({slopex(0,0)},{slopey(0,0)});
  \coordinate (slopeA1) at (.48,.54);
  \coordinate (slopeA2) at (.48,-.54);
  \foreach \slopePoint in {slopeA1,slopeA2,slopeA3,slopeA4}
   \fill (\slopePoint) circle[radius=1.5pt];
  \node[slope point label,right=3pt] at (slopeA1) {$a_1$};
  \node[slope point label,right=3pt] at (slopeA2) {$a_2$};
  \node[slope point label,above=9pt,xshift=2pt] at (slopeA3) {$a_3$};
  \node[slope point label,right=7pt] at (slopeA4) {$a_4$};
 \end{scope}
}
\node at (12,1.05) {$\eta_{1/2}$};
\node at (9.45,-2.35) {$e_{1/2}$};
\node at (3.75,-2.35) {$\beta_{e_{1/2}}=\partial N(\eta_{1/2})$};
\draw[slope arrow,->] (3.3,2.8)--(4.7,2.8)
 node[midway,above=4pt] {$\pi$};
\draw[slope arrow,->] (8.6,2.8)--(10.35,2.8)
 node[midway,above=4pt] {$\varpi$};
\draw[slope arrow,<->] (11.25,1.70)--(10.25,.50);
\draw[slope arrow,<->] (7.9,-.7)--(5.3,-.7);
\end{tikzpicture}
\caption{Reduced fractions and closed curves ($t=\tfrac12$)}
\label{fig:slope-return-construction}
\end{figure}
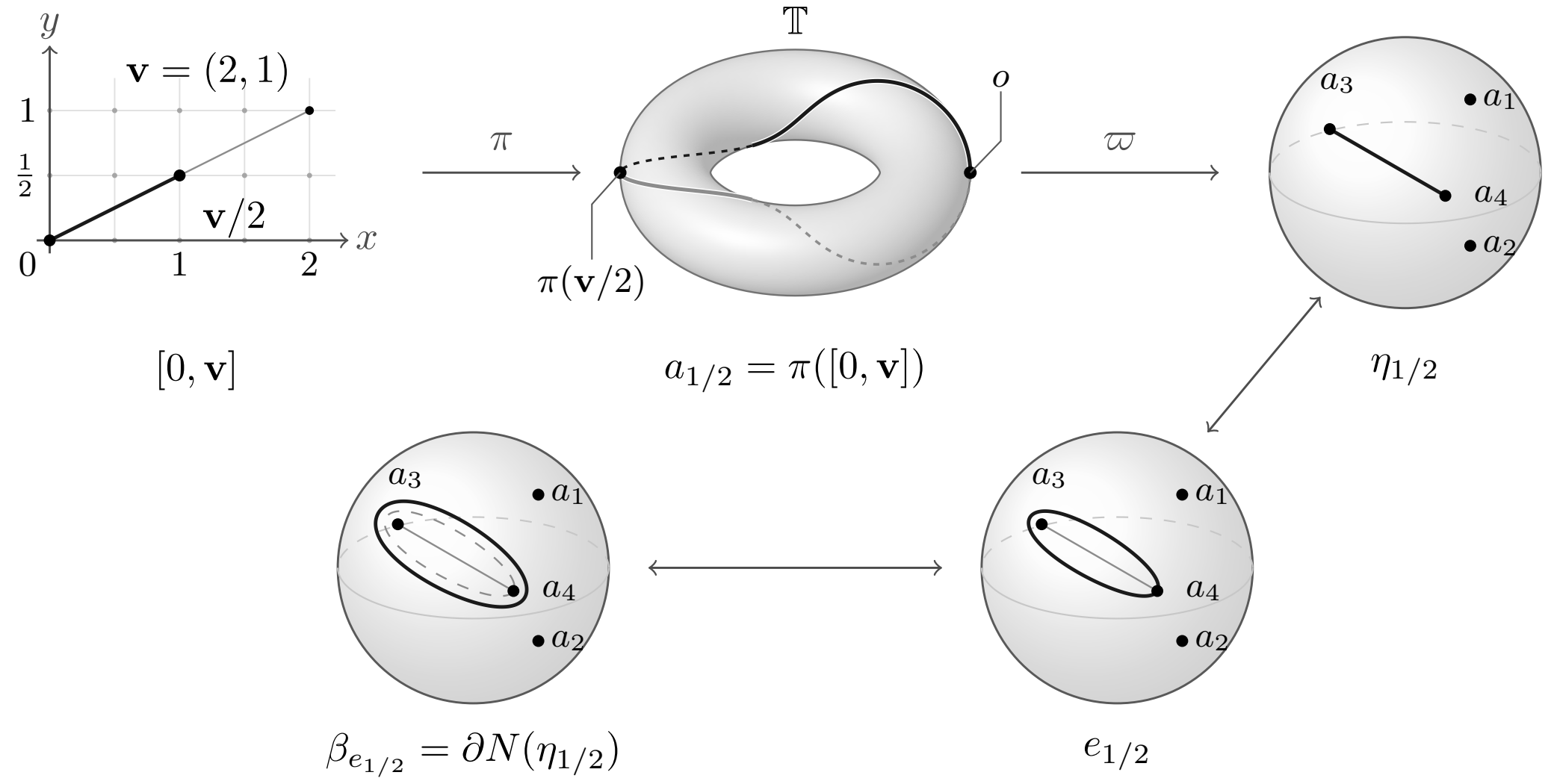

For a reduced fraction $t=p/q$, put $\mathbf v=(q,p)$. The standard arc $a_t$ in the lemma is preserved by $\iota$, which fixes $o$ and the midpoint $\pi(\mathbf v/2)$. The image of its first half on the sphere,
\[
 \eta_t=\Pi([0,\mathbf v/2])
\]
is a simple arc joining $a_4$ to one of the $a_i$. Let $e_t$ be the class of the ideal arc that follows $\eta_t$ from $a_4$, goes once around $a_i$, and returns. Take the two paths slightly apart on opposite sides of $\eta_t$, so that they are disjoint except at $a_4$. The corresponding closed curve is $\beta_{e_t}=\partial N(\eta_t)$ (Figure~\ref{fig:slope-return-construction}).

\begin{proposition}\label{prop:slope-bijection}
The correspondence
\begin{equation}\label{eq:slope-curve}
 \frac pq\ \longleftrightarrow\ [\beta_{e_{p/q}}]
\end{equation}
is a bijection between $\widehat{\mathbb Q}$ and the unoriented isotopy classes of essential nonperipheral simple closed curves on $\Sigma_{0,4}$.
\end{proposition}
\begin{proof}
We construct a bijection $J:\EE\to\mathcal A$ from isotopy classes of ideal arcs to isotopy classes of arcs on the torus, to compose with the bijection in Lemma~\ref{lem:torus-arc-slopes}.
First, construct the correspondence. For $e\in\EE$, choose a simple arc $\eta_e$ joining $a_4$ to the single marked point $a_i$ on the side of $e$ containing only that point. Separating the two paths traveling out and back along $\eta_e$ gives $e$, and $\partial N(\eta_e)$ is $\beta_e$. These operations can be reversed up to isotopy, so the isotopy classes of $\eta_e,e,\beta_e$ correspond bijectively.

$C_e=\varpi^{-1}(\eta_e)$ is a simple closed curve through $o$. Let $D_+,D_-$ be the two disks bounded by $\beta_e$, each containing two marked points; take $D_+$ to contain $a_4,a_i$ and $D_-$ to contain the other marked points $a_j,a_\ell$. Cut each disk along an arc joining its two marked points and glue two copies crosswise along the cuts. Each inverse image $\widehat D_\pm$ is then an annulus. Since $\widehat D_+$ is a regular neighborhood of $C_e$ and its complement $\widehat D_-$ is also an annulus, $C_e$ is essential. Thus, regarding $o$ as both endpoints of $C_e$ defines an element of $\mathcal A$. An isotopy of $\eta_e$ fixing the marked points lifts to the torus, so $J(e)=[C_e]$ is independent of the representative.

$\varpi^{-1}(\beta_e)=\alpha_1\sqcup\alpha_2$ consists of the two boundary components of $\widehat D_+$. The involution $\iota([\mathbf u])=[-\mathbf u]$ of the cover $\varpi$ exchanges the two copies covering each disk $D_\pm$ and hence exchanges $\alpha_1,\alpha_2$, obtained from their outer boundaries.
Figure~\ref{fig:return-arc-cover-construction} shows the construction of $J(e)$ from an ideal arc.

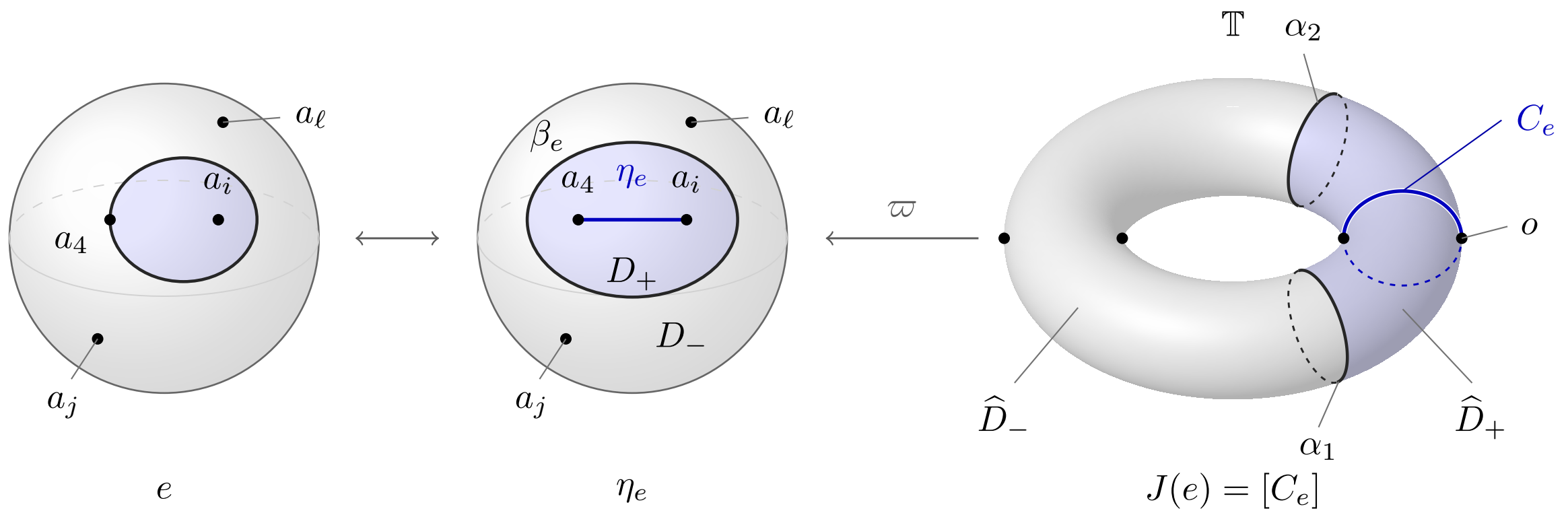
\begin{figure}[!htbp]
\centering
\begin{tikzpicture}[x=1cm,y=1cm,
 every node/.style={font=\small,inner sep=1.8pt},
 return cover curve/.style={black!85,line width=.85pt},
 return cover bridge/.style={blue!75!black,line width=1.0pt},
 return cover hidden/.style={line width=.55pt,dash pattern=on 1.5pt off 1.8pt},
 return cover leader/.style={black!55,line width=.4pt},
 return cover label/.style={font=\small}]

\foreach \retPanel/\retShift in {1/1.65,2/6.25} {
 \begin{scope}[shift={(\retShift,.8)},x=1.52cm,y=1.52cm]
  \ifnum\retPanel=1\relax
   \fill[blue!12] (.125,.12) ellipse[x radius=.475,y radius=.40];
  \else
   \fill[blue!12] (0,.12) ellipse[x radius=.68,y radius=.50];
  \fi
  \shade[ball color=gray!10,opacity=.24] (0,0) circle[radius=1];
  \draw[black!60,line width=.5pt] (0,0) circle[radius=1];
  \draw[black!18,dashed,line width=.35pt]
   plot[domain=0:180,samples=41,variable=\retAngle]
    ({cos(\retAngle)},{cos(68)*sin(\retAngle)});
  \draw[black!18,line width=.35pt]
   plot[domain=180:360,samples=41,variable=\retAngle]
    ({cos(\retAngle)},{cos(68)*sin(\retAngle)});
  \ifnum\retPanel=1\relax
   \draw[return cover curve] (.125,.12) ellipse[x radius=.475,y radius=.40];
   \node[return cover label,anchor=east] at (-.44,-.03) {$a_4$};
   \node[return cover label,anchor=south] at (.35,.24) {$a_i$};
  \else
   \draw[return cover curve] (0,.12) ellipse[x radius=.68,y radius=.50];
   \draw[return cover bridge] (-.35,.12)--(.35,.12);
   \node[blue!75!black,return cover label] at (0,.40) {$\eta_e$};
   \node[return cover label] at (-.55,.66) {$\beta_e$};
   \node[return cover label] at (0,-.23) {$D_+$};
   \node[return cover label] at (.32,-.65) {$D_-$};
   \node[return cover label,anchor=south] at (-.35,.25) {$a_4$};
   \node[return cover label,anchor=south] at (.35,.25) {$a_i$};
  \fi
  \foreach \retPoint in {(-.35,.12),(.35,.12),(-.43,-.65),(.38,.75)}
   \fill \retPoint circle[radius=1.6pt];
  \draw[return cover leader] (-.43,-.65)--(-.60,-.91);
  \node[return cover label,anchor=north] at (-.65,-.95) {$a_j$};
  \draw[return cover leader] (.38,.75)--(.75,.78);
  \node[return cover label,anchor=west] at (.80,.78) {$a_\ell$};
 \end{scope}
}

\pgfmathdeclarefunction{retcoverx}{2}{\pgfmathparse{(1.15+.40*cos(#1))*cos(#2)}}
\pgfmathdeclarefunction{retcovery}{2}{\pgfmathparse{.60*(1.15+.40*cos(#1))*sin(#2)+.32*sin(#1)}}
\begin{scope}[shift={(12.15,.8)},x=1.45cm,y=1.45cm]
 \node[inner sep=0pt,outer sep=0pt] at (0,0)
  {\sourcestorusimage{band}{4.785cm}{3.625cm}};
 \foreach \retPhi/\retStart/\retEnd in
  {-60/-49.106605/130.893395,60/49.106605/229.106605} {
  \draw[return cover curve,return cover hidden]
   plot[domain=\retEnd:{\retStart+360},samples=91,variable=\retTheta]
    ({retcoverx(\retTheta,\retPhi)},{retcovery(\retTheta,\retPhi)});
 }
 \draw[return cover bridge,return cover hidden]
  plot[domain=180:360,samples=91,variable=\retTheta]
   ({retcoverx(\retTheta,0)},{retcovery(\retTheta,0)});
 \foreach \retPhi/\retStart/\retEnd in
  {-60/-49.106605/130.893395,60/49.106605/229.106605} {
  \draw[return cover curve]
   plot[domain=\retStart:\retEnd,samples=91,variable=\retTheta]
    ({retcoverx(\retTheta,\retPhi)},{retcovery(\retTheta,\retPhi)});
 }
 \draw[return cover bridge,preaction={draw=white,line width=1.6pt}]
  plot[domain=0:180,samples=91,variable=\retTheta]
   ({retcoverx(\retTheta,0)},{retcovery(\retTheta,0)});
 \foreach \retFixed in {-1.55,-.75,.75,1.55}
  \fill (\retFixed,0) circle[radius=1.6pt];
 \draw[blue!65!black,line width=.4pt] (1.15,.32)--(1.82,.80);
 \node[blue!75!black,return cover label,anchor=west] at (1.87,.80) {$C_e$};
 \draw[return cover leader] (1.55,0)--(1.86,.08);
 \node[return cover label,anchor=west] at (1.90,.08) {$o$};
 \draw[return cover leader] (.71642,-.97077)--(.58,-1.28);
 \node[return cover label,anchor=north] at (.58,-1.31) {$\alpha_1$};
 \draw[return cover leader] (.575,.91756)--(.48,1.25);
 \node[return cover label,anchor=south] at (.48,1.29) {$\alpha_2$};
 \draw[return cover leader] (1.15,-.48)--(1.62,-.98);
 \node[return cover label,anchor=north] at (1.68,-1.02) {$\widehat D_+$};
 \draw[return cover leader] (-1.05,-.47)--(-1.48,-.98);
 \node[return cover label,anchor=north] at (-1.55,-1.02) {$\widehat D_-$};
\end{scope}

\node at (12.15,2.91) {$\mathbb T$};
\node at (1.65,-1.68) {$e$};
\node at (6.25,-1.68) {$\eta_e$};
\node at (12.15,-1.68) {$J(e)=[C_e]$};
\draw[<->,black!65,line width=.55pt] (3.53,.8)--(4.35,.8);
\draw[->,black!65,line width=.55pt] (9.65,.8)--(8.15,.8)
 node[midway,above=4pt] {$\varpi$};
\end{tikzpicture}
\caption{Ideal arcs and arcs on the torus}
\label{fig:return-arc-cover-construction}
\end{figure}

Next, we prove surjectivity of $J$. By Lemma~\ref{lem:torus-arc-slopes}, every element of $\mathcal A$ is represented by a standard arc $a_t=\pi([0,(q,p)])$. For the ideal arc $e_t$ constructed above,
\[
 \varpi^{-1}(\eta_t)=a_t,\qquad J(e_t)=[a_t],
\]
so $J$ is surjective.

Finally, we prove injectivity of $J$. Suppose $J(e)=J(e')$. Put $\beta=\beta_e,\beta'=\beta_{e'}$ in transverse position and let $A=\varpi^{-1}(\beta)=\alpha_1\sqcup\alpha_2$ and $A'=\varpi^{-1}(\beta')=\alpha_1'\sqcup\alpha_2'$ be their inverse images under the branched double cover $\varpi:\mathbb T\to S^2$. Here $A$ is the union of the two boundary components in Figure~\ref{fig:return-arc-cover-construction}, and $A'$ is the analogous pair for $e'$. The curves $C_e,C_{e'}$ are isotopic fixing $o$, and $A,A'$ bound annular neighborhoods of these curves. Thus all the components of $A,A'$ are isotopic on the closed torus $\mathbb T$, and their algebraic intersection numbers are $0$.

If there are intersections, an annulus obtained by cutting along $A$ contains an arc of $A'$ returning to the same boundary component. Otherwise, every intersecting component would cross $\alpha_1$ always with the same sign, contradicting algebraic intersection number $0$. Among the disks bounded by such an arc and a boundary interval, choose one, denoted by $D$, that is innermost among all arcs of $A'$. By innermostness and transversality, there are no other intersections in its interior or in the interiors of its sides, and no essential closed component can lie inside the disk.

To descend this intersection removal to an isotopy of the sphere fixing all four marked points, perform corresponding deformations simultaneously in $D$ and $\iota(D)$. Their boundaries lie on distinct components and are disjoint. If the disks overlapped, one would contain the other, but applying $\iota$ would reverse this inclusion, which is impossible. Thus the two disks are disjoint and contain no fixed points. In small neighborhoods of these disks, push the $A'$ arcs across $A$, with the two deformations related by $\iota$. This decreases the number of intersections by an isotopy of the sphere fixing the four marked points (Figure~\ref{fig:equivariant-bigon-removal}). The black dots in the figure denote fixed points of the involution on the torus and the four marked points on the sphere.

\begin{figure}[!htbp]
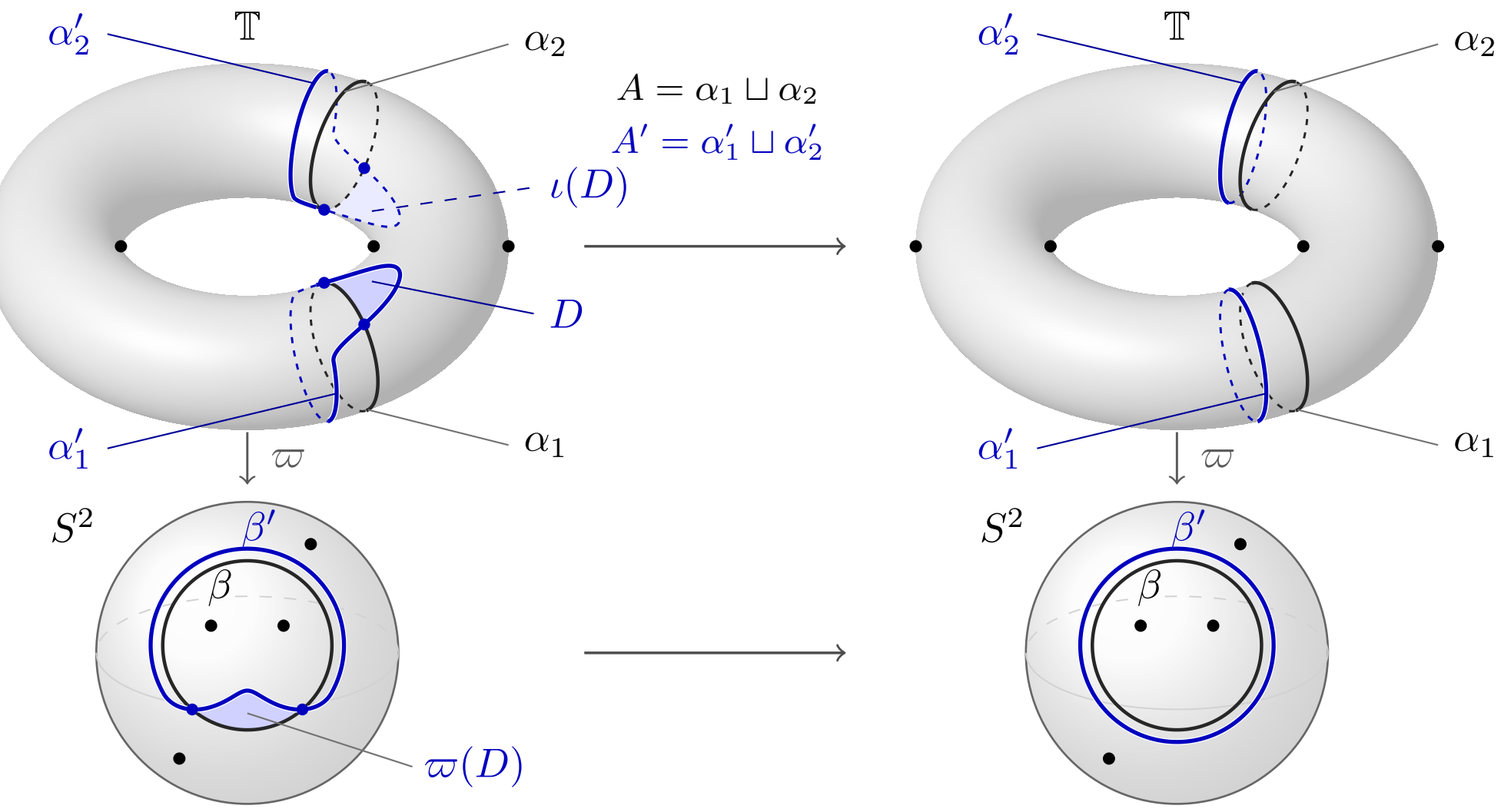

\centering

\caption{Removing intersections}
\label{fig:equivariant-bigon-removal}
\end{figure}

Repeating this operation makes $\beta,\beta'$ disjoint. They cobound an annulus containing no marked points and hence are isotopic. The correspondence between ideal arcs and closed curves then gives $e=e'$, proving injectivity of $J$.

Thus $J$ is bijective. Composing $\frac pq\leftrightarrow[a_{p/q}]$ from Lemma~\ref{lem:torus-arc-slopes}, $J^{-1}$, and $e\leftrightarrow[\beta_e]$, and using $J(e_{p/q})=[a_{p/q}]$, gives the bijection~\eqref{eq:slope-curve}.
\end{proof}

\Needspace{8\baselineskip}
The endpoint $\Pi(\mathbf v/2)$ is determined by the parity of $\mathbf v$. Its correspondence with peripheral labels is therefore
\begingroup
\setlength{\topsep}{10pt}
\renewcommand{\arraystretch}{1.2}
\begin{center}
\begin{tabular}{ccc}
\toprule
$(q,p)\bmod2$&Enclosed label $i(e_{p/q})$&Root fraction\\
\midrule
\addlinespace[4pt]
$(1,0)$&$1$&$\dfrac01$\\
\addlinespace[6pt]
$(1,1)$&$2$&$\dfrac11$\\
\addlinespace[6pt]
$(0,1)$&$3$&$\dfrac10$\\
\addlinespace[4pt]
\bottomrule
\end{tabular}
\end{center}
\endgroup
For example, $1/2$ has $(q,p)=(2,1)$, so $e_{1/2}$ is an ideal arc enclosing the third peripheral component. Each Farey triple contains one of each of the three parity types. Since $\mathbf v+\mathbf w\equiv\mathbf v-\mathbf w\pmod2$, the exchange~\eqref{eq:farey-exchange} preserves the peripheral label of the exchanged component.
Thus the peripheral labels given by parity agree with the labels propagated from the root along exchange edges.

\begin{figure}[!htbp]
\centering
\begin{tikzpicture}[x=2.15cm,y=2.15cm,
 every node/.style={font=\small,inner sep=2pt},
 root beta zero/.style={black!85,line width=1.0pt},
 root beta one/.style={blue!75!black,line width=1.0pt},
 root beta infinity/.style={red!70!black,line width=1.0pt},
 root beta point leader/.style={black!55,line width=.4pt}]
 \shade[ball color=gray!10,opacity=.28] (0,0) circle[radius=1];
 \draw[black!60,line width=.5pt] (0,0) circle[radius=1];
 \draw[black!16,dashed,line width=.35pt]
  plot[domain=0:180,samples=61,variable=\rootAngle]
   ({cos(\rootAngle)},{cos(68)*sin(\rootAngle)});
 \draw[black!16,line width=.35pt]
  plot[domain=180:360,samples=61,variable=\rootAngle]
   ({cos(\rootAngle)},{cos(68)*sin(\rootAngle)});

 \draw[root beta zero] (0,.28) circle[radius=.56];
 \draw[root beta one] (-.25,-.16) circle[radius=.56];
 \draw[root beta infinity] (.25,-.16) circle[radius=.56];

 \foreach \rootPoint in {(0,0),(0,.66),(-.58,-.36),(.58,-.36)}
  \fill \rootPoint circle[radius=1.8pt];
 \node[anchor=west] at (.08,-.02) {$a_4$};
 \draw[root beta point leader] (0,.66)--(0,1.00);
 \node[anchor=south] at (0,1.05) {$a_1$};
 \draw[root beta point leader] (-.58,-.36)--(-1.00,-.62);
 \node[anchor=east] at (-1.05,-.64) {$a_2$};
 \draw[root beta point leader] (.58,-.36)--(1.00,-.62);
 \node[anchor=west] at (1.05,-.64) {$a_3$};

 \draw[black!75,line width=.45pt] (.36,.708952)--(.95,.95);
 \node[black!85,anchor=west] at (1.00,.98) {$\beta_{e_{0/1}}$};
 \draw[blue!65!black,line width=.45pt] (-.735,.12)--(-1.10,.25);
 \node[blue!75!black,anchor=east] at (-1.15,.25) {$\beta_{e_{1/1}}$};
 \draw[red!65!black,line width=.45pt] (.735,.12)--(1.10,.25);
 \node[red!70!black,anchor=west] at (1.15,.25) {$\beta_{e_{1/0}}$};
\end{tikzpicture}
\caption{Closed curves corresponding to the root Farey triple}
\label{fig:farey-root-beta}
\end{figure}
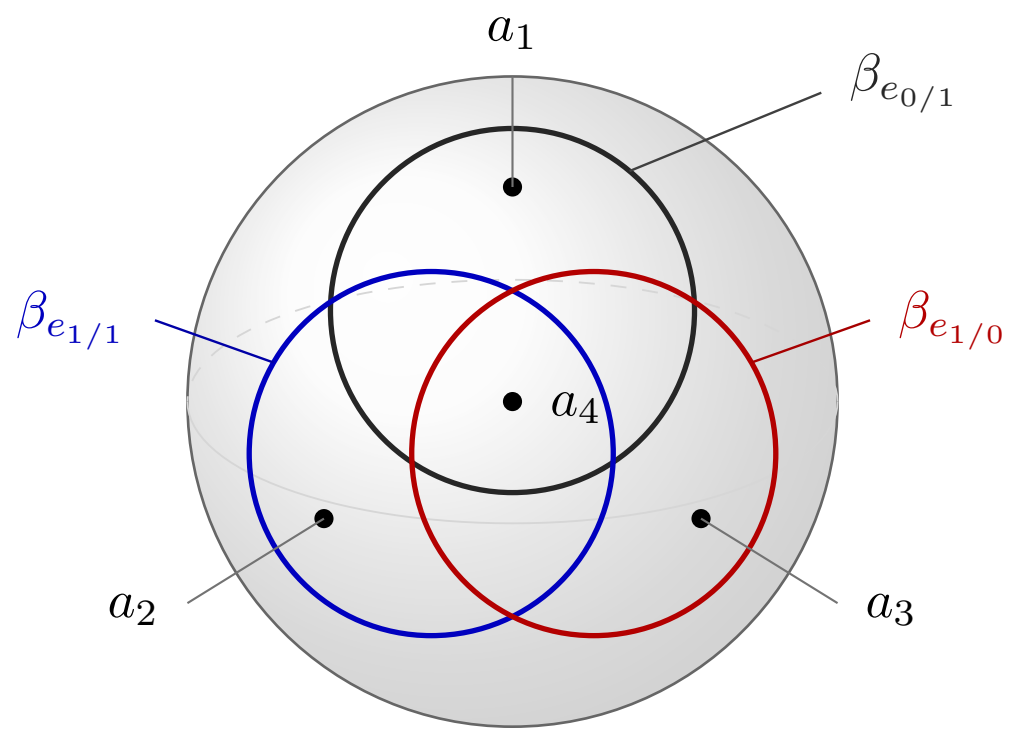

The minimal intersection number of the corresponding closed curves is
\begin{equation}\label{eq:slope-intersection}
 i(\beta_{e_{p/q}},\beta_{e_{r/s}})=2|ps-qr|.
\end{equation}
Indeed, the determinant of the period vectors on the torus gives a lower bound $|ps-qr|$ for the intersection number of each pair of lifts. Choosing parallel straight-line representatives makes each of the four pairs of lifts intersect exactly $|ps-qr|$ times, and the double cover to the sphere halves the total. The Farey relation is therefore exactly the relation that the corresponding closed curves have minimal intersection number $2$. A vertex of the Farey tree corresponds to a configuration of three ideal arcs, and an edge corresponds to exchanging one ideal arc. Figure~\ref{fig:farey-root-beta} shows the three curves $\beta_{e_{0/1}},\beta_{e_{1/1}},\beta_{e_{1/0}}$ corresponding to the root $\tau_0=\{\tfrac01,\tfrac11,\tfrac10\}$.

\subsection{Peripheral monodromy and closed geodesic lengths}

Henceforth, the same notation $\beta_e$ also denotes the element of $\pi_1(\Sigma_{0,4},b)$ obtained by choosing an orientation of the closed curve and a connecting path from the basepoint $b$, as well as its image in $\Gamma_{\mathrm{orb}}$. Changing these choices changes the element only by conjugation or inversion; the traces and group-theoretic properties used below are invariant under both operations.

Along the connecting paths of the peripheral loops corresponding to the initial triple, choose ideal arcs $e_i^0$ ($i=1,2,3$) enclosing only the $i$th peripheral component, pairwise disjoint except at the cusp $a_4$. Let $\beta_i^0=\beta_{e_i^0}=[\gamma_i\gamma_4]$ be the corresponding separating curves. Their holonomies are $-A_iB$, where $A_1=L_0,A_2=M_0,A_3=N_0$. Indeed,
\[
 (\beta_1^0)^{-1}\sim\gamma_2\gamma_3,\quad
 (\beta_3^0)^{-1}=\gamma_1\gamma_2,\quad
 (\beta_2^0)^{-1}=\gamma_1\gamma_2\gamma_3\gamma_2^{-1}.
\]
Let $\eta_i^0$ be the arc between marked points corresponding to $e_i^0$. Choose these arcs disjoint except at $a_4$, and take $\beta_i^0$ as the boundaries of thin regular neighborhoods. Any two of these boundaries then meet in two points. Disjoint essential nonperipheral simple closed curves on a sphere with four marked points are isotopic, whereas the curves $\beta_i^0$ give different partitions into two pairs. Since their intersection numbers are even, their minimal intersection numbers are $2$. By~\eqref{eq:slope-intersection}, the reduced fractions corresponding to the three curves form a Farey triple. Choose the marking of the planar model so that $e_1^0,e_2^0,e_3^0$ correspond to $\tfrac01,\tfrac11,\tfrac10$, respectively. Thus $e_1^0=e_{0/1}$, $e_2^0=e_{1/1}$, and $e_3^0=e_{1/0}$, corresponding to the root $\tau_0$ of the Farey tree.

The half-twists in~\eqref{eq:half-twist-loops} carry the connecting paths and peripheral loops together. Since they fix the fourth peripheral loop, for any loop $\alpha$ we have
\[
 h_*(\alpha\gamma_4)=h_*(\alpha)\gamma_4.
\]
Thus the corresponding separating curves and ideal arcs are carried along as well. If $\alpha_1,\alpha_2,\alpha_3$ are the three based peripheral loops in the current display order, the cyclic change of representative~\eqref{eq:cyclic-Hurwitz} corresponds to replacing them by
\[
 (\alpha_1,\alpha_2,\alpha_3)
 \longmapsto(\alpha_2,\alpha_3,\gamma_4\alpha_1\gamma_4^{-1}).
\]
The separating curve of the moved component is
\[
 (\gamma_4\alpha_1\gamma_4^{-1})\gamma_4
 =\gamma_4(\alpha_1\gamma_4)\gamma_4^{-1}
\]
and represents the same free homotopy class as the original curve. Since these are essential simple closed curves, their isotopy classes agree, so the bijection $e\leftrightarrow\beta_e$ also leaves the ideal arc unchanged. The other two components are unchanged. Hence a cyclic change of representative preserves the triple of labeled ideal arcs.

Carry the triple of ideal arcs by half-twists along the unique path from the root in the tree of Proposition~\ref{prop:Hurwitz-tree}. By Artin's action formula~\eqref{eq:half-twist-loops} and~\eqref{eq:three-Hurwitz}, $\mathcal H_r$ preserves the two loops with labels other than $r$, up to a cyclic change of representative. It therefore preserves the isotopy classes of their two closed curves. Half-twists also preserve minimal intersection numbers, so the three transformed curves still have pairwise minimal intersection number $2$, and the corresponding reduced fractions form a Farey triple.

Let $A,A'$ be the matrices with label $r$ before and after the exchange. Then $\tr A=\tr A'=-k_r$, and the holonomies of the corresponding closed curves are $-AB,-A'B$. If the third curve were also isotopic to its original curve, these holonomies would have equal traces under the fixed representation $\rho_0$. The identity $\tr(-AB)=k_r-KA_{21}$ would give $A_{21}=A'_{21}$. This contradicts Proposition~\ref{prop:Hurwitz} and the GM tree of Theorem~\ref{thm:GM-enumeration}, where the value of the exchanged component always changes.

For primitive representatives $\mathbf v,\mathbf w$ of the two shared fractions, the only possibilities for the third fraction are $[\mathbf v+\mathbf w]$ and $[\mathbf v-\mathbf w]$ in~\eqref{eq:farey-exchange}. Thus $\mathcal H_r$ crosses the Farey edge labeled $r$ and exchanges only that ideal arc. This root- and edge-label-preserving correspondence agrees with the isomorphism $\mathcal T_{\mathrm H}\cong\mathcal T_{\mathrm F}$ above. Hence every Farey vertex is realized by such a triple of ideal arcs.

For a general ideal arc $e\in\EE$, we can likewise take a matrix $A_e$ corresponding to a positively oriented peripheral loop such that
\begin{equation}\label{eq:curve-holonomy}
 \tr A_e=-k_e,\qquad
 \rho_0(\beta_e)\sim -A_e B.
\end{equation}
Here $\sim$ denotes identification up to reversal of orientation and conjugation. A lift of the geodesic realization of the marked ideal arc $f(e)$ is the geodesic joining $\infty$ to $A_e(\infty)$. If an ideal arc encloses an order-$2$ cone point, its geodesic realization folds back at that point; the ideal arc is counted as one class. The set $\EE$ itself is defined by isotopy classes on the sphere with four marked points before taking geodesic representatives.

\Needspace{5\baselineskip}
\begin{proposition}\label{prop:curve-number}
Set $m_e=(A_e)_{21}$. Then
\begin{equation}\label{eq:trace-number}
 m_e\in\ZZ_{>0},\qquad
 -\tr\rho_0(\beta_e)=\tr(A_e B)=Km_e-k_e\ge3.
\end{equation}
In particular, $m_e$ is independent of the auxiliary connecting paths. Order each Farey-tree vertex as $\{t_1,t_2,t_3\}$ with $i(e_{t_r})=r$. The three values $(x,y,z)=(m_{e_{t_1}},m_{e_{t_2}},m_{e_{t_3}})$ form a GM triple, and every GM triple occurs.
\end{proposition}
\begin{proof}
At the root, the three values for $e_i^0$ are $(1,1,1)$. Along the root- and edge-label-preserving isomorphism $\mathcal T_{\mathrm H}\cong\mathcal T_{\mathrm F}$, realize the exchange of the arc labeled $r$ by $\mathcal H_r$. Propositions~\ref{prop:Hurwitz} and~\ref{prop:Hurwitz-tree} show that the $(2,1)$ entries at every tree vertex form a positive integer GM triple, and that every GM triple is obtained.

For any $A=\left(\begin{smallmatrix}a&b\\c&d\end{smallmatrix}\right)$, we have $\tr(AB)=\tr A+Kc$. Equation~\eqref{eq:curve-holonomy} gives~\eqref{eq:trace-number}. Since this formula determines $m_e$ from the trace of the closed curve $\beta_e$, independence of the choices follows as well. Moreover,
\[
 Km_e-k_e\ge K-k_e=3+\sum_{j\ne i(e)}k_j\ge3.
\]
This proves all the assertions of the proposition.
\end{proof}

\begin{remark}
When $k_1=k_2=k_3=k$, the matrix $A_eB$ is a $k$-generalized Cohn matrix introduced by the author and Maruyama~\cite[Definition~3.1]{GyM}, under our convention using the $(2,1)$ entry. In the convention of the original paper, this is the conjugate $P(A_eB)P$, with $P$ as in Remark~\ref{rem:BG}. For general $k_1,k_2,k_3$, the matrices $A_eB$ extend that construction. These generalized Cohn matrices also play an important direct role in describing Lagrange constants; see~\cite{GyoSurvey} for the author's survey of this relationship.
\end{remark}

\Needspace{5\baselineskip}
We next relate the length of the geodesic realization of $\beta_e$ to the GM number $m_e$ (Proposition~\ref{prop:geodesic-length}). We first define generalized simple closed geodesics.

\begin{definition}[{\cite[Definition~1.7]{TWZ}}]\label{def:generalized-geodesic}
On a hyperbolic surface allowing cone points of angles at most $\pi$, cusps, and geodesic boundary, a \emph{generalized simple closed geodesic} is one of the following.
\begin{enumerate}
 \item An ordinary simple closed geodesic avoiding cone points, cusps, and boundary.
 \item A simple geodesic segment joining two cone points of angle $\pi$, traversed out and back, whose interior avoids cone points, cusps, and boundary.
 \item A cusp, a cone point, or a geodesic boundary component itself.
\end{enumerate}
The first two types are called \emph{interior generalized simple closed geodesics}. The length of a segment traversed out and back is twice the length of the segment.
\end{definition}

\begin{lemma}[{\cite[Theorem~4.1(i), Addendum~4.2(i)]{TWZ}}]\label{lem:geodesic-realization}
Let $M$ be a connected orientable hyperbolic cone surface with finitely many cone points and cusps. Assume that each cone angle lies in $(0,\pi]$, that each boundary component is a smooth simple closed geodesic, and that filling the cusps gives a compact surface. Write $M^\circ$ for the complement of the cone points, cusps, and boundary. Let $c$ be a simple closed curve in $M^\circ$ that is essential and nonperipheral there. It has a unique realization as an interior generalized simple closed geodesic, as follows.

Suppose $c$ bounds an embedded disk $D$ containing exactly two cone points, both of angle $\pi$, and no cusps or boundary components. The realization traverses a simple geodesic segment between the two cone points out and back. This segment is homotopic to a simple arc in $D$ joining the two points, with endpoints fixed and arc interiors in $M^\circ$ throughout. In all other cases the realization is an ordinary simple closed geodesic freely homotopic to $c$ in $M^\circ$.
\end{lemma}

\begin{lemma}[{see~\cite{Bea}}]\label{lem:half-turns}
Let $\sigma_P,\sigma_Q$ be rotations through angle $\pi$ about distinct points $P,Q\in\HH$. Then $\sigma_Q\sigma_P$ is hyperbolic, with axis the geodesic joining $P,Q$, and
\[
 \ell(\sigma_Q\sigma_P)=2d_{\HH}(P,Q).
\]
\end{lemma}

\begin{proposition}\label{prop:geodesic-length}
For each $e\in\EE$, the element $\beta_e$ of the orbifold group has infinite order, is nonperipheral, and is not a proper power (a power of another element with exponent at least $2$). Moreover, the marked curve $f(\beta_e)$ has a unique realization in $\OO_{\kk}$ as an interior generalized simple closed geodesic, whose length $\ell_{\beta_e}$ equals the translation length of $\bar\rho_0(\beta_e)$. Consequently,
\[
 2\cosh(\ell_{\beta_e}/2)=|\tr\rho_0(\beta_e)|=Km_e-k_e.
\]
The set $\EE$ indexes the nonperipheral interior generalized simple closed geodesics on this sphere without repetition.
\end{proposition}
\begin{proof}
Let $\{i(e),j,\ell\}=\{1,2,3\}$. Choose an orientation-preserving homeomorphism $h$ fixing the basepoint and marked points such that its induced automorphism $h_*$ of $\Gamma_{\mathrm{orb}}$ satisfies
\[
 h_*(\beta_e)\sim\gamma_j\gamma_\ell.
\]
Here $\sim$ means equality up to conjugation and reversal of orientation. We examine the properties of $w=\gamma_j\gamma_\ell$ and then apply $h_*^{-1}$ to return to the original $\beta_e$.

Write the free product in Proposition~\ref{prop:orbifold-group} as
\[
 \Gamma_{\mathrm{orb}}=G_1*G_2*G_3,\qquad
 G_r=\langle\gamma_r\rangle\quad(r=1,2,3).
\]
Each cyclic group $G_r$ making up this free product is called a free factor, and each of its nontrivial elements is counted as one syllable. For example, $\gamma_r^n$ is one syllable whenever it is nontrivial. An expression whose consecutive syllables belong to different factors is called a reduced word, and this expression is unique. If a reduced word with at least two syllables begins and ends in the same factor, conjugating to move its first syllable to the end combines these two syllables and decreases their number. Repeating this operation gives either a word of zero or one syllable, or a reduced word whose first and last syllables also belong to different factors. For a word of the latter form with at least two syllables, the number of syllables is minimal in its conjugacy class. Conjugate words of this form differ by a cyclic permutation and therefore have the same number of syllables.

First, we show that $w$ has infinite order. Since $j\ne\ell$, it is a reduced word with two syllables belonging to different factors, so its minimal syllable length in its conjugacy class is $2$. For each positive integer $n$, there is no cancellation at the joins in $w^n=(\gamma_j\gamma_\ell)^n$, which is a reduced word of $2n$ syllables. Thus $w$ has infinite order.

Next, we show that $w$ is not a proper power. Suppose $w=u^q$ with $q\ge2$. Write $u=gvg^{-1}$, where $v$ is a reduced word with the fewest syllables in its conjugacy class. Then
\[
 w=gv^qg^{-1}.
\]
If $v$ has zero or one syllable, then $v^q$ is either the identity or an element of one factor, so its minimal syllable length after conjugation is $0$ or $1$. If $v$ has $s\ge2$ syllables, its first and last syllables belong to different factors. Hence $v^q$ has no cancellation between successive copies, and its minimal syllable length in its conjugacy class is $qs\ge4$. Neither case agrees with the syllable length $2$ of $w$, a contradiction.

We also show that $w$ is not conjugate to a power of a peripheral element. A power $\gamma_r^n$ of one of the first three peripheral elements lies in one factor $G_r$, so its minimal syllable length after conjugation is $0$ or $1$. On the other hand,
\[
 \gamma_4=(\gamma_1\gamma_2\gamma_3)^{-1},
\]
so for $n\ne0$ there is no cancellation at the joins in $\gamma_4^n$, and its minimal syllable length in its conjugacy class is $3|n|$. None of these lengths equals $2$, so $w$ is not conjugate to a power of any peripheral element.

We have proved the three properties of $w$. Now use the automorphism $h_*^{-1}$ to return to the original curve. Automorphisms preserve the order of an element and whether it is a proper power. Since $h$ fixes each marked point, $h_*^{-1}$ also preserves each peripheral conjugacy class. Thus $h_*^{-1}(w)$ has infinite order, is not a proper power, and is not conjugate to a power of a peripheral element. Moreover,
\[
 \beta_e\sim h_*^{-1}(w),
\]
and these properties are also invariant under conjugation and inversion. The same conclusions therefore hold for the original $\beta_e$.

If $\OO_{\kk}$ has funnels, cut along the closed geodesic forming the entrance of each funnel and remove the funnel side. The resulting surface is orientable and becomes compact after filling its cusps. Its cone points have angle $\pi$ or $2\pi/3$ (Theorem~\ref{thm:orbifold}). The marked curve $f(\beta_e)$ separates the four peripheral components into two pairs, so it is essential and nonperipheral in the geometric interior. Lemma~\ref{lem:geodesic-realization} therefore applies.

The holonomy is hyperbolic by Proposition~\ref{prop:curve-number}. A lift of an ordinary closed geodesic is the axis of its holonomy. Since $\beta_e$ is not a proper power in the orbifold group, the length of one circuit equals the translation length of its holonomy. For a segment joining two cone points and traversed out and back, the holonomy is the product of two half-turns about the endpoints of a lift, so Lemma~\ref{lem:half-turns} gives the same conclusion. Theorem~\ref{thm:trace-length} and~\eqref{eq:trace-number} give the asserted equalities.

An ordinary geodesic determines its closed-curve class. In the degenerate case, the original closed-curve class is recovered from the boundary of a thin regular neighborhood of the segment. Since $e\leftrightarrow\beta_e$ is bijective, the indexing by $\EE$ has no repetitions.
\end{proof}

\subsection{Horocycle normalization and \texorpdfstring{$\lambda$}{lambda}-lengths}
\begin{figure}[!htbp]
\centering
\begin{tikzpicture}[x=2.2cm,y=2.2cm,
 every node/.style={font=\small,inner sep=1.5pt},
 horocycle boundary/.style={black!80,line width=.75pt},
 horocycle guide/.style={black!45,line width=.45pt,densely dashed},
 horocycle cutoff/.style={black,line width=1.5pt},
 horocycle leader/.style={black!55,line width=.4pt}]
\foreach \horopanel/\horoshift/\horodiameter in {1/0/.25,2/3.55/1.4} {
 \begin{scope}[shift={(\horoshift,0)}]
  \fill[gray!12] (-.1,1) rectangle (2.75,1.62);
  \fill[gray!12] (1.15,\horodiameter/2)
   circle[radius=\horodiameter/2];
  \ifnum\horopanel=2\relax
   \begin{scope}
    \clip (-.1,1) rectangle (2.75,1.62);
    \fill[gray!25] (1.15,\horodiameter/2)
     circle[radius=\horodiameter/2];
   \end{scope}
  \fi
  \draw[black!55,line width=.45pt,->] (-.1,0)--(2.85,0)
   node[right,font=\footnotesize] {$\mathbb R$};
  \draw[horocycle boundary] (-.1,1)--(2.75,1);
  \node[below left=3pt] at (2.75,1) {$H_1$};
  \node[left=3pt,font=\footnotesize] at (-.1,1) {$1$};
  \draw[horocycle boundary] (1.15,\horodiameter/2)
   circle[radius=\horodiameter/2];
  \draw[horocycle guide,->] (1.15,0)--(1.15,1.72)
   node[above] {$\infty$};
  \draw[horocycle cutoff] (1.15,1)--(1.15,\horodiameter);
  \foreach \horoy in {1,\horodiameter}
   \filldraw[fill=white,draw=black,line width=.65pt]
    (1.15,\horoy) circle[radius=1.5pt];
  \fill (1.15,0) circle[radius=1.4pt];
  \ifnum\horopanel=1\relax
   \node[right=7pt] at (1.15,.625) {$\delta>0$};
   \draw[horocycle leader] (1.275,.125)--(1.65,.30);
   \node[right=3pt] at (1.65,.30) {$A(H_1)$};
   \node[below=5pt] at (1.15,0) {$A(\infty)$};
  \else
   \draw[horocycle leader] (1.18,1.20)--(2.03,1.20);
   \node[right=3pt] at (2.03,1.20) {$\delta<0$};
   \node[right=3pt] at (1.85,.70) {$A(H_1)$};
   \node[below=5pt] at (1.15,0) {$A(\infty)$};
  \fi
 \end{scope}
}
\end{tikzpicture}
\caption{Horocycles and signed truncated length}
\label{fig:horocycle-truncation}
\end{figure}
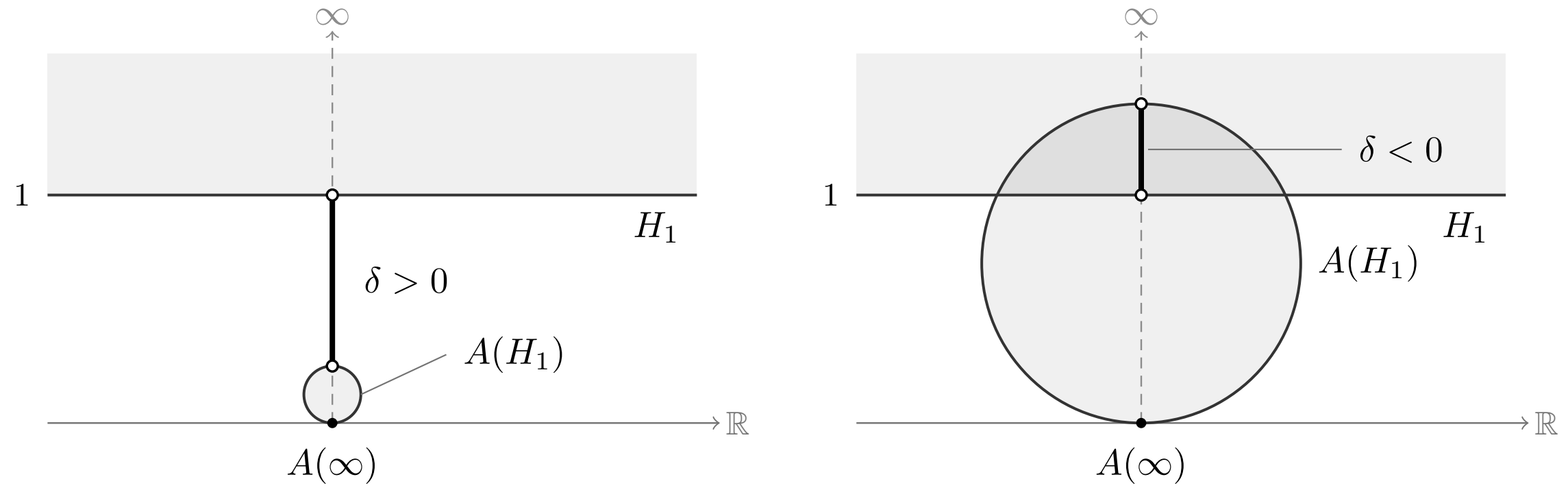

To decorate the cusp of $\OO_{\kk}$ corresponding under the marking to $a_4$, choose a horocycle at each of its lifts equivariantly under $\Lambda_{\kk}$. Here we choose $H_1=\{\operatorname{Im}z=1\}$ at $\infty$ and its $\Lambda_{\kk}$-translates at the other lifts. We do not require its image in the quotient to be embedded. One circuit corresponding to $B$ has length $K$. Let $\delta(e)$ be the signed truncated length of the geodesic in $\HH$ joining $\infty$ to $A_e(\infty)$, cut off at $H_1$ and $A_e(H_1)$. Define the $\lambda$-length (see~\cite{Pen12}) by
\[
 \lambda(e)=\exp\bigl(\delta(e)/2\bigr).
\]
When the horodisks overlap, the truncated length is taken to be negative (Figure~\ref{fig:horocycle-truncation}).

\begin{proposition}\label{prop:lambda}
For the decoration of length $K$ above, $\lambda(e)=m_e$ for every ideal arc $e$.
\end{proposition}
\begin{proof}
Let $A=\left(\begin{smallmatrix}a&b\\c&d\end{smallmatrix}\right)\in\SL_2(\RR)$ with $c\ne0$. Applying the imaginary-part transformation formula $\operatorname{Im}(A(z))=\operatorname{Im}z/|cz+d|^2$ to the inverse matrix $A^{-1}=\left(\begin{smallmatrix}d&-b\\-c&a\end{smallmatrix}\right)$ gives, for $w=u+iv\in\HH$,
\[
 \operatorname{Im}(A^{-1}(w))
 =\frac{v}{|-cw+a|^2}
 =\frac{v}{(a-cu)^2+c^2v^2}.
\]
Since $w\in A(H_1)$ is equivalent to $\operatorname{Im}(A^{-1}(w))=1$, completing the square gives
\[
 \begin{aligned}
 w\in A(H_1)
 &\Longleftrightarrow (a-cu)^2+c^2v^2=v\\
 &\Longleftrightarrow \left(u-\frac ac\right)^2+v^2=\frac{v}{c^2}\\
 &\Longleftrightarrow \left(u-\frac ac\right)^2
   +\left(v-\frac1{2c^2}\right)^2=\frac1{4c^4}.
 \end{aligned}
\]
This is the equation of a circle with center $\frac ac+\frac{i}{2c^2}$ and radius $\frac1{2c^2}$. Hence $A(H_1)$ is a horocycle of Euclidean diameter $1/c^2$, tangent to the real axis at $A(\infty)=a/c$.

The geodesic joining $\infty$ to $A(\infty)$ is the vertical line $x=a/c$. Its intersection with $H_1$ has height $1$, and its intersection with $A(H_1)$ is the top of the circle, of height $1/c^2$. On this line, the length element of the hyperbolic metric $ds^2=(dx^2+dy^2)/y^2$ is $|dy|/y$. Thus the truncated length, with the sign convention above, is
\[
 \delta=\int_{1/c^2}^{1}\frac{dy}{y}
 =\log\frac1{1/c^2}=2\log|c|.
\]
If $1/c^2>1$, the horodisks overlap and this integral is negative, in agreement with the sign convention. Taking $A=A_e$ and substituting $c=m_e>0$ from Proposition~\ref{prop:curve-number} gives $\lambda(e)=\exp(\delta(e)/2)=m_e$.
\end{proof}

\subsection{Exchanges of ideal arcs and GM exchange relations}
Fix a Farey triple, and let $(x,y,z)$ be the $\lambda$-lengths of the corresponding ideal arcs, ordered by their peripheral labels $1,2,3$. By Propositions~\ref{prop:lambda} and~\ref{prop:Hurwitz}, for example, $\mathcal H_1$ gives
\[
 (x,y,z)\longmapsto(x',y,z),\qquad
 x'=Kyz-k_2z-k_3y-x,
 \qquad xx'=y^2+k_1yz+z^2.
\]
Similarly, $\mathcal H_2,\mathcal H_3$ exchange only the second and third components, respectively. Write $\mathcal M_{\kk}$ for the set of cyclic classes of Markov--monodromy triples. For the Vieta involutions $\mu_r$ of Section~\ref{sec:Vieta}, the diagram
\begingroup
\postdisplaypenalty=10000
\[
 \begin{tikzpicture}[every node/.style={inner sep=2pt},line width=.5pt]
  \node (commMMleft) at (0,0) {$\mathcal M_{\kk}$};
  \node (commMMright) at (3.0,0) {$\mathcal M_{\kk}$};
  \node (commGMleft) at (0,-1.25) {$\mathcal G_{\kk}$};
  \node (commGMright) at (3.0,-1.25) {$\mathcal G_{\kk}$};
  \draw[->] (commMMleft)--(commMMright)
   node[midway,above=3pt] {$\mathcal H_r$};
  \draw[->] (commGMleft)--(commGMright)
   node[midway,above=3pt] {$\mu_r$};
  \draw[->] (commMMleft)--(commGMleft)
   node[midway,left=4pt] {$\Phi=\lambda$};
  \draw[->] (commMMright)--(commGMright)
   node[midway,right=4pt] {$\Phi=\lambda$};
 \end{tikzpicture}
\]
commutes. Here $\lambda$ is the map that reads the $\lambda$-lengths of the three arcs in label order.
\par\endgroup

\Needspace{5\baselineskip}
Combining the results above gives the following correspondence.

\begin{theorem}\label{thm:lambda-GM}
Choose a horocycle of length $K$ at the cusp $a_4$ of Theorem~\ref{thm:orbifold}. Each ideal arc $e\in\EE$ has an associated positive integer $m_e$ satisfying
\begin{equation}\label{eq:lambda-main}
 \lambda(e)=m_e,\qquad
 2\cosh\frac{\ell_{\beta_e}}{2}=Km_e-k_e.
\end{equation}
The three values $(x,y,z)$ corresponding to each vertex of the Farey tree form a GM triple, and every GM triple occurs. This correspondence intertwines the involutive Hurwitz move $\mathcal H_r$ associated with an edge labeled $r$ of the trivalent tree with the GM exchanges
\begin{equation}\label{eq:exchange-main}
 xx'=y^2+k_1yz+z^2,\qquad
 yy'=z^2+k_2zx+x^2,\qquad
 zz'=x^2+k_3xy+y^2.
\end{equation}
\end{theorem}
\begin{proof}
This follows from Propositions~\ref{prop:curve-number}, \ref{prop:geodesic-length}, \ref{prop:lambda}, and~\ref{prop:Hurwitz}.
\end{proof}
\section{McShane's identity and the summation formula with multiplicities}\label{sec:McShane}
\subsection{Applying the Tan--Wong--Zhang identity}
If $\OO_{\kk}$ has funnels, cut along the closed geodesic at the entrance of each funnel and remove the funnel side. We apply McShane's identity to the resulting hyperbolic orbifold with geodesic boundary. The geometric pairs of pants adjacent to the cusp $a_4$ are precisely those corresponding to ideal arcs $e\in\EE$, with the $i(e)$th peripheral component and $\beta_e$ as their other two geometric boundary components. Here cone points and cusps are allowed as boundary components. Conversely, each such pair of pants contains a unique unoriented ideal arc based at $a_4$. By Proposition~\ref{prop:geodesic-length}, indexing by $\EE$ counts each term exactly once. Since the sphere has four peripheral components, no such pair of pants has two interior curves as its boundary components other than $a_4$.

The peripheral parameter in TWZ is $\cos(\theta_i/2)$ for a cone point of angle $\theta_i$, $1$ for a cusp, and $\cosh(P_i/2)$ for a geodesic boundary of length $P_i$. In this paper all these parameters equal $k_i/2$. Below we multiply the original $\operatorname{Gap}'$ by $2$, normalizing the sum to $1$.

\begin{lemma}[{\cite[Theorem~1.8, (10), (21)--(23)]{TWZ}}]\label{thm:TWZ}
For $i=1,2,3$ and $\ell\ge0$, define
\[
 g_i(\ell)=1-\frac{\sinh(\ell/2)}{k_i/2+\cosh(\ell/2)}.
\]
Then for the orbifold considered in this paper,
\begin{equation}\label{eq:McShane}
 \sum_{e\in\EE}g_{i(e)}(\ell_{\beta_e})=1.
\end{equation}
\end{lemma}

Using Proposition~\ref{prop:geodesic-length} to rewrite Lemma~\ref{thm:TWZ} in terms of GM numbers gives the following summation formula.

\begin{theorem}\label{thm:weighted-sum}
For the values $m_e$ of the ideal arcs in Theorem~\ref{thm:lambda-GM}, the following identity holds:
\begin{equation}\label{eq:weighted-main}
 \sum_{e\in\EE}
 \left(K-\frac{\sqrt{(Km_e-k_e)^2-4}}{m_e}\right)=K.
\end{equation}
This is an absolutely convergent series of positive terms.
\end{theorem}

\begin{proof}[Proof of Theorem~\ref{thm:weighted-sum}]
By Proposition~\ref{prop:geodesic-length},
\[
 \cosh\frac{\ell_{\beta_e}}{2}=\frac{Km_e-k_e}{2},\qquad
 \sinh\frac{\ell_{\beta_e}}{2}=\frac{\sqrt{(Km_e-k_e)^2-4}}2.
\]
Hence
\begin{equation}\label{eq:gap-number}
 g_{i(e)}(\ell_{\beta_e})=1-\frac{\sqrt{(Km_e-k_e)^2-4}}{Km_e}
 =\frac{K-\LL_{k_e}(m_e)}K.
\end{equation}
Substituting this into equation~\eqref{eq:McShane} of Lemma~\ref{thm:TWZ} and multiplying by $K$ gives equation~\eqref{eq:weighted-main}. Since every term is positive, absolute convergence follows as well.
\end{proof}

Theorem~\ref{thm:weighted-sum} counts each ideal arc once. To obtain Theorem~\ref{thm:arithmetic-sum}, we now count the arcs with a prescribed peripheral label $i$ and GM number $n$.

\begin{proposition}\label{cor:multiplicity}
The multiplicity $r_i(n)$ defined in equation~\eqref{eq:arithmetic-multiplicity} satisfies
\[
 r_i(n)=\#\{e\in\EE\mid i(e)=i,\ m_e=n\}.
\]
More precisely, the three root ideal arcs correspond to the three components of the initial solution $(1,1,1)$, and all other ideal arcs correspond bijectively to ordered GM triples other than the initial solution. In the latter correspondence, the label $i(e)$ and value $m_e$ of the arc agree with the position and value of the unique largest component of the triple. This correspondence commutes with a permutation $\sigma$ of the root labels satisfying $k_{\sigma(i)}=k_i$ ($i=1,2,3$) and the action that rearranges the components of a GM solution by the same permutation.
\end{proposition}
\begin{proof}
The isomorphisms of the three trees preserving edge labels, together with Theorem~\ref{thm:GM-enumeration}, give a bijection between the vertices of the Farey tree and ordered GM triples. At every nonroot vertex, the component newly introduced by the exchange from its parent is the unique largest component.
By Lemma~\ref{lem:farey-fan}, the vertices containing a given fraction form a connected path. Thus each fraction other than the three root fractions is introduced exactly once, at the vertex containing it that is closest to the root. Conversely, each nonroot vertex introduces one new fraction through the exchange from its parent. By Proposition~\ref{prop:slope-bijection}, assigning this fraction to the closed curve $\beta_e$ and then assigning its ideal arc $e$ to the GM triple at that vertex gives a bijection preserving the position and value of the largest component. The three root ideal arcs all have value $1$, whereas the largest value of every nonroot solution is greater than $1$. This also proves the expression for $r_i(n)$.

A permutation $\sigma$ of the labels satisfying $k_{\sigma(i)}=k_i$ ($i=1,2,3$) leaves the three root values $(1,1,1)$ unchanged. The exchange relation~\eqref{eq:exchange-main} is also invariant under simultaneously rearranging $(k_1,k_2,k_3)$ and the components of a GM solution by the same permutation. Induction on the distance from the root proves that the above correspondence commutes with these permutations.
\end{proof}

\begin{proof}[Proof of Theorem~\ref{thm:arithmetic-sum}]
Grouping the sum in Theorem~\ref{thm:weighted-sum} by $(i,n)$ according to Proposition~\ref{cor:multiplicity} gives equation~\eqref{eq:regrouped}. The original series is absolutely convergent with positive terms, so this regrouping is valid and the resulting series also converges absolutely.
\end{proof}

\section{The summation inequality and its equality conditions}\label{sec:uniqueness}
Equation~\eqref{eq:weighted-main} counts each ideal arc as one term. In this section we organize the multiplicities to prove Theorem~\ref{thm:GL}, and state the equality criterion for general $\kk$ and an injectivity conjecture.

\subsection{Peripheral permutation symmetries and the orbit sum}
Let
\[
 \PP_{\kk}=\{\sigma\in S_3\mid k_{\sigma(i)}=k_i\ (1\le i\le3)\},
 \qquad p_{\kk}=|\PP_{\kk}|
\]
be the group of permutations exchanging peripheral components with the same value of $k_i$. Here $S_3$ is the permutation group of the peripheral labels $\{1,2,3\}$. For example, if $\kk=(0,0,2)$, peripheral components $1$ and $2$ can be exchanged, so $\PP_{\kk}=\{\mathrm{id},(12)\}$. If the three $k_i$ are all distinct, only the identity is allowed; if exactly two are equal, their transposition is also allowed; and if all are equal, every permutation is allowed. Thus the corresponding values of $p_{\kk}$ are $1,2,6$.

To let these permutations act on ideal arcs, write the three root fractions as $t_1=\tfrac01,t_2=\tfrac11,t_3=\tfrac10$. For each $\sigma\in S_3$, there is a unique fractional linear transformation $g_\sigma\in\mathrm{PGL}_2(\ZZ)$ satisfying
\[
 g_\sigma(t_i)=t_{\sigma(i)}\qquad(1\le i\le3).
\]
These transformations are
\[
 \begin{aligned}
 g_{\mathrm{id}}(t)&=t,&
 g_{(123)}(t)&=\frac1{1-t},&
 g_{(132)}(t)&=\frac{t-1}{t},\\[4pt]
 g_{(12)}(t)&=1-t,&
 g_{(13)}(t)&=\frac1t,&
 g_{(23)}(t)&=\frac{t}{t-1}.
 \end{aligned}
\]
They preserve $\{\tfrac01,\tfrac11,\tfrac10\}$ as a set and permute its three elements. Acting by the same transformation on every fraction, define the action on ideal arcs by $\sigma\cdot e_t=e_{g_\sigma(t)}$. This action permutes the root labels by $\sigma$, and hence agrees with the root symmetry used in Proposition~\ref{cor:multiplicity}. For $\sigma\in\PP_{\kk}$ it preserves $k_e$ and $m_e$, so the following map is defined. We prove this preservation in Lemma~\ref{lem:symmetry}.
\begin{equation}\label{eq:spectral-map}
 F_{\kk}:\EE/\PP_{\kk}\longrightarrow\Spec_{\kk},
 \qquad [e]\longmapsto\LL_{k_e}(m_e).
\end{equation}

\Needspace{7\baselineskip}
Recall the constant $C_{\kk}$ defined in equation~\eqref{eq:Cdef}:
\[
 C_{\kk}=\begin{cases}
 K,&k_1,k_2,k_3\text{ distinct},\\[7pt]
 \displaystyle\frac{3K}{2}-
 \frac{\LL_k(1)+\LL_k(k+2)}2,
 &\text{exactly two }k_i\text{ equal; remaining value }k,\\[7pt]
 \displaystyle\frac{7K}{6}-
 \frac{\LL_k(1)+\LL_k(k+2)}2,
 &k_1=k_2=k_3=k.
 \end{cases}
\]

\begin{proposition}\label{prop:equality}
For general $\kk$, the sum over distinct values equals the upper bound $C_{\kk}$, that is,
\[
 \sum_{s\in\Spec_{\kk}}(K-s)=C_{\kk},
\]
if and only if the map $F_{\kk}$ is injective.
\end{proposition}
We will prove Proposition~\ref{prop:equality} below using Proposition~\ref{prop:orbit-sum}.

Each $g_\sigma$ preserves Farey adjacency and the root $\tau_0$, so it defines an automorphism of $\mathcal T_{\mathrm F}$ permuting the root labels by $\sigma$.

\begin{lemma}\label{lem:symmetry}
If $\sigma\in\PP_{\kk}$, then for $e\in\EE$,
\[
 i(\sigma e)=\sigma(i(e)),\qquad
 k_{\sigma e}=k_e,\qquad m_{\sigma e}=m_e.
\]
Consequently, $K-\LL_{k_e}(m_e)$ is constant on each $\PP_{\kk}$-orbit, and equation~\eqref{eq:spectral-map} defines a well-defined surjection.
\end{lemma}
\begin{proof}
The correspondence in Proposition~\ref{cor:multiplicity} commutes with permutations of the root labels and preserves the position and value of the largest component. Since the three root values are also all $1$, this gives the formulas for $i$ and $m$. The equality $k_{\sigma(i)}=k_i$ gives preservation of $k$ as well. Surjectivity follows from the enumeration in Proposition~\ref{prop:curve-number}.
\end{proof}

\begin{lemma}\label{lem:exceptions}
Suppose that $p_{\kk}>1$. Exactly two orbits in $\EE/\PP_{\kk}$ have size different from $p_{\kk}$. Both have size $p_{\kk}/2$, and the pairs $(k_e,m_e)$ for their representatives $e$ are
\[
 (k_i,1),\qquad(k_i,k_i+2).
\]
Here choose $i$ such that $\{i,j,\ell\}=\{1,2,3\}$ and $k_j=k_\ell$.
\end{lemma}
\begin{proof}
The transformations $g_{(123)},g_{(132)}$ have no fixed points on the real projective line. The fixed points of $g_{(13)}$ are $\tfrac11,\tfrac{-1}1$: one is a root fraction, and the other is the fraction obtained by the first exchange of that fraction. The transformations corresponding to the other transpositions are conjugate to this one and have the same property.

If $p_{\kk}=2$, the nontrivial permutation exchanges the two labels $j,\ell$ with $k_j=k_\ell$ and fixes two fractions of label $i$. If $p_{\kk}=6$, the fixed points of the three transpositions form two orbits: the three root fractions, and the three fractions obtained by their first exchanges. Each of these points has a stabilizer of order $2$, and all other stabilizers are trivial. The root value is $1$, and the exchanged value is $k_i+2$ by equation~\eqref{eq:exchange-main}.
\end{proof}

\begin{proposition}\label{prop:orbit-sum}
The following identity holds as a series of positive terms:
\begin{equation}\label{eq:orbit-sum}
 \sum_{[e]\in\EE/\PP_{\kk}}\bigl(K-\LL_{k_e}(m_e)\bigr)=C_{\kk}.
\end{equation}
\end{proposition}
\begin{proof}
\postdisplaypenalty=10000
For $p_{\kk}=1$, this is equation~\eqref{eq:weighted-main}. Let $p=p_{\kk}>1$, and choose $i$ such that $\{i,j,\ell\}=\{1,2,3\}$ and $k_j=k_\ell$. The spectral values on the two exceptional orbits are $\LL_{k_i}(1)$ and $\LL_{k_i}(k_i+2)$. Thus Lemma~\ref{lem:exceptions} and equation~\eqref{eq:weighted-main} give
\[
 K=p\sum_{[e]}\bigl(K-\LL_{k_e}(m_e)\bigr)
 -\frac p2\bigl(2K-\LL_{k_i}(1)-\LL_{k_i}(k_i+2)\bigr).
\]
Solving this identity gives equations~\eqref{eq:Cdef} and~\eqref{eq:orbit-sum}.
\end{proof}

For $s\in\Spec_{\kk}$, let $\nu_{\kk}(s)=\#F_{\kk}^{-1}(s)$ be the number of orbits giving the value $s$.

\begin{proof}[Proof of Proposition~\ref{prop:equality}]
Lemma~\ref{lem:symmetry} and Proposition~\ref{prop:orbit-sum} imply $1\le\nu_{\kk}(s)<\infty$. Grouping the series of positive terms by spectral value gives
\[
 C_{\kk}=\sum_{s\in\Spec_{\kk}}\nu_{\kk}(s)(K-s).
\]
Since every $K-s$ is positive, this sum equals $\sum_{s\in\Spec_{\kk}}(K-s)$ if and only if every $\nu_{\kk}(s)$ equals $1$, which is precisely the injectivity of $F_{\kk}$.
\end{proof}

\subsection{Coincident spectral values and the injectivity conjecture}
Consider the set
\[
 \mathscr A_{\kk}=\{(k_i,n)\mid n\text{ is a GM number with label }i\},
\]
in which identical pairs $(k_i,n)$ are counted only once. Proposition~\ref{prop:orbit-sum} and the positivity of every term give
\begin{equation}\label{eq:pair-sum}
 \sum_{s\in\Spec_{\kk}}(K-s)
 \le\sum_{(a,n)\in\mathscr A_{\kk}}\bigl(K-\LL_a(n)\bigr)\le C_{\kk}.
\end{equation}
Equality on the right is equivalent to injectivity of $\EE/\PP_{\kk}\to\mathscr A_{\kk}$. Equality on the left is equivalent to injectivity of $(a,n)\mapsto\LL_a(n)$ on $\mathscr A_{\kk}$. For fixed $a$,
\[
 \frac{d}{dn}\LL_a(n)^2
 =\frac{2}{n^3}\bigl(Ka n-a^2+4\bigr)>0\qquad(n\ge1),
\]
so $\LL_a(n)$ is strictly increasing. However, different values of $a$ may yield the same value. Such coincidences are described by the following proposition.

\begin{proposition}\label{prop:cross-spectral-collision}
For integers $0\le a<b$, $K\ge a+b+3$, and $n,m\ge1$,
\[
 \frac{\sqrt{(Kn-a)^2-4}}n
 =\frac{\sqrt{(Km-b)^2-4}}m
 \quad\Longleftrightarrow\quad
 a=0,\quad b=2,\quad m=Kn^2.
\]
\end{proposition}
\begin{proof}
In this proof write $\LL_c(t)=\sqrt{(Kt-c)^2-4}/t$. For $0\le c\le b$ and $t\ge1$,
\[
 \frac{\partial}{\partial t}\LL_c(t)^2
 =\frac{2}{t^3}(Kct-c^2+4)>0,\qquad
 \frac{\partial}{\partial c}\LL_c(t)^2
 =-\frac{2(Kt-c)}{t^2}<0.
\]
Thus $\LL_c(t)$ is strictly increasing in $t$ and strictly decreasing in $c$.

Suppose that equality holds, and denote the common value by $s$. The above monotonicity gives $n<m$, and
\begin{equation}\label{eq:cross-collision-s-bound}
 s\ge\LL_a(1)=\sqrt{(K-a)^2-4}
 \ge\sqrt{(b+3)^2-4}>b+2.
\end{equation}
Let $g=\gcd(n,m)$, and write $n=gu$, $m=gv$, with $\gcd(u,v)=1$ and $u<v$. Then
\[
 Kgu-a=\sqrt{s^2g^2u^2+4},\qquad
 Kgv-b=\sqrt{s^2g^2v^2+4}.
\]
Setting $\delta=bu-av$ and rationalizing gives
\[
 \delta
 =\frac{4(v^2-u^2)}
 {v\sqrt{s^2g^2u^2+4}+u\sqrt{s^2g^2v^2+4}}.
\]
The numerator is positive and the denominator is greater than $2sguv$, so $\delta$ is a positive integer and
\[
 1\le\delta<\frac{2v}{sgu}.
\]
Moreover, squaring the original equality and clearing denominators gives
\begin{equation}\label{eq:cross-collision-cleared}
 2Kguv\delta=(b^2-4)u^2-(a^2-4)v^2.
\end{equation}

First exclude $a>0$. Since $\delta=bu-av>0$, we have $v/u<b/a$, and equation~\eqref{eq:cross-collision-s-bound} gives $s>b$. Hence
\[
 1\le\delta<\frac{2b}{sga}<\frac2{ga}.
\]
This rules out $a\ge2$. If $a=1$, the inequality gives $g=\delta=1$, and hence $v=bu-1$. However, the right-hand side of equation~\eqref{eq:cross-collision-cleared} satisfies
\[
 (b^2-4)u^2+3v^2
 \equiv(bu)^2+(bu-1)^2\equiv1\pmod2,
\]
contradicting the evenness of its left-hand side. Thus $a=0$.

In this case equation~\eqref{eq:cross-collision-cleared} becomes
\[
 4v^2=u^2\{2Kbgv-(b^2-4)\}.
\]
Since $\gcd(u,v)=1$, we have $u^2\mid4$, so $u=1$ or $2$.

First exclude $b=1$. In this case, if $u=1$, the preceding equation is impossible by parity. If $u=2$, then $v^2=2Kgv+3$, so $v\mid3$. Since $v>u=2$, we obtain $v=3$. This gives $Kg=1$, contradicting $K\ge4$.

Next suppose that $b\ge3$. Viewing the preceding equation modulo $v$ gives $v\mid b^2-4$. Taking the positive integer $q=(b^2-4)/v$, we have
\[
 Kg=\frac{q+4v/u^2}{2b},\qquad qv=b^2-4.
\]
If $u=2$, the inequality $(q-1)(v-1)\ge0$ gives $q+v\le qv+1=b^2-3$, and hence
\[
 Kg\le\frac{b^2-3}{2b}<b<K,
\]
a contradiction. If $u=1$ and $q\ge2$, then $(q-2)(2v-1)\ge0$ gives $q+4v\le2qv+2=2b^2-6$, again implying $Kg<b<K$, a contradiction. In the remaining case $u=1$, $q=1$,
\[
 Kg=\frac{1+4(b^2-4)}{2b}=2b-\frac{15}{2b}.
\]
But the even integer $2b$ cannot divide $15$, so the right-hand side is not an integer.

Thus only $b=2$ remains. In this case
\[
 \LL_0(n)^2=K^2-\frac4{n^2},\qquad
 \LL_2(m)^2=K^2-\frac{4K}{m},
\]
so equality of the two values is equivalent to $m=Kn^2$. The converse also follows immediately from these formulas.
\end{proof}

Taking $a,b$ to be two distinct values among the $k_i$ in this paper gives $K\ge a+b+3$. Consequently, GM numbers associated with two distinct values of $k_i$ can give the same Lagrange constant only when one value of $k_i$ is $0$ and the other is $2$. Such examples do occur; see Example~\ref{ex:spectral-collision}.

If $k_1=k_2=k_3=k$, the strict monotonicity of $\LL_k$ implies that the first inequality in equation~\eqref{eq:pair-sum} is always an equality. In this case the condition for the sum to equal $C_{\kk}$ can be expressed by uniqueness of a GM triple with a prescribed largest value, as follows.

\begin{lemma}\label{cor:symmetric-uniqueness}
For $\kk=(k,k,k)$, equality in equation~\eqref{eq:GL-main} holds if and only if, for each GM number $n$, there is exactly one GM triple in nondecreasing order whose largest value is $n$.
\end{lemma}
\begin{proof}
For $\kk=(k,k,k)$, the correspondence in Proposition~\ref{cor:multiplicity} commutes with the $S_3$-action. The three root ideal arcs form a single orbit, and quotienting the nonroot ordered GM triples by $S_3$ gives GM triples in nondecreasing order. Thus there is a bijection
\[
 \EE/S_3\longleftrightarrow\{\text{nondecreasing GM triples}\}
\]
preserving the value of the ideal arc and the largest value of the triple. Since $\LL_k(n)$ is strictly increasing for fixed $k$, injectivity of $F_{\kk}$ is equivalent to uniqueness of a GM triple in nondecreasing order with a prescribed largest value. Proposition~\ref{prop:equality} gives the conclusion.
\end{proof}

\begin{proof}[Proof of Theorem~\ref{thm:GL}]
Grouping the sum in Proposition~\ref{prop:orbit-sum} by spectral value gives
\begin{equation}\label{eq:defect}
 D_{\kk}=C_{\kk}-\sum_{s\in\Spec_{\kk}}(K-s)
 =\sum_{s\in\Spec_{\kk}}(\nu_{\kk}(s)-1)(K-s)\ge0.
\end{equation}
Thus equation~\eqref{eq:GL-main} holds, and the series of positive terms converges absolutely.

If $p_{\kk}=1$, equation~\eqref{eq:Cdef} gives $C_{\kk}=K$. If $p_{\kk}>1$, at least two of the root ideal arcs belong to the same orbit. Since every term is positive, the orbit sum in Proposition~\ref{prop:orbit-sum} is strictly smaller than the sum $K$ in equation~\eqref{eq:weighted-main}. When $k_1=k_2=k_3$, the equality criterion follows from Proposition~\ref{prop:equality} and Lemma~\ref{cor:symmetric-uniqueness}.

\nopagebreak[4]
For $\kk=\mathbf0$, we have $K=3$, and equation~\eqref{eq:Cdef} gives $C_{\mathbf0}=\frac72-\frac{\sqrt5+\sqrt8}{2}=4-\varphi-\sqrt2$.
\end{proof}

Equation~\eqref{eq:defect} gives the interpretation of $D_{\kk}$ introduced in Section~\ref{sec:intro-main}: for each spectral value $s$, it counts the $\nu_{\kk}(s)-1$ excess orbits with weight $K-s$. By Proposition~\ref{prop:equality}, the equality $D_{\kk}=0$ is equivalent to injectivity of $F_{\kk}$.

\Needspace{5\baselineskip}
To formulate the arithmetic uniqueness in Lemma~\ref{cor:symmetric-uniqueness} in terms of orbits of ideal arcs, we introduce the following injectivity conjecture.

\begin{conjecture}\label{conj:injectivity}
For every $k\in\ZZ_{\ge0}$, the map in equation~\eqref{eq:spectral-map},
\[
 F_{(k,k,k)}:\EE/S_3\longrightarrow\Spec_{(k,k,k)},
\]
is injective. Equivalently, if $e,e'\in\EE$ have the same spectral value, then $e'=\sigma e$ for some $\sigma\in S_3$. Since $\LL_k$ is strictly increasing for fixed $k$, this can also be written as
\[
 m_e=m_{e'}\quad\Longrightarrow\quad
 e'=\sigma e\quad\text{for some }\sigma\in S_3.
\]
\end{conjecture}

\begin{corollary}\label{cor:injectivity-arithmetic}
For each $k\in\ZZ_{\ge0}$, injectivity of $F_{(k,k,k)}$ is equivalent to the following assertion: for every $(k,k,k)$-GM number $n$, there is exactly one $(k,k,k)$-GM triple $(x,y,z)$ satisfying \mbox{$0<x\le y\le z=n$}. Consequently, Conjecture~\ref{conj:injectivity} and Conjecture~\ref{conj:GM-uniqueness} are equivalent.
\end{corollary}
\begin{proof}
This follows from Proposition~\ref{prop:equality} and Lemma~\ref{cor:symmetric-uniqueness}.
\end{proof}

The author and Maruyama~\cite[Corollary~4.2]{GyM} proved that arithmetic uniqueness is equivalent to injectivity of their labeling $t\mapsto m_t$ for reduced fractions $t\in[0,1]\cap\mathbb Q$. By Theorem~\ref{thm:GL}, Conjecture~\ref{conj:injectivity} is also equivalent to equality in equation~\eqref{eq:GL-main} for every $\kk=(k,k,k)$. For $k=0$, this is the Markov uniqueness conjecture, Conjecture~\ref{conj:classical-uniqueness}. Since $F_{\kk}$ need not be injective for arbitrary $\kk$, the conjecture is restricted to $k_1=k_2=k_3$; see Example~\ref{ex:spectral-collision}.

\Needspace{4\baselineskip}
\begin{example}
For $\kk=(2,2,2)$, we have $K=9$, the exceptional values are $1,4$, and
\[
 C_{(2,2,2)}=\frac96+\frac12\bigl((9-3\sqrt5)+(9-6\sqrt2)\bigr)
 =3(4-\varphi-\sqrt2).
\]
Squaring the components of a classical Markov triple gives a $(2,2,2)$-GM triple~\cite[Theorem~11]{GM}. The associated Lagrange constants are $3$ times the classical ones, in agreement with this value of the constant.
\end{example}

\subsection{Improving the upper bound using coincident spectral values}\label{sec:improvement}
Equation~\eqref{eq:defect} shows that if different orbits still give the same spectral value after passing to $\PP_{\kk}$-orbits, the upper bound can be improved by the resulting excess. More precisely, if two distinct orbits share a value $s_0$, then $\nu_{\kk}(s_0)\ge2$, so
\begin{equation}\label{eq:additional-bound}
 \sum_{s\in\Spec_{\kk}}(K-s)
 \le C_{\kk}-(K-s_0)<C_{\kk}.
\end{equation}
If such repetitions are identified for several distinct values, the corresponding positive differences $K-s_0$ can each be subtracted. Below, $\mu_i$ denotes the Vieta involution in the $i$th component.

\Needspace{5\baselineskip}
\begin{example}\label{ex:spectral-collision}
Let $\kk=(0,2,4)$ and $K=9$. The exchange path
\[
 (1,1,1)\xrightarrow{\mu_1}(2,1,1)
 \xrightarrow{\mu_2}(2,9,1)
\]
gives the value $9$ in the second component. Together with the initial value $1$ in the first component, this gives
\[
 \LL_0(1)=\sqrt{77},\qquad
 \LL_2(9)=\frac{\sqrt{79^2-4}}9=\sqrt{77}.
\]
The two ideal arcs have different peripheral labels, and $p_{\kk}=1$, so they belong to distinct orbits. Removing just this one repetition gives
\[
 \sum_{s\in\Spec_{(0,2,4)}}(9-s)
 \le9-(9-\sqrt{77})=\sqrt{77}<9.
\]
Thus, even when $k_1,k_2,k_3$ are pairwise distinct, repetitions of spectral values can sometimes improve the bound $C_{\kk}=K$ in Theorem~\ref{thm:GL}. For general $\kk$, the equality criterion is injectivity of the map in equation~\eqref{eq:spectral-map}.
\end{example}

\Needspace{5\baselineskip}
\begin{example}\label{ex:partial-symmetry-improvement}
Let $\kk=(2,2,10)$ and $K=17$. Then $\PP_{\kk}=\{\mathrm{id},(12)\}$. Both exchange paths
\[
 \begin{aligned}
 (1,1,1)&\xrightarrow{\mu_3}(1,1,12)
 \xrightarrow{\mu_1}(169,1,12),\\
 (1,1,1)&\xrightarrow{\mu_1}(4,1,1)
 \xrightarrow{\mu_2}(4,25,1)
 \xrightarrow{\mu_1}(169,25,1)
 \end{aligned}
\]
introduce the largest value $169$ in the first component at the final step. The correspondence in Proposition~\ref{cor:multiplicity} commutes with permutations of labels. The third components of the two triples above are $12$ and $1$, respectively, so exchanging the first and second components cannot take one triple to the other. Thus the corresponding ideal arcs belong to distinct $\PP_{\kk}$-orbits. Their common spectral value is
\[
 s_0=\LL_2(169)
 =\frac{\sqrt{2871^2-4}}{169}
 =\frac{\sqrt{48773}}{13}.
\]
The constant in Theorem~\ref{thm:GL} is
\[
 C_{(2,2,10)}
 =\frac{51}{2}-\frac{3\sqrt5}{2}-\frac{14\sqrt3}{3}.
\]
\par\Needspace{3\baselineskip}\noindent
Hence equation~\eqref{eq:additional-bound} gives
\[
 \sum_{s\in\Spec_{(2,2,10)}}(17-s)
 \le C_{(2,2,10)}-
 \left(17-\frac{\sqrt{48773}}{13}\right)
 <C_{(2,2,10)}.
\]
This example improves the upper bound in Theorem~\ref{thm:GL} even after repetitions arising from exchanging the first and second labels have already been removed.
\end{example}

\bibliographystyle{amsalpha}
\bibliography{generalized_markov_mcshane}

@article{Artin47,
  author  = {Artin, E.},
  title   = {Theory of braids},
  journal = {Ann. of Math. (2)},
  volume  = {48},
  number  = {1},
  year    = {1947},
  pages   = {101--126},
  doi     = {10.2307/1969218},
  url     = {https://webhomes.maths.ed.ac.uk/~v1ranick/papers/artinbraids.pdf},
  note    = {\href{https://doi.org/10.2307/1969218}{doi:10.2307/1969218}}
}

@misc{BG,
  author        = {Banaian, E. and Gyoda, Y.},
  title         = {Cluster algebraic interpretation of generalized {Markov} numbers and their matrixizations},
  year          = {2025},
  eprint        = {2507.06900v3},
  archivePrefix = {arXiv},
  url           = {https://arxiv.org/abs/2507.06900v3},
  note          = {\href{https://arxiv.org/abs/2507.06900v3}{arXiv:2507.06900v3}}
}

@book{Bea,
  author    = {Beardon, A. F.},
  title     = {The Geometry of Discrete Groups},
  series    = {Graduate Texts in Mathematics},
  volume    = {91},
  publisher = {Springer-Verlag},
  year      = {1983},
  doi       = {10.1007/978-1-4612-1146-4},
  url       = {https://doi.org/10.1007/978-1-4612-1146-4},
  note      = {\href{https://doi.org/10.1007/978-1-4612-1146-4}{doi:10.1007/978-1-4612-1146-4}}
}

@article{BIW,
  author  = {Burger, M. and Iozzi, A. and Wienhard, A.},
  title   = {Surface group representations with maximal {Toledo} invariant},
  journal = {Ann. of Math. (2)},
  volume  = {172},
  number  = {1},
  year    = {2010},
  pages   = {517--566},
  doi     = {10.4007/annals.2010.172.517},
  url     = {https://annals.math.princeton.edu/2010/172-1/p10},
  note    = {\href{https://doi.org/10.4007/annals.2010.172.517}{doi:10.4007/annals.2010.172.517}}
}

@article{Cohn71,
  author  = {Cohn, Harvey},
  title   = {Representation of {Markoff}'s binary quadratic forms by geodesics on a perforated torus},
  journal = {Acta Arith.},
  volume  = {18},
  year    = {1971},
  pages   = {125--136},
  doi     = {10.4064/aa-18-1-125-136},
  url     = {https://www.impan.pl/en/publishing-house/journals-and-series/acta-arithmetica/all/18/0/97677/representation-of-markoff-s-binary-quadratic-forms-by-geodesics-on-a-perforated-torus},
  note    = {\href{https://doi.org/10.4064/aa-18-1-125-136}{doi:10.4064/aa-18-1-125-136}}
}

@article{CS,
  author        = {Chekhov, L. and Shapiro, M.},
  title         = {Teichm\"uller spaces of {Riemann} surfaces with orbifold points of arbitrary order and cluster variables},
  journal       = {Int. Math. Res. Not. IMRN},
  volume        = {2014},
  number        = {10},
  year          = {2014},
  pages         = {2746--2772},
  eprint        = {1111.3963v3},
  archivePrefix = {arXiv},
  url           = {https://arxiv.org/abs/1111.3963v3},
  note          = {\href{https://arxiv.org/abs/1111.3963v3}{arXiv:1111.3963v3}}
}

@article{FZ02,
  author        = {Fomin, S. and Zelevinsky, A.},
  title         = {The {Laurent} phenomenon},
  journal       = {Adv. in Appl. Math.},
  volume        = {28},
  number        = {2},
  year          = {2002},
  pages         = {119--144},
  doi           = {10.1006/aama.2001.0770},
  eprint        = {math/0104241},
  archivePrefix = {arXiv},
  url           = {https://doi.org/10.1006/aama.2001.0770},
  note          = {\href{https://doi.org/10.1006/aama.2001.0770}{doi:10.1006/aama.2001.0770}}
}

@article{Frob,
  author  = {Frobenius, G.},
  title   = {{\"Uber die Markoffschen Zahlen}},
  journal = {Sitzungsberichte der K\"oniglich Preussischen Akademie der Wissenschaften zu Berlin},
  year    = {1913},
  pages   = {458--487},
  url     = {https://edoc.bbaw.de/frontdoor/index/index/docId/5095},
  archiveurl = {https://archive.org/details/dli.ministry.22787},
}

@article{GL,
  author        = {Gaster, J. and Loustau, B.},
  title         = {The sum of {Lagrange} numbers},
  journal       = {Proc. Amer. Math. Soc.},
  volume        = {149},
  number        = {12},
  year          = {2021},
  pages         = {5385--5391},
  eprint        = {2008.07659},
  archivePrefix = {arXiv},
  url           = {https://arxiv.org/abs/2008.07659},
  note          = {\href{https://arxiv.org/abs/2008.07659}{arXiv:2008.07659}}
}

@misc{GMS,
  author        = {Gyoda, Y. and Maruyama, S. and Sato, Y.},
  title         = {{SL(2,Z)}-matrixizations of generalized {Markov} numbers},
  year          = {2024},
  eprint        = {2407.08203v3},
  archivePrefix = {arXiv},
  url           = {https://arxiv.org/abs/2407.08203v3},
  note          = {\href{https://arxiv.org/abs/2407.08203v3}{arXiv:2407.08203v3}}
}

@phdthesis{Gol80,
  author = {Goldman, W. M.},
  title  = {Discontinuous groups and the {Euler} class},
  school = {University of California, Berkeley},
  year   = {1980},
  url    = {https://www.math.umd.edu/~wmg/PhDthesis.pdf},
  note   = {\url{https://www.math.umd.edu/~wmg/PhDthesis.pdf}}
}

@article{GM,
  author  = {Gyoda, Y. and Matsushita, K.},
  title   = {Generalization of {Markov} {Diophantine} equation via generalized cluster algebra},
  journal = {Electron. J. Combin.},
  volume  = {30},
  number  = {4},
  year    = {2023},
  pages   = {Paper P4.10},
  doi     = {10.37236/11420},
  url     = {https://doi.org/10.37236/11420},
  note    = {\href{https://doi.org/10.37236/11420}{doi:10.37236/11420}}
}

@article{GyM,
  author        = {Gyoda, Y. and Maruyama, S.},
  title         = {Uniqueness theorem of generalized {Markov} numbers that are prime powers},
  journal       = {Integers},
  volume        = {26},
  year          = {2026},
  pages         = {Paper A91},
  eprint        = {2312.07329},
  archivePrefix = {arXiv},
  url           = {https://math.colgate.edu/~integers/aa91/aa91.pdf},
  note          = {\href{https://arxiv.org/abs/2312.07329}{arXiv:2312.07329}}
}

@misc{Gyo,
  author        = {Gyoda, Y.},
  title         = {Generalized discrete {Markov} spectra},
  year          = {2025},
  eprint        = {2512.04547v5},
  archivePrefix = {arXiv},
  url           = {https://arxiv.org/abs/2512.04547v5},
  note          = {\href{https://arxiv.org/abs/2512.04547v5}{arXiv:2512.04547v5}}
}

@misc{GyoSurvey,
  author        = {Gyoda, Y.},
  title         = {An introduction to the {Lagrange} and {Markov} spectra through the lens of generalized {Markov} numbers},
  year          = {2026},
  eprint        = {2609.28097v1},
  archivePrefix = {arXiv},
  url           = {https://arxiv.org/abs/2609.28097v1},
  note          = {\href{https://arxiv.org/abs/2609.28097v1}{arXiv:2609.28097v1}}
}

@article{FH82,
  author  = {Floyd, W. and Hatcher, A.},
  title   = {Incompressible surfaces in punctured-torus bundles},
  journal = {Topology Appl.},
  volume  = {13},
  number  = {3},
  year    = {1982},
  pages   = {263--282},
  doi     = {10.1016/0166-8641(82)90035-9},
  url     = {https://pi.math.cornell.edu/~hatcher/Papers/PuncturedTorusBundles.pdf},
  note    = {\href{https://doi.org/10.1016/0166-8641(82)90035-9}{doi:10.1016/0166-8641(82)90035-9}}
}

@article{HP10,
  author  = {Hacking, P. and Prokhorov, Y.},
  title   = {Smoothable del {Pezzo} surfaces with quotient singularities},
  journal = {Compos. Math.},
  volume  = {146},
  number  = {1},
  year    = {2010},
  pages   = {169--192},
  doi     = {10.1112/S0010437X09004370},
  url     = {https://doi.org/10.1112/S0010437X09004370},
  note    = {\href{https://doi.org/10.1112/S0010437X09004370}{doi:10.1112/S0010437X09004370}}
}

@article{HT,
  author  = {Hatcher, A. and Thurston, W.},
  title   = {Incompressible surfaces in 2-bridge knot complements},
  journal = {Invent. Math.},
  volume  = {79},
  year    = {1985},
  pages   = {225--246},
  url     = {https://pi.math.cornell.edu/~hatcher/Papers/2-bridgeKnots.pdf},
}

@book{Kat,
  author    = {Katok, S.},
  title     = {Fuchsian Groups},
  series    = {Chicago Lectures in Mathematics},
  publisher = {University of Chicago Press},
  address   = {Chicago},
  year      = {1992},
  isbn      = {978-0-226-42583-2},
  url       = {https://press.uchicago.edu/ucp/books/book/chicago/F/bo3621076.html}
}

@article{Mar79,
  author  = {Markoff, A.},
  title   = {Sur les formes quadratiques binaires ind\'efinies},
  journal = {Math. Ann.},
  volume  = {15},
  year    = {1879},
  pages   = {381--406},
  doi     = {10.1007/BF02086269},
  url     = {https://gdz.sub.uni-goettingen.de/download/pdf/PPN235181684_0015/LOG_0031.pdf},
  note    = {\href{https://doi.org/10.1007/BF02086269}{doi:10.1007/BF02086269}}
}

@article{Mar80,
  author  = {Markoff, A.},
  title   = {Sur les formes quadratiques binaires ind\'efinies (second m\'emoire)},
  journal = {Math. Ann.},
  volume  = {17},
  year    = {1880},
  pages   = {379--399},
  doi     = {10.1007/BF01446234},
  url     = {https://archive.org/details/sim_mathematische-annalen_1880_17},
  note    = {\href{https://doi.org/10.1007/BF01446234}{doi:10.1007/BF01446234}}
}

@article{McS,
  author  = {McShane, G.},
  title   = {Simple geodesics and a series constant over {Teichm\"uller} space},
  journal = {Invent. Math.},
  volume  = {132},
  number  = {3},
  year    = {1998},
  pages   = {607--632},
  doi     = {10.1007/s002220050235},
  url     = {https://doi.org/10.1007/s002220050235},
  note    = {\href{https://doi.org/10.1007/s002220050235}{doi:10.1007/s002220050235}}
}

@book{Pen12,
  author    = {Penner, R. C.},
  title     = {Decorated {Teichm\"uller} Theory},
  series    = {QGM Master Class Series},
  publisher = {European Mathematical Society},
  address   = {Z\"urich},
  year      = {2012},
  doi       = {10.4171/075},
  url       = {https://ems.press/books/qgm/106},
  note      = {\href{https://doi.org/10.4171/075}{doi:10.4171/075}}
}

@article{Propp,
  author        = {Propp, J.},
  title         = {The combinatorics of frieze patterns and {Markoff} numbers},
  journal       = {Integers},
  volume        = {20},
  year          = {2020},
  pages         = {Paper A12},
  eprint        = {math/0511633},
  archivePrefix = {arXiv},
  url           = {https://math.colgate.edu/~integers/u12/u12.pdf},
  note          = {\href{https://arxiv.org/abs/math/0511633}{arXiv:math/0511633}}
}

@book{Rat,
  author    = {Ratcliffe, J. G.},
  title     = {Foundations of Hyperbolic Manifolds},
  edition   = {2nd},
  series    = {Graduate Texts in Mathematics},
  volume    = {149},
  publisher = {Springer},
  address   = {New York},
  year      = {2006},
  doi       = {10.1007/978-0-387-47322-2},
  url       = {https://doi.org/10.1007/978-0-387-47322-2},
  note      = {\href{https://doi.org/10.1007/978-0-387-47322-2}{doi:10.1007/978-0-387-47322-2}}
}

@article{Rud,
  author  = {Rudakov, A. N.},
  title   = {The {Markov} numbers and exceptional bundles on {$\mathbf P^2$}},
  journal = {Izv. Akad. Nauk SSSR Ser. Mat.},
  volume  = {52},
  number  = {1},
  year    = {1988},
  pages   = {100--112},
  url     = {https://www.mathnet.ru/im1170},
  note    = {English transl., Math. USSR-Izv. 32 (1989), no.~1, 99--112, \href{https://doi.org/10.1070/IM1989v032n01ABEH000738}{doi:10.1070/IM1989v032n01ABEH000738}}
}

@article{Via16,
  author        = {Vianna, R.},
  title         = {Infinitely many exotic monotone {Lagrangian} tori in {$\mathbb{CP}^{2}$}},
  journal       = {J. Topol.},
  volume        = {9},
  number        = {2},
  year          = {2016},
  pages         = {535--551},
  doi           = {10.1112/jtopol/jtw002},
  eprint        = {1409.2850},
  archivePrefix = {arXiv},
  url           = {https://doi.org/10.1112/jtopol/jtw002},
  note          = {\href{https://doi.org/10.1112/jtopol/jtw002}{doi:10.1112/jtopol/jtw002}}
}

@article{TWZ,
  author        = {Tan, S. P. and Wong, Y. L. and Zhang, Y.},
  title         = {Generalizations of {McShane}'s identity to hyperbolic cone-surfaces},
  journal       = {J. Differential Geom.},
  volume        = {72},
  number        = {1},
  year          = {2006},
  pages         = {73--112},
  eprint        = {math/0404226v1},
  archivePrefix = {arXiv},
  url           = {https://arxiv.org/abs/math/0404226v1},
  note          = {\href{https://arxiv.org/abs/math/0404226v1}{arXiv:math/0404226v1}}
}
\end{document}